\documentclass[a4paper,11pt]{amsart}
\usepackage{amsmath,amsthm,amsfonts,amssymb,graphicx,psfrag,dsfont,amscd,mathrsfs,mathabx,tikz-cd}
\usepackage{color}

\usepackage{hyperref}
\usepackage[alphabetic,initials]{amsrefs}

\theoremstyle{plain}
\newtheorem{prop}{Proposition}[section]
\newtheorem{thm}[prop]{Theorem}

\newtheorem{cor}[prop]{Corollary}

\newtheorem{lemma}[prop]{Lemma}
\newtheorem{lemmq}[prop]{Lemma(?)}
\newtheorem{lemmf}[prop]{(False) Lemma}

\newtheorem{theorem}{Theorem}

\theoremstyle{definition}
\newtheorem{example}[prop]{Example}

\newtheorem{defn}[prop]{Definition}

\newtheorem{assumption}[prop]{Assumption}
\newtheorem{remark}[prop]{Remark}

\theoremstyle{remark}

\newcommand{\Aut}{\operatorname{Aut}}

\newcommand{\codim}{\operatorname{codim}}
\newcommand{\coker}{\operatorname{coker}}

\newcommand{\CRR}{{\mathcal{CR}}_\RR}

\newcommand{\defin}{\textbf}

\newcommand{\diff}{\mathscr{D}}
\newcommand{\diffaff}{\widehat{\mathscr{D}}}
\newcommand{\Diff}{\operatorname{Diff}}
\newcommand{\dist}{\operatorname{dist}}
\newcommand{\End}{\operatorname{End}}
\newcommand{\ev}{\operatorname{ev}}

\newcommand{\fix}{_{\operatorname{fix}}}

\newcommand{\Hom}{\operatorname{Hom}}

\newcommand{\Id}{\operatorname{Id}}

\newcommand{\im}{\operatorname{im}}
\newcommand{\ind}{\operatorname{ind}}

\newcommand{\Jethol}{\operatorname{Jet}_J}
\newcommand{\Jfix}{J_{\operatorname{fix}}}
\newcommand{\Jref}{J_{\operatorname{ref}}}
\newcommand{\Lin}{\mathscr L}

\newcommand{\met}{{\mathfrak{m}}}
\newcommand{\muCZ}{\mu_{\operatorname{CZ}}}
\newcommand{\Obst}{{\mathcal{O}b}}
\newcommand{\ord}{\operatorname{ord}}

\newcommand{\Petrimod}{{\mathscr{P}}}

\newcommand{\rank}{\operatorname{rank}}
\newcommand{\rk}{\operatorname{rank}}

\newcommand{\reg}{{\operatorname{reg}}}

\newcommand{\seqs}{\boldsymbol{\mathcal E}}

\newcommand{\Span}{\operatorname{Span}}

\newcommand{\transpose}{{\operatorname{T}}}

\newcommand{\univ}{{\mathscr U}}
\newcommand{\virdim}{\operatorname{vir-dim}}

\newcommand{\volume}{\mathfrak{v}}
\newcommand{\wind}{\operatorname{wind}}

\newcommand{\GL}{\operatorname{GL}}

\newcommand{\CC}{{\mathbb C}}
\newcommand{\DD}{{\mathbb D}}

\newcommand{\HH}{{\mathbb H}}
\newcommand{\KK}{{\mathbb K}}
\newcommand{\NN}{{\mathbb N}}
\newcommand{\QQ}{{\mathbb Q}}
\newcommand{\RR}{{\mathbb R}}
\newcommand{\TT}{{\mathbb T}}
\newcommand{\ZZ}{{\mathbb Z}}

\newcommand{\bB}{{\mathcal B}}

\newcommand{\eE}{{\mathcal E}}

\newcommand{\jJ}{{\mathcal J}}

\newcommand{\mM}{{\mathcal M}}
\newcommand{\nN}{{\mathcal N}}
\newcommand{\oO}{{\mathcal O}}

\newcommand{\tT}{{\mathcal T}}
\newcommand{\uU}{{\mathcal U}}
\newcommand{\vV}{{\mathcal V}}
\newcommand{\yY}{{\mathcal Y}}
\newcommand{\1}{\mathds{1}}
\newcommand{\p}{\partial}
\renewcommand{\dbar}{\bar{\partial}}

\numberwithin{equation}{section}

\definecolor{Chris}{rgb}{0.5,0.25,1}
\definecolor{blue}{rgb}{0,0,1}
\definecolor{red}{rgb}{1,0,0}
\definecolor{green}{rgb}{0,.7,0}

\title[Transversality and super-rigidity for holomorphic curves]{Transversality and super-rigidity for multiply covered holomorphic curves}

\author{Chris Wendl}
\address{Institut f\"ur Mathematik \\ 
Humboldt-Universit\"at zu Berlin \\
Unter den Linden 6 \\
10099 Berlin \\ 
Germany}
\email{wendl@math.hu-berlin.de}

\dedicatory{To Elisabeth}

\thanks{Research partially supported by a Royal Society University Research 
Fellowship, a Leverhulme Research Project Grant and the ERC grant
\textsf{TRANSHOLOMORPHIC}}

\begin{document}

\begin{abstract}
We develop new techniques to study regularity questions for moduli spaces of
pseudoholomorphic curves that are multiply covered.  Among the main results,
we show that unbranched multiple covers of closed holomorphic curves are 
generically regular, and simple index~$0$ curves in dimensions greater than
four are generically super-rigid, implying e.g.~that the Gromov-Witten invariants
of Calabi-Yau $3$-folds reduce to sums of local invariants for finite sets of
embedded curves.  We also establish partial results on super-rigidity in dimension four
and regularity of branched covers, and briefly discuss the
outlook for bifurcation analysis.  The proofs are based on a general
stratification result for moduli spaces of multiple covers, framed in terms
of a representation-theoretic
splitting of Cauchy-Riemann operators with symmetries.
\end{abstract}

\maketitle

\tableofcontents

\section{Introduction}
\label{sec:intro}

\subsubsection*{Motivation}
The issue of transversality in Gromov's theory of pseudoholomorphic curves 
\cite{Gromov} has always been problematic, and has attracted renewed
interest in recent years.  While many powerful symplectic invariants such
as Gromov-Witten theory, Hamiltonian Floer homology and symplectic field
theory are based on holomorphic curves, most of them run into severe
technical complications unless multiply covered curves can be excluded,
thus necessitating rather sophisticated techniques that typically replace
the standard nonlinear Cauchy-Riemann equation by an abstract perturbation, see
e.g.~\cites{LiTian,FukayaOno,Ruan:virtual,
Siebert:GW,CieliebakMohnke:transversality,HWZ:GW,Pardon:virtual}.
Aside from the technical challenges that these methods pose, they are
non-ideal for many applications: for instance abstract perturbations destroy
intersection theory in symplectic $4$-manifolds, and in Calabi-Yau $3$-folds
they obscure information that one might hope to find in the geometric
relationship between simple curves and their multiple covers,
as exemplified by the Gopakumar-Vafa formula 
\cites{GopakumarVafa,BryanPandharipande:BPS,PandharipandeThomas:blackbird,IonelParker:GV,DoanIonelWalpuski}.

The motivating principle of this paper is in some sense orthogonal to that of
abstract perturbations: our aim will be to extend the transversality theory for the
\emph{standard} pseudoholomorphic curve equation as far as it can reasonably be
pushed, i.e.~to prove transversality when it is possible, and in other cases
to isolate the precise phenomena which make it impossible and explain what is
true instead.  Let us start by singling out two situations in which this
program is not obviously hopeless.

\begin{example}
\label{ex:unbranched}
If $u : (\Sigma,j) \to (M,J)$ is a closed $J$-holomorphic curve and
$\varphi : (\widetilde{\Sigma},\tilde{\jmath}) \to (\Sigma,j)$ is an \emph{unbranched} cover of
closed connected Riemann surfaces with degree $d \in \NN$, then the virtual
dimensions of the moduli spaces containing $u$ and
$u \circ \varphi : (\widetilde{\Sigma},\tilde{\jmath}) \to (M,J)$, also known as the
\emph{indices} of these two curves, are related by
$$
\ind(u \circ \varphi) = d \cdot \ind(u).
$$
Since $\ind(u \circ \varphi)$ is then nonnegative whenever $\ind(u) \ge 0$,
there is no obvious reason why $u \circ \varphi$ could not achieve transversality
generically, but traditional methods in the theory of $J$-holomorphic
curves do not prove this except when $u \circ \varphi$ is simply covered,
or in certain $4$-dimensional cases
\cite{HoferLizanSikorav}, or more recently, when $\ind(u) = 0$ if a
sufficiently large space of perturbed almost complex structures is allowed
\cite{GerigWendl}.
\end{example}

\begin{example}
\label{ex:super}
Suppose $u : (\Sigma,j) \to (M,J)$ is a closed simply covered curve with index~$0$ and 
$\varphi : (\widetilde{\Sigma},\tilde{\jmath}) \to (\Sigma,j)$ is a branched cover of closed connected
Riemann surfaces with degree~$d \in \NN$ and
$Z(d\varphi) \ge 0$ as the algebraic count of branch points.  Then combining
the Riemann-Hurwitz formula
\begin{equation}
\label{eqn:RiemannHurwitz}
-\chi(\widetilde{\Sigma}) + d \cdot \chi(\Sigma) = Z(d\varphi)
\end{equation}
with the standard index formula for closed holomorphic curves gives the relation
\begin{equation}
\label{eqn:indexCover}
\ind(u \circ \varphi) = d \cdot \ind(u) - (n-3) Z(d\varphi) = -(n-3) Z(d\varphi),
\end{equation}
where $\dim_\RR M = 2n$.  This shows that $u \circ \varphi$ lives in a space of 
nonpositive virtual dimension when $\dim M \ge 6$ and thus cannot 
achieve transversality if $\varphi$ has branch points, as the space of holomorphic
branched covers then has dimension $2 Z(d\varphi) > 0$.
It is interesting however to observe that $u$ must be immersed if $J$ is generic, so
it has a well-defined normal bundle $N_u \to \Sigma$, and restricting the
linearized Cauchy-Riemann operators for $u$ and $u \circ \varphi$ to the
normal bundle and its pullback gives operators $\mathbf{D}_u^N$ and
$\mathbf{D}_{u\circ \varphi}^N$ with indices related by
$$
\ind(\mathbf{D}_{u \circ \varphi}^N) = d \cdot \ind(\mathbf{D}_u^N) -
(n-1) Z(d\varphi) = - (n-1) Z(d\varphi).
$$
The latter is always nonpositive, so $\mathbf{D}_{u \circ \varphi}^N$
can be injective, and this condition has a geometric meaning: it implies 
that $u \circ \varphi$ can never be the
limit of a sequence of somewhere injective curves
(see Proposition~\ref{prop:isolated}).  In fact, the only other
curves near $u \circ \varphi$ are other branched covers of the form
$u \circ \varphi'$ for $\varphi'$ near~$\varphi$, and the cokernels of the
operators $\mathbf{D}_{u \circ \varphi}^N$ define an obstruction bundle
over the space of branched covers which can be used to compute
Gromov-Witten invariants.  This phenomenon is known as \emph{super-rigidity},
see Definition~\ref{defn:super}.
\end{example}

Considerable interest in super-rigidity has been motivated by the study
of Gromov-Witten invariants in Calabi-Yau $3$-folds, where all moduli
spaces of holomorphic curves without marked points have virtual dimension zero.  In this case
it can be interpreted as a Morse-Bott condition for families
of ``degenerate'' (i.e.~multiply covered) curves, so that the Gromov-Witten
counts of these curves are expressed by integrating Euler classes of obstruction
bundles over finitely many such families---these integrals define the
so-called ``multiple cover contributions,'' also known as the 
\emph{local Gromov-Witten invariants} of the underlying embedded curves.  
A substantial body of results has emerged during the past two decades 
on local Gromov-Witten
invariants and their consequences for Calabi-Yau $3$-folds in the presence of the super-rigidity hypothesis,
using both algebro-geometric
\cites{Pandharipande:degenerate,BryanKatzLeung,BryanPandharipande:BPS,BryanPandharipande:TQFT,BryanPandharipande:local}
and symplectic methods \cites{LiZinger:rigidity,Zinger:comparison,DoanWalpuski:Castelnuovo}.
In spite of these developments, a general result establishing the super-rigidity 
hypothesis itself has thus far been unavailable.  In the algebraic
category it is known to hold in some cases and not in others \cite{BryanPandharipande:rigidity},
and while it was conjectured in \cite{BryanPandharipande:BPS} to hold generically in symplectic manifolds,
proofs have been found only in very special settings (e.g.~\cites{LeeParker:structure,LeeParker:obstruction}
for certain K\"ahler surfaces),
and a strategy was even outlined in \cite{LiZinger:rigidity} to \emph{disprove} the 
conjecture for higher genus curves.

\subsubsection*{Results}
The first of the main results stated in \S\ref{sec:mainTheorems} below
settles the super-rigidity question for symplectic manifolds
of dimension at least six: by Theorem~\ref{thm:super}, super-rigidity
does hold in this setting for all simple closed $J$-holomorphic curves of index~$0$ if $J$ is generic,
and it also holds in dimension four for curves of low genus.
Complementary to this, we will see in Theorem~\ref{thm:unbranched} that transversality
holds for the unbranched multiple covers in Example~\ref{ex:unbranched},
and we will also be able to prove some transversality results for
branched covers (Theorem~\ref{thm:branched}).  The actual main result of this
paper is Theorem~\ref{thm:submanifolds0}, which implies the aforementioned
results by stratifying the space of all multiply covered $J$-holomorphic
curves into smooth submanifolds, with precise formulas for their dimensions.
The dimensions are determined by a
general picture of Cauchy-Riemann type operators with symmetries described in
\S\ref{sec:theBigIdea}, which has its origins in Taubes's work on the
Gromov invariant of symplectic $4$-manifolds \cite{Taubes:counting}.  As in 
Taubes's paper, the approach adopted here also lends itself to the study of
bifurcations and wall crossing for multiple covers, 
on which we will make some brief remarks in \S\ref{sec:bifurcations}
but save the detailed examination for future work.

\subsubsection*{The difficulty}
As with any transversality result, the proof of our main theorem boils down
to establishing that a certain bounded linear operator is surjective.  The type
of operator that arises has appeared before, e.g.~in the context of
wall-crossing arguments \cites{Taubes:counting,IonelParker:GV} (see also
\cite{Eftekhary:superrigidity}), and it has previously been dealt with by
various ad hoc methods that suffice for certain specific applications, but would
not be general enough for the problems studied here.  The solution to this
difficulty is probably the most technically novel element in the present paper:
it is reduced to a local property of Cauchy-Riemann type operators known as
\emph{Petri's condition}, which involves a ``decoupling'' between the
pointwise linear dependence relations for local solutions of a 
linear Cauchy-Riemann type equation and of its formal adjoint equation.
Section~\ref{sec:unique} of this paper 
proves that Petri's condition holds generically for Cauchy-Riemann type operators,
and this should be regarded as the main step that makes all of our other results
possible.

\subsubsection*{Outlook}
While the results in this paper focus specifically on closed holomorphic curves,
there is no obvious obstruction to applying the same techniques to study
punctured curves in symplectic cobordisms.  As with \cite{Taubes:counting}
and the Gromov invariant, this can be expected to have important applications
to the foundations of Embedded Contact Homology \cite{Hutchings:lectures}, 
e.g.~for defining cobordism maps and proving invariance without reliance on
Seiberg-Witten theory.  It also raises the intriguing possibility of 
localizing (in the sense of Corollary~\ref{cor:GW} below)
and/or proving integrality results for invariants in symplectic field 
theory \cite{SFT}.
A few special cases of super-rigidity in the punctured case have
previously been observed in \cites{Wendl:automatic,Fabert:local}; those examples were restricted to
dimension four, but the results of the present article suggest
that super-rigidity is likely to be a considerably more general phenomenon.

Since the first version of this paper appeared, A.~Doan and T.~Walpuski have
initiated a program extending the equivariant
transversality methods introduced here to more general
classes of elliptic problems; see \cite{DoanWalpuski:BrillNoether}.
More recently, Bai and Swaminathan \cite{BaiSwaminathan} have
also carried out the first step in the bifurcation analysis
proposed in \S\ref{sec:bifurcations}, and applied it toward
defining an extension of Taubes's Gromov invariant to
Calabi-Yau $3$-folds.

\subsection{Super-rigidity and transversality theorems}
\label{sec:mainTheorems}

To state the main results,
assume $(M,\omega)$ is a symplectic manifold with
$$
\dim M = 2n \ge 4,
$$
and $\Jfix$ is a smooth almost complex structure that is \defin{compatible}
with~$\omega$, meaning that $\omega(\cdot,\Jfix\cdot)$ defines a Riemannian
metric on~$M$.  We fix also an open subset
$\uU \subset M$ with compact closure, and consider the space
$$
\jJ(M,\omega\,;\,\uU,\Jfix)
$$
of smooth $\omega$-compatible almost complex structures on 
$M$ that match $\Jfix$ outside of~$\uU$, with its natural $C^\infty$-topology.

\begin{remark}
The existence of a symplectic form on $M$ is not required for
any of the arguments in this paper, but we are including it in the setup
since it is important in applications---all results could alternatively be
stated and proved for the larger space of $\omega$-tame
almost complex structures, or for arbitrary almost complex structures on a
smooth (not necessarily symplectic) manifold.
\end{remark}

Following the usual convention among symplectic topologists, we will say that
a subset of a topological space is a \defin{Baire subset} if it is comeager,
i.e.~it is a countable intersection of open and dense subsets.  The intersection
of a countable sequence of Baire subsets is again a Baire subset, and by the
Baire category theorem, any Baire subset of a complete metric space is dense.
We will say that a given property is true \defin{generically} (e.g.~for
generic~$J$) whenever there exists a Baire subset of the space of all
admissible data (e.g.~in $\jJ(M,\omega\,;\,\uU,\Jfix)$) such that the
property holds for all choices of data in that subset.

Given $J \in \jJ(M,\omega\,;\,\uU,\Jfix)$, a closed connected Riemann surface
$(\Sigma,j)$ and a $J$-holomorphic curve
$u : (\Sigma,j) \to (M,J)$, the \defin{index} of $u$ is the integer
\begin{equation}
\label{eqn:index}
\ind(u) = (n-3) \chi(\Sigma) + 2 c_1(u),
\end{equation}
where we abbreviate $c_1(u) := \langle c_1(TM,J), [u] \rangle$,
$[u] := u_*[\Sigma] \in H_2(M)$.  A closed and connected $J$-holomorphic curve
$\tilde{u} : (\widetilde{\Sigma},\tilde{\jmath}) \to (M,J)$ is said to be
a ($d$-fold) \defin{multiple cover} of $u$ if $\tilde{u} = u \circ \varphi$
for some holomorphic map $\varphi : (\widetilde{\Sigma},\tilde{\jmath}) \to
(\Sigma,j)$ of degree $d \ge 2$, and $u$ is called \defin{simple} if it is
nonconstant and is not a multiple cover of any other curve.  

The notion of super-rigidity was outlined already in Example~\ref{ex:super};
see Definition~\ref{defn:super} for a more precise formulation.  We will also use
the term \emph{Fredholm regular} to refer to the standard notion of transversality
for moduli spaces of unparametrized $J$-holomorphic curves,
cf.~Proposition~\ref{prop:regular} below.
In each of the following theorems, $(M,\omega)$ is a symplectic manifold
of dimension $2n$ with a compatible almost complex structure $\Jfix$, and
$\uU \subset M$ is an open subset with compact closure.

\begin{theorem}[super-rigidity]
\label{thm:super}
If $\dim M \ge 6$, then
there exists a Baire subset $\jJ^\reg$ of the space $\jJ(M,\omega\,;\,\uU,\Jfix)$
such that for all $J \in \jJ^\reg$, every simple $J$-holomorphic curve
of index~$0$ that intersects $\uU$ is super-rigid.  Moreover, this result
also holds when $\dim M = 4$ for all simple index~$0$ curves of genus~$0$ or~$1$.
\end{theorem}

Super-rigidity has a number of well-known consequences, which are especially 
important in the case $\dim M = 6$.  These are based partly on the observation that
the space of all covers of super-rigid curves is an open and closed subset of
the ambient moduli space of $J$-holomorphic curves, see Proposition~\ref{prop:isolated}
in Appendix~\ref{sec:isolated}.
Applying Gromov compactness and the standard implicit function theorem
for simple curves, plus the fact that simple $J$-holomorphic curves of
index~$0$ are
generically embedded and disjoint from each other in dimensions greater
than four, this implies:

\begin{cor}
\label{cor:dimension6}
For generic compatible $J$ in a closed symplectic $6$-manifold $(M,\omega)$, 
there exist for each integer $g \ge 0$ and real number $E > 0$ 
at most finitely many distinct simple $J$-holomorphic curves $u$ of genus~$g$ 
in homology classes $[u] = A \in H_2(M)$ with $c_1(A) = 0$ and $\omega(A) \le E$.
Moreover, these curves are embedded and pairwise disjoint.
\qed
\end{cor}

\begin{remark}
Doan and Walpuski \cite{DoanWalpuski:Castelnuovo} have recently shown that
if one fixes the class $A \in H_2(M)$ in Corollary~\ref{cor:dimension6},
then it is not actually necessary to fix the genus~$g$,
i.e.~for generic~$J$, there exist at most
finitely-many simple curves \emph{of any genus} homologous to~$A$.
Their proof uses techniques from geometric measure theory.
\end{remark}

Using results of Zinger \cite{Zinger:comparison} (see also Lee-Parker \cite{LeeParker:obstruction}),
Theorem~\ref{thm:super} also implies that for generic~$J$, the space of branched covers 
of an embedded index~$0$ curve admits a well-defined obstruction
bundle which can be used to compute Gromov-Witten invariants.
In particular, if $\dim M \ge 6$ and $u : (\Sigma,j) \to (M,J)$ is an embedded $J$-holomorphic curve
of genus $g$ with $c_1(u)=0$, one can apply \cite{Zinger:comparison}*{Theorem~1.2} with no marked
point constraints to study the space of $J$-holomorphic curves with
image in $u(\Sigma)$,
so that Theorem~\ref{thm:super} establishes hypothesis~(b) in Zinger's result,
implying that the cokernels of the normal operators $\mathbf{D}_{u \circ \varphi}^N$ for 
$\varphi$ varying in the space 
$\widebar{\mM}_h(d[\Sigma],j)$ of degree~$d$ nodal holomorphic curves
in $(\Sigma,j)$ with arithmetic genus~$h$ 
form a well-defined and oriented orbibundle
$$
\Obst^u \to \widebar{\mM}_h(d[\Sigma],j)
$$
with $\rank_\RR \Obst^u = (n-1) (2h - 2 + d(2-2g))$.  Note that by the Riemann-Hurwitz
formula, the term $2h - 2 + d(2-2g)$ is simply the algebraic count of branch
points $Z(d\varphi)$ for any map $\varphi$ in the non-nodal stratum 
of $\widebar{\mM}_h(d[\Sigma],j)$.\footnote{One must keep in mind however
that the non-nodal stratum of $\widebar{\mM}_h(d[\Sigma])$ may be empty
even if $\widebar{\mM}_h(d[\Sigma])$ itself is not, e.g.~this is the case 
whenever $d=1$ and $h > g$.}
The obstruction bundle is interesting mainly in the $6$-dimensional case,
since $n=3$ means that $\rank_\RR \Obst^u$ matches the real virtual dimension
of $\widebar{\mM}_h(d[\Sigma],j)$,
and the count of solutions to an abstract perturbation of the holomorphic curve
equation can then be computed by integrating the Euler class $e(\Obst^u)$ over the
virtual fundamental cycle of~$\widebar{\mM}_h(d[\Sigma],j)$ in the sense
of \cites{LiTian:algebraic,LiTian,FukayaOno}.  This produces a formula for the \emph{local}
Gromov-Witten invariants of the curve~$u$,
$$
N_d^h(u) = \int_{[\widebar{\mM}_h(d[\Sigma],j)]^{\operatorname{vir}}}
e(\Obst^u) \in \QQ,
$$
defined for every $d \in \NN$ and $h \ge g$.  These numbers depend only on the
germ of the almost complex manifold $(M,J)$ at~$u(\Sigma)$.  Note that
$N_1^g(u) = \pm 1$, with the sign depending on the canonically oriented
determinant line of~$\mathbf{D}_u^N$.

Combining the obstruction bundle discussion with Corollary~\ref{cor:dimension6},
let
$$
N_A^g(M,\omega) \in \QQ
$$
denote the $0$-point Gromov-Witten invariant of $(M,\omega)$ for genus~$g$
curves in a class $A \in H_2(M)$ with $c_1(A)=0$.

\begin{cor}[via \cite{Zinger:comparison}*{Theorem~1.2}]
\label{cor:GW}
Suppose $(M,\omega)$ is a closed symplectic $6$-manifold, $g \ge 0$ is an
integer and $A \in H_2(M)$ satisfies $c_1(A)=0$.  Then for generic
$\omega$-compatible almost complex structures~$J$, 
$$
N_A^g(M,\omega) = \sum_{i=1}^N N_{d_i}^g(u_i),
$$
where the sum ranges over the (by Corollary~\ref{cor:dimension6}) finite set 
of pairwise disjoint embedded
$J$-holomorphic curves $u_1,\ldots,u_N$ that have genera at most~$g$ and
homology classes satisfying $d_i [u_i] = A$
for some $d_1,\ldots,d_N \in \NN$.
\qed
\end{cor}
In particular in the Calabi-Yau case, with $c_1(TM,\omega) = 0$,
this corollary localizes all of the Gromov-Witten invariants of $(M,\omega)$.

We next state two results on transversality for multiple covers.

\begin{theorem}[transversality, unbranched]
\label{thm:unbranched}
There exists a Baire subset $\jJ^\reg \subset \jJ(M,\omega\,;\,\uU,\Jfix)$
such that for all $J \in \jJ^\reg$, for every simple $J$-holomorphic curve
$u : (\Sigma,j) \to (M,J)$ intersecting $\uU$ and every unbranched holomorphic cover
$\varphi : (\widetilde{\Sigma},\tilde{\jmath}) \to (\Sigma,j)$ of closed Riemann surfaces,
the $J$-holomorphic curve $u \circ \varphi : (\widetilde{\Sigma},\tilde{\jmath}) \to (M,J)$ is
Fredholm regular.
\end{theorem}

\begin{remark}
The case $\ind(u)=0$ of Theorem~\ref{thm:unbranched} has been proved previously
in \cite{GerigWendl}, though with stronger assumptions: for technical reasons, it
was necessary in that paper to assume that $u(\Sigma)$ is \emph{contained entirely} in~$\uU$,
and in dimension four also to allow perturbations of $J$ that are $\omega$-tame but
not necessarily $\omega$-compatible.  The present paper uses a completely different
approach to the transversality problem and is thus able to remove these restrictions.
As explained in \cite{GerigWendl}, the theorem implies an integrality result for the
Gromov-Witten invariants in dimension four.
\end{remark}

It is generally harder to achieve transversality for covers $u \circ \varphi$ 
with branch points, e.g.~the index relation \eqref{eqn:indexCover} shows that 
$\ind(u\circ \varphi)$ can easily become negative
in dimensions greater than six.  More seriously, if $u$ is Fredholm regular, then one
can always find a smooth family of other multiple covers near $u \circ \varphi$
obtained by varying both $u$ and $\varphi$ in their respective moduli spaces;
since the latter lives in a space of real dimension $2 Z(d\varphi)$, the condition
$$
\ind(u \circ \varphi) \ge \ind(u) + 2 Z(d\varphi)
$$
is evidently necessary in order for $u \circ \varphi$ to be Fredholm regular.
Observe that if $\varphi$ has $r \ge 0$ critical values, then this condition
is satisfied whenever $\ind(u) \ge (n-1) r$: indeed, each critical value
is the image of at most $d-1$ branch points (counted algebraically), 
so we have $Z(d\varphi) \le (d-1)r$ and
\eqref{eqn:indexCover} implies
\begin{equation*}
\begin{split}
\ind(u \circ \varphi) &= \ind(u) + (d-1) \ind(u) - (n-3) Z(d\varphi) \\
&\ge \ind(u) + (n-1) Z(d\varphi) - (n-3) Z(d\varphi) = \ind(u) + 2 Z(d\varphi).
\end{split}
\end{equation*}
The next result states that the condition $\ind(u) \ge (n-1)r$ is also, in some
sense, sufficient.

\begin{theorem}[transversality, branched]
\label{thm:branched}
There exists a Baire subset $\jJ^\reg \subset \jJ(M,\omega\,;\, \uU,\Jfix)$ such that
the following holds for all $J \in \jJ^\reg$.  Suppose $u : (\Sigma,j) \to (M,J)$ is a
simple $J$-holomorphic curve intersecting $\uU$ and satisfying 
$$
\ind(u) \ge (n-1)r
$$
for some integer $r \ge 0$, and $\varphi : (\widetilde{\Sigma},\tilde{\jmath}) \to (\Sigma,j)$ is a
holomorphic branched cover of closed connected Riemann surfaces with $r$ distinct
critical values.  Then there exists a $J$-holomorphic curve and a holomorphic
branched cover
$$
u_\epsilon : (\Sigma,j_\epsilon) \to (M,J) \quad \text{ and } \quad
\varphi_\epsilon : (\widetilde{\Sigma},\tilde{\jmath}_\epsilon) \to (\Sigma,j_\epsilon)
$$ 
such that $u_\epsilon$, $\varphi_\epsilon$, $j_\epsilon$ and $\tilde{\jmath}_\epsilon$ are
arbitrarily $C^\infty$-close to $u$, $\varphi$, $j$ and $\tilde{\jmath}$ respectively,
and $u_\epsilon \circ \varphi_\epsilon : (\widetilde{\Sigma},\tilde{\jmath}_\epsilon) \to (M,J)$
is Fredholm regular.
\end{theorem}

Note that whenever $\ind(u \circ \varphi)$ is also \emph{strictly} greater than
$\ind(u) + 2 Z(d\varphi)$, one can combine this result with the implicit function
theorem to deduce the existence of simple
$J$-holomorphic curves that are small perturbations of multiple covers of~$u$.

The proofs of these theorems are inspired by the work of Taubes \cite{Taubes:counting},
whose definition of the Gromov invariant for symplectic $4$-manifolds required
a special case of Theorem~\ref{thm:unbranched} along with related 
bifurcation-theoretic results (cf.~\S\ref{sec:bifurcations}) for multiply 
covered holomorphic tori.  Roughly speaking, the idea is to study the local
structure of spaces of the form
\begin{equation}
\label{eqn:theBigIdea}
\mM(k,c) := \left\{ \tilde{u} = u \circ \varphi\ \big|\ 
\text{$\dim \ker\mathbf{D}_{\tilde{u}}^N = k$ and $\dim \coker \mathbf{D}_{\tilde{u}}^N = c$}
\right\},
\end{equation}
where $k,c \ge 0$ are fixed integers, $u$ varies in the moduli space of simple 
$J$-holomorphic curves and $\varphi$ varies in the moduli space of holomorphic
branched covers.  Ideally, one would like to show that these spaces are
smooth manifolds for generic $J$, and to compute their codimensions in the
space of pairs $(u,\varphi)$.  This turns Theorems~\ref{thm:super}
and~\ref{thm:unbranched} into ``dimension counting'' problems, as whenever one
can show that the codimension of $\mM(k,c)$ is larger than the
dimension of the ambient space for suitable values of $k$ and~$c$, one may conclude that
either $\ker \mathbf{D}_{\tilde{u}}^N$ or $\coker \mathbf{D}_{\tilde{u}}^N$
must be trivial.  This discussion is oversimplified in at least three
respects: first, we will not be able to find any nice structure on
$\mM(k,c)$ if $\varphi$ varies in the space of \emph{all}
branched covers, but it will help to confine it to certain substrata of that
space in which all branch points have prescribed branching orders.  For similar
reasons, it will also help to confine $u$ to substrata in which its number of
critical points and their orders are constrained, and this is easily done.
More seriously, the space $\mM(k,c)$ as sketched above can have different 
codimensions on different components, as its codimension depends intricately on
symmetry information which is ignored in \eqref{eqn:theBigIdea}.  We will therefore
need to define a more elaborate version of $\mM(k,c)$ which depends on a
splitting of the operator $\mathbf{D}_{\tilde{u}}^N$ into summands
corresponding to irreducible representations of the (generalized) symmetry
group of the cover.  This idea is borrowed directly from \cite{Taubes:counting}, 
though the details are somewhat more involved since, in contrast to the case
of unbranched covers of tori, we cannot assume that all covers are regular
or that their symmetry groups are abelian.  We will see that once the
formalism is developed in sufficient generality,
it ``breaks the symmetry'' of $\mathbf{D}_{\tilde{u}}^N$ enough to make 
dimension counting arguments much more effective.

\begin{remark}
A slightly different variation on the ideas in \cite{Taubes:counting} has
been implemented by Eftekhary to prove a partial result toward
super-rigidity in dimension six, see \cite{Eftekhary:superrigidity}.
\end{remark}

Here is an outline of the rest of the paper.

After establishing some standard definitions and notation, \S\ref{sec:mainIdea}
will further elucidate the ideas sketched above and formulate a precise version
of the statement that $\mM(k,c)$ from \eqref{eqn:theBigIdea} is a smooth submanifold,
Theorem~\ref{thm:submanifolds0}.  This will then be used as a black box to prove
Theorems~\ref{thm:super}, \ref{thm:unbranched} and~\ref{thm:branched} in \S\ref{sec:everything},
followed in \S\ref{sec:bifurcations} by a brief informal discussion of
bifurcation theory.
The remainder of the paper is then devoted to the proof of Theorem~\ref{thm:submanifolds0}.
In \S\ref{sec:prep}, we explain the splitting construction for Cauchy-Riemann
operators with symmetries and prove some lemmas based on a mixture of 
elliptic regularity for punctured Cauchy-Riemann operators,
topology of covering spaces, and representation theory of finite groups.
The summands in the splitting are also Cauchy-Riemann operators, whose
indices are a somewhat delicate computation, carried out in \S\ref{sec:index}.
In \S\ref{sec:unique} we prove a local genericity result for Cauchy-Riemann
operators that takes on the role usually played by
unique continuation in applications of the
Sard-Smale theorem, and the latter will be used in \S\ref{sec:SardSmale}
to complete the proof of Theorem~\ref{thm:submanifolds0}.
Finally,
\S\ref{sec:dimensionFour} deals with super-rigidity in the four-dimensional
case, which is something of an anomaly and requires different techniques based
on intersection theory.  The appendices provide various
results that may be considered ``standard'' and yet, in this
author's experience, seem to cause sufficient confusion among experts to
warrant some discussion; their proofs require a few ideas that will in any
case be useful elsewhere in the paper.

\subsection{Apologies and acknowledgements}

The super-rigidity problem has a slightly troubled history, and as the
author of a new paper on the subject, it would behoove me at this point to 
apologize for having caused some of that trouble: I am aware of three previous 
attempts to prove some version of Theorem~\ref{thm:super} which were
later either withdrawn or revised to prove much weaker statements, and I was an
author of one of them (the original version of \cite{GerigWendl}).
To make matters worse, earlier versions of the present paper also contained a major
error in \S\ref{sec:unique} on which the main results were crucially dependent,
causing the paper to be withdrawn for several months while
the offending section underwent an extensive rewrite.
(For more on the history of failed super-rigidity proofs, see Appendix~\ref{sec:catastrophe}.)
With all this in mind, I would sympathize with any reader's inclination to
greet this paper with a dose of skepticism, though it seems worth
pointing out that rather than being an attempt to rescue the (probably
unrescuable) proof originally attempted in \cite{GerigWendl},
the approach taken here has almost nothing in common with the
previous one, other than the considerable debt that both of them owe to the ideas of
Taubes \cites{Taubes:SWtoGr,Taubes:counting}.

I would like to thank Dan Cristofaro-Gardiner, Chris Gerig, Michael Hutchings,
Eleny Ionel, Mihai Munteanu, Tom Parker, Cliff Taubes
and Aleksey Zinger for conversations and correspondence which helped to
improve my understanding of the problems studied in this paper.
Special thanks are due to Aleksander Doan and Thomas Walpuski for having uncovered
a few minor errors and one major error in the original version; my
discussions with them were invaluable in the effort toward fixing those errors.
Finally, many thanks to the anonymous referees for their impressively
careful reading of earlier drafts, which has induced measurable
improvements in the exposition.

\section{The main idea}
\label{sec:mainIdea}

\subsection{Some definitions}
\label{sec:defn}

Let us now fix some notation
and definitions that will be essential in the rest of the paper.

Given integers $g,m \ge 0$ and a class $A \in H_2(M)$, 
the moduli space of \defin{unparametrized $J$-holomorphic curves}
$\mM_{g,m}(A,J)$ can be defined as the set
of equivalence classes of tuples
$(\Sigma,j,\Theta,u)$ where $(\Sigma,j)$ is a closed connected Riemann
surface of genus~$g$, $\Theta \subset \Sigma$ is an ordered set of
$m$ distinct points (the \defin{marked points}), 
and $u : (\Sigma,j) \to (M,J)$ is a $J$-holomorphic
map satisfying $[u] := u_*[\Sigma] = A$, with equivalence defined by $(\Sigma,j,\Theta,u) \sim
(\Sigma',\psi^*j,\psi^{-1}(\Theta),u \circ \psi)$ for diffeomorphisms
$\psi : \Sigma' \to \Sigma$.  
The \defin{Gromov compactification} of
$\mM_{g,m}(A,J)$ is the space $\widebar{\mM}_{g,m}(A,J)$ of 
(equivalence classes of) \defin{stable nodal
curves} $(S,j,\Theta,\Delta,u)$, where now $S$ may be disconnected, and
the original data are augmented by an unordered
set of distinct points in $S \setminus \Theta$, arranged into unordered pairs
$$
\Delta = \left\{ \{\widehat{z}_1,\widecheck{z}_1\},\ldots,\{\widehat{z}_r,\widecheck{z}_r\} \right\},
$$
such that $u(\widehat{z}_i) = u(\widecheck{z}_i)$ for each $i=1,\ldots,r$.  
We call the pairs $\{\widehat{z}_i,\widecheck{z}_i\}$ \defin{nodes}, and each individual
$\widehat{z}_i$ or $\widecheck{z}_i \in S$ a \defin{nodal point}.
The curves in
$\widebar{\mM}_{g,m}(A,J)$ are required to have \defin{arithmetic genus}~$g$,
which means that the surface obtained from $S$ by performing connected sums
at all matched pairs of nodal points is a closed connected
surface of genus~$g$.  The stability condition requires that any component of
$S \setminus (\Theta \cup \Delta)$ on which $u$ is constant should have negative
Euler characteristic.  With this condition, $\widebar{\mM}_{g,m}(A,J)$ can be
given a natural topology as a metrizable Hausdorff space, and it is compact
whenever $J$ is tamed by a symplectic form.  A definition of the topology may be
found e.g.~in \cite{SFTcompactness}; for convergent
sequences in $\mM_{g,m}(A,J)$, it amounts to the notion of $C^\infty$-convergence
for $j$ and $u$ after a choice of parametrization for which all domains and
marked point sets are identified.  {Curves $[(S,j,\Theta,\Delta,u)] \in
\widebar{\mM}_{g,m}(A,J)$ with $\Delta = \emptyset$ can equivalently be
regarded as elements of $\mM_{g,m}(A,J)$, and are thus called \defin{smooth}
curves to distinguish them from nodal curves.}

\begin{remark}
In this paper, the word ``curve'' always means ``smooth curve'' (i.e.~without nodes)
unless the word ``nodal'' is explicitly included.  Similarly, all dimensions
and Fredholm indices in this paper are \emph{real} (not complex) unless otherwise
specified.  This usage differs somewhat 
from the algebraic geometry literature.
\end{remark}

When there is no danger of confusion, we shall sometimes 
abuse notation by writing equivalence classes $[(\Sigma,j,\Theta,u)] \in \mM_{g,m}(A,J)$
or $[({S},j,\Theta,\Delta,u)] \in \widebar{\mM}_{g,m}(A,J)$ via the abbreviations
$u \in \mM_{g,m}(A,J)$ or $u \in \widebar{\mM}_{g,m}(A,J)$ respectively, 
and we will refer to the
restriction of a nodal curve {$[(S,j,\Theta,\Delta,u)]$ to any connected component of 
its domain $S$} as a
\defin{smooth component} of~$u$.  We shall also abbreviate 
$$
\mM_g(A,J) := \mM_{g,0}(A,J),
\quad \text{ and } \quad
\widebar{\mM}_g(A,J) := \widebar{\mM}_{g,0}(A,J).
$$
Recall that $\mM_g(A,J)$ has
\defin{virtual dimension} equal to the index of any curve $u \in \mM_g(A,J)$
as written in \eqref{eqn:index},
while the virtual dimension of the moduli space with marked points is
$$
\virdim \mM_{g,m}(A,J) = \virdim \mM_g(A,J) + 2m.
$$

The multiply covered curves form a distinguished closed subset of
$\widebar{\mM}_g(A,J)$.  Given any $u \in \mM_g(A,J)$ with domain
$(\Sigma,j)$, and integers $h \ge 0$, $d \ge 1$, define the space of
stable \defin{nodal $d$-fold covers} of~$u$,
$$
\widebar{\mM}_{h}(d ; u) = \left\{ [(S,\tilde{\jmath},\Delta,u \circ \varphi)]
\in \widebar{\mM}_{h}(dA,J)\ \big|\ [(S,\tilde{\jmath},\Delta,\varphi)] \in
\widebar{\mM}_{h}(d[\Sigma],j) \right\},
$$
so in particular,
each {smooth} component $\tilde{u}_i$ of $\tilde{u} \in \widebar{\mM}_{h}(d ; u)$
belongs to a space $\mM_{g_i}(d_i ; u)$ of smooth branched covers $u \circ \varphi_i$
of some degreee $d_i \ge 0$, such that $\sum_i d_i = d$.  
Note that $\widebar{\mM}_h(d ; u)$ may in general be strictly larger
than the closure of $\mM_h(d ; u)$ in the Gromov topology---to cite one well-known
example, {the space $\mM_1([S^2],i)$ of smooth degree~$1$ holomorphic tori in
$(S^2,i)$ is empty}, but $\widebar{\mM}_1([S^2],i)$ contains
a nodal curve with a constant component of genus~$1$.

Recall next that every $J$-holomorphic curve $u : (\Sigma,j) \to (M,J)$
gives rise to a \defin{linearized Cauchy-Riemann operator}
$$
\mathbf{D}_u : \Gamma(u^*TM) \to \Omega^{0,1}(\Sigma,u^*TM),
$$
i.e.~the linearization at $u$ of the nonlinear Cauchy-Riemann operator
$\dbar_J(u) := Tu + J \circ Tu \circ j \in \Omega^{0,1}(\Sigma,u^*TM)$, 
whose zero-set is the space of all $J$-holomorphic maps with domain~$(\Sigma,j)$.  The operator $\mathbf{D}_u$
takes vector fields along~$u$
to $(0,1)$-forms valued in the complex vector bundle $(u^*TM,J)$,
and can be written explicitly as
$$
\mathbf{D}_u \eta = \nabla \eta + J(u) \circ \nabla\eta \circ j +
(\nabla_\eta J) \circ Tu \circ j
$$
for any choice of symmetric connection~$\nabla$
(cf.~\cite{Wendl:lecturesV2}*{\S 2.4}).
Recall moreover that whenever $u$ is
nonconstant, its critical points are isolated and one can find a
smooth splitting of complex vector bundles
\begin{equation}
\label{eqn:splitting}
u^*TM = T_u \oplus N_u
\end{equation}
such that $T_u$ matches the image of $du$ at regular points;
see e.g.~\cite{Wendl:automatic}*{\S 3.3} for details.  We shall refer to
$N_u$ as the \defin{generalized normal bundle} of~$u$.
In many cases
of interest in this paper, $u$ will be a cover of
an immersed $J$-holomorphic curve $v$, so $N_u$ is then simply the
pullback of the normal bundle of~$v$ via the cover.
We define the \defin{normal Cauchy-Riemann operator} at $u$ as the restriction 
of $\mathbf{D}_u$ to sections of $N_u$, composed with
the projection $\pi_N : u^*TM \to N_u$ along~$T_u$, hence
$$
\mathbf{D}_u^N = \pi_N \circ \mathbf{D}_u|_{\Gamma(N_u)} : \Gamma(N_u) \to 
\Omega^{0,1}(\Sigma,N_u).
$$

In general, a neighborhood of any element in $\mM_{g,m}(A,J)$ can be identified
with the zero-set of a smooth Fredholm section of a Banach space bundle, modulo
a finite group action if there are nontrivial automorphisms.
We say that $u \in \mM_g(A,J)$ is \defin{Fredholm regular} whenever it is a
transverse intersection of this section with the zero-section.  Note that
whenever this condition holds, it automatically also holds after adding any
finite collection of marked points and viewing $u$ as an element 
of $\mM_{g,m}(A,J)$.  The implicit function theorem gives the open set of
regular curves in $\mM_{g,m}(A,J)$ the structure of a smooth orbifold with 
dimension equal to its virtual dimension, and local isotropy groups
determined by the automorphism groups of the curves---in particular, the set
of regular simple curves forms a manifold, though orbifold singularities can 
appear when multiple covers are included.
The following convenient repackaging of the regularity
condition comes from \cite{Wendl:automatic}*{Corollary~3.13}.

\begin{prop}
\label{prop:regular}
A closed and connected $J$-holomorphic curve $u : (\Sigma,j) \to (M,J)$
is Fredholm regular if and only if its normal
operator $\mathbf{D}_u^N : W^{k,p}(N_u) \to W^{k-1,p}(\overline{\Hom}_\CC(T\Sigma,N_u))$
is surjective for some (and therefore all) $k \in \NN$ and $p \in (1,\infty)$.
\qed
\end{prop}

\begin{defn}
\label{defn:super}
A closed, connected, simple $J$-holomorphic curve 
$u : (\Sigma,j) \to (M,J)$ is called
\defin{super-rigid} if it satisfies the following:
\begin{enumerate}
\item $\ind(u) = 0$;
\item $u : \Sigma \to M$ is an immersion;
\item For all closed connected Riemann surfaces $(\widetilde{\Sigma},\tilde{\jmath})$
and holomorphic maps $\varphi : (\widetilde{\Sigma},\tilde{\jmath})
\to (\Sigma,j)$ of positive degree, the curve $\tilde{u} := u \circ \varphi : 
(\widetilde{\Sigma},\tilde{\jmath}) \to (M,J)$ admits no nontrivial solutions
to the normal linearized equation $\mathbf{D}_{\tilde{u}}^N \eta = 0$.
\end{enumerate}
\end{defn}

Proposition~\ref{prop:isolated} in Appendix~\ref{sec:isolated} proves that
if $u$ is a super-rigid curve, then the only possible sequences that
converge to a nodal branched cover of $u$ consist of other covers of~$u$.
In the language of the present section, this means:

\begin{cor}[of Proposition~\ref{prop:isolated}]
\label{cor:isolated}
Suppose $(M,J)$ is an almost complex manifold and 
$u \in \mM_{g}(A,J)$ is a super-rigid curve in~$M$.  Then for
every $h \ge 0$ and $d \ge 1$,
$\widebar{\mM}_{h}(d ; u)$ is an open and closed subset of
$\widebar{\mM}_{h}(dA,J)$.
\qed
\end{cor}

\subsection{A stratification theorem}
\label{sec:theBigIdea}

We now explain in precise terms the stratification result that underlies
the main theorems of \S\ref{sec:mainTheorems}.

\subsubsection{Splitting the linearization at a doubly covered curve}

Suppose $v : (\Sigma,j) \to (M,J)$ is a simple $J$-holomorphic curve with genus $g \ge 0$, and
$\varphi : (\Sigma',j') \to (\Sigma,j)$ is a holomorphic
branched cover with degree $d \ge 1$, giving rise to the multiply covered
curve $u = v \circ \varphi : (\Sigma',j') \to
(M,J)$ of genus $h \ge 0$.  We assume as always that $\Sigma$ and $\Sigma'$
are both closed and connected, and for the sake of intuition, we begin 
in this subsection with the special case $d=2$.
The automorphism group
$$
\Aut(u) = \Aut(\varphi) := \left\{ \psi : (\Sigma',j') \stackrel{\cong}{\longrightarrow}
(\Sigma',j')\ \Big|\ \varphi = \varphi \circ \psi \right\}
$$
then contains a unique nontrivial element $\psi$, and the space of sections
$\Gamma(N_u)$ has a natural splitting
$$
\Gamma(N_u) = \Gamma_+(N_u) \oplus \Gamma_-(N_u)
$$
where $\Gamma_\pm(N_u) := \{ \eta \in \Gamma(N_u)\ |\ \eta = \pm \eta \circ \psi \}$.
Splitting $\Omega^{0,1}(\Sigma',N_u) = \Gamma(\overline{\Hom}_\CC(T\Sigma',N_u))$
in the same way, one obtains a splitting of the normal Cauchy-Riemann operator
\begin{equation}
\label{eqn:firstSplitting}
\mathbf{D}_u^N = \mathbf{D}_{u,+}^N \oplus \mathbf{D}_{u,-}^N
\end{equation}
into two operators $\mathbf{D}_{u,\pm}^N : \Gamma_\pm(N_u) \to
\Gamma_\pm(\overline{\Hom}_\CC(T\Sigma',N_u))$.  It is not hard to see
that $\mathbf{D}_{u,+}^N$ is in some sense equivalent to~$\mathbf{D}_v^N$, as its
domain and target both consist of sections that are pullbacks via $\varphi$
of sections over~$\Sigma$.  The operators $\mathbf{D}_{u,+}^N$ and
$\mathbf{D}_{u,-}^N$ have unique extensions over the spaces
of symmetric/antisymmetric sections of Sobolev class $W^{k,p}$ for
$k \in \NN$ and $p \in (1,\infty)$, giving bounded linear operators
$$
\mathbf{D}_{u,\pm}^N : W^{k,p}_\pm(N_u) \to
W^{k,p}_\pm(\overline{\Hom}_\CC(T\Sigma',N_u)),
$$
and the standard transversality theory for simple
curves then implies that $\mathbf{D}_{u,_+}^N$ can be assumed surjective
(and also injective if $v$ is immersed with index~$0$) if $J$ is chosen generically.  
We will see that the problem of proving surjectivity or injectivity for
$\mathbf{D}_u^N$ becomes more tractable when viewed as two independent problems
for the operators $\mathbf{D}_{u,+}^N$ and~$\mathbf{D}_{u,-}^N$.

In order to generalize this discussion beyond the degree~$2$ case, it helps 
to adopt an alternative perspective based on representation theory.
Let $\Theta \subset \Sigma$ denote a finite subset that contains all critical 
values of~$\varphi$, and set 
\begin{equation}
\label{eqn:punctures}
\Theta' := \varphi^{-1}(\Theta), \qquad
\dot{\Sigma} := \Sigma \setminus \Theta, \qquad
\dot{\Sigma}' := \Sigma' \setminus \Theta',
\end{equation}
so that $\dot{\Sigma}' \stackrel{\varphi}{\longrightarrow} \dot{\Sigma}$
is a smooth covering map with $G := \Aut(\varphi) \cong \ZZ_2$ as its
group of deck transformations.  Define
$$
\rho : G \to S_2 : g \mapsto \rho_g
$$
as the isomorphism to the symmetric group on $\{1,2\}$.
We can then identify the covering map $\dot{\Sigma}' \stackrel{\varphi}{\longrightarrow} \dot{\Sigma}$
with
$$
\left( \dot{\Sigma}' \times \{1,2\} \right) \Big/ G \to \dot{\Sigma} :
[(z,i)] \mapsto \varphi(z),
$$
where $G$ acts on $\dot{\Sigma}'$ by deck transformations and on $\{1,2\}$
via~$\rho$.  Now if $(e_1,e_2)$ denotes the standard basis of $\RR^2$, 
then $\rho$ also gives rise to a real permutation representation
$$
\boldsymbol{\rho} : G \to \GL(2,\RR), \qquad \boldsymbol{\rho}(g) e_i :=
e_{\rho_g(i)},
$$
and a corresponding real vector bundle $V^{\boldsymbol{\rho}} \to \dot{\Sigma}$ defined 
as the $\ZZ_2$-quotient of a trivial bundle over $\dot{\Sigma}'$,
$$
V^{\boldsymbol{\rho}} := \left( \dot{\Sigma}' \times \RR^2 \right) \Big/ G.
$$
The space of sections of the twisted normal bundle
$$
N_v^{\boldsymbol{\rho}} := N_v \otimes_\RR V^{\boldsymbol{\rho}} \to \dot{\Sigma}
$$
then has a natural identification with the space of sections of
$N_u = \varphi^*N_v$: indeed, we can represent sections of $N_v^{\boldsymbol{\rho}}$
as $\ZZ_2$-equivariant sections $\eta = \sum_{i=1}^2 \eta^i \otimes e_i$ of 
$\varphi^*N_v \otimes_\RR \RR^2$, which satisfy the relation $\eta^i \circ \psi = \eta^{\rho_\psi(i)}$,
thus a corresponding section $\widehat{\eta} \in \Gamma(\varphi^*N_v)$ can be
defined under the identification of $\dot{\Sigma}'$ with 
$(\dot{\Sigma}' \times \{1,2\}) / G$ by
$$
\widehat{\eta}([(z,i)]) = \eta^i(z).
$$
Under this identification, $\mathbf{D}_u^N$ becomes a Cauchy-Riemann type
operator on the twisted bundle~$N_v^{\boldsymbol{\rho}}$, defined locally
by $\mathbf{D}_u^N (\eta \otimes s) = (\mathbf{D}_v^N \eta)\otimes s$
whenever $s$ is a local section of $V^{\boldsymbol{\rho}}$ that has a
constant lift to the trivial bundle $\dot{\Sigma}' \times \RR^2$.

The above construction appears cumbersome at first glance, but it has the
following advantage: the decomposition $\Gamma(N_u) = \Gamma_+(N_u) \oplus
\Gamma_-(N_u)$ now corresponds to a splitting of the twisted bundle
$N_v^{\boldsymbol{\rho}}$ into subbundles
$$
N_v^{\boldsymbol{\rho}} = N_v^{\boldsymbol{\theta}_+} \oplus N_v^{\boldsymbol{\theta}_-} :=
(N_v \otimes_\RR V^{\boldsymbol{\theta}_+}) \oplus
(N_v \otimes_\RR V^{\boldsymbol{\theta}_-})
$$
where $V^{\boldsymbol{\theta}_\pm} := (\dot{\Sigma}' \times W_\pm) / G$
are defined in terms of the natural splitting of $\RR^2 = W_+ \oplus W_-$
into irreducible $G$-invariant subspaces
$$
W_\pm = \RR\begin{pmatrix} 1 \\ \pm 1 \end{pmatrix} \subset \RR^2.
$$
This is the simplest nontrivial example of what turns out to be a general
principle: splittings of Cauchy-Riemann operators for multiply covered
curves arise from decompositions of permutation representations into
irreducible summands.
To turn $\boldsymbol{\rho} = \boldsymbol{\theta}_+ \oplus \boldsymbol{\theta}_-$
into a splitting of Cauchy-Riemann operators, we still have a small
analytical issue to cope with since the bundles $N_v^{\boldsymbol{\theta}_\pm}$
are defined over $\dot{\Sigma}$ and do not both extend over the punctures.
In place of \eqref{eqn:firstSplitting}, we therefore obtain a splitting
$$
\dot{\mathbf{D}}_u^N = \dot{\mathbf{D}}_{u,\boldsymbol{\theta}_+}^N \oplus
\dot{\mathbf{D}}_{u,\boldsymbol{\theta}_-}^N,
$$
where the dots over the operators indicate that we are restricting them to
the punctured domain~$\dot{\Sigma}'$.  We will see in \S\ref{sec:CRpunctured}
how to define suitable weighted Sobolev spaces over $\dot{\Sigma}$ and
$\dot{\Sigma}'$ so that the punctured operators have the same indices, kernels
and cokernels as their unpunctured counterparts.

\begin{remark}
A slightly different approach to defining twisted Cauchy-Riemann operators
is taken by Doan and Walpuski \cite{DoanWalpuski:BrillNoether}, who express
it in the elegant language of local systems.
\end{remark}

\subsubsection{The codimension of a multiply covered curve}

We return now to the general case of a closed connected $J$-holomorphic curve
$u = v \circ \varphi : (\Sigma',j') \to (M,J)$ of genus~$h$, where $v : (\Sigma,j) \to (M,J)$ is
simple with genus~$g$ and $\varphi : (\Sigma',j') \to (\Sigma,j)$ has degree $d \in \NN$.
We continue using the notation
$\dot{\Sigma}' \stackrel{\varphi}{\longrightarrow} \dot{\Sigma}$ for the
$d$-fold covering map obtained by deleting some finite subsets that include the
critical values and their preimages.
Recall that $\varphi$ is called \defin{regular} if $|\Aut(\varphi)| =
\deg(\varphi) = d$.  This condition was secretly important in the
above discussion of the $d=2$ case, as the definition of the twisted bundle
$N_v^{\boldsymbol{\rho}}$ required identifying $\dot{\Sigma}$ with the
quotient of $\dot{\Sigma}'$ by deck transformations.
In general, $\Aut(\varphi)$ can have order smaller than $d$ and 
may even be trivial, but we can use some notions from elementary covering space 
theory to get around this.

\begin{defn}
\label{defn:genAut}
The \defin{generalized automorphism group} of a $d$-fold branched cover
$\varphi : \Sigma' \to \Sigma$ is the quotient $G := \pi_1(\dot{\Sigma}) / H$,
where $H$ is the normal core\footnote{Recall that the \defin{normal core} of
a subgroup $H$ in a group $\Gamma$ is the largest normal subgroup of $\Gamma$
that is contained in~$H$.} of $\varphi_*(\pi_1(\dot{\Sigma}'))$, 
and $\dot{\Sigma}$ and $\dot{\Sigma}'$ are defined by \eqref{eqn:punctures}
with $\Theta$ as the set of critical values of~$\varphi$.
\end{defn}

\begin{remark}
Like fundamental groups, the generalized automorphism group $G$ of
$\varphi : \Sigma' \to \Sigma$ depends on choices of base points in
$\dot{\Sigma}$ and $\dot{\Sigma}'$, but its isomorphism class is
independent of these choices.  We will see below that $G$ is a finite
group of order at most~$d!$ that is isomorphic
to $\Aut(\varphi)$ if and only if $\varphi : \Sigma' \to \Sigma$ is regular,
and more generally, $G$ has a natural identification with the automorphism
group of a certain regular branched cover of $\Sigma$ that is determined
by $\varphi$ and a choice of base points, and factors through~$\varphi$.
\end{remark}

\begin{defn}
\label{defn:presentation}
A \defin{regular presentation} of the holomorphic $d$-fold branched cover 
$\varphi : (\Sigma',j') \to (\Sigma,j)$ 
is a tuple $(\Theta,\dot{\Sigma}'',\pi,G,\rho,I,f)$ consisting of:
\begin{itemize}
\setlength{\itemsep}{0cm}
\item A finite subset $\Theta \subset \Sigma$ containing the critical values
of $\varphi$ and defining the punctured surfaces $\dot{\Sigma}$ and $\dot{\Sigma}'$
via \eqref{eqn:punctures};
\item A connected surface $\dot{\Sigma}''$ and regular covering map
$\pi : \dot{\Sigma}'' \to
\dot{\Sigma}$ with finite automorphism group $G := \Aut(\pi)$;
\item A set~$I$ with $d$ elements;
\item A transitive action of $G$ on~$I$, defined via a homomorphism
$\rho : G \to S(I)$ from $G$ to the symmetric group on~$I$;
\item A diffeomorphism $f : \dot{\Sigma}' \to (\dot{\Sigma}'' \times I) / G$,
where $G$ acts on $\dot{\Sigma}''$ by deck transformations and on $I$
via~$\rho$, such that $\varphi \circ f^{-1}$ takes the form
$$
\left( \dot{\Sigma}'' \times I \right) \Big/ G \to \dot{\Sigma} :
[(z,i)] \mapsto \pi(z).
$$
\end{itemize}
We say that $(\Theta,\dot{\Sigma}'',\pi,G,\rho,I,f)$ is \defin{minimal} if 
$\Theta \subset \Sigma$ is the set of critical values of $\varphi$ and
$\rho : G \to S(I)$ is injective.
Two regular presentations $(\Theta_j,\dot{\Sigma}''_j,\pi_j,G_j,\rho_j,I_j,f_j)$
of $\varphi : \Sigma' \to \Sigma$ for $j=1,2$ are \defin{isomorphic}
if $\Theta_1 = \Theta_2$ and there exists a diffeomorphism $\Psi : \dot{\Sigma}''_1 \to \dot{\Sigma}''_2$,
a bijection $\beta : I_1 \to I_2$, and a group isomorphism $\Phi : G_1 \to G_2$
such that:
\begin{enumerate}
\setlength{\itemsep}{0cm}
\item $\pi_2 \circ \Psi = \pi_1$ and for all $g \in G_1$, 
$\Psi \circ g = \Phi(g) \circ \Psi$;
\item For all $g \in G_1$, $\beta \circ \rho_1(g) = \rho_2(\Phi(g)) \circ \beta$;
\item $f_2 \circ f_1^{-1}$ takes the form
$$
\left( \dot{\Sigma}''_1 \times I_1 \right) \Big/ G_1 \to
\left( \dot{\Sigma}''_2 \times I_2 \right) \Big/ G_2 :
[(z,i)] \mapsto [(\Psi(z),\beta(i))].
$$
\end{enumerate}
\end{defn}

Most of the regular presentations we encounter in this paper will be
minimal, though an important example that is not (in particular where $\Theta$
may contain more than just the critical values) will arise in Example~\ref{ex:factorization}.
Standard results about Riemann surfaces (see \S\ref{sec:regPres}) imply
that the regular cover
$\pi : \dot{\Sigma}'' \to \dot{\Sigma}$ in any regular presentation
can be extended to a holomorphic
branched cover of closed connected Riemann surfaces $(\Sigma'',j'') \to
(\Sigma,j)$ such that $\dot{\Sigma}'' = \Sigma'' \setminus \pi^{-1}(\Theta)$.
Observe that if $i \in I$ and $G_i \subset G$ denotes the stabilizer of $i$
under the $G$-action defined by~$\rho$, then
$$
\dot{\Sigma}'' / G_i \to \left(\dot{\Sigma}'' \times I\right) \Big/ G : [z] \mapsto [(z,i)]
$$
is a diffeomorphism identifying $\varphi \circ f^{-1}$ with the natural
projection $\dot{\Sigma}'' / G_i \to \dot{\Sigma}'' / G = \dot{\Sigma}$.
Thus one can associate to any regular presentation a (non-unique) factorization
of $\pi : \dot{\Sigma}'' \to \dot{\Sigma}$ by covering maps
$\dot{\Sigma}'' \to \dot{\Sigma}' \stackrel{\varphi}{\longrightarrow} \dot{\Sigma}$,
which extends over the punctures to a factorization of 
$\pi : (\Sigma'',j'') \to (\Sigma,j)$ by holomorphic branched covers
$$
(\Sigma'',j'') \to (\Sigma',j') \stackrel{\varphi}{\longrightarrow} (\Sigma,j).
$$
We will also show in Lemma~\ref{lemma:regPres} that $\varphi : \Sigma' \to \Sigma$
always admits a unique isomorphism class of minimal regular presentations
$(\Theta,\dot{\Sigma}'',\pi,G,\rho,I,f)$, for which $G$ is isomorphic to the
generalized automorphism group of~$\varphi$, and in this case
$\pi : \dot{\Sigma}'' \to \dot{\Sigma}$ is isomorphic to
$\varphi : \dot{\Sigma}' \to \dot{\Sigma}$ whenever the latter happens to be
already regular (cf.~Example~\ref{ex:regular}).  

Given a choice of regular presentation $(\Theta,\dot{\Sigma}'',\pi,G,\rho,I,f)$,
the discussion of the degree~$2$ case can be generalized as follows.
The transitive action $\rho : G \to S(I)$ induces a
permutation representation $\boldsymbol{\rho} : G \to
\Aut_\RR(\RR^I)$ on the real vector space $\RR^I$ with basis labeled by the
elements of~$I$, and a twisted bundle $N_v^{\boldsymbol{\rho}} =
N_v \otimes_\RR V^{\boldsymbol{\rho}} \to \dot{\Sigma}$, where
$$
V^{\boldsymbol{\rho}} := (\dot{\Sigma}'' \times \RR^I) / G,
$$
with a natural isomorphism
$$
\Gamma(N_v^{\boldsymbol{\rho}}) = \Gamma(\varphi^*N_v|_{\dot{\Sigma}'})
= \Gamma(N_u|_{\dot{\Sigma}'})
$$
that identifies $\mathbf{D}_u^N$ with a Cauchy-Riemann operator
$$
\dot{\mathbf{D}}_{u,\boldsymbol{\rho}}^N : 
\Gamma(N_v^{\boldsymbol{\rho}}) \to
\Omega^{0,1}(\dot{\Sigma},N_v^{\boldsymbol{\rho}}),
$$
defined on suitable exponentially weighted
Sobolev spaces of sections of~$N_v^{\boldsymbol{\rho}}$.
(The appropriate functional-analytic setting for this operator
will be specified precisely in \S\ref{sec:CRpunctured}.)
Any representation $\boldsymbol{\theta} : G \to \Aut_\RR(W)$ on a real
finite-dimensional vector space $W$
similarly gives rise to a twisted bundle $N_v^{\boldsymbol{\theta}} =
N_v \otimes_\RR V^{\boldsymbol{\theta}} \to \dot{\Sigma}$, with
$V^{\boldsymbol{\theta}} := (\dot{\Sigma}'' \times W) / G$,
and a twisted Cauchy-Riemann operator 
$$
\dot{\mathbf{D}}_{u,\boldsymbol{\theta}}^N : 
\Gamma(N_v^{\boldsymbol{\theta}}) \to
\Omega^{0,1}(\dot{\Sigma},N_v^{\boldsymbol{\theta}}),
$$
which (up to conjugacy) depends only on $\dot{\mathbf{D}}_v^N$ and the isomorphism 
classes of the regular presentation and the representation $\boldsymbol{\theta}$.
Now any representation-theoretic decomposition
$\boldsymbol{\rho} = \boldsymbol{\theta}_1^{\oplus m_1} \oplus \ldots \oplus
\boldsymbol{\theta}_p^{\oplus m_p}$
induces a splitting of the punctured Cauchy-Riemann operator
\begin{equation}
\label{eqn:CRsplitting}
\mathbf{D}_u^N \cong \dot{\mathbf{D}}_{u,\boldsymbol{\rho}}^N = 
(\dot{\mathbf{D}}_{u,\boldsymbol{\theta}_1}^N)^{\oplus m_1}
\oplus \ldots \oplus
(\dot{\mathbf{D}}_{u,\boldsymbol{\theta}_p}^N)^{\oplus m_p},
\end{equation}
with the following useful property:

\begin{lemma}
\label{lemma:splitting}
The normal Cauchy-Riemann operator $\mathbf{D}_u^N$ for a multiple cover is 
surjective or injective if and only if
the same holds for all of the summands 
$\dot{\mathbf{D}}_{u,\boldsymbol{\theta}_j}^N$ in \eqref{eqn:CRsplitting}
with $m_j > 0$.
\end{lemma}

\begin{remark}
We will see below that the splitting \eqref{eqn:CRsplitting} for a
multiply covered curve $u = v \circ \varphi$ can be arranged
to vary smoothly as $v$ and $\varphi$ move about in their respective (suitably constrained) moduli
spaces, so the indices of the summands $\dot{\mathbf{D}}_{u,\boldsymbol{\theta}_j}^N$
are constant under such variations.  This immediately gives rise to ``no-go'' results about
transversality and super-rigidity: the former is impossible on components
of the moduli space where the
$\dot{\mathbf{D}}_{u,\boldsymbol{\theta}_j}^N$ do not all have nonnegative index,
and the latter requires them instead to have nonpositive index.  Conversely, whenever
either of these index conditions holds for all summands given by
irreducible representations, Theorem~\ref{thm:submanifolds0} below will 
imply that the desired transversality or super-rigidity result holds for
all pairs $(v,\varphi)$ lying in some open and dense subset.  This is the
main idea behind Theorem~\ref{thm:branched}, and it similarly can be used
to determine the feasibility of obstruction bundle arguments in general
situations.
\end{remark}

It should be emphasized that the representations of $G$ in this discussion
are \emph{real}, not complex.  We will
need to use the standard fact (see \S\ref{sec:representations})
that for any finite group~$G$, real irreducible representations 
$\boldsymbol{\theta} : G \to \Aut_\RR(W)$ come in three types, characterized
via the algebra $\KK := \End_G(W)$ of $G$-equivariant real-linear maps $W \to W$:
\begin{itemize}
\item \defin{Real type}: $\KK \cong \RR$;
\item \defin{Complex type}: $\KK \cong \CC$;
\item \defin{Quaternionic type}: $\KK \cong \HH$.
\end{itemize}
The endomorphism algebra $\KK = \End_G(W)$ endows the domain and target of the
operator $\dot{\mathbf{D}}_{u,\boldsymbol{\theta}}^N$ with $\KK$-module structures,
for which $\dot{\mathbf{D}}_{u,\boldsymbol{\theta}}^N$ is
$\KK$-linear.\footnote{In cases where $\mathbf{D}_u^N$ is already 
complex linear with respect to the natural complex structure on~$N_u$,
it is important to keep in mind that this natural complex structure has
nothing to do with the one induced on $\dot{\mathbf{D}}_{u,\boldsymbol{\theta}}^N$
when $\KK = \CC$.  In fact, these are two distinct complex structures that
commute with each other, and $\dot{\mathbf{D}}_{u,\boldsymbol{\theta}}^N$
is then complex linear with respect to both of them.}

The purpose of the following definition
will become clear in the statement of Theorem~\ref{thm:submanifolds0} below;
it is independent of choices due to the uniqueness of minimal regular presentations.

\begin{defn}
\label{defn:codimCurve}
The \defin{codimension} $\codim(u) \ge 0$ of the closed, connected,
$d$-fold covered $J$-holomorphic curve $u = v \circ \varphi$ is
a nonnegative integer defined as follows.  Choose a minimal regular presentation
$(\Theta,\dot{\Sigma}'',\pi,G,\rho,I,f)$ of $\varphi$
and a complete list of pairwise non-isomorphic irreducible real
representations
$\left\{\boldsymbol{\theta}_i : G \to \Aut_\RR(W_i)\right\}_{i=1,\ldots,p}$ of~$G$,
whose equivariant endomorphism algebras we denote by
$$
\KK_i := \End_G(W_i) \in \{\RR,\CC,\HH\}, \qquad i=1,\ldots,p.
$$
Then
$$
\codim(u) := \sum_{i=1}^p t_i k_i c_i,
$$
where $t_i := \dim_\RR \KK_i \in \{1,2,4\}$,
$k_i := \dim_{\KK_i} \ker \dot{\mathbf{D}}_{u,\boldsymbol{\theta}_i}^N$
and $c_i := \dim_{\KK_i} \coker \dot{\mathbf{D}}_{u,\boldsymbol{\theta}_i}^N$
for $i=1,\ldots,p$.
\end{defn}

\begin{example}
When $d = 1$, $u$ is a simple curve and its generalized automorphism
group $G$ is trivial, so there is only the trivial representation
$\boldsymbol{\theta} : G \to \Aut_\RR(\RR)$ to consider in
Definition~\ref{defn:codimCurve}, with $\End_G(\RR) = \RR$
and $\dot{\mathbf{D}}_{u,\boldsymbol{\theta}}^N \cong \mathbf{D}_u^N$.
So in this case,
$\codim(u) = \dim(\ker \mathbf{D}_u^N) \cdot \dim(\coker \mathbf{D}_u^N)$
can be interpreted as a measurement of the failure of transversality at~$u$, and the
standard transversality results imply that all simple curves have
codimension~$0$ for generic~$J$.  One of the consequences of
Theorem~\ref{thm:submanifolds0} will be that generically, this is also
true for \emph{generic} curves in the space of multiple covers,
though not necessarily for all of them.
\end{example}

\subsubsection{Isosymmetric strata}

In order to discuss what happens to the splitting of Cauchy-Riemann operators
\eqref{eqn:CRsplitting} as $v$ and
$\varphi$ move in their respective moduli spaces, we observe that the
construction depends quite heavily on the branching structure of
$\varphi : \Sigma' \to \Sigma$, i.e.~the number of punctures $\Theta' \subset
\Sigma'$ and the topological behavior of $\varphi$ in their vicinity.
This necessitates decomposing the space of all degree~$d$ branched covers
into strata
$$
\bigcup_{h \ge 0} \mM_h(d[\Sigma],j) = \bigcup_{\mathbf{b}} \mM^d_{\mathbf{b}}(j)
$$
labeled by their so-called \emph{branching data}~$\mathbf{b}$.  
Choose an integer $r \ge 0$, and associate to each of the numbers 
$i=1,\ldots,r$ a nonempty finite ordered set of natural numbers
$$
\mathbf{b}_i = ( b_i^1,\ldots,b_i^{q_i} )
$$
such that
$$
b_i^1 + \ldots + b_i^{q_i} = d
$$
and at least one of the numbers $b_i^1,\ldots,b_i^{q_i}$ is strictly greater
than~$1$.  We denote the totality of this data by $\mathbf{b} = 
(\mathbf{b}_1,\ldots,\mathbf{b}_r)$ and call it \defin{branching data of
degree~$d$ with $r$ critical values}.  
Given this, let $\widetilde{\mM}_{\mathbf{b}}^d(j)$ denote the moduli
space of all closed and connected unparametrized $j$-holomorphic 
curves $\varphi$ of degree~$d$ mapping into $(\Sigma,j)$ with $q_1 + \ldots + q_r$ marked 
points
$$
\zeta_1^1,\ldots,\zeta_1^{q_1},\zeta_2^1,\ldots,\zeta_2^{q_2},\ldots,
\zeta_r^1,\ldots,\zeta_r^{q_r}
$$
such that
\begin{enumerate}
\item there are distinct points $w_1,\ldots,w_r \in \Sigma$ such that
$\varphi^{-1}(w_i) = \{ \zeta_i^1,\ldots,\zeta_i^{q_i} \}$ for each
$i=1,\ldots,r$;
\item for each $i=1,\ldots,r$ and $j=1,\ldots,q_i$, $\varphi$ is
$b_i^j$-to-$1$ on a punctured neighborhood of~$\zeta_i^j$;
\item $\varphi$ has no critical points outside of the marked points.
\end{enumerate}
Note that we do not require \emph{every} marked point of $\varphi$ to be
a critical point, but we are assuming $\{w_1,\ldots,w_r\}$ is the set of
critical values, whose preimages are marked points and
may include both critical and regular points.
For any $\varphi \in \widetilde{\mM}_{\mathbf{b}}(j)$, we have
$$
Z(d\varphi) = \sum_{i=1}^r \sum_{j=1}^{q_i} ( b_i^j - 1),
$$
thus $d$ and $\mathbf{b}$ determine the genus $h$ of $\varphi$ via the
Riemann-Hurwitz formula, and we shall denote by
$$
\mM^d_{\mathbf{b}}(j) \subset \mM_h(d[\Sigma],j)
$$
the image of the natural map $\widetilde{\mM}^d_{\mathbf{b}}(j) \to \mM_h(d[\Sigma],j)$
defined by forgetting the marked points.
Note that in some cases, the Riemann-Hurwitz calculation may produce a negative genus,
which just means that $\mM_{\mathbf{b}}^d(j)$ is empty.  If 
$\mathbf{b}$ is empty, i.e.~$r=0$, 
it means every $\varphi \in \mM_{\mathbf{b}}^d(j)$ is unbranched.  

It is a classical fact that $\mM_{\mathbf{b}}^d(j)$ is a smooth manifold of
real dimension $2r$, as it can be parametrized locally by the positions of
the critical values $w_1,\ldots,w_r \in \Sigma$ (cf.~Example~\ref{ex:Teichmueller}).  
Moreover, it depends smoothly
on $j$ in the sense that if $P$ is any smooth finite-dimensional family of
complex structures on~$\Sigma$, then
$$
\bigcup_{j \in P} \mM_{\mathbf{b}}^d(j) \to P
$$
defines a smooth fiber bundle.  
We will show in \S\ref{sec:regPres} that regular presentations 
of $\varphi : \Sigma' \to \Sigma$ can also be arranged to vary
smoothly as $\varphi$ varies with fixed branching data.

Constraints must also be imposed on the simple $J$-holomorphic curve $v$ so that
the normal Cauchy-Riemann operators $\mathbf{D}_v^N$ and $\mathbf{D}_u^N$ vary smoothly as
$v$ moves in its moduli space.
Given integers $m \ge 0$ and $\ell_1,\ldots,\ell_m \ge 1$, let
$$
\mM_{g,m}(A,J\,;\,\ell_1,\ldots,\ell_m) \subset \mM_{g,m}(A,J)
$$
denote the subset consisting of curves that have critical points of critical
order $\ell_i$ at the $i$th marked point for $i=1,\ldots,m$ and are 
immersed everywhere else.
As explained in Appendix~\ref{sec:critical}, the simple
curves in this space form a smooth submanifold
for generic~$J$, with codimension $2n \sum_i \ell_i$ in $\mM_{g,m}(A,J)$.
Moreover, the generalized normal bundles $N_v$ of curves
$v \in \mM_{g,m}(A,J\,;\,\ell_1,\ldots,\ell_m)$ can be regarded as a smooth
family (cf.~Lemma~\ref{lemma:contin}).
This is not generally true if $v$ is allowed to move freely in
$\mM_{g,m}(A,J)$, as the topology of $N_v$ changes when critical
points of $v$ appear, disappear or change order.

Given an integer $d \in \NN$ and branching data
$\mathbf{b}$ of degree~$d$ with $r \ge 0$ critical values, define
$$
\mM^d_{\mathbf{b}}(\mM_{g,m}(A,J\,;\,\ell_1,\ldots,\ell_m)) \subset
\mM_h(dA,J)
$$
to be the set of all curves admitting representatives of the form
$u = v \circ \varphi : (\Sigma',j') \to (M,J)$, where
$\varphi : (\Sigma',j') \to (\Sigma,j)$ parametrizes an element in
$\mM^d_{\mathbf{b}}(j)$
and $v : (\Sigma,j) \to (M,J)$ is a simple curve that intersects $\uU$ and
(after labeling its critical points as marked points in a suitable order)
parametrizes an element of $\mM_{g,m}(A,J\,;\,\ell_1,\ldots,\ell_m)$.
If $J$ is generic on~$\uU$, then standard results give
$\mM^d_{\mathbf{b}}(\mM_{g,m}(A,J\,;\,\ell_1,\ldots,\ell_m))$ the structure
of a smooth manifold with
$$
\dim \mM^d_{\mathbf{b}}(\mM_{g,m}(A,J\,;\,\ell_1,\ldots,\ell_m)) =
2r + (n-3) (2-2g) + 2 c_1(A) - 2 \sum_{i=1}^m (n \ell_i - 1).
$$
Since every closed connected $J$-holomorphic curve belongs to
such a space for a unique (up to ordering) choice of branching data
$\mathbf{b}$ and critical orders $\ell_1,\ldots,\ell_m$, these spaces
form a smooth stratification of the moduli space of all $J$-holomorphic
curves.  They are sometimes called \defin{isosymmetric strata}, as they
have the property that all curves in the same connected component
of $\mM^d_{\mathbf{b}}(\mM_{g,m}(A,J\,;\,\ell_1,\ldots,\ell_m))$ have
isomorphic generalized automorphism groups.  More importantly,
each isosymmetric stratum admits a smooth family of normal
Cauchy-Riemann operators $\mathbf{D}_u^N$ with a smooth family of
splittings as in \eqref{eqn:CRsplitting} with respect to the irreducible
representations of their generalized automorphism groups.

\subsubsection{Walls}

Here is the main stratification result.

\begin{theorem}[stratification]
\label{thm:submanifolds0}
There exists a Baire subset 
$$
\jJ^\reg \subset \jJ(M,\omega\,;\,\uU,\Jfix)
$$
such that the following holds for all $J \in \jJ^\reg$.  For all choices of
integers $g , m \ge 0$, $d,\ell_1,\ldots,\ell_m \ge 1$,
branching data $\mathbf{b}$ of degree~$d$
and homology classes $A \in H_2(M)$,
the smooth isosymmetric stratum
$\mM^d_{\mathbf{b}}(\mM_{g,m}(A,J\,;\,\ell_1,\ldots,\ell_m))$ is a union
of countably many pairwise disjoint connected smooth submanifolds, referred to in
the following as \defin{walls}, which have the following properties:
\begin{enumerate}
\item For $u \in \mM^d_{\mathbf{b}}(\mM_{g,m}(A,J\,;\,\ell_1,\ldots,\ell_m))$,
the vector spaces $\ker\mathbf{D}_u^N$ and $\coker\mathbf{D}_u^N$
form the fibers of smooth vector bundles over each wall;
\item The codimension in $\mM^d_{\mathbf{b}}(\mM_{g,m}(A,J\,;\,\ell_1,\ldots,\ell_m))$
of the wall containing any given curve $u$ is $\codim(u)$.
\end{enumerate}
\end{theorem}

\begin{remark}
\label{remark:bifurcation}
The statement of Theorem~\ref{thm:submanifolds0} is specifically geared toward the applications treated
in this paper, but for different purposes one could formulate various other
versions, e.g.~one could add more marked points to $\mM_{g,m}(A,J\,;\,\ell_1,\ldots,\ell_m)$
and impose intersection constraints on them, or one could consider generic
finite-dimensional families $\{J_s\}_{s \in P}$ of almost complex structures 
and thus replace $\mM_{g,m}(A,J\,;\,\ell_1,\ldots,\ell_m)$ with a parametric 
moduli space of pairs $(u,s)$ where $s \in P$ and $u$ is $J_s$-holomorphic.
Either would require no serious modifications to the proof, other than more
cumbersome notation (cf.~Remark~\ref{remark:families}).
\end{remark}

\begin{remark}
\label{remark:whatisawall}
A natural guess for the precise definition of the walls mentioned in
Theorem~\ref{thm:submanifolds0} would be that they are maximal connected
subsets of $\mM^d_{\mathbf{b}}(\mM_{g,m}(A,J\,;\,\ell_1,\ldots,\ell_m))$
satisfying the constraint that $\dim \ker\mathbf{D}_u^N$ and
$\dim \coker\mathbf{D}_u^N$ are constant.  In fact, smooth walls can be
defined in that way using the methods of \cite{DoanWalpuski:BrillNoether},
but the actual definition used in this paper is slightly more complicated:
it requires a choice of a smooth family of minimal regular presentations,
and the constraint to impose is then that for every finite-dimensional
representation $\boldsymbol{\theta}$ of the resulting generalized automorphism
group, the kernels and cokernels of the twisted
Cauchy-Riemann operators $\dot{\mathbf{D}}_{u,\boldsymbol{\theta}}^N$
should have constant dimension as $u$ varies in the wall.
This would give the same result as
the simpler definition if one could guarantee that every summand 
in the splitting \eqref{eqn:CRsplitting} of $\mathbf{D}_u^N$ appears with positive multiplicity,
i.e.~that $m_i > 0$ for each of the irreducible representations~$\boldsymbol{\theta}_i$, 
but the latter is not always true.  As a consequence,
a maximal connected subset on which $\ker\mathbf{D}_u^N$ and $\coker\mathbf{D}_u^N$
have constant dimension may in general contain multiple walls of varying codimensions,
distinguished from each other by twisted Cauchy-Riemann operators corresponding
to representations that play no role in the splitting of~$\mathbf{D}_u^N$.
This phenomenon is harmless: the important detail for our purposes is that
whenever transversality or super-rigidity fails for a particular curve~$u$, it implies
that $u$ belongs to a wall whose codimension is positive and satisfies certain estimates.  
The converse is neither true nor necessary.
\end{remark}

We need two further ingredients in order to turn Theorem~\ref{thm:submanifolds0} 
into a powerful enough
tool for proving the theorems of \S\ref{sec:mainTheorems}.  The first is an
index calculation for the twisted operators $\dot{\mathbf{D}}_{u,\boldsymbol{\theta}}^N$.
The precise result is stated and proved in \S\ref{sec:index}, but for the
main applications we only need the following estimate, which is a corollary:

\begin{lemma}
\label{lemma:index}
Given a $J$-holomorphic curve $v : (\Sigma,j) \to (M,J)$ with normal
Cauchy-Riemann operator $\mathbf{D}_v^N$,
a $d$-fold branched cover $\varphi : (\Sigma',j') \to (\Sigma,j)$ with
$r \ge 0$ critical values, a regular presentation $(\Theta,\dot{\Sigma}'',\pi,G,\rho,I,f)$ 
for $\varphi$ and a
representation $\boldsymbol{\theta} : G \to \Aut_\RR(W)$, the resulting
twisted Cauchy-Riemann operator $\dot{\mathbf{D}}_{u,\boldsymbol{\theta}}^N$
for $u = v \circ \varphi$ satisfies
$$
\dim W \cdot \left[ \ind (\mathbf{D}_v^N) - (n-1) r \right]
\le \ind(\dot{\mathbf{D}}_{u,\boldsymbol{\theta}}^N)  \le
\dim W \cdot \ind (\mathbf{D}_v^N).
$$
Moreover, if the regular presentation is minimal and
$\boldsymbol{\theta}$ is a faithful irreducible representation with
$\End_G(W) \cong \KK \in \{\RR,\CC,\HH\}$, then the second estimate can be improved to
$$
\ind_\KK(\dot{\mathbf{D}}_{u,\boldsymbol{\theta}}^N)  \le
\dim_\KK W \cdot \ind_\RR (\mathbf{D}_v^N) - (n-1) r,
$$
and this estimate is strict in the case $\KK = \RR$
unless all branch points of $\varphi$ have branching order~$2$.
\end{lemma}

For the proof of super-rigidity, we will need the next result as a means of 
improving the upper bound in Lemma~\ref{lemma:index} for representations that 
are not faithful.

\begin{lemma}[see \S\ref{sec:nonfaithfulRevisit}]
\label{lemma:nonfaithful}
Under the assumptions of Lemma~\ref{lemma:index}, suppose the regular
presentation is minimal, and
the splitting \eqref{eqn:CRsplitting} of $\mathbf{D}_u^N$ includes a summand
$\dot{\mathbf{D}}_{u,\boldsymbol{\theta}}^N$ for which the representation
$\boldsymbol{\theta} : G \to \Aut_\RR(W)$ is not faithful.  
Then $\varphi : (\Sigma',j') \to (\Sigma,j)$ admits a factorization by 
holomorphic branched covers
$$
(\Sigma',j') \to (\Sigma'_0,j'_0) \stackrel{\varphi_0}{\longrightarrow} (\Sigma,j)
$$
with $\deg(\varphi_0) < d$, and
$\dot{\mathbf{D}}_{u,\boldsymbol{\theta}}^N$ is conjugate to
an operator $\dot{\mathbf{D}}_{u_0,\boldsymbol{\theta}_0}^N$ defined with
respect to a regular presentation $(\Theta,\dot{\Sigma}''_0,\pi_0,G_0,\rho_0,I_0,f_0)$
for $\varphi_0$, where $u_0 := v \circ \varphi_0 : (\Sigma'_0,j'_0) \to (M,J)$, 
$G_0 := G / \ker \boldsymbol{\theta}$, and 
$$
\boldsymbol{\theta}_0 : G / \ker\boldsymbol{\theta} \to \Aut_\RR(W)
$$
is the faithful representation of $G_0$ determined by~$\boldsymbol{\theta}$.
Moreover, $\mathbf{D}_{u_0}^N$ also admits a splitting in the form
\eqref{eqn:CRsplitting} which has
$\dot{\mathbf{D}}_{u_0,\boldsymbol{\theta}_0}^N$ as a summand.
\end{lemma}

\subsection{Proof of the main theorems modulo stratification}
\label{sec:everything}

Let us now take the results of the previous section as black boxes and
prove the main theorems from \S\ref{sec:mainTheorems}.

\begin{proof}[Proof of Theorem~\ref{thm:super} (super-rigidity) in dimension greater than four]
We argue by induction on the degrees $d \in \NN$ of branched covers.
For $d=1$, we only need to know that generic perturbations of $J$ suffice
to make all simple index~$0$ curves through $\uU$ regular and immersed; this 
is standard (see Appendix~\ref{sec:critical} for the immersion property).  
Thus for $d \ge 2$, assume we
have already found a Baire subset in $\jJ(M,\omega\,;\,\uU,\Jfix)$ for which
all branched covers $u := v \circ \varphi$ with $v : (\Sigma,j) \to (M,J)$
a simple curve of index~$0$ and $\deg(\varphi) \le d-1$
have $\mathbf{D}_u^N$ injective.  Suppose $\varphi \in \mM^d_{\mathbf{b}}(j)$
has $r \ge 0$ critical values and $\deg(\varphi) = d$ and 
$\mathbf{D}_u^N$ is not injective for $u := v \circ \varphi$.  Then
picking the minimal regular presentation $(\Theta,\dot{\Sigma}'',\pi,G,\rho,I,f)$
for~$\varphi$ and decomposing $\boldsymbol{\rho}$ into irreducible
representations $\boldsymbol{\theta}_1^{\oplus \ell_1} \oplus \ldots \oplus
\boldsymbol{\theta}_p^{\oplus \ell_p}$ of $G$ splits $\mathbf{D}_u^N$ into
twisted Cauchy-Riemann operators
$\dot{\mathbf{D}}_{u,\boldsymbol{\theta}_i}^N$ for $i=1,\ldots,p$ with
$$
k_i := \dim_{\KK_i} \ker \dot{\mathbf{D}}_{u,\boldsymbol{\theta}_i}^N,
$$
and at least one of the $k_i$ must be strictly positive by
Lemma~\ref{lemma:splitting}.  If $k_i > 0$ and $\boldsymbol{\theta}_i$ is
non-faithful, then Lemma~\ref{lemma:nonfaithful} identifies
$\dot{\mathbf{D}}_{u,\boldsymbol{\theta}_i}^N$ with a summand of
$\mathbf{D}_{u_0}^N$ for some other cover $u_0$ of $v$ with strictly
smaller degree, implying $\dim \ker \mathbf{D}_{u_0}^N > 0$ and thus
violating the inductive hypothesis.  We can therefore assume $k_i > 0$
for some faithful representation $\boldsymbol{\theta}_i$.  But then
Theorem~\ref{thm:submanifolds0} and Lemma~\ref{lemma:index} imply that $u$ lives in a submanifold
of the $2r$-dimensional space of branched covers of $v$ with branching
data~$\mathbf{b}$, having dimension at most
$$
2r - t_i k_i \left[k_i - \ind_{\KK_i}(\dot{\mathbf{D}}_{u,\boldsymbol{\theta}_i}^N)\right]
\le 2r - t_i k_i[k_i + (n-1)r] =
r [2 - t_i k_i (n-1)] - t_i k_i^2 < 0
$$
since we are assuming $n \ge 3$.  This gives a contradiction and thus completes
the induction.
\end{proof}

In dimension four, the above argument fails to exclude the possibility of
$\dim \ker \dot{\mathbf{D}}_{u,\boldsymbol{\theta}_i}^N = 1$ for some
real-type repesentation~$\boldsymbol{\theta}_i$, and this is why we do not know whether super-rigidity
always holds in dimension four.  We will prove in \S\ref{sec:dimensionFour}
that it does hold for covers of genus zero and one curves, using 
different techniques based on intersection theory.

\begin{proof}[Proof of Theorem~\ref{thm:unbranched} (transversality, unbranched)]
Suppose $v : (\Sigma,j) \to (M,J)$ is a simple curve intersecting $\uU$
and $\varphi : (\Sigma',j') \to (\Sigma,j)$ is a $d$-fold unbranched cover
for which $u := v \circ \varphi$ is not Fredholm regular, hence by
Prop.~\ref{prop:regular}, $\mathbf{D}_u^N$ is not surjective.
Fixing the minimal regular presentation of $\varphi$ and considering the
splitting \eqref{eqn:CRsplitting}, we find a twisted Cauchy-Riemann
operator $\dot{\mathbf{D}}_{u,\boldsymbol{\theta}_i}^N$ with 
$$
c_i := \dim_{\KK_i} \coker \dot{\mathbf{D}}_{u,\boldsymbol{\theta}_i}^N > 0
$$
for some irreducible representation $\boldsymbol{\theta}_i : G \to
\Aut_\RR(W_i)$ of the
generalized automorphism group $G$ of~$\varphi$, with $\End_G(W_i) \cong
\KK_i \in \{\RR,\CC,\HH\}$.  Suppose $v$ has
exactly $m \ge 0$ critical points, with critical orders
$\ell_1,\ldots,\ell_m$, so viewing these as marked points allows us to consider
$v$ as an element in the space $\mM_{g,m}(A,J\,;\,\ell_1,\ldots,\ell_m)$,
which has dimension 
$$
\dim \mM_{g,m}(A,J\,;\,\ell_1,\ldots,\ell_m) = 
\ind(v) + 2m - 2 n Z(dv) \ge 0.
$$
The count of critical points $Z(dv)$ also appears in the relation between
$\ind(v)$ and $\ind \mathbf{D}_v^N$: indeed, writing $v^*TM = T_v \oplus N_v$,
we can view $dv$ as a holomorphic section of $\Hom_\CC(T\Sigma,T_v)$,
hence
$$
Z(dv) = c_1\big(\Hom_\CC(T\Sigma,T_v)\big) = -c_1(T\Sigma) + c_1(T_v) =
-\chi(\Sigma) + c_1(T_v),
$$
implying $c_1(N_v) = c_1(v^*TM) - c_1(T_v) = c_1(v^*TM) -\chi(\Sigma) - Z(dv)$.
Plugging in this into the Riemann-Roch formula then gives
\begin{equation*}
\begin{split}
\ind \mathbf{D}_v^N &= (n-1) \chi(\Sigma) + 2 c_1(N_v) =
(n-3) \chi(\Sigma) + 2 c_1(v^*TM) - 2 Z(dv) \\
&= \ind(v) - 2 Z(dv).
\end{split}
\end{equation*}
Meanwhile, $\varphi$ lives in a discrete stratum of the space of branched 
covers since it has no branch points, and Lemma~\ref{lemma:index} reduces to
an equality
$$
\ind_{\KK_i} \dot{\mathbf{D}}_{u,\boldsymbol{\theta}_i}^N =
\dim_{\KK_i} W_i \cdot \ind_\RR(\mathbf{D}_v^N).
$$
Now using Theorem~\ref{thm:submanifolds0}, we find that if $J$ is generic, 
$u$ lives in a manifold of dimension at most
\begin{equation*}
\begin{split}
\dim \mM_{g,m}(&A,J\,;\,\ell_1,\ldots,\ell_m) - t_i c_i (c_i + \ind_{\KK_i} \dot{\mathbf{D}}_{u,\boldsymbol{\theta}_i}^N) \\
&= \ind(v) + 2m - 2 n Z(dv) - t_i c_i (c_i + \dim_{\KK_i} W_i \cdot \ind \mathbf{D}_v^N) \\
&= \ind(v) + 2m - 2 n Z(dv) - t_i c_i (c_i + \dim_{\KK_i} W_i \cdot \left[ \ind(v) - 2 Z(dv) \right] ) \\
&= (1 - t_i c_i \dim_{\KK_i} W_i) \left[ \ind(v) + 2m - 2 n Z(dv)\right] \\
& \qquad - 2 t_i c_i \dim_{\KK_i} W_i \cdot \left[ (n-1) Z(dv) - m \right] - t_i c_i^2 < 0,
\end{split}
\end{equation*}
where we note that $(n-1) Z(dv) - m \ge 0$ since $n \ge 2$ and every critical
point has order at least~$1$.
\end{proof}

\begin{proof}[Proof of Theorem~\ref{thm:branched} (transversality, branched)]
Assume $v : (\Sigma,j) \to (M,J)$ is simple and satisfies
$\ind(v) \ge (n-1) r$, while $\varphi : (\Sigma',j') \to (\Sigma,j)$ has
degree $d \in \NN$ and $r$ critical values.  If $J$ is generic, then by 
Proposition~\ref{prop:critical} the moduli space containing $v$
has an open and dense subset consisting of immersed curves, so we are free 
to assume $v$ is immersed and thus $\ind(v) = \ind \mathbf{D}_v^N$.
The key observation is then that by Lemma~\ref{lemma:index}, the twisted
operators $\dot{\mathbf{D}}_{u,\boldsymbol{\theta}}^N$ all have nonnegative
index, hence Theorem~\ref{thm:submanifolds0} implies that all of them are surjective
unless $(v,\varphi)$ lies in a countable union of submanifolds with positive
codimension.
\end{proof}

\subsection{Some remarks on wall crossing}
\label{sec:bifurcations}

Part of the point of Taubes's twisted bundle setup in \cite{Taubes:counting}
was to understand bifurcations of isolated $J$-ho\-lo\-mor\-phic tori under
generic $1$-parameter deformations in~$J$.  While bifurcation theory is
not the main topic of this article, it should be clear that such a theory
could be developed based on Theorem~\ref{thm:submanifolds0}, thus we take
this opportunity to make a few observations about it.

\begin{remark}
In the time since the present article first appeared in preprint form, some interesting cases of the bifurcation
analysis proposed below have been worked out in detail by Bai and Swaminathan,
see \cite{BaiSwaminathan}.
\end{remark}

If $\{J_s\}_{s \in [0,1]}$ is a generic homotopy of compatible almost
complex structures whose endpoints are generic, then as mentioned in
Remark~\ref{remark:bifurcation}, one can modify Theorem~\ref{thm:submanifolds0}
to the statement that the parametric moduli space
$$
\mM^d_{\mathbf{b}}(\mM_{g,m}(A,\{J_s\}\,;\,\ell_1,\ldots,\ell_m))
$$
consisting of pairs $(u,s)$ where $s \in [0,1]$ and
$u \in \mM^d_{\mathbf{b}}(\mM_{g,m}(A,J_s\,;\,\ell_1,\ldots,\ell_m))$
is stratified by smooth submanifolds characterized by the dimensions of the
kernels and cokernels of twisted Cauchy-Riemann operators, and their
codimensions are given by the same formula.  In this setting,
suppose $\{v_\tau\}$ is a smooth $1$-parameter family of 
simple $J_{s(\tau)}$-holomorphic curves with
index~$0$ for some function $s(\tau) \in [0,1]$, and $\{u_\tau = v_\tau \circ \varphi_\tau\}$
defines a corresponding $1$-parameter family of unbranched covers.  The latter
have index~$0$ and will be regular for almost every~$\tau$, but a bifurcation
or ``wall crossing'' phenomenon
occurs at any parameter value $\tau_0$ for which the family $\{u_\tau\}$ passes 
(necessarily transversely) through one of the codimension~$1$ walls given
by Theorem~\ref{thm:submanifolds0}.
When this happens, most of the twisted operators $\dot{\mathbf{D}}_{u_{\tau_0},\boldsymbol{\theta}}^N$
remain both injective and surjective, but there will be
exactly one irreducible representation $\boldsymbol{\theta}$ for which
$$
\dim \ker \dot{\mathbf{D}}_{u_{\tau_0},\boldsymbol{\theta}}^N =
\dim \coker \dot{\mathbf{D}}_{u_{\tau_0},\boldsymbol{\theta}}^N = 1,
$$
and $\boldsymbol{\theta}$ is necessarily of real type.
Whenever $\boldsymbol{\theta}$ is not faithful, one can factor $\varphi_\tau$
through a cover $\widehat{\varphi}_\tau$ of smaller degree and instead examine
$\widehat{u}_\tau := v_\tau \circ \widehat{\varphi}_\tau$, so that $\boldsymbol{\theta}$
becomes faithful without loss of generality (cf.~Lemma~\ref{lemma:nonfaithful}).  
For the trivial representation, this means replacing $u_\tau$ with $v_\tau$
itself, so regularity fails for the underlying simple curve at $\tau=\tau_0$: as shown in
\cite{Taubes:counting}, this is the case where the family $\{v_\tau\}$
undergoes a \emph{birth-death} bifurcation.  The other interesting phenomenon
examined by Taubes was the \emph{degree-doubling} bifurcation, in which
$v_\tau$ remains regular but it has a double cover $u_\tau = v_\tau \circ \varphi_\tau$
which loses regularity at $\tau=\tau_0$, causing an additional $1$-parameter
family of simple curves $\{w_\tau\}$ to collide with $\{u_\tau\}$ at $\tau=\tau_0$.
This is what happens when $\dot{\mathbf{D}}_{u_\tau,\boldsymbol{\theta}}^N$
remains an isomorphism for the trivial representation but acquires
$1$-dimensional kernel and cokernel for the nontrivial irreducible
representation of~$\ZZ_2$.  

In \cite{Taubes:counting}, no further bifurcations
beyond these two types are possible: this can be attributed to the fact that
since Taubes only considers unbranched covers of tori, all covers are regular 
and abelian.  As a consequence, all the complex irreducible representations in 
the picture are $1$-dimensional, implying that the only faithful real-type 
irreducible representations one needs to consider are the trivial 
representation of the trivial group and the nontrivial representation of~$\ZZ_2$.
We should not expect this fortunate situation to hold more generally:
for unbranched covers with higher genus, one certainly encounters generalized
automorphism groups that are non-abelian and thus have faithful real-type
representations of dimension greater than one.  These should presumably 
give rise to bifurcation phenomena involving covers of arbitrarily high degree.

In the context of super-rigidity, it is also important to consider bifurcations
that involve branched covers of index~$0$ curves under generic
homotopies of~$J$.  Inspecting the proof of Theorem~\ref{thm:super}, one should
expect to see interesting phenomena whenever the dimension that was estimated
at the end of the proof turns out to be at least~$-1$, i.e.
$$
2r - t_i k_i \left[ k_i - \ind_{\KK_i}(\dot{\mathbf{D}}_{u,\boldsymbol{\theta}_i}^N)\right] \ge -1.
$$
Assuming we're in dimension at least six, this can only mean $t_i = k_i = 1$
and either $r=0$ or $n=3$.  The case $r=0$ means the cover is unbranched, so
this is what we discussed in the previous paragraphs.  Bifurcations involving
branched covers can evidently also occur in dimension six, 
and in this case the improved index bound from Lemma~\ref{lemma:index} must
be an equality.  The scenario is therefore that the rank of the
obstruction bundle over the space of covers $\{v_\tau \circ \varphi_\tau\}$ 
jumps at a particular parameter value $\tau=\tau_0$ and for some 
isolated element $\varphi_{\tau_0}$ in the space of branched covers with
only simple (i.e.~two-to-one) branch points: this can presumably cause
both a change in the Euler class of the obstruction bundle and the
breaking off of a new family of simple curves from $v_{\tau_0} \circ \varphi_{\tau_0}$.
Once again the irreducible representation involved must be of real type
but can have arbitrary dimension, meaning we
should not expect any limitation on the degree of $\varphi_{\tau_0}$,
contrary to the situation in \cite{Taubes:counting}.

\section{Splitting Cauchy-Riemann operators with symmetries}
\label{sec:prep}

In this section we give a detailed account of the twisted bundle formalism
behind Theorem~\ref{thm:submanifolds0} and prove several lemmas required
for its proof, as well as Lemma~\ref{lemma:nonfaithful}.
Instead of talking directly about $J$-holomorphic curves, we shall work
in the context of abstract Cauchy-Riemann operators on vector bundles and
their pullbacks.

\subsection{Regular presentations of branched covers}
\label{sec:regPres}

The notion of a regular presentation was introduced in Definition~\ref{defn:presentation}.
The following standard result from the theory of Riemann surfaces
(see e.g.~\cite{Donaldson:RS}*{Chapter~4, Theorem~2})
allows us to move freely back and forth between
talking about holomorphic branched covers of closed Riemann surfaces and
honest covering maps of punctured surfaces.

\begin{lemma}
\label{lemma:closedUp}
Suppose $(\dot{\Sigma},j)$ is the complement of a finite set of points
$\Theta$ in a closed connected Riemann surface $(\Sigma,j)$,
$(\dot{\Sigma}',j')$ is a connected noncompact Riemann surface, and
$$
\varphi : (\dot{\Sigma}',j') \to (\dot{\Sigma},j)
$$
is a holomorphic covering map of finite degree.  
Then there exists a closed connected Riemann surface
$(\Sigma',j')$ with a finite set of points $\Theta' \subset \Sigma'$ such
that $(\dot{\Sigma}',j')$ admits a biholomorphic identification with
$(\Sigma' \setminus \Theta',j')$ and $\varphi$ extends over the punctures
to a holomorphic branched cover $\varphi : (\Sigma',j') \to (\Sigma,j)$ with
$\varphi^{-1}(\Theta) = \Theta'$.
\qed
\end{lemma}

Assume $\varphi : (\Sigma',j') \to (\Sigma,j)$ is a $d$-fold holomorphic
branched cover of closed connected Riemann surfaces with branching
data $\mathbf{b}$ as defined in \S\ref{sec:theBigIdea}, having $r \ge 0$ 
distinct critical values.  Recall from 
Definition~\ref{defn:presentation} that for a regular presentation
$(\Theta,\dot{\Sigma}'',\pi,G,\rho,I,f)$
of~$\varphi$, $\Theta \subset \Sigma$ is a finite set containing the critical
values of~$\varphi$, giving rise to the punctured surfaces
$$
\dot{\Sigma} := \Sigma \setminus \Theta, \qquad
\dot{\Sigma}' := \Sigma' \setminus \Theta',
$$
where $\Theta' := \varphi^{-1}(\Theta)$.

\begin{lemma}
\label{lemma:regPres}
There exists a natural bijection between the set of isomorphism classes of 
regular presentations of $\varphi$ and the set of pairs $(\Theta,H)$ where
$\Theta \subset \Sigma$ is a finite subset containing the critical values
of~$\varphi$ and $H$ is a finite-index normal subgroup
$H \subset \pi_1(\dot{\Sigma})$ that is contained in
$\varphi_*(\pi_1(\dot{\Sigma}'))$.  This bijection matches any
minimal regular presentation to the smallest possible choice of $\Theta$ and
largest possible choice of~$H$, i.e.~the normal core of $\varphi_*(\pi_1(\dot{\Sigma}'))$.  
Moreover, if $\varphi$ is regular and $(\Theta,\dot{\Sigma}'',\pi,G,\rho,I,f)$ is a
minimal regular presentation, then there exists a diffeomorphism
$g : \dot{\Sigma}' \to \dot{\Sigma}''$ such that $\pi \circ g = \varphi$.
\end{lemma}
\begin{proof}
Given a finite set $\Theta \subset \Sigma$ containing the critical values
of~$\varphi$, pick a base point $w \in \dot{\Sigma}$ and let $\tilde{\pi} : \univ \to
\dot{\Sigma}$ denote the universal cover, with $\univ$ defined as a
space of homotopy classes of paths beginning at $w$, so that 
$\pi_1(\dot{\Sigma}) := \pi_1(\dot{\Sigma},w)$ acts naturally on $\univ$
as the group of deck transformations for~$\tilde{\pi}$.
Lifting loops based
at $w$ to paths in $\dot{\Sigma}'$ then defines a homomorphism
$$
\tilde{\rho} : \pi_1(\dot{\Sigma}) \to S(\varphi^{-1}(w)) :
\gamma \mapsto \tilde{\rho}_\gamma
$$
so that the covering map $\dot{\Sigma}' \stackrel{\varphi}{\longrightarrow}
\dot{\Sigma}$ can be identified with
$$
\dot{\Sigma}' = \left( \univ \times \varphi^{-1}(w) \right) \Big/ \pi_1(\dot{\Sigma})
\to \dot{\Sigma} : [(z,\zeta)] \mapsto \tilde{\pi}(z),
$$
where $\pi_1(\dot{\Sigma})$ acts on $\univ$ by deck transformations and
on $\varphi^{-1}(w)$ via~$\tilde{\rho}$.  We claim that
$$
\ker \tilde{\rho} \subset \pi_1(\dot{\Sigma})
$$
is the normal core of $\varphi_*(\pi_1(\dot{\Sigma}'))$.
Indeed, selecting a base point $w' \in \varphi^{-1}(w) \subset \dot{\Sigma}'$
to define $\pi_1(\dot{\Sigma}') := \pi_1(\dot{\Sigma}',w')$, we have
$$
\varphi_*(\pi_1(\dot{\Sigma}')) = \left\{ \gamma \in \pi_1(\dot{\Sigma})\ \Big|\ 
\tilde{\rho}_\gamma(w') = w' \right\},
$$
which obviously contains $\ker \tilde{\rho}$.  Changing the base
point $w' \in \varphi^{-1}(w)$ changes the subgroup
$\varphi_*(\pi_1(\dot{\Sigma}'))$ by conjugation with arbitrary elements
of $\pi_1(\dot{\Sigma}')$, and the normal core is the intersection of
all these conjugates, which we can now recognize as the intersection of all
the stabilizers of the permutation action on $\varphi^{-1}(w)$, and that
is $\ker \tilde{\rho}$.

Suppose $H \subset \pi_1(\dot{\Sigma})$ is a finite-index normal subgroup
contained in $\varphi_*(\pi_1(\dot{\Sigma}'))$, and therefore also in
$\ker \tilde{\rho}$.  Then $\tilde{\rho}$ descends to the finite group
$G := \pi_1(\dot{\Sigma}) / H$, giving a homomorphism
$$
\rho : G \to S(\varphi^{-1}(w)),
$$
which is injective if and only if $H = \ker \tilde{\rho}$.  It is now possible to define a
regular presentation $(\Theta,\dot{\Sigma}'',\pi,G,\rho,\varphi^{-1}(w),f)$ of $\varphi$
with $\pi$ as the natural quotient projection
$$
\dot{\Sigma}'' := \univ / H \stackrel{\pi}{\longrightarrow}
\univ / \pi_1(\dot{\Sigma}) = \dot{\Sigma}
$$
and
$$
\dot{\Sigma}' = \left(\univ \times \varphi^{-1}(w) \right) \Big/ \pi_1(\dot{\Sigma})
\stackrel{f}{\longrightarrow} \left( \dot{\Sigma}'' \times \varphi^{-1}(w) \right) \Big/ G
$$
defined via the quotient projection $\univ \to \univ / H = \dot{\Sigma}''$.
Observe that if we choose $H = \ker\tilde{\rho}$ and $\varphi$ is regular,
then $\varphi_*(\pi_1(\dot{\Sigma}')) \subset \pi_1(\dot{\Sigma})$ is normal
and is therefore identical to~$H$, so the natural identification of
$\dot{\Sigma}'$ with $\univ \big/ \varphi_*(\pi_1(\dot{\Sigma}')) =
\univ / H = \dot{\Sigma}''$ gives an isomorphism between the covering maps
$\varphi$ and~$\pi$.

Finally, suppose $(\Theta,\dot{\Sigma}'',\pi,G,\rho,I,f)$ is a
regular presentation of~$\varphi$, and define the subgroup
$H := \pi_*(\pi_1(\dot{\Sigma}''))$, which is normal since
$\pi : \dot{\Sigma}'' \to \dot{\Sigma}$ is regular and has finite index
since $\Aut(\pi) = G = \pi_1(\dot{\Sigma}) / H$ is finite.  We claim
$H \subset \varphi_*(\pi_1(\dot{\Sigma}'))$: indeed, any 
$\gamma \in H$ is represented by a loop $\dot{\Sigma}$ based at $w$ that 
lifts to a loop $\gamma''$ in 
$\dot{\Sigma}''$ and thus has $d$ lifts to $\dot{\Sigma}' \cong (\Sigma'' \times I)/G$
in the form $\gamma \times \{i\}$ for $i \in I$.  We can therefore use
$H$ to define the regular presentation from the previous paragraph, with
$G = \pi_1(\dot{\Sigma}'') / H$ acting on $\varphi^{-1}(w)$ via $\tilde{\rho}$,
and we claim that this is isomorphic to $(\Theta,\dot{\Sigma}'',\pi,G,\rho,I,f)$.
Indeed, choosing a base point $w'' \in \pi^{-1}(w) \subset \dot{\Sigma}''$,
the identification
$f : \dot{\Sigma}' \to (\dot{\Sigma}'' \times I) / G$ provides a bijection
$$
\beta : \varphi^{-1}(w) \to I \quad \text{ such that } \quad
f(w') = [(w'',\beta(w'))] \text{ for $w' \in \varphi^{-1}(w)$},
$$
and combining this with the natural identification of $\dot{\Sigma}''$
with $\univ / H$ gives an isomorphism of regular presentations.
\end{proof}

\begin{lemma}
\label{lemma:branchPoints}
Suppose $(\Theta,\dot{\Sigma}'',\pi,G,\rho,I,f)$ is a minimal regular presentation
of $\varphi : (\Sigma',j') \to (\Sigma,j)$, and let $\pi : (\Sigma'',j'') \to (\Sigma,j)$
denote the branched cover of closed Riemann surfaces provided by
Lemma~\ref{lemma:closedUp} such that $\dot{\Sigma}'' = \Sigma'' \setminus \pi^{-1}(\Theta)$.
Then for each $w \in \Theta$ and $\zeta \in \pi^{-1}(w) \subset \Sigma''$,
the branching order of $\pi$ at $\zeta$ is the least common multiple of the
branching orders of $\varphi$ at all $z \in \varphi^{-1}(w)$.  In particular,
$\pi$ and $\varphi$ have the same sets of critical values.
\end{lemma}
\begin{proof}
If $k \in \NN$ is the branching order of $\pi$ at~$\zeta$, we can find punctured
neighborhoods $\uU_w \subset \dot{\Sigma}$ of $w$ and
$\uU_\zeta \subset \dot{\Sigma}''$ of $\zeta$ and identify both with the
half-cylinder $[0,\infty) \times S^1$ with coordinates $(s,t)$ such that
$\pi(s,t) = (ks,kt)$.  Let $G_\zeta \subset G$ denote the group of automorphisms
of $\pi$ that fix~$\zeta$; since $\pi : \dot{\Sigma}'' \to \dot{\Sigma}$ 
is a regular cover, $G_\zeta$ is necessarily a cyclic group of order~$k$,
with a generator $g \in G_\zeta$ that acts on $\uU_\zeta \cong [0,\infty) \times S^1$
as the rotation $(s,t) \mapsto (s,t + 1/k)$.
Appealing again to regularity, we can then restrict the identification 
$\dot{\Sigma}' = (\dot{\Sigma}'' \times I) / G$ to $\uU_\zeta$ and obtain an
identification
$$
\varphi^{-1}(\uU_w) = \left( \uU_\zeta \times I \right) \Big/ G_\zeta.
$$
The connected components
of $\varphi^{-1}(\uU_w)$ are then in bijective correspondence to the
orbits of the $G_\zeta$-action on~$I$ defined by~$\rho : G \to S(I)$, with
the branching order $k_z \in \NN$ of each corresponding point $z \in \varphi^{-1}(w)$
given by the number of points in its respective orbit in~$I$.  By the
orbit-stabilizer theorem, all of these numbers $k_z$ must divide~$k = |G_\zeta|$.
If $\ell$ is their least common multiple, we conclude that $g^\ell \in G_\zeta$
acts trivially on~$I$, which means $g^\ell$ is the identity since
$\rho : G \to S(I)$ is injective for the minimal regular presentation, 
hence $\ell = k$.
\end{proof}

\begin{example}
\label{ex:regular}
If $\varphi$ is regular with $\Aut(\varphi) = G$, then it admits a canonical
minimal regular presentation $(\Theta,\dot{\Sigma}'',\pi,G,\rho,I,f)$ where
$\dot{\Sigma}'' := \dot{\Sigma}'$, $\pi := \varphi$, $I := G$, and the action
$\rho : G \to S(G)$ of $G$ on itself is defined by left multiplication
$$
\rho_g(h) := gh.
$$
Here the identification $\dot{\Sigma}' \stackrel{f}{\longrightarrow} (\dot{\Sigma}'' \times G) / G$ 
sends $z \in \dot{\Sigma}'$ to $[(z,e)]$, where $e \in G$ is the identity element.
The action of
$G$ on $\dot{\Sigma}' = (\dot{\Sigma}'' \times G) / G$ by deck transformations can now be 
presented as the action via right multiplication
$$
G \times \dot{\Sigma}' \to \dot{\Sigma}' : (g,[(z,h)]) \mapsto
[(z,h g^{-1})].
$$
Notice that any regular presentation in which $\rho : G \to S(I)$ acts on $I$ both
transitively and without fixed points is isomorphic to one of this form, since for
any $i \in I$, the map $G \to I : g \mapsto \rho_g(i)$ defines a bijection that
transforms the action by left multiplication into~$\rho$.
\end{example}

\begin{example}
\label{ex:factorization}
The following construction underlies Lemma~\ref{lemma:nonfaithful}:
any proper normal subgroup $H \subset G$
gives rise to a factorization of $\varphi : (\Sigma',j') \to (\Sigma,j)$
in the following way.  Let $I / H$ denote the set of orbits for the 
action $\rho|_H : H \to S(I)$.
Then $G/H$ is a finite group and $\rho$ descends to a homomorphism
$$
\rho_H : G/H \to S(I/H),
$$
which acts transitively on~$I/H$.  The regular cover $\pi : \dot{\Sigma}''
\to \dot{\Sigma} = \dot{\Sigma}'' / G$ now factors through the obvious
projections
$$
\dot{\Sigma}'' \to \dot{\Sigma}''_H := \dot{\Sigma}'' / H 
\stackrel{\pi_H}{\longrightarrow} \dot{\Sigma} = \dot{\Sigma}'' / G,
$$
and $\pi_H : \dot{\Sigma}''_H \to \dot{\Sigma}$ is a regular holomorphic 
cover with automorphism group~$G/H$.  We can thus define
$$
\dot{\Sigma}'_H := \left( \dot{\Sigma}''_H \times (I/H) \right) \Big/ (G/H)
\stackrel{\varphi_H}{\longrightarrow} \dot{\Sigma} : [(z,i)] \mapsto
\pi_H(z),
$$
as well as a factorization of $\varphi : \dot{\Sigma}' \to \dot{\Sigma}$
by covering maps
$$
\dot{\Sigma}' = \left( \dot{\Sigma}'' \times I \right)\Big/ G
\longrightarrow
\dot{\Sigma}'_H \stackrel{\varphi_H}{\longrightarrow}
\dot{\Sigma},
$$
where the first map is also defined via the obvious quotient projections.  
It follows from Lemma~\ref{lemma:closedUp} that $\dot{\Sigma}'_H$ and
$\dot{\Sigma}''_H$ each arise by
puncturing closed connected Riemann surfaces $(\Sigma'_H,j'_H)$ and
$(\Sigma''_H,j''_H)$ respectively, and in particular we obtain a factorization
of $\varphi$ via holomorphic branched covers
$$
(\Sigma',j') \to (\Sigma'_H,j'_H) \stackrel{\varphi_H}{\longrightarrow} (\Sigma,j)
$$
with $\deg(\varphi_H) \le d$ equal to the number of distinct orbits of the
$H$-action on~$I$, hence
$$
\deg(\varphi_H) < d
$$
holds whenever the action of $H$ on $I$ is nontrivial.  Note that
$\varphi_H$ inherits from this construction a regular presentation
$(\Theta,\dot{\Sigma}''_H,\pi_H,G/H,\rho_H,I/H,f_H)$, though it need not
be minimal and $\Theta$ may contain points that are not critical values
of~$\varphi_H$, even if $(\Theta,\dot{\Sigma}'',\pi,G,\rho,I,f)$ is minimal.
This is the main reason why non-minimal regular presentations have been
included in the discussion.
\end{example}

It will be important to understand how the various objects constructed out of
a regular presentation vary smoothly under changes in $\varphi$ and~$j$.
To this end, we shall fix the following data for the remainder of 
\S\ref{sec:prep}:
\begin{itemize}
\item $\varphi : (\Sigma',j') \to (\Sigma,j)$ is a holomorphic branched cover
of degree $d \in \NN$ with branching data $\mathbf{b}$;
\item $(\Theta,\dot{\Sigma}'',\pi,G,\rho,I,f)$ is a regular presentation of~$\varphi$;
\item $P$ is a connected smooth Banach manifold;
\item $\vV \subset \dot{\Sigma}$ is an open subset with compact closure;
\item $\{j_\tau\}_{\tau \in P}$ is a smooth family of complex structures 
on $\Sigma$ that match $j$ outside of~$\vV$;
\item $\{\psi_\tau\}_{\tau \in P}$ is a smooth family of diffeomorphisms
$\psi_\tau : \Sigma \to \Sigma$ which restrict to the identity on~$\vV$
and are $j$-holomorphic near~$\Theta$.
\end{itemize}
We shall abbreviate the family of closed Riemann surfaces determined by $j_\tau$ as
$$
\Sigma_\tau := (\Sigma,j_\tau),
$$
and denote by
$$
\pi : (\Sigma'',j'') \to (\Sigma,j), \qquad \Theta'' = \pi^{-1}(\Theta) \subset \Sigma''
$$
the holomorphic branched cover of closed surfaces provided by
Lemma~\ref{lemma:closedUp} such that $\dot{\Sigma}'' = \Sigma'' \setminus \Theta''$.
These choices produce a family of punctured Riemann surfaces
$$
\dot{\Sigma}_\tau := (\Sigma \setminus \Theta_\tau,j_\tau) \quad \text{ where } \quad
\Theta_\tau := \psi_\tau(\Theta) \subset \Sigma,
$$
and we define
$$
\varphi_\tau := \psi_\tau \circ \varphi : \Sigma' \to \Sigma, \qquad
j_\tau' := \varphi_\tau^*j_\tau \text{ on $\Sigma'$},
$$
where we observe that $j_\tau'$ is always well defined and matches $j'$
near $\Theta'$ since $\psi_\tau$ is holomorphic near~$\Theta$.  This makes
$$
\varphi_\tau : \Sigma'_\tau \to \Sigma_\tau
$$
a smooth family of holomorphic branched covers, where
$$
\Sigma'_\tau := (\Sigma',j'_\tau),
$$
and they restrict to
holomorphic covering maps of punctured surfaces
$\dot{\Sigma}'_\tau \stackrel{\varphi}{\longrightarrow} \dot{\Sigma}_\tau$, where
$$
\dot{\Sigma}'_\tau := (\dot{\Sigma}',j'_\tau).
$$

\begin{example}
\label{ex:Teichmueller}
Suppose $\Theta$ is the set of critical values of~$\varphi$, $r := |\Theta|$,
$P$ is the $2r$-dimensional open ball~$B^{2r}$, $j_\tau := j$ for all~$\tau$,
and $\psi_\tau : \Sigma \to \Sigma$ is chosen to be any smooth family of
diffeomorphisms supported near~$\Theta$ that are holomorphic in a smaller
neighborhood of $\Theta$ and such that $\psi_0 = \Id$ and
$$
B^{2r} \to \Sigma^{\times r} : \tau \mapsto (\psi_\tau(w_1),\ldots,\psi_\tau(w_r))
$$
is an embedding onto an open subset, where $\Theta = \{w_1,\ldots,w_r\}$.
Then the branched covers $\varphi_\tau : (\Sigma',j'_\tau) \to (\Sigma,j)$
parametrize a neighborhood of $\varphi$ in $\mM^d_{\mathbf{b}}(j)$.
\end{example}
\begin{example}
\label{ex:Teichmueller2}
If $v_0 : (\Sigma,j_0) \to (M,J_0)$ represents a simple element of the moduli
space $\mM_{g,m}(A,J_0\,;\,\ell_1,\ldots,\ell_m)$ defined in Appendix~\ref{sec:critical}
and $J_0$ is generic, then one can enhance the previous example as follows to parametrize
a neighborhood of $u_0 := v_0 \circ \varphi$ in the space
$\mM^d_{\mathbf{b}}(\mM_{g,m}(A,J_0\,;\,\ell_1,\ldots,\ell_m))$.
A neighborhood of $v_0$ in $\mM_{g,m}(A,J_0\,;\,\ell_1,\ldots,\ell_m)$ can be
identified with a smooth submanifold $X$ of $\dbar_{J_0}^{-1}(0)$, where
$\dbar_{J_0} : \tT \times \bB \to \eE$ is the nonlinear Cauchy-Riemann operator
defined on the product of $\bB := W^{k,p}(\Sigma,M)$ with a Teichm\"uller
slice $\tT$ through~$j_0$, cf.~Appendix~\ref{sec:critical}.  Here $\tT$ is 
a finite-dimensional smooth family of complex structures on $\Sigma$, which
can all be arranged to match $j_0$ near~$\Theta$.  A neighborhood in
$\mM^d_{\mathbf{b}}(\mM_{g,m}(A,J_0\,;\,\ell_1,\ldots,\ell_m))$ is now
parametrized by
$$
P := B^{2r} \times X,
$$
namely via the curves $v \circ (\psi_\sigma \circ \varphi) : (\Sigma',\varphi^*\psi_\sigma^*j)
\to (M,J_0)$ for each $\tau := (\sigma,(j,v)) \in P$, and we associate to
these parameters the families $j_\tau := j$ and 
$\psi_\tau := \psi_\sigma$.
\end{example}
\begin{example}
\label{ex:Teichmueller3}
Enhancing the previous example one step further, suppose $\jJ_\varepsilon$
is an infinite-dimensional Banach manifold consisting of smooth almost complex structures and we consider a
neighborhood of $(v_0,J_0)$ in the universal moduli space
$$
\univ^*(\jJ_\varepsilon\,;\,\ell_1,\ldots,\ell_m) = \left\{ (v,J) \ |\ 
J \in \jJ_\varepsilon,\ 
v \in \mM_{g,m}(A,J\,;\,\ell_1,\ldots,\ell_m) \right\}.
$$
Such a neighborhood can be identified with a finite-codimensional submanifold
$X$ in the infinite-dimensional Banach manifold $\dbar^{-1}(0) \subset
\tT \times \bB \times \jJ_\varepsilon$, where $\dbar(j,u,J) := \dbar_J(j,u)$.
Defining $P := B^{2r} \times X$ and the families $\{j_\tau\}$ and
$\{\psi_\tau\}$ as in Example~\ref{ex:Teichmueller2}, the parameter space
$P$ is now infinite dimensional.
\end{example}

Observe that the branched covers in the family $\varphi_\tau$ all have 
essentially the same topological properties, e.g.~their branch points and
automorphism groups are identical.  It is therefore trivial to extend
$(\Theta,\dot{\Sigma}'',\pi,G,\rho,I,f)$ to a smooth
family of regular presentations 
$$
(\Theta_\tau,\dot{\Sigma}'',\pi_\tau,G,\rho,I,f)
$$
for~$\varphi_\tau$, where $\pi_\tau := \psi_\tau \circ \pi$.
By the same reasoning as above, we can define on $\Sigma''$ a smooth family of
complex structures $j''_\tau := \pi_\tau^*j_\tau$ such that
$$
\pi_\tau : \Sigma''_\tau \to \Sigma_\tau, \qquad \Sigma''_\tau := (\Sigma'',j''_\tau)
$$
becomes a smooth family of holomorphic branched covers, restricting to a smooth family of holomorphic
covering maps $\dot{\Sigma}''_\tau \stackrel{\pi_\tau}{\longrightarrow} \dot{\Sigma}_\tau$,
defined on the family of punctured Riemann surfaces 
$$
\dot{\Sigma}''_\tau := (\dot{\Sigma}'',j''_\tau).
$$

\subsection{Cauchy-Riemann operators on closed and punctured domains}
\label{sec:CRpunctured}

Fix a complex vector bundle
$$
(E,J) \to (\Sigma,j)
$$
of rank $m \ge 1$, and define the rank~$m$ bundle of complex-antilinear maps
$$
F = \overline{\Hom}_\CC(T\Sigma,E) = \Lambda^{0,1}T^*\Sigma \otimes E.
$$
Recall that a first-order real-linear partial differential operator
$\mathbf{D} : \Gamma(E) \to \Gamma(F) = \Omega^{0,1}(\Sigma,E)$ is then called a
\defin{Cauchy-Riemann type operator} on $E$ if it
satisfies the Leibniz rule
$$
\mathbf{D} (f\eta) = (\dbar f) \eta + f \mathbf{D} \eta
$$
for all $\eta \in \Gamma(E)$ and $f \in C^\infty(\Sigma,\RR)$, where
$\dbar f = df + i\, df \circ j \in \Omega^{0,1}(\Sigma)$.  The space
$$
\CRR(E)
$$
of all such operators is an affine space modelled on the space of smooth
real-linear bundle maps $\Gamma(\Hom_\RR(E,F)) = 
\Omega^{0,1}(\Sigma,\End_\RR(E,J))$.  The \defin{pullback} of
$\mathbf{D} \in \CRR(E)$ via $\varphi : (\Sigma',j') \to (\Sigma,j)$ defines
a Cauchy-Riemann operator
$$
\varphi^*\mathbf{D} : \Gamma(E^\varphi) \to \Gamma(F^\varphi),
$$
where we define two bundles over $\Sigma'$ by
$$
E^\varphi := \varphi^*E, \qquad F^\varphi := \overline{\Hom}_\CC(T\Sigma',\varphi^*E)
= \Lambda^{0,1}T^*\Sigma' \otimes E^\varphi
$$
and characterize $\varphi^*\mathbf{D}$ via the relation
$$
(\varphi^*\mathbf{D})(\eta \circ \varphi) = \varphi^*\left( \mathbf{D}\eta \right)
\quad \text{ for all } \quad \eta \in \Gamma(E).
$$

\begin{example}
If $v : (\Sigma,j) \to (M,J)$ is a $J$-holomorphic curve with generalized
normal bundle $N_v \to \Sigma$, its normal Cauchy-Riemann operator
$\mathbf{D}_v^N$ belongs to $\CRR(N_v)$, and if $u = v \circ \varphi :
(\Sigma',j') \to (M,J)$, then $N_u = \varphi^*N_v$ and
$\mathbf{D}_u^N = \varphi^*\mathbf{D}_v^N \in \CRR(N_u)$.
\end{example}

\begin{remark}
\label{remark:dbar}
Note that the operator $\dbar : C^\infty(\Sigma,\CC) \to \Omega^{0,1}(\Sigma)$
used in our definition of Cauchy-Riemann type operators makes $\dbar f$
\emph{twice} the complex-antilinear part of the differential~$d f$.  
This is a common convention
in $J$-holomorphic curve theory, but differs from the standard convention in
complex analysis.  We will also often use the symbol $\dbar$ to mean
the coordinate-based differential operator
$$
\dbar := \p_s + i\p_t,
$$
acting on functions valued in a complex vector space and defined
on open domains in $\CC$ with complex coordinate $s + it$.  The
alternative convention would be to write $\dbar = \frac{1}{2} (\p_s + i\p_t)$.
\end{remark}

Fixing Hermitian bundle metrics 
$\langle\ ,\ \rangle_E$ and $\langle\ ,\ \rangle_\Sigma$ on $E$ and $T\Sigma$
respectively, we can integrate real parts of bundle metrics
to define real-valued $L^2$-pairings
$\langle\ ,\ \rangle_{L^2}$ on $\Gamma(E)$ and $\Gamma(F)$, which determines 
a \defin{formal adjoint} operator $\mathbf{D}^* : \Gamma(F) \to \Gamma(E)$ 
via the relation
$$
\langle \alpha ,\mathbf{D} \eta \rangle_{L^2} = 
\langle \mathbf{D}^*\alpha , \eta \rangle_{L^2}
$$
for all smooth sections $\alpha \in \Gamma(F)$ and $\eta \in \Gamma(E)$
with compact support.\footnote{The compact support condition is vacuous
in the present context since $\Sigma$ is compact, but the same definition
is also valid on punctured domains.}
Viewing $\mathbf{D}$ as a Fredholm operator on Sobolev spaces
$W^{k,p}(E) \to W^{k-1,p}(F)$ for some $k \in \NN$ and $p \in (1,\infty)$,
we can then identify $\coker \mathbf{D}$ with $\ker \mathbf{D}^* \subset
\Gamma(F)$, which is the $L^2$-orthogonal complement of
$\im \mathbf{D} \subset W^{k-1,p}(F)$ and is a finite-dimensional
space of smooth sections by elliptic regularity.  
Using the Riemann-Roch formula $\ind(\mathbf{D}) = m\chi(\Sigma) + 2 c_1(E)$
and computing the algebraic count of branch points $Z(d\varphi)$ from
the Riemann-Hurwitz formula, the (real) Fredholm indices of 
$\mathbf{D}$ and $\varphi^*\mathbf{D}$ are related by
$$
\ind (\varphi^*\mathbf{D}) = d \cdot \ind \mathbf{D} - m Z(d\varphi).
$$

In order to exploit the topological constructions in the previous section,
we will need to work with Cauchy-Riemann type operators on punctured surfaces
instead of closed surfaces.  We shall now show that this can be done without
loss of generality by choosing suitable weighted Sobolev spaces.  Assume 
$$
E_\tau \to \Sigma_\tau
$$
is a smooth family of rank~$m$ complex vector bundles with complex structures~$J_\tau$,
equipped with a smooth family
of Cauchy-Riemann operators $\mathbf{D}_\tau \in \CRR(E_\tau)$.  Denote
the restrictions of the bundles $E_\tau$ and 
$$
F_\tau := \overline{\Hom}_\CC(T\Sigma_\tau,E_\tau)
$$
to the punctured surfaces $\dot{\Sigma}_\tau$ by
$$
\dot{E}_\tau := E_\tau|_{\dot{\Sigma}_\tau}, \qquad
\dot{F}_\tau := F_\tau|_{\dot{\Sigma}_\tau} = \overline{\Hom}_\CC(T\dot{\Sigma}_\tau,\dot{E}_\tau).
$$
Restricting $\mathbf{D}_\tau$ to $\dot{\Sigma}_\tau$ then defines a family
of Cauchy-Riemann type operators
$$
\dot{\mathbf{D}}_\tau \in \CRR(\dot{E}_\tau).
$$
In order to understand the functional-analytic properties of
$\dot{\mathbf{D}}_\tau$, we must examine its asymptotic behavior fairly 
carefully.  Fix local holomorphic coordinate charts to identify a neighborhood of
each $w \in \Theta$ in $\Sigma$ with the closed unit disk $\DD \subset \CC$,
with $w$ corresponding to $0 \in \DD$, and use the maps 
$\psi_\tau$ introduced at the end of \S\ref{sec:regPres} to produce from
these a smooth family of holomorphic charts on neighborhoods of 
$\psi_\tau(w) \in \Theta_\tau$ for $\tau \in P$.
In these coordinates, use the biholomorphic map
$$
[0,\infty) \times S^1 \to \DD \setminus \{0\} : (s,t) \mapsto
e^{-2\pi(s+it)}
$$
to define cylindrical ends of $\dot{\Sigma}_\tau$ with holomorphic
coordinates $(s,t) \in [0,\infty) \times S^1$.
Choose also a smooth family of trivializations
of $E_\tau$ near $\Theta_\tau$ and denote the resulting trivialization 
of $\dot{E}_\tau$ over the cylindrical ends by~$\Phi$.  The relative
first Chern number\footnote{Recall that for any complex line bundle $E$
over a surface $\Sigma$ with a trivialization $\Phi$ specified outside of
some open subset in $\Sigma$ with compact closure, the relative first
Chern number $c_1^\Phi(E) \in \ZZ$ is defined by algebraically counting
the zeroes of a generic section that is constant with respect to $\Phi$
wherever the latter is defined.  This definition extends uniquely to
higher rank bundles via the relation $c_1^{\Phi_1 \oplus \Phi_2}(E_1 \oplus E_2)
= c_1^{\Phi_1}(E_1) + c_1^{\Phi_2}(E_2)$.}
of $\dot{E}_\tau$ is then given by
\begin{equation}
\label{eqn:c1}
c_1^{\Phi}(\dot{E}_\tau) = c_1(E_\tau) \in \ZZ.
\end{equation}
For any tuple of real numbers 
$$
\boldsymbol{\delta} = \{ \delta_w \in \RR \}_{w \in \Theta},
$$
we can use the chosen coordinates and trivializations
over the cylindrical ends of $\dot{\Sigma}_\tau$ to define the Sobolev
space with \defin{exponential weights}
$$
W^{k,p,\boldsymbol{\delta}}(\dot{E}_\tau) := \Big\{ \eta \in 
W^{k,p}_{\operatorname{loc}}(\dot{E}_\tau)\ \Big|
\ \text{$e^{\delta_w s} \eta \in W^{k,p}([0,\infty) \times S^1)$
on the end near $\psi_\tau(w) \in \Theta_\tau$} \Big\}.
$$
We will also write
$$
L^{p,\boldsymbol{\delta}}(\dot{E}_\tau) := W^{0,p,\boldsymbol{\delta}}(\dot{E}_\tau).
$$
Note that sections $\eta \in W^{k,p,\boldsymbol{\delta}}(\dot{E}_\tau)$ have 
exponential decay
at any end where $\delta_w > 0$, but one can also take $\delta_w < 0$, in which
case $\eta$ may be unbounded with exponential \emph{growth} near~$w$.
In order to emphasize when we are using negative exponential weights, we
associate to $\boldsymbol{\delta} = \{\delta_w\}_{w \in \Theta}$ the
inverse set of weights
$$
-\boldsymbol{\delta} := \{ -\delta_w \}_{w \in \Theta}.
$$
The asymptotic coordinates and trivializations also naturally give rise to
asymptotic trivializations of 
$\dot{F}_\tau = \overline{\Hom}_\CC(T\dot{\Sigma}_\tau,\dot{E}_\tau)$,
so we can similarly define the Banach space
$W^{k-1,p,\boldsymbol{\delta}}(\dot{F}_\tau)$,
which is a completion of some subset of 
$\Omega^{0,1}(\dot{\Sigma}_\tau,\dot{E}_\tau)$
determined by the asymptotic conditions.

Choose a smooth $\tau$-parametrized family of Hermitian bundle metrics and 
connections on $E_\tau$ which match the trivial
metric and connection in our chosen family of trivializations 
near~$\Theta_\tau$.  Any 
Cauchy-Riemann type operator on $E_\tau$ can then be written as
$\mathbf{D}_\tau = \dbar_\nabla + A$ for some
$A \in \Omega^{0,1}(\Sigma_\tau,\End_\RR(E_\tau))$, where $\dbar_\nabla :=
\nabla + J_\tau \circ \nabla \circ j_\tau : \Gamma(E_\tau) \to
\Omega^{0,1}(\Sigma_\tau,E_\tau)$.  In the chosen coordinates
and trivialization near a point $w \in \Theta_\tau$, the $(0,1)$-form 
$A$ can be written as
$$
A = A^{(w)}_\tau(z)\, d\bar{z}
$$
for some smooth function $A^{(w)}_\tau : \DD \to \End_\RR(\CC^m)$.  The restriction of $A$
to an $\End_\RR(\dot{E}_\tau)$-valued $(0,1)$-form $\dot{A}_\tau \in 
\Omega^{0,1}(\dot{\Sigma}_\tau,\End_\RR(\dot{E}_\tau))$ can then be written 
on the corresponding cylindrical end as
$$
\dot{A}_\tau = \dot{A}^{(w)}_\tau(s,t)\, (-ds + i\, dt)
$$
where
\begin{equation}
\label{eqn:dotAw}
\dot{A}^{(w)}_\tau(s,t) := 2\pi e^{-2\pi (s - it)} A^{(w)}_\tau\big( e^{-2\pi (s+it)}\big),
\end{equation}
and given a section $\eta \in \Gamma(\dot{E}_\tau)$ expressed as a function
$\eta(s,t) \in \CC^m$ with respect to the trivialization on the same end, 
$\dot{\mathbf{D}}_\tau \eta$ on this end takes the form
\begin{equation}
\label{eqn:asymptoticCR}
\dot{\mathbf{D}}_\tau \eta = 
\left( \p_s \eta + i \p_t \eta + \dot{A}^{(w)}_\tau \eta \right) (-ds + i\, dt) =:
\left( \dbar \eta + \dot{A}^{(w)}_\tau \eta \right) (-ds + i\, dt).
\end{equation}
(Here and in further local expressions below, we are using the abbreviation $\dbar := \p_s + i\p_t$ as mentioned
in Remark~\ref{remark:dbar}.)
Observe that $\dot{A}^{(w)}_\tau(s,\cdot) \to 0$ with all derivatives as $s \to \infty$.  
This expression shows that $\dot{\mathbf{D}}_\tau$ extends to a bounded linear
operator
$$
\dot{\mathbf{D}}_\tau : W^{k,p,\boldsymbol{\delta}}(\dot{E}_\tau) \to
W^{k-1,p,\boldsymbol{\delta}}(\dot{F}_\tau)
$$
for any choices of $k \in \NN$, $p \in (1,\infty)$ and
exponential weights~$\boldsymbol{\delta} = \{\delta_w \in \RR\}_{w \in \Theta}$.
Operators of this type are standard in Floer-type theories, and especially in
symplectic field theory.  Appealing to the Fredholm theory on punctured
surfaces developed in \cite{Schwarz}, the asymptotic decay of 
$\dot{A}^{(w)}_\tau(s,\cdot)$ means that
$\dot{\mathbf{D}}_\tau : W^{k,p}(\dot{E}_\tau) \to W^{k-1,p}(\dot{F}_\tau)$ 
is controlled at every puncture by the so-called \emph{trivial}
asymptotic operator~$-i\p_t : H^1(S^1,\CC^m) \to L^2(S^1,\CC^m)$, for which $0$
is an eigenvalue of maximal
multiplicity.  In this sense, the asymptotics are degenerate, i.e.~in the
SFT setting, such an operator can arise as the linearized Cauchy-Riemann operator
of a holomorphic curve asymptotic
to periodic orbits that live in Morse-Bott families foliating an open set.
In particular, $\dot{\mathbf{D}}_\tau : W^{k,p} \to W^{k-1,p}$ 
is \emph{not} Fredholm, but it 
becomes Fredholm when we introduce suitable weights:
conjugating $\dot{\mathbf{D}}_\tau : W^{k,p,\boldsymbol{\delta}} \to
W^{k-1,p,\boldsymbol{\delta}}$ with a map of the form $\Psi(\eta) = e^f \eta$
for a suitable function $f : \dot{\Sigma}_\tau \to \RR$ 
(cf.~\cite{HWZ:props3}*{\S 6} or \cite{Wendl:automatic}*{\S 2.1})
produces a commutative diagram
\begin{equation}
\label{eqn:weights}
\begin{CD}
W^{k,p,\boldsymbol{\delta}}(\dot{E}_\tau)
@>\dot{\mathbf{D}}_\tau>>
W^{k-1,p,\boldsymbol{\delta}}(\dot{F}_\tau) \\
@VV\Psi V     @VV\Psi V \\
W^{k,p}(\dot{E}_\tau) 
@>\widehat{\mathbf{D}}_\tau>>
W^{k-1,p}(\dot{F}_\tau), \\
\end{CD}
\end{equation}
where $\widehat{\mathbf{D}} : W^{k,p} \to W^{k-1,p}$ is another Cauchy-Riemann
type operator whose asymptotic operators are offset by constants
depending on the weights~$\boldsymbol{\delta}$, and thus is Fredholm for suitable 
choices.
In particular, the computation in \eqref{eqn:conjugate} and~\eqref{eqn:Adelta}
below will show that imposing the exponential growth condition 
$e^{-\delta s} \eta \in
W^{k,p}([0,\infty) \times S^1)$ on each cylindrical end for sufficiently small
$\delta > 0$ adjusts the asymptotic operators of $\widehat{\mathbf{D}}_\tau$ 
so that each acquires an effective Conley-Zehnder index $m$ relative to 
the trivialization~$\Phi$.

We need to be a bit cautious with the weights when discussing elliptic
regularity and formal adjoints: as a rule, the Sobolev constants $k \in \NN$ 
and $p \in (1,\infty)$ can be changed freely, but the weights cannot.  
The following are immediate consequences
of \eqref{eqn:weights} after applying standard regularity arguments
to~$\widehat{\mathbf{D}}_\tau$, plus (in the case of Lemma~\ref{lemma:regularity2}) 
the fact that Cauchy-Riemann operators with
nondegenerate asymptotics automatically impose exponential decay conditions
on their kernels (cf.~\cite{Schwarz}*{Prop.~3.1.26}):

\begin{lemma}
\label{lemma:regularity}
Suppose $k \in \NN$, $1 < p < \infty$, and $\boldsymbol{\delta} =
\{ \delta_w \in \RR\}_{w \in \Theta}$ is any choice of exponential
weights.  If $\eta \in  L^{p,\boldsymbol{\delta}}(\dot{E}_\tau)$ is a weak solution to
$\dot{\mathbf{D}}_\tau \eta = \xi$ for $\xi \in W^{k-1,p,\boldsymbol{\delta}}(\dot{F}_\tau)$,
then $\eta \in W^{k,p,\boldsymbol{\delta}}(\dot{E}_\tau)$.
\qed
\end{lemma}
\begin{lemma}
\label{lemma:regularity2}
Suppose $1 < p < \infty$ and the weights $\boldsymbol{\delta}$ are chosen
such that $\dot{\mathbf{D}}_\tau : W^{k,p,\boldsymbol{\delta}}(\dot{E}_\tau) \to
W^{k-1,p,\boldsymbol{\delta}}(\dot{F}_\tau)$ is Fredholm.  If 
$\eta \in L^{p,\boldsymbol{\delta}}(\dot{E}_\tau)$ is a weak solution to
$\dot{\mathbf{D}}_\tau \eta = 0$, then
$\eta \in W^{k,q,\boldsymbol{\delta}}(\dot{E}_\tau)$ for all $k \in \NN$
and $q \in (1,\infty)$.
\qed
\end{lemma}

To discuss the formal adjoint on punctured domains, 
one should define real $L^2$-products for
$\Gamma(\dot{E}_\tau)$ and $\Gamma(\dot{F}_\tau)$ in terms of a family of Hermitian
bundle metrics on $E_\tau$ and Riemannian metrics on $\dot{\Sigma}_\tau$ that 
are compatible with the conformal structure and standard on the cylindrical 
ends; in particular, the right metric to use on the cylindrical ends
is the Euclidean metric in the coordinates $(s,t) \in [0,\infty) \times S^1$,
so that ends have infinite area and the metric does not extend over the punctues.
The key technical point is then the
following: there are well-defined $L^2$-pairings
\begin{equation}
\label{eqn:L2pairing}
L^{p,\boldsymbol{\delta}} \otimes L^{q,-\boldsymbol{\delta}} \to \RR :
\eta \otimes \xi \mapsto \langle \eta,\xi \rangle_{L^2}
\end{equation}
whenever $1/p + 1/q = 1$, and using the density of $C_0^\infty$, the
usual relation
\begin{equation}
\label{eqn:usualRelation}
\langle \alpha , \dot{\mathbf{D}}_\tau \eta \rangle_{L^2} =
\langle \dot{\mathbf{D}}_\tau^*\alpha , \eta \rangle_{L^2}
\end{equation}
for smooth compactly supported sections $\eta$ and $\alpha$ remains valid whenever
$\eta \in W^{1,p,-\boldsymbol{\delta}}(\dot{E}_\tau)$ and
$\alpha \in W^{1,q,\boldsymbol{\delta}}(\dot{F}_\tau)$ for $1/p + 1/q = 1$.  
Using \eqref{eqn:weights}, one finds $\dot{\mathbf{D}}_\tau^* =
\Psi \widehat{\mathbf{D}}_\tau^* \Psi^{-1}$, from which one can check that
$\dot{\mathbf{D}}_\tau^* : W^{k,p,\boldsymbol{\delta}}(\dot{F}_\tau) \to
W^{k-1,p,\boldsymbol{\delta}}(\dot{E}_\tau)$ satisfies the Fredholm property
and Lemmas~\ref{lemma:regularity} and~\ref{lemma:regularity2}
under the same conditions on $\boldsymbol{\delta}$ as
$\dot{\mathbf{D}}_\tau : W^{k,p,-\boldsymbol{\delta}}(\dot{E}_\tau) \to
W^{k-1,p,-\boldsymbol{\delta}}(\dot{F}_\tau)$.
The next result appears standard at first glance, but the reader should be 
cautioned that it depends on inclusions $W^{k,p,\boldsymbol{\delta}} \hookrightarrow
W^{k,p,-\boldsymbol{\delta}}$ which hold only when all the weights are 
nonnegative, so e.g.~one does not obtain any similar result with the roles of
$\dot{\mathbf{D}}_\tau$ and $\dot{\mathbf{D}}_\tau^*$ reversed.

\begin{prop}
\label{prop:formalAdjoint}
Assume $k \in \NN$, $1 < p < \infty$, and $\boldsymbol{\delta} =
\{ \delta_w \ge 0 \}_{w \in \Theta}$ is a set of nonnegative exponential
weights such that
$$
\dot{\mathbf{D}}_\tau : W^{k,p,-\boldsymbol{\delta}}(\dot{E}_\tau) \to
W^{k-1,p,-\boldsymbol{\delta}}(\dot{F}_\tau)
$$
is Fredholm.  Defining its formal adjoint as a bounded linear map
$$
\dot{\mathbf{D}}_\tau^* : W^{k,p,\boldsymbol{\delta}}(\dot{F}_\tau) \to
W^{k-1,p,\boldsymbol{\delta}}(\dot{E}_\tau)
$$
and using the obvious inclusions $W^{k,p,\boldsymbol{\delta}}(\dot{F}_\tau)
\hookrightarrow W^{k-1,p,\boldsymbol{\delta}}(\dot{F}_\tau)
\hookrightarrow W^{k-1,p,-\boldsymbol{\delta}}(\dot{F}_\tau)$, we have
$$
W^{k-1,p,-\boldsymbol{\delta}}(\dot{F}_\tau) = \im \dot{\mathbf{D}}_\tau
\oplus \ker \dot{\mathbf{D}}_\tau^*.
$$
In particular, $\coker \dot{\mathbf{D}}_\tau$ is isomorphic to the space
of all sections in $L^{q,\boldsymbol{\delta}}(\dot{F}_\tau)$ for $1/p + 1/q = 1$
that are $L^2$-orthogonal to $\im \dot{\mathbf{D}}_\tau \subset L^{p,-\boldsymbol{\delta}}(\dot{F}_\tau)$
under the pairing \eqref{eqn:L2pairing}.
\end{prop}
\begin{proof}
If $\alpha \in \im \dot{\mathbf{D}}_\tau \cap \ker \dot{\mathbf{D}}_\tau^*$,
then $\alpha = \dot{\mathbf{D}}_\tau\eta$ for some $\eta \in W^{k,p,-\boldsymbol{\delta}}(\dot{E}_\tau) \subset
W^{1,p,-\boldsymbol{\delta}}(\dot{E}_\tau)$, while $\alpha$ also belongs to
$W^{1,q,\boldsymbol{\delta}}(\dot{F}_\tau)$ for $1/p + 1/q = 1$ by
Lemma~\ref{lemma:regularity2}.  Thus $\alpha$ has a well-defined $L^2$-pairing
with itself and \eqref{eqn:usualRelation} gives
$$
\| \alpha \|_{L^2}^2 = \langle \alpha , \dot{\mathbf{D}}_\tau \eta \rangle_{L^2} =
\langle \dot{\mathbf{D}}_\tau^*\alpha , \eta \rangle_{L^2} = 0.
$$
To show that $\im \dot{\mathbf{D}}_\tau + \ker \dot{\mathbf{D}}_\tau^*$
is $W^{k-1,p,-\boldsymbol{\delta}}(\dot{F}_\tau)$, note first that it is a 
closed subspace since $\dot{\mathbf{D}}_\tau$ is Fredholm.  Then in the case
$k=1$, the contrary would mean there exists a nontrivial 
$\lambda \in (L^{p,-\boldsymbol{\delta}}(\dot{F}_\tau))^* =
L^{q,\boldsymbol{\delta}}(\dot{F}_\tau)$ for $1/p + 1/q = 1$ such that
$\langle \dot{\mathbf{D}}_\tau \eta , \lambda \rangle_{L^2} = 0$ for all
$\eta \in W^{1,p,-\boldsymbol{\delta}}(\dot{E}_\tau)$ and
$\langle \alpha , \lambda \rangle_{L^2} = 0$ for all $\alpha \in 
\ker\dot{\mathbf{D}}_\tau^*$.  The first condition means 
$\lambda \in \ker\dot{\mathbf{D}}_\tau^*$ by Lemma~\ref{lemma:regularity2}
and thus contradicts the second unless $\lambda = 0$.  To extend this result
to all $k \in \NN$, note that if $\lambda \in 
W^{k-1,p,-\boldsymbol{\delta}}(\dot{F}_\tau) \subset L^{p,-\boldsymbol{\delta}}(\dot{F}_\tau)$
then the $k=1$ case gives $\eta \in W^{1,p,-\boldsymbol{\delta}}(\dot{E}_\tau)$
and $\alpha \in \ker \dot{\mathbf{D}}_\tau^*$ such that
$\dot{\mathbf{D}}_\tau \eta + \alpha = \lambda$.  Then Lemma~\ref{lemma:regularity2}
implies $\alpha \in W^{k-1,p,\boldsymbol{\delta}}(\dot{F}_\tau) \subset
W^{k-1,p,-\boldsymbol{\delta}}(\dot{F}_\tau)$, implying that
$\dot{\mathbf{D}}_\tau \eta$ is also in $W^{k-1,p,-\boldsymbol{\delta}}(\dot{F}_\tau)$,
so Lemma~\ref{lemma:regularity} implies
$\eta \in W^{k,p,-\boldsymbol{\delta}}(\dot{E}_\tau)$ and we are done.
\end{proof}

This discussion extends easily to the pulled back operators
$$
\varphi_\tau^*\mathbf{D}_\tau \in \CRR(\varphi_\tau^*E_\tau) \quad \text{ and } \quad
\varphi_\tau^*\dot{\mathbf{D}}_\tau \in \CRR(\varphi_\tau^*\dot{E}_\tau)
$$
on bundles over $\Sigma'_\tau$ and $\dot{\Sigma}'_\tau$ respectively.
Observe that since $\dot{\Sigma}'_\tau \stackrel{\varphi_\tau}{\longrightarrow}
\dot{\Sigma}_\tau$ has no branch points, $d\varphi_\tau$ gives a bundle isomorphism
$T\dot{\Sigma}'_\tau \to \varphi_\tau^*T\dot{\Sigma}_\tau$ and we can thus
identify
$$
F_\tau^{\varphi_\tau}|_{\dot{\Sigma}'_\tau} = 
\overline{\Hom}_\CC(T\dot{\Sigma}'_\tau,\varphi_\tau^*\dot{E}_\tau) =
\overline{\Hom}_\CC(\varphi_\tau^*T\dot{\Sigma}_\tau,\varphi_\tau^*\dot{E}_\tau) =
\varphi_\tau^*\dot{F}_\tau,
$$
so that $\varphi_\tau^*\dot{\mathbf{D}}_\tau $ can be viewed as a map
$\Gamma(\varphi_\tau^*\dot{E}_\tau) \to \Gamma(\varphi_\tau^*\dot{F}_\tau)$.
We can now define fixed holomorphic cylindrical coordinate systems
$(s,t) \in [0,\infty) \times S^1$ on punctured neighborhoods of each point
$\zeta \in \Theta' = \varphi_\tau^{-1}(\Theta_\tau)$ such that $\varphi_\tau$
takes the form
\begin{equation*}
\begin{split}
\dot{\Sigma}_\tau' \supset [0,\infty) \times S^1 &\stackrel{\varphi_\tau}{\longrightarrow}
[0,\infty) \times S^1 \subset \dot{\Sigma}_\tau,\\  
(s,t) &\mapsto (k_\zeta s, k_\zeta t),
\end{split}
\end{equation*}
where $k_\zeta \in \NN$ is the branching order of $\varphi$ at~$\zeta$.  
Pulling back the trivializations $\Phi$ on $E_\tau$ 
near $\Theta_\tau$ to define corresponding trivializations of 
$\varphi_\tau^*E_\tau$ near~$\Theta'$, we obtain asymptotic trivializations of
$\varphi_\tau^*\dot{E}_\tau$ and $\varphi_\tau^*\dot{F}_\tau$ on the 
cylindrical ends and can thus define weighted Sobolev norms for sections of
these bundles, producing a bounded linear operator
$$
\varphi_\tau^*\dot{\mathbf{D}}_\tau : 
W^{k,p,\boldsymbol{\delta}}(\varphi_\tau^*\dot{E}_\tau) \to 
W^{k-1,p,\boldsymbol{\delta}}(\varphi_\tau^*\dot{F}_\tau)
$$
for all choices of $k \in \NN$, $p \in (1,\infty)$ and exponential weights
$\boldsymbol{\delta} = \{ \delta_\zeta \in \RR \}_{\zeta \in \Theta'}$.
If $\boldsymbol{\delta} = \{ \delta_w \}_{w \in \Theta}$ is a choice of
weights for $\dot{\mathbf{D}}_\tau$, there is an induced set of weights for
$\varphi_\tau^*\dot{\mathbf{D}}_\tau$ defined by
$$
\varphi^*\boldsymbol{\delta} := \left\{ k_\zeta \delta_{\varphi(\zeta)} \right\}_{\zeta \in \Theta'},
$$
where $k_\zeta \in \{1,\ldots,d\}$ again denotes the branching order of 
$\varphi$ at~$\zeta$.  

\begin{prop}
\label{prop:punctured}
Suppose $k \in \NN$, $p \in (1,\infty)$, and the exponential weights
$\boldsymbol{\delta} = \{ \delta_w \}_{w \in \Theta}$ are chosen to satisfy
$$
0 < \delta_w < \frac{2\pi}{d}
$$
for every $w \in \Theta$.  Then for any
$\mathbf{D}_\tau \in \CRR(E_\tau)$, the operators
\begin{equation*}
\begin{split}
\dot{\mathbf{D}}_\tau : W^{k,p,-\boldsymbol{\delta}}(\dot{E}_\tau) &\to
W^{k-1,p,-\boldsymbol{\delta}}(\dot{F}_\tau), \\
\varphi_\tau^*\dot{\mathbf{D}}_\tau : 
W^{k,p,-\varphi^*\boldsymbol{\delta}}(\varphi_\tau^*\dot{E}_\tau) &\to 
W^{k-1,p,-\varphi^*\boldsymbol{\delta}}(\varphi_\tau^*\dot{F}_\tau)
\end{split}
\end{equation*}
are Fredholm and satisfy
$$
\ind(\dot{\mathbf{D}}_\tau) = \ind(\mathbf{D}_\tau), \quad \text{ and }\quad
\ind(\varphi_\tau^*\dot{\mathbf{D}}_\tau) =
\ind(\varphi_\tau^*\mathbf{D}_\tau).
$$
Moreover, the maps 
$\Gamma(E_\tau) \to \Gamma(\dot{E}_\tau)$ and
$\Gamma(\varphi_\tau^*E_\tau) \to \Gamma(\varphi_\tau^*\dot{E}_\tau)$ defined
by restricting smooth sections to the corresponding punctured domains
define isomorphisms
$$
\ker \mathbf{D}_\tau \stackrel{\cong}{\longrightarrow} \ker \dot{\mathbf{D}}_\tau
\quad\text{ and }\quad
\ker(\varphi_\tau^*\mathbf{D}_\tau) \stackrel{\cong}{\longrightarrow}
\ker(\varphi_\tau^*\dot{\mathbf{D}}_\tau).
$$
\end{prop}
\begin{proof}
We will prove the correspondence between $\mathbf{D}_\tau$ and $\dot{\mathbf{D}}_\tau$,
as the result for the pulled back operators follows by the same argument
simply replacing the bundles $E_\tau \to \Sigma$ and $\dot{E}_\tau \to \dot{\Sigma}_\tau$
with $\varphi_\tau^*E_\tau \to \Sigma'$ and $\varphi_\tau\dot{E}_\tau \to
\dot{\Sigma}'_\tau$ respectively.

The Fredholm property for $\dot{\mathbf{D}}_\tau$ and the index calculation
follow from  the usual index formula for Cauchy-Riemann operators on
Riemann surfaces with cylindrical ends, proved in \cite{Schwarz}
(see also \cite{Wendl:SFT}*{Lecture~5}),
supplemented by the transformation \eqref{eqn:weights} to handle the
exponential weights (cf.~\cite{HWZ:props3}*{\S 6}).
In particular, the condition
$-2\pi < -\delta_w < 0$ for each $w \in \Theta_\tau$ guarantees that
$\dot{\mathbf{D}}_\tau$ is conjugate (cf.~\eqref{eqn:conjugate} and~\eqref{eqn:Adelta} below)
to a Cauchy-Riemann type operator
$W^{k,p}(\dot{E}_\tau) \to W^{k-1,p}(\dot{F}_\tau)$
with nondegenerate asymptotic operators at every puncture whose Conley-Zehnder
indices with respect to the trivialization $\Phi$ are $m = \rank_\CC E_\tau$.
In light of \eqref{eqn:c1}, the index formula from \cite{Schwarz} thus gives
$$
\ind(\dot{\mathbf{D}}_\tau) = m\chi(\dot{\Sigma}_\tau) + 2 c_1^{\Phi}(\dot{E}_\tau) + 
m \cdot |\Theta_\tau| = m\chi(\Sigma) + 2 c_1(E_\tau) = \ind(\mathbf{D}_\tau).
$$
Note that doing the same computation for the pulled back operators requires
the stronger condition $-2\pi/d < -\delta_w < 0$ in order to ensure that
all of the pulled back weights in the set $-\varphi^*\boldsymbol{\delta}$
lie in the interval $(-2\pi,0)$.

To understand the kernels, observe that
since any $\eta \in \ker \mathbf{D}_\tau$ is smooth, its restriction
to $\dot{\Sigma}_\tau$ belongs to $W^{k,p,-\boldsymbol{\delta}}(\dot{E}_\tau)$
and is thus in $\ker \dot{\mathbf{D}}_\tau$.\footnote{Note 
that $\eta|_{\dot{\Sigma}_\tau}$
would not belong to $W^{k,p,-\boldsymbol{\delta}}(\dot{E}_\tau)$ in general if
$\eta$ were an arbitrary (not necessarily smooth) 
section of class $W^{k,p}$ on~$E_\tau$, nor if any of the
exponential weights were nonnegative---the latter in particular permits
sections in $W^{k,p,-\boldsymbol{\delta}}(\dot{E}_\tau)$ that do not decay to
zero at infinity, which is crucial since arbitrary smooth sections
$\eta \in \ker \mathbf{D}_\tau$ may indeed be nonzero at points in~$\Theta_\tau$.}
Conversely, we need to show that any section
$\eta \in W^{k,p,-\boldsymbol{\delta}}(\dot{E}_\tau)$ annihilated by
$\dot{\mathbf{D}}_\tau$ can be extended over the punctures to a section in
$W^{k,p}(E_\tau)$, which is then automatically annihilated by
$\mathbf{D}_\tau$.  This will follow from the asymptotic elliptic theory of
the equation $\dot{\mathbf{D}}_\tau \eta = 0$.  Indeed,
recall from \eqref{eqn:asymptoticCR} that on the cylindrical end
near any puncture $w \in \Theta_\tau$, the function
$\eta(s,t) \in \CC^m$ representing $\eta \in \ker \dot{\mathbf{D}}_\tau$
in some trivialization satisfies
$$
\dbar \eta + \dot{A}^{(w)}_\tau \eta \equiv 0,
$$
and
$$
\eta = e^{\delta s} f \quad \text{ for some }
f \in W^{k,p}([0,\infty) \times S^1, \CC^m),
$$
where $\delta := \delta_w \in (0,2\pi)$.
Then $f = e^{-\delta s} \eta$ satisfies the Cauchy-Riemann type equation
\begin{equation}
\label{eqn:conjugate}
\dbar f + (\delta + \dot{A}^{(w)}_\tau) f = \p_s f - 
[ - i \p_t - (\delta + \dot{A}^{(w)}_\tau)] f 
= 0.
\end{equation}
Since $\dot{A}^{(w)}_\tau(s,\cdot) \to 0$ as $s \to \infty$, this equation is
asymptotic to the equation $(\p_s - \mathbf{A}_\delta) f = 0$ for the
asymptotic operator
\begin{equation}
\label{eqn:Adelta}
\mathbf{A}_\delta := - i \p_t - \delta : H^1(S^1,\CC^m) \to L^2(S^1,\CC^m),
\end{equation}
which can be regarded as a densely defined unbounded self-adjoint operator on
$L^2(S^1,\CC^m)$.  The function $A^{(w)}_\tau : \DD \to \End_\RR(\CC^m)$ is
smooth by assumption, and \eqref{eqn:dotAw} then implies that the derivatives
$\p^\alpha\dot{A}^{(w)}_\tau(s,t)$ of $\dot{A}^{(w)}_\tau$ for arbitrary multi-indices $\alpha$
satisfy exponential decay conditions
$$
| \p^\alpha \dot{A}^{(w)}_\tau(s,t) | \le M_\alpha e^{-2\pi s}
$$
for suitable constants $M_\alpha > 0$.  Applying
\cite{Siefring:asymptotics}*{Theorem~A.1}, $f$ therefore satisfies
$$
f(s,t) = e^{\lambda s} \left[ e(t) + r(s,t) \right],
$$
where $e : S^1 \to \CC^m$ is a nontrivial eigenfunction of $\mathbf{A}_\delta$
with eigenvalue $\lambda < 0$, and the remainder $r(s,t) \in \CC^m$ 
decays to zero with all its derivatives uniformly in~$t$ as $s \to \infty$.
The spectrum
of $\mathbf{A}_\delta$ is $\{ 2\pi k - \delta\ |\ k \in \ZZ \} \subset \RR$,
hence the assumption $\delta \in (0,2\pi)$ implies $\lambda \le - \delta$, 
and we conclude that
$$
\eta(s,t) = e^{(\delta + \lambda)s} \left[ e(t) + r(s,t) \right]
$$
is bounded on the cylindrical end; in fact, one can use this to show
that the smooth
function $\DD \setminus \{0\} \to \CC^m : z \mapsto \eta(z)$ defined via
the transformation
$z = e^{-2\pi (s+it)}$ has finite $W^{1,p}$-norm on $\DD \setminus \{0\}$.
Moreover, $\eta(z)$ has a continuous extension to $z=0$: indeed, the
extension is obviously $\eta(0)=0$ if $\lambda < - \delta$, while in
the case $\lambda = -\delta$, the eigenfunction $e(t)$ is necessarily
constant, so that $\eta(s,\cdot)$ converges to this constant value
as $s \to \infty$.  All these conditions together imply that the
continuous extension of $\eta$ over the punctures is of class $W^{k,p}$,
e.g.~the case $k=1$ is a standard exercise using the definition of weak 
derivatives (cf.~\cite{Wendl:lecturesV2}*{Exercise~2.118}), and the general
case follows from this by elliptic regularity.
\end{proof}

\begin{remark}
Since sections in $W^{k,p,-\boldsymbol{\delta}}(\dot{E}_\tau)$ and its
pulled back counterpart need not be bounded when the weights
$-\boldsymbol{\delta}$ are negative,
the punctured operators in Proposition~\ref{prop:punctured} cannot
be interpreted in any reasonable way as linearizations of nonlinear
Cauchy-Riemann operators, e.g.~$W^{k,p,-\boldsymbol{\delta}}(\dot{E}_\tau)$ in 
this case is not a subspace of
a tangent space in any reasonable Banach manifold.  For our purposes,
the exponential growth condition is merely a technical convenience so
that we can consider operators with the right index and the right kernel
and cokernel while
dealing with honest covering maps instead of branched covers.
The geometrically meaningful operators are still $\mathbf{D}_\tau$ and
$\varphi_\tau^*\mathbf{D}_\tau$, on unpunctured domains.
\end{remark}

\begin{remark}
Suppose $E_\tau$, $\Sigma_\tau$ and $\mathbf{D}_\tau$ are independent of $\tau$
but $\varphi_\tau$ moves in $\mM^d_{\mathbf{b}}(j)$ as $\tau$ varies,
e.g.~this is the relevant situation for the proof of
super-rigidity.  There is then a subtle but important difference between what Proposition~\ref{prop:punctured}
says about $\dot{\mathbf{D}}_\tau$ and what it says about
$\varphi_\tau^*\dot{\mathbf{D}}_\tau$.  The former is a family of operators
whose relationship to each other for different values of $\tau$ is not
obvious from the definitions, but the proposition implies that they
are all in some sense equivalent to a single operator $\mathbf{D}$ on the
closed domain, so they all have isomorphic kernels.  No such thing can be assumed for
the pulled back operators: while $\varphi_\tau^*\dot{\mathbf{D}}_\tau$ must
have the same index for all~$\tau$, there is nothing in this setup to stop
the dimension of its kernel from varying wildly with~$\tau$.
\end{remark}

\subsection{A digression on representation theory}
\label{sec:representations}

In preparation for the twisted bundle construction in the next section,
we now collect some general facts from representation theory.

\subsubsection{Real permutation representations and subrepresentations}
\label{sec:realComplex}

Given a finite set $I$ with $d := |I| \in \NN$ elements and a finite group
with a homomorphism
$$
\rho : G \to S(I) : g \mapsto \rho_g
$$
defining a transitive group action on~$I$,
we denote by $\RR^I$ the real vector space spanned by basis vectors
$\{e_i\}_{i \in I}$, with an inner product such that this basis is
orthonormal.  We shall use the boldface symbol $\boldsymbol{\rho}$ to
denote the corresponding real $d$-dimensional representation of~$G$,
\begin{equation}
\label{eqn:permRep}
\boldsymbol{\rho} : G \to \Aut_\RR(\RR^I) \quad\text{ such that } \quad
\boldsymbol{\rho}(g) e_i := e_{\rho_g(i)}.
\end{equation}
We will be interested in the decomposition of $\RR^I$ into irreducible
$G$-invariant summands.  This can be understood in terms of its
complexification
$$
\boldsymbol{\rho}_\CC : G \to \Aut_\CC(\CC^I),
$$
defined by viewing $\{e_i\}_{i \in I}$ as a complex basis of~$\CC^I$.  
In general, we say that a complex representation $\boldsymbol{\lambda} : G \to \Aut_\CC(V)$
is the \defin{complexification} of a real representation
$\boldsymbol{\theta} : G \to \Aut_\RR(W)$ if $V$ is isomorphic to $W \oplus iW$ such that
$G$ acts on the latter by the complex-linear extension of its action on~$W$.
Recall from \cite{Serre}*{\S 13.2} that irreducible complex representations
${\boldsymbol{\lambda}} : G \to \Aut_\CC(V)$ come in three mutually exclusive types:
\begin{itemize}
\setlength{\itemsep}{0cm}
\item \defin{Real type}: $V$ admits a complex-antilinear $G$-invariant involution.
Then ${\boldsymbol{\lambda}}$ is the complexification of a real irreducible representation
$\boldsymbol{\theta} : G \to \Aut_\RR(W)$.
It follows that ${\boldsymbol{\lambda}}$ is isomorphic to its dual 
representation~$\boldsymbol{\lambda}^* : G \to \Aut_\CC(V^*)$, and
all $G$-equivariant linear maps $W \to W$ are given by scalar multiplication:
$$
\End_G(W) \cong \RR.
$$
\item \defin{Complex type}: ${\boldsymbol{\lambda}}$ is not isomorphic to its dual 
representation ${\boldsymbol{\lambda}}^* : G \to \Aut_\CC(V^*)$.  Then
${\boldsymbol{\lambda}} \oplus {\boldsymbol{\lambda}}^* : G \to \Aut_\CC(V \oplus V^*)$ is the
complexification of a real irreducible representation
$\boldsymbol{\theta} : G \to \Aut_\RR(W)$ obtained from ${\boldsymbol{\lambda}} : G \to \Aut_\CC(V)$
by setting $W := V$ and using the obvious inclusion 
$\Aut_\CC(V) \subset \Aut_\RR(W)$.  The algebra of $G$-equivariant real-linear
maps on $W$ is then
$$
\End_G(W) \cong \CC.
$$
\item \defin{Quaternionic type}: ${\boldsymbol{\lambda}}$ is not of real type but is 
nonetheless isomorphic to its dual representation.  Then
${\boldsymbol{\lambda}} \oplus {\boldsymbol{\lambda}} : G \to \Aut_\CC(V \oplus V)$ is the
complexification of a real irreducible representation
$\boldsymbol{\theta} : G \to \Aut_\RR(W)$ obtained from ${\boldsymbol{\lambda}} : G \to \Aut_\CC(V)$
by setting $W := V$ and using the obvious inclusion
$\Aut_\CC(V) \subset \Aut_\RR(W)$, and the algebra of $G$-equivariant real-linear
maps on $W$ is isomorphic to the quaternions:
$$
\End_G(W) \cong \HH.
$$
\end{itemize}

We shall also refer to a real irreducible representation as 
``of \defin{real} / \defin{complex} / \defin{quaternionic type}''
according to which of these three constructions it comes from.
With this classification in mind, we denote the various complex 
irreducible representations of $G$ by
$$
{\boldsymbol{\lambda}}_{j,\KK} : G \to \Aut_\CC(V_{j,\KK}),
$$
where $\KK$ stands for $\RR$, $\CC$ or $\HH$ depending on the type, and arrange
a complete list of pairwise non-isomorphic irreducible representations in the 
form
$$
{\boldsymbol{\lambda}}_{1,\RR},\ldots,{\boldsymbol{\lambda}}_{p,\RR},\ 
{\boldsymbol{\lambda}}_{1,\CC},{\boldsymbol{\lambda}}^*_{1,\CC},\ldots,{\boldsymbol{\lambda}}_{q,\CC},{\boldsymbol{\lambda}}^*_{q,\CC},\ 
{\boldsymbol{\lambda}}_{1,\HH},\ldots,{\boldsymbol{\lambda}}_{n,\HH}.
$$
This gives rise to a corresponding complete list
$$
\boldsymbol{\theta}_{1,\RR},\ldots,\boldsymbol{\theta}_{p,\RR},\ 
\boldsymbol{\theta}_{1,\CC},\ldots,\boldsymbol{\theta}_{q,\CC},\ 
\boldsymbol{\theta}_{1,\HH},\ldots,\boldsymbol{\theta}_{n,\HH}
$$
of pairwise non-isomorphic real irreducible representations
$$
\boldsymbol{\theta}_{j,\KK} : G \to \Aut_\RR(W_{j,\KK}) \quad \text{ satisfying } \quad
\End_G(W_{j,\KK}) \cong \KK,
$$
where for each~$j$, the complexification of $\boldsymbol{\theta}_{j,\KK}$ is
${\boldsymbol{\lambda}}_{j,\RR}$ for $\KK=\RR$, ${\boldsymbol{\lambda}}_{j,\CC} \oplus {\boldsymbol{\lambda}}^*_{j,\CC}$
for $\KK=\CC$, and ${\boldsymbol{\lambda}}_{j,\HH} \oplus {\boldsymbol{\lambda}}_{j,\HH}$ for $\KK=\HH$.
Note that the $G$-equivariant endomorphisms endow each $W_{j,\KK}$
with the structure of a left $\KK$-module such that the representation
$\boldsymbol{\theta}_{j,\KK}$ is $\KK$-linear.

We recall a standard fact from representation theory:

\begin{prop}
Every finite-dimensional representation $\boldsymbol{\theta} : G \to \Aut(W)$
of a finite group $G$ has a unique \defin{isotypic decomposition}, meaning
a splitting $W = W_1 \oplus \ldots \oplus W_N$ such that:
\begin{enumerate}
\item For each $i=1,\ldots,N$, $W_i \subset W$ is a
$G$-invariant subspace on which $\boldsymbol{\theta}$ is isomorphic to a
direct sum of copies of a single irreducible representation;
\item The irreducible representations corresponding any two distinct subspaces 
in the splitting are not isomorphic.
\end{enumerate}
\qed
\end{prop}

Since $\boldsymbol{\rho}_\CC$ itself is a complexification of a real representation,
every subspace in the resulting isotypic decomposition of $\CC^I$ is either
identical or orthogonal to its complex conjugate, where the conjugate always
carries the dual representation.  Thus we can uniquely
decompose $\CC^I$ into pairwise orthogonal $G$-invariant complex subspaces
\begin{equation}
\label{eqn:isotypic}
\CC^I = X_{1,\RR} \oplus \ldots \oplus X_{p,\RR} \oplus X_{1,\CC} \oplus 
\widebar{X}_{1,\CC} \oplus \ldots \oplus X_{q,\CC} \oplus \widebar{X}_{q,\CC}
\oplus X_{1,\HH} \oplus \ldots \oplus X_{n,\HH},
\end{equation}
where each $X_{j,\RR}$ and $X_{j,\HH}$ is of the form
$Y_{j,\KK} \oplus i Y_{j,\KK}$ for some real subspace $Y_{j,\KK} \subset \RR^I$,
and each $X_{j,\CC}$ has trivial intersection with~$\RR^I$.
Next, observe that every irreducible $G$-invariant subspace in $\CC^I$ is
either identical to its complex conjugate or intersects it trivially:
indeed, any other option would produce an intersection which is a
nontrivial but smaller $G$-invariant subspace.  We can thus further decompose
$X_{j,\RR}$ and $X_{j,\CC}$ into irreducible $G$-invariant subspaces
$$
X_{j,\RR} \cong V_{j,\RR}^{\oplus k_j}, \qquad
X_{j,\CC} \cong V_{j,\CC}^{\oplus m_j}
$$
for some integers $k_j, m_j \ge 0$, where each $V_{j,\RR}$ summand in
$X_{j,\RR}$ can be assumed of the form $W_{j,\RR} \oplus i W_{j,\RR}$ for
some irreducible $G$-invariant real subspace $W_{j,\RR} \subset Y_{j,\RR}$.
In $X_{j,\HH}$, the irreducible $G$-invariant subspaces cannot be
complexifications since the corresponding representation is not realizable
over~$\RR$, thus these subspaces have trivial intersection with $\RR^I$ and
can instead be arranged in conjugate pairs:
$$
X_{j,\HH} \cong V_{j,\HH}^{\oplus \ell_j} \oplus \widebar{V_{j,\HH}^{\oplus \ell_j}}
$$
for some integers $\ell_j \ge 0$.  From this decomposition of $\boldsymbol{\rho}_\CC$ we can
immediately read off a corresponding decomposition of $\boldsymbol{\rho}$: we have
\begin{equation}
\label{eqn:isotypicReal}
\RR^I = Y_{1,\RR} \oplus \ldots \oplus Y_{p,\RR} \oplus
Y_{1,\CC} \oplus \ldots \oplus Y_{q,\CC} \oplus
Y_{1,\HH} \oplus \ldots \oplus Y_{n,\HH},
\end{equation}
where the summands are all $G$-invariant and pairwise orthogonal,
$Y_{j,\KK} = X_{j,\KK} \cap \RR^I$ for $\KK = \RR,\HH$, and
$Y_{j,\CC} = (X_{j,\CC} \oplus \widebar{X}_{j,\CC}) \cap \RR^I$,
hence,
$$
\dim_\RR Y_{j,\KK} = \begin{cases}
\dim_\CC X_{j,\KK} & \text{ if $\KK = \RR$ or $\HH$},\\
2 \dim_\CC X_{j,\KK} & \text{ if $\KK = \CC$}.
\end{cases}
$$
These summands admit further (non-unique) decompositions into real 
irreducible $G$-invariant subspaces
$$
Y_{j,\RR} \cong W_{j,\RR}^{\oplus k_j}, \qquad
Y_{j,\CC} \cong W_{j,\CC}^{\oplus m_j}, \qquad
Y_{j,\KK} \cong W_{j,\HH}^{\oplus \ell_j}.
$$

\subsubsection{The regular case}
\label{sec:regular}

We now specialize the above discussion to the case 
$$
I := G, \qquad \rho_g(h) := gh,
$$
in which case $\boldsymbol{\rho} : G \to \Aut_\RR(\RR^G)$ is the so-called
\defin{regular representation} of~$G$.  By a standard theorem in complex
representation theory, the complexification $\boldsymbol{\rho}_\CC : G \to \Aut_\CC(\CC^G)$ then
contains every irreducible complex representation 
$\boldsymbol{\lambda}_{j,\KK} : G \to \Aut_\CC(V_{j,\KK})$
as a subrepresentation with multiplicity equal to $\dim_\CC V_{j,\KK}$.
This implies a similar fact about $\boldsymbol{\rho}$ that we will make use
of in \S\ref{sec:SardSmale} for proving Theorem~\ref{thm:submanifolds0}:

\begin{lemma}
\label{lemma:regular}
The real regular representation $\boldsymbol{\rho} : G \to \Aut_\RR(\RR^G)$ contains
every irreducible representation $\boldsymbol{\theta}_{j,\KK} : G \to
\Aut_\RR(W_{j,\KK})$ of $G$ as a subrepresentation with multiplicity
equal to $\dim_\KK W_{j,\KK}$.
\qed
\end{lemma}

Next, recall that the action of $G$ on itself by right multiplication
$$
G \to S(G) : g \mapsto \rho'_g, \qquad \rho'_g h := h g^{-1}
$$
commutes with $\rho$ and thus defines a second permutation representation
$$
\boldsymbol{\rho}' : G \to \Aut_\RR(\RR^G), \qquad
\boldsymbol{\rho}'(g) e_h = e_{h g^{-1}}
$$
which commutes with~$\boldsymbol{\rho}$, giving rise to a representation
\begin{equation}
\label{eqn:GtimesG}
G \times G \to \Aut_\RR(\RR^G) : (g,h) \mapsto \boldsymbol{\rho}(g)
\boldsymbol{\rho}'(h).
\end{equation}
By another standard theorem of complex representation theory, the summands in the
isotypic decomposition \eqref{eqn:isotypic} of $\CC^G$ are then invariant
under the complexification of the $(G \times G)$-action \eqref{eqn:GtimesG},
and they define irreducible complex representations of $G \times G$.
In particular, $\boldsymbol{\rho}'$ therefore preserves each isotypic component
for $\boldsymbol{\rho}$ but does not preserve any further decomposition of that component into
irreducible $G$-invariant subspaces.  For future use, we note one additional 
fact from complex representation theory: the action of $G \times G$ on an 
isotypic component in $\CC^G$ corresponding to a given irreducible
representation ${\boldsymbol{\lambda}} : G \to \Aut_\CC(V)$ is isomorphic to
$V \otimes V^*$, with $G \times G$ acting by
$$
(G \times G) \times (V \otimes V^*) \to V \otimes V^* :
\left( (g,h) , v \otimes \alpha \right) \mapsto {\boldsymbol{\lambda}}(g)v \otimes {\boldsymbol{\lambda}}^*(h)\alpha,
$$
cf.~\cite{Serre}*{\S 6.2}.

\subsubsection{Non-faithful representations}
\label{sec:nonfaithful}

An important special case of the factorization construction in
Example~\ref{ex:factorization} arises when
$$
\boldsymbol{\theta} : G \to \Aut_\RR(W)
$$
is an irreducible representation that is not faithful.  Choosing
$H$ to be any nontrivial normal subgroup of $G$ contained in its kernel
$$
H \subset \ker \boldsymbol{\theta} \subset G,
$$
$G/H$ then inherits an irreducible representation
$$
\boldsymbol{\theta}_H : G/H \to \Aut_\RR(W).
$$
For example one can take $H = \ker \boldsymbol{\theta}$, in which case $\boldsymbol{\theta}_H$ becomes
faithful.  Now if $\rho : G \to S(I)$ is a transitive action on the set $I$
of $d$ elements, let 
$$
\rho_H : G/H \to S(I/H)
$$
denote the induced action on the set $I / H$ of $H$-orbits, and
consider the corresponding permutation representations
$$
\boldsymbol{\rho} : G \to \Aut_\RR(\RR^I), \qquad
\boldsymbol{\rho}_H : G/H \to \Aut_\RR(\RR^{I/H}).
$$

\begin{lemma}
\label{lemma:characters}
Under the assumptions described above,
the multiplicity of $\boldsymbol{\theta} : G \to \Aut_\RR(W)$ as a
subrepresentation of
$\boldsymbol{\rho} : G \to \Aut_\RR(\RR^I)$ matches the multiplicity of
$\boldsymbol{\theta}_H : G/H \to \Aut_\RR(W)$ as a
subrepresentation of $\boldsymbol{\rho}_H : G/H \to \Aut_\RR(\RR^{I/H})$.
\end{lemma}
\begin{proof}
Observe that in terms of the real/complex/quaternionic distinction described in 
\S\ref{sec:realComplex}, $\boldsymbol{\theta}$ and $\boldsymbol{\theta}_H$ 
are necessarily of the same type: indeed, the spaces of linear maps on $W$
that are $G$-equivariant or $(G/H)$-equivariant are the same since
$H$ acts trivially on~$W$.  The multiplicities
of both are therefore determined in the same way by the multiplicities of
the corresponding \emph{complex} irreducible representations in the
complexifications of $\boldsymbol{\rho}$ and $\boldsymbol{\rho}_H$ respectively, thus it will suffice
to prove a similar statement about complex representations.
Namely, assume ${\boldsymbol{\lambda}} : G \to \Aut_\CC(V)$ is
complex irreducible, $H \subset \ker {\boldsymbol{\lambda}} \subset G$ is a normal subgroup
and ${\boldsymbol{\lambda}}_H : G/H \to \Aut_\CC(V)$ is the resulting irreducible
representation of~$G/H$.  By orthonormality of characters, it will suffice
to prove
$$
\langle \chi_{\boldsymbol{\rho}},\chi_{\boldsymbol{\lambda}} \rangle = \langle \chi_{\boldsymbol{\rho}_H},\chi_{{\boldsymbol{\lambda}}_H} \rangle,
$$
where the inner product of characters $\chi_{\boldsymbol{\lambda}} : G \to \CC$ is given 
in general by
$$
\langle \chi_{\boldsymbol{\lambda}},\chi_{{\boldsymbol{\lambda}}'} \rangle := \frac{1}{|G|} \sum_{g \in G}
\widebar{\chi_{\boldsymbol{\lambda}}(g)} \chi_{{\boldsymbol{\lambda}}'}(g) \in \CC.
$$
For each $i \in I$, let $G_i \subset G$ denote the stabilizer subgroup
for $i$ under the $G$-action on~$I$ via~$\boldsymbol{\rho}$.  Since the
action is transitive, the orbit-stabilizer theorem implies $|G_i| = |G| / d$.
The trace of
a permutation matrix is the number of elements that it fixes, in other words
the number of stabilizer subgroups that it belongs to, hence for each
$g \in G$,
$$
\chi_{\boldsymbol{\rho}}(g) = \left| \left\{ i \in I  \ \big|\ 
g \in G_i \right\} \right|.
$$
This implies
\begin{equation}
\label{eqn:GkThing}
\langle \chi_{\boldsymbol{\rho}}, \chi_{\boldsymbol{\lambda}} \rangle =
\frac{1}{|G|} \sum_{i \in I} \sum_{g \in G_i} \chi_{\boldsymbol{\lambda}}(g).
\end{equation}
This can be simplified since $G$ acts transitively on~$I$, so
the subgroups $G_i$ for distinct $i \in I$ are all
conjugate.  By the conjugation-invariance of
characters, this implies that all $d$ of the sums over $G_i$
in \eqref{eqn:GkThing} are identical, so plugging in $|G_i| = |G| / d$, we have
$$
\langle \chi_{\boldsymbol{\rho}},\chi_{\boldsymbol{\lambda}} \rangle =
\frac{1}{|G_i|} \sum_{g \in G_i} \chi_{\boldsymbol{\lambda}}(g),
$$
where $i \in I$ in this expression can be chosen arbitrarily.

To write down a similar expression for $\langle \chi_{{\boldsymbol{\rho}}_H} , \chi_{{\boldsymbol{\lambda}}_H} \rangle$,
define for each $i \in I$
$$
H_i := H \cap G_i \subset G,
$$
which is a subgroup of both $H$ and $G_i$ and is normal in the latter.
There is then a natural inclusion of $G_i / H_i$ as a subgroup of $G/H$,
and it is the stabilizer subgroup of $[i] \in I/H$ for the permutation
action of $G/H$ on~$I/H$.  The same computation thus gives
$$
\langle \chi_{{\boldsymbol{\rho}}_H}, \chi_{{\boldsymbol{\lambda}}_H} \rangle =
\frac{1}{| G_i / H_i|} \sum_{[g] \in G_i/H_i} \chi_{{\boldsymbol{\lambda}}_H}([g])
= \frac{|H_i|}{|G_i|} \sum_{[g] \in G_i/H_i} \chi_{{\boldsymbol{\lambda}}_H}([g]).
$$
Finally, observe that $\chi_{\boldsymbol{\lambda}}(g) = \chi_{{\boldsymbol{\lambda}}_H}([g])$ for each
$g \in G$ since both are traces of the same linear operator acting on~$V$, so one can 
replace the last expression with a sum over $g \in G_i$, giving
$$
\langle \chi_{{\boldsymbol{\rho}}_H}, \chi_{{\boldsymbol{\lambda}}_H} \rangle =
\frac{1}{|G_i|} \sum_{g \in G_i} \chi_{\boldsymbol{\lambda}}(g) = \langle \chi_{\boldsymbol{\rho}}, \chi_{\boldsymbol{\lambda}} \rangle.
$$
\end{proof}

\subsection{Twisted bundles and splittings of operators}
\label{sec:twisted}

We can now make precise the splitting of pulled back Cauchy-Riemann type 
operators that was sketched in \S\ref{sec:theBigIdea}.

\subsubsection{Twisted bundles from representations}

We associate to any representation $\boldsymbol{\theta} : G \to \Aut_\RR(W)$ the
family of real vector bundles $W_\tau^{\boldsymbol{\theta}} \to \dot{\Sigma}_\tau$ 
defined by
$$
W_\tau^{\boldsymbol{\theta}} = \left( \dot{\Sigma}''_\tau \times W \right) \Big/ G,
$$
where $G$ acts on $W$ via~$\boldsymbol{\theta}$ and on $\dot{\Sigma}''_\tau$ 
by deck transformations, so that $\pi_\tau : \dot{\Sigma}''_\tau \to \dot{\Sigma}_\tau$
identifies $\dot{\Sigma}_\tau$ with $\dot{\Sigma}''_\tau / G$.  
This gives rise to complex vector bundles
$\dot{E}^{\boldsymbol{\theta}}_\tau, \dot{F}^{\boldsymbol{\theta}}_\tau \to 
\dot{\Sigma}_\tau$ of rank $m \cdot \dim_\RR W$, defined by
$$
\dot{E}^{\boldsymbol{\theta}}_\tau = \dot{E}_\tau \otimes_\RR W^{\boldsymbol{\theta}}_\tau, \qquad
\dot{F}^{\boldsymbol{\theta}}_\tau = \dot{F}_\tau \otimes_\RR W^{\boldsymbol{\theta}}_\tau = 
\overline{\Hom}_\CC(T\dot{\Sigma}_\tau,\dot{E}^{\boldsymbol{\theta}}_\tau).
$$
Each of the bundles $W^{\boldsymbol{\theta}}_\tau$ has a canonical flat 
structure, i.e.~it comes with a
well-defined notion of constant local sections, thus $\mathbf{D}_\tau \in
\CRR(E_\tau)$ determines a family of Cauchy-Riemann type operators
$$
\dot{\mathbf{D}}^{\boldsymbol{\theta}}_\tau : \Gamma(\dot{E}^{\boldsymbol{\theta}}_\tau) 
\to \Gamma(\dot{F}^{\boldsymbol{\theta}}_\tau) =
\Omega^{0,1}(\dot{\Sigma}_\tau,\dot{E}^{\boldsymbol{\theta}}_\tau)
$$
such that $\dot{\mathbf{D}}^{\boldsymbol{\theta}}_\tau (\eta \otimes v) = 
\dot{\mathbf{D}}_\tau \eta 
\otimes v$ whenever $v$ is a constant local section of~$W^{\boldsymbol{\theta}}_\tau$.
Since $\dot{\mathbf{D}}^{\boldsymbol{\theta}}_\tau \in \CRR(\dot{E}_\tau^{\boldsymbol{\theta}})$,
it is Fredholm in suitable Banach space settings, in particular as a bounded linear
operator
$$
\dot{\mathbf{D}}^{\boldsymbol{\theta}}_\tau :
W^{k,p,-\boldsymbol{\delta}}(\dot{E}^{\boldsymbol{\theta}}_\tau) \to
W^{k-1,p,-\boldsymbol{\delta}}(\dot{F}^{\boldsymbol{\theta}}_\tau)
$$
for any $k \in \NN$, $p \in (1,\infty)$, and negative exponential weights
$-\boldsymbol{\delta} = \{-\delta_w\}_{w \in \Theta}$ with all
$\delta_w > 0$ sufficiently small.  We will formulate a precise
version of this statement and compute the index in \S\ref{sec:index}.
Observe that aside from its obvious dependence on $\mathbf{D}_\tau$,
$\dot{\mathbf{D}}^{\boldsymbol{\theta}}_\tau$ depends on our choice of
regular presentation for $\varphi$ and on the representation~$\boldsymbol{\theta}$,
but both of them only up to isomorphism.

If $\boldsymbol{\theta}$ is irreducible with $\End_G(W) = \KK \in \{\CC,\HH\}$,
then the resulting left $\KK$-module structure of $W$ induces a left $\KK$-module
structure on each fiber of the twisted bundles $\dot{E}^{\boldsymbol{\theta}}_\tau$
and~$\dot{F}^{\boldsymbol{\theta}}_\tau$, for which the twisted operator
$\dot{\mathbf{D}}^{\boldsymbol{\theta}}_\tau$ commutes with the action of~$\KK$,
thus its kernel and cokernels are also left $\KK$-modules.
Note that if $\KK=\CC$, the resulting complex structure on 
$\dot{E}^{\boldsymbol{\theta}}_\tau$ and~$\dot{F}^{\boldsymbol{\theta}}_\tau$
is different from the one defined by~$J$; the latter does not commute with
$\dot{\mathbf{D}}^{\boldsymbol{\theta}}_\tau$ unless $\mathbf{D}_\tau$ is a
$J$-linear operator to start with.

The most important special case of the above construction 
is $\dot{E}^{\boldsymbol{\rho}}_\tau \to \dot{\Sigma}_\tau$,
where ${\boldsymbol{\rho}} : G \to \Aut_\RR(\RR^I)$ is the permutation 
representation associated to our regular presentation of~$\varphi$.
We define
$\dot{E}^{\boldsymbol{\rho}}_\tau = \dot{E}_\tau \otimes (\RR^I)^{\boldsymbol{\rho}}_\tau \to \dot{\Sigma}_\tau$ as above
and can identify it canonically with
$$
\dot{E}^{\boldsymbol{\rho}}_\tau = \left( \pi_\tau^*\dot{E}_\tau \otimes \RR^I \right) \Big/ G,
$$
so that sections of
$\dot{E}^{\boldsymbol{\rho}}_\tau$ are written as $G$-equivariant sections of
$\pi_\tau^*\dot{E}_\tau \otimes \RR^I$, hence
$$
\eta = \sum_{i \in I} \eta^i \otimes e_i
$$
for $\eta^i \in \Gamma(\pi_\tau^*\dot{E}_\tau)$.  Here $G$-equivariance means that
for all $z \in \dot{\Sigma}''_\tau$ and $g \in G$,
$$
\eta(gz) = (\1 \otimes {\boldsymbol{\rho}}(g)) \eta(z) = \sum_{i \in I} \eta^i(z) \otimes
e_{\rho_g(i)},
$$
hence
\begin{equation}
\label{eqn:etaEquivariance}
\eta^i(z) = \eta^{\rho_g(i)}(gz) \quad \text{ for all } \quad
\text{$z \in \dot{\Sigma}''_\tau$, $g \in G$ and $i \in I$}.
\end{equation}
Writing $\dot{\Sigma}_\tau' = (\dot{\Sigma}''_\tau \times I ) / G$,
this relation gives rise to a bijective correspondence
\begin{equation}
\label{eqn:twistedAndPullback}
\begin{split}
\Gamma(\dot{E}^{\boldsymbol{\rho}}_\tau) \to \Gamma(\varphi_\tau^*\dot{E}_\tau) 
& : \eta \mapsto \widehat{\eta} \\
\widehat{\eta}([(z,i)]) = \eta^i(z)
\end{split}
\end{equation}
and thus natural isomorphisms
\begin{equation}
\label{eqn:twistedPullbackWeights}
W^{k,p,-\boldsymbol{\delta}}(\dot{E}^{\boldsymbol{\rho}}_\tau) \to
W^{k,p,-\varphi^*\boldsymbol{\delta}}(\varphi_\tau^*\dot{E}_\tau)
\end{equation}
for every $k \ge 0$ and $p \in (1,\infty)$, where we recall from
\S\ref{sec:CRpunctured} that the pulled back exponential weights
are defined by
$$
\varphi^*\boldsymbol{\delta} := \left\{ k_\zeta \delta_{\varphi(\zeta)} \right\}_{\zeta \in \Theta'},
$$
with $k_\zeta \in \{1,\ldots,d\}$ denoting the branching order of
$\varphi : \Sigma' \to \Sigma$ at $\zeta \in \Theta'$.
The reason for using these particular weights in the isomorphism
\eqref{eqn:twistedPullbackWeights} is as follows.
We observe first that if $\varphi : [0,\infty) \times S^1 \to [0,\infty) \times S^1$
is a holomorphic covering map of the form $(s,t) \mapsto (ms,mt)$
and $\ZZ_m$ is defined to act on $[0,\infty) \times S^1$ via the
transformation $(s,t) \mapsto (s,t + 1/m)$ and its iterates,
then the map $f \mapsto f \circ \varphi$ defines for each integer $k \ge 0$
and $p \in (1,\infty)$ an isomorphism from $W^{k,p}([0,\infty) \times S^1)$
to the closed subspace of $W^{k,p}([0,\infty) \times S^1)$ consisting
of $\ZZ_m$-invariant functions.  It follows that for any exponential
weight $\delta$, a function $f$ on $[0,\infty) \times S^1$ is of
class $W^{k,p,\delta}$ if and only if $f \circ \varphi$ 
is of class $W^{k,p,m\delta}$.  The global consequence of these
observations is that for $\eta \in \Gamma(\dot{E}_\tau^{\boldsymbol{\rho}})$
and the corresponding section $\widehat{\eta} \in \Gamma(\varphi_\tau^*\dot{E}_\tau)$,
the $W^{k,p,-\varphi^*\boldsymbol{\delta}}$-norm of $\widehat{\eta}$ can be
bounded in terms of the $W^{k,p,-\boldsymbol{\delta}}$-norm of $\eta$, and
vice versa.

Observe that $(\RR^I)^{\boldsymbol{\rho}}_\tau \to \dot{\Sigma}_\tau$ also
has a well-defined real bundle metric since $\boldsymbol{\rho}$ acts on
$\RR^I$ by orthogonal transformations, so endowing $E_\tau$ with a Hermitian
bundle metric induces a Hermitian bundle metric on
$\dot{E}^{\boldsymbol{\rho}}_\tau = \dot{E}_\tau \otimes (\RR^I)^{\boldsymbol{\rho}}_\tau$
such that the correspondence \eqref{eqn:twistedAndPullback} 
also preserves $L^2$-products.
After writing down a similar correspondence for the bundles $\dot{F}^{\boldsymbol{\rho}}_\tau$ and
$\varphi_\tau^*\dot{F}_\tau$, we obtain an identification between the
Cauchy-Riemann operators $\varphi_\tau\dot{\mathbf{D}}_\tau$ and 
$\dot{\mathbf{D}}^{\boldsymbol{\rho}}_\tau$:
\begin{equation}
\label{eqn:CD}
\begin{CD}
W^{k,p,-\boldsymbol{\delta}}(\dot{E}^{\boldsymbol{\rho}}_\tau) 
@>\dot{\mathbf{D}^{\boldsymbol{\rho}}_\tau}>>
W^{k-1,p,-\boldsymbol{\delta}}(\dot{F}^{\boldsymbol{\rho}}_\tau) \\
@VV\cong V     @VV\cong V \\
W^{k,p,-\varphi^*\boldsymbol{\delta}}(\varphi_\tau^*\dot{E}_\tau) 
@>\varphi_\tau^*\dot{\mathbf{D}}_\tau>>
W^{k-1,p,-\varphi^*\boldsymbol{\delta}}(\varphi_\tau^*\dot{F}_\tau), \\
\end{CD}
\end{equation}

\subsubsection{Splitting the twisted Cauchy-Riemann operator}

If $W \subset \RR^I$ is any $G$-invariant
subspace and $\boldsymbol{\theta} : G \to \Aut_\RR(W)$ denotes the resulting
subrepresentation, then we obtain corresponding subbundles
$$
\dot{E}^{\boldsymbol{\theta}}_\tau \subset \dot{E}^{\boldsymbol{\rho}}_\tau, 
\qquad \dot{F}^{\boldsymbol{\theta}}_\tau \subset \dot{F}^{\boldsymbol{\rho}}_\tau
$$
such that $\dot{\mathbf{D}}^{\boldsymbol{\rho}}_\tau$ takes sections of 
$\dot{E}^{\boldsymbol{\theta}}_\tau$ to
sections of $\dot{F}^{\boldsymbol{\theta}}_\tau$, acting as the operator 
$\dot{\mathbf{D}}^{\boldsymbol{\theta}}_\tau$.
Under the correspondence \eqref{eqn:twistedAndPullback}, one can understand this
as identifying $\Gamma(\dot{E}^{\boldsymbol{\theta}}_\tau)$ and 
$\Gamma(\dot{F}^{\boldsymbol{\theta}}_\tau)$ with closed subspaces
$$
\Gamma_{\boldsymbol{\theta}}(\varphi_\tau^*\dot{E}_\tau) \subset \Gamma(\varphi_\tau^*\dot{E}_\tau), \qquad
\Gamma_{\boldsymbol{\theta}}(\varphi_\tau^*\dot{F}_\tau) \subset \Gamma(\varphi_\tau^*\dot{F}_\tau),
$$
with a similar definition for closed subspaces of the relevant weighted 
Sobolev spaces, such that $\varphi_\tau^*\dot{\mathbf{D}}_\tau$ restricts to a bounded 
linear operator
$$
W^{k,p,-\varphi^*\boldsymbol{\delta}}_{\boldsymbol{\theta}}(\varphi_\tau^*\dot{E}_\tau)
\stackrel{\varphi_\tau^*\dot{\mathbf{D}}_\tau}{\longrightarrow}  W^{k-1,p,-\varphi^*\boldsymbol{\delta}}_{\boldsymbol{\theta}}(\varphi_\tau^*\dot{F}_\tau),
$$
which is conjugate to $\dot{\mathbf{D}}_\tau^{\boldsymbol{\theta}}
: W^{k,p,-\boldsymbol{\delta}}(\dot{E}_\tau^{\boldsymbol{\theta}}) \to
W^{k-1,p,-\boldsymbol{\delta}}(\dot{F}_\tau^{\boldsymbol{\theta}})$ and will 
thus be Fredholm with any negative exponential weights that are close enough to~$0$.  Now if
$$
\RR^I = W_1 \oplus \ldots \oplus W_N
$$
is a decomposition of ${\boldsymbol{\rho}}$ into subrepresentations
$\boldsymbol{\theta}_j : G \to \Aut_\RR(W_j)$ for $j=1,\ldots,N$, we obtain a direct
sum decomposition 
$$
\dot{\mathbf{D}}^{\boldsymbol{\rho}}_\tau = 
\dot{\mathbf{D}}^{\boldsymbol{\theta}_1}_\tau \oplus 
\ldots \oplus \dot{\mathbf{D}}^{\boldsymbol{\theta}_N}_\tau,
$$
which is equivalent via \eqref{eqn:CD} to a decomposition of $\varphi_\tau^*\dot{\mathbf{D}}_\tau$ 
over a splitting of Banach spaces
$$
W^{k,p,-\varphi^*\boldsymbol{\delta}}(\varphi_\tau^*E_\tau) =
\bigoplus_{j=1}^N W^{k,p,-\varphi^*\boldsymbol{\delta}}_{\boldsymbol{\theta}_j}(\varphi_\tau^*E_\tau)
$$
and the corresponding decomposition of 
$W^{k-1,p,-\varphi^*\boldsymbol{\delta}}(\varphi_\tau^*F_\tau)$.
Observe that if the subspaces $W_1,\ldots,W_N \subset \RR^I$ are pairwise orthogonal,
then the corresponding spaces of sections of $\varphi_\tau^*\dot{E}_\tau$ and
$\varphi_\tau^*\dot{F}_\tau$ are $L^2$-orthogonal as a consequence.
It is useful to note that whenever two of the representations
$\boldsymbol{\theta}_i : G \to \Aut_\RR(W_i)$ and
$\boldsymbol{\theta}_j : G \to \Aut_\RR(W_j)$ are isomorphic, the 
$G$-equivariant isomorphism $W_i \to W_j$ 
induces bundle isomorphisms $\dot{E}_\tau^{\boldsymbol{\theta}_i} \to
\dot{E}_\tau^{\boldsymbol{\theta}_j}$ and
$\dot{F}_\tau^{\boldsymbol{\theta}_i} \to \dot{F}_\tau^{\boldsymbol{\theta}_j}$
that identify $\dot{\mathbf{D}}_\tau^{\boldsymbol{\theta}_i}$ with
$\dot{\mathbf{D}}_\tau^{\boldsymbol{\theta}_j}$, so these two operators have
isomorphic kernels and cokernels.  This implies:

\begin{lemma}
\label{lemma:kernelDecomp}
Suppose $\boldsymbol{\theta}_j : G \to \Aut_\RR(W_j)$ for $j=1,\ldots,N$ 
is a collection of representations of~$G$, and
$\boldsymbol{\theta} : G \to \Aut_\RR(W)$ is another representation such that
$$
\boldsymbol{\theta} \cong \bigoplus_{j=1}^N \boldsymbol{\theta}_j^{\oplus k_j}
$$
for some integers $k_1,\ldots,k_N \ge 0$.
Then there exist isomorphisms
$$
\ker \dot{\mathbf{D}}_\tau^{\boldsymbol{\theta}} \cong
\bigoplus_{j=1}^N \left(\ker \dot{\mathbf{D}}_\tau^{\boldsymbol{\theta}_j}\right)^{\oplus k_j}
\quad \text{ and } \quad
\coker \dot{\mathbf{D}}_\tau^{\boldsymbol{\theta}} \cong
\bigoplus_{j=1}^N \left(\coker \dot{\mathbf{D}}_\tau^{\boldsymbol{\theta}_j}\right)^{\oplus k_j}.
$$
In particular, if $\boldsymbol{\theta}$ is the permutation representation 
$\boldsymbol{\rho} : G \to \Aut_\RR(\RR^I)$, this gives isomorphisms
$$
\ker (\varphi_\tau^*\dot{\mathbf{D}}_\tau)  \cong
\bigoplus_{j=1}^N \left(\ker \dot{\mathbf{D}}_\tau^{\boldsymbol{\theta}_j}\right)^{\oplus k_j}
\quad \text{ and } \quad
\coker (\varphi_\tau^*\dot{\mathbf{D}}_\tau) \cong
\bigoplus_{j=1}^N \left(\coker \dot{\mathbf{D}}_\tau^{\boldsymbol{\theta}_j}\right)^{\oplus k_j}.
$$
\qed
\end{lemma}

\subsubsection{Non-faithful representations revisited}
\label{sec:nonfaithfulRevisit}

Here is a proof of Lemma~\ref{lemma:nonfaithful}.  For the present discussion
we drop the parameter $\tau$ from the notation since it does not play
any important role.

Suppose $\boldsymbol{\theta} : G \to \Aut_\RR(W)$ is a representation and
$H \subset \ker \boldsymbol{\theta} \subset G$ is a nontrivial normal subgroup of~$G$, giving 
rise to a representation 
$$
\boldsymbol{\theta}_H : G/H \to \Aut_\RR(W),
$$ 
and (following Example~\ref{ex:factorization}) a factorization of 
$\varphi : \Sigma' \to \Sigma$ as
$$
\Sigma' \to \Sigma'_H \stackrel{\varphi_H}{\longrightarrow} 
\Sigma.
$$
By assumption we are using a minimal regular presentation and thus 
$\rho : G \to S(I)$ is injective, so $H$ acts nontrivially on $I$, implying
$\deg(\varphi_H) < d$.
Writing $\dot{\Sigma}''_H = \dot{\Sigma}'' \big/ H$, the obvious
projection map
$$
\left(\dot{\Sigma}'' \times W\right) \Big/ G \to 
\left(\dot{\Sigma}''_H \times W\right) \Big/ (G/H)
$$
is then an isomorphism of real vector bundles over $\dot{\Sigma}$
and thus gives rise to a canonical identification between the 
twisted bundles 
$\dot{E}^{\boldsymbol{\theta}}$ and $\dot{E}^{\boldsymbol{\theta}_H}$ with their Cauchy-Riemann operators
$\dot{\mathbf{D}}^{\boldsymbol{\theta}}$ and $\dot{\mathbf{D}}^{\boldsymbol{\theta}_H}$.  
To prove the lemma, we now just need to observe that Lemma~\ref{lemma:characters}
implies $\boldsymbol{\theta}$ is a subrepresentation of $\boldsymbol{\rho}$
if and only if $\boldsymbol{\theta}_H$ is a subrepresentation of $\boldsymbol{\rho}_H$,
hence the corresponding twisted operators appear simultaneously as summands
in the decompositions of $\varphi^*\dot{\mathbf{D}}$ and
$\varphi_H^*\dot{\mathbf{D}}$ from Lemma~\ref{lemma:kernelDecomp}.

\begin{remark}
In the situation above, one should interpret $\ker \dot{\mathbf{D}}^{\boldsymbol{\theta}}$
as the set of all sections in $\ker (\varphi^*\dot{\mathbf{D}})$ that are
pullbacks of sections in $\ker \dot{\mathbf{D}}^{\boldsymbol{\theta}_H}$
(interpreted as a subspace of $\ker (\varphi_H^*\dot{\mathbf{D}})$) via
the branched cover $\Sigma' \to \Sigma'_H$.
\end{remark}

\subsubsection{The regular case revisited}

Now consider the special case where $\boldsymbol{\rho}$ is the regular 
representation $G \to \Aut_\RR(\RR^G)$, defined via
$$
\rho : G \to S(G), \qquad \rho_g(h) = gh.
$$
We saw in Example~\ref{ex:regular}
that this means $\varphi_\tau : \dot{\Sigma}_\tau'
\to \dot{\Sigma}_\tau$ are all regular covers isomorphic to
$\pi : \dot{\Sigma}''_\tau \to \dot{\Sigma}_\tau$, and the action of $G$
on $\dot{\Sigma}'_\tau = (\dot{\Sigma}''_\tau \times G) / G$ by deck 
transformations takes the form
$$
g [(z,h)] := [(z,\rho'_g(h))]
$$
where $\rho' : G \to S(G)$ is the action of $G$ on itself by right multiplication,
$\rho'_g(h) = h g^{-1}$.  The induced $G$-action on spaces of sections
$\eta$ of $\varphi_\tau^*\dot{E}_\tau$ is defined by
$$
(g \eta)([(z,h)]) := \eta(g^{-1} [(z,h)]) = \eta([(z,hg)]).
$$
Recall now from \S\ref{sec:regular} that the
permutation representation $\boldsymbol{\rho}' : G \to \Aut_\RR(\RR^G)$
arising from $\rho'$ commutes with $\boldsymbol{\rho}$ and preserves the
isotypic components of~$\boldsymbol{\rho}$.  It therefore defines an
action on $\dot{E}^{\boldsymbol{\rho}}_\tau$ by fiber-preserving bundle
isomorphisms, and these isomorphisms preserve each of the subbundles in the 
splitting
\begin{equation}
\label{eqn:isoSubbundles}
\dot{E}^{{\boldsymbol{\rho}}}_\tau = \bigoplus_{j=1}^p (\dot{E}^{{\boldsymbol{\rho}}}_\tau)_{j,\RR} \oplus
\bigoplus_{j=1}^q (\dot{E}^{{\boldsymbol{\rho}}}_\tau)_{j,\CC} \oplus
\bigoplus_{j=1}^n (\dot{E}^{{\boldsymbol{\rho}}}_\tau)_{j,\HH}
\end{equation}
corresponding to the isotypic decomposition \eqref{eqn:isotypicReal}
of~${\boldsymbol{\rho}}$.  In particular, this $G$-action by bundle isomorphisms 
gives a linear $G$-action on each of the subspaces 
$\Gamma((\dot{E}^{\boldsymbol{\rho}}_\tau)_{j,\KK}) \subset 
\Gamma(\dot{E}^{\boldsymbol{\rho}}_\tau)$, and there is a 
similar action on sections of $\dot{F}^{\boldsymbol{\rho}}_\tau$ such that 
the restriction of $\dot{\mathbf{D}}^{\boldsymbol{\rho}}_\tau$ to each of 
these subspaces is $G$-equivariant.  Its kernel and cokernel thus inherit natural
$G$-actions.  Under the correspondence \eqref{eqn:twistedAndPullback},
this action on sections of $\dot{E}_\tau^{\boldsymbol{\rho}}$ matches the 
action by deck transformations on $\Gamma(\varphi_\tau^*\dot{E}_\tau)$.

\begin{lemma}
\label{lemma:ndimension}
Suppose $\rho : G \to S(G)$ is defined by left multiplication,
$\boldsymbol{\theta}_0 : G \to \Aut_\RR(W)$ is an irreducible representation of $G$,
and $\boldsymbol{\theta} : G \to \Aut_\RR(Y)$ denotes 
the corresponding summand in the isotypic decomposition 
\eqref{eqn:isotypicReal} of the regular representation~$\boldsymbol{\rho} :
G \to \Aut_\RR(\RR^G)$.  Then every irreducible subrepresentation for the
natural $G$-action on
$\ker \dot{\mathbf{D}}^{\boldsymbol{\theta}}_\tau$ or
$\coker \dot{\mathbf{D}}^{\boldsymbol{\theta}}_\tau$ is 
isomorphic to~$\boldsymbol{\theta}_0$.
\end{lemma}
\begin{proof}
Suppose first that $\boldsymbol{\theta}_0$ is of either real or quaternionic
type, in which case the complexification $X := Y \oplus iY \subset \CC^G$
of $Y \subset \RR^G$ is also an isotypic component for the complexified
regular representation $\boldsymbol{\rho}_\CC : G \to \Aut_\CC(\CC^G)$.
We shall denote the restriction of $\boldsymbol{\rho}_\CC$ to $X$ by
$$
\boldsymbol{\lambda} : G \to \Aut_\CC(X),
$$
and let $\boldsymbol{\lambda}_0 : G \to \Aut_\CC(V)$ denote the underlying
complex irreducible representation.  Regarding these complex representations
as real representations on $X$ and $V$ respectively gives rise to
corresponding twisted bundles and Cauchy-Riemann operators on them, along
with a natural linear inclusion of vector bundles
$$
\dot{E}_\tau^{\boldsymbol{\theta}} \hookrightarrow \dot{E}_\tau^{\boldsymbol{\lambda}} \qquad
\text{ such that } \qquad
\ker \dot{\mathbf{D}}_\tau^{\boldsymbol{\theta}} =
\ker \dot{\mathbf{D}}_\tau^{\boldsymbol{\lambda}} \cap 
\Gamma(\dot{E}_\tau^{\boldsymbol{\theta}}).
$$
It will be useful to think of $\dot{E}_\tau^{\boldsymbol{\lambda}}$ as a
\emph{complexification} of $\dot{E}_\tau^{\boldsymbol{\theta}}$, in the 
following sense.  While $\dot{E}_\tau^{\boldsymbol{\theta}}$ is already a
complex vector bundle, $\dot{E}_\tau^{\boldsymbol{\lambda}} =
\dot{E}_\tau \otimes_\RR X^{\boldsymbol{\lambda}}_\tau$ naturally carries 
\emph{two} complex structures $J_\tau$ and~$i$, which commute with each other:
the former acts on $\eta \otimes v \in \dot{E}_\tau \otimes_\RR X^{\boldsymbol{\lambda}}_\tau$
by $J_\tau\eta \otimes v$ and the latter by $\eta \otimes iv$, using the fact that
$\boldsymbol{\lambda}$ is a complex representation and $X^{\boldsymbol{\lambda}}_\tau$
is therefore naturally a complex vector bundle.
From this perspective, $\dot{\mathbf{D}}_\tau^{\boldsymbol{\lambda}}$ is
the natural $i$-complex-linear extension of
$\dot{\mathbf{D}}_\tau^{\boldsymbol{\theta}}$ to its complexified domain,
and the representations defined by the $G$-action on 
$\ker \dot{\mathbf{D}}_\tau^{\boldsymbol{\lambda}}$ and
$\coker \dot{\mathbf{D}}_\tau^{\boldsymbol{\lambda}}$ will be the 
complexifications of the real representations it defines on
$\ker \dot{\mathbf{D}}_\tau^{\boldsymbol{\theta}}$ and
$\coker \dot{\mathbf{D}}_\tau^{\boldsymbol{\theta}}$ respectively.
In the following we shall use the symbol ``$\otimes_i$'' to denote complex
tensor products of vector spaces and bundles with $i$ (instead of~$J_\tau$)
as the complex structure.

Recall now that as an isotypic component of the complex regular representation,
$X$ admits a complex-linear 
isomorphism to $V \otimes_i V^*$ such that for all $g \in G$, 
$\boldsymbol{\rho}(g)$ acts on $V \otimes_i V^*$ as $\boldsymbol{\lambda}_0 \otimes \1$,
while $\boldsymbol{\rho}'(g)$ acts as $\1 \otimes \boldsymbol{\lambda}_0^*$.
The isomorphism $X \to V \otimes_i V^*$ thus gives rise to $i$-complex bundle 
isomorphisms
$$
\dot{E}_\tau^{\boldsymbol{\lambda}} \to
\dot{E}_\tau^{\boldsymbol{\lambda}_0} \otimes_i V^*, \qquad
\dot{F}_\tau^{\boldsymbol{\lambda}} \to
\dot{F}_\tau^{\boldsymbol{\lambda}_0} \otimes_i V^*,
$$
where we are abusing notation to let $V^*$ denote the trivial bundle over
$\dot{\Sigma}_\tau$ with fiber~$V^*$, and this identifies
$\dot{\mathbf{D}}_\tau^{\boldsymbol{\lambda}}$ with 
$\dot{\mathbf{D}}_\tau^{\boldsymbol{\lambda}_0} \otimes \1$.
We therefore have
$$
\ker \dot{\mathbf{D}}_\tau^{\boldsymbol{\lambda}} \cong
\ker \dot{\mathbf{D}}_\tau^{\boldsymbol{\lambda}_0} \otimes_i V^*, \qquad
\coker \dot{\mathbf{D}}_\tau^{\boldsymbol{\lambda}} \cong
\coker \dot{\mathbf{D}}_\tau^{\boldsymbol{\lambda}_0} \otimes_i V^*,
$$
with $G$ acting on both by $\1 \otimes \boldsymbol{\lambda}^*_0$,
hence all irreducible subrepresentations in these spaces are isomorphic
to~$\boldsymbol{\lambda}^*_0$, which is isomorphic to $\boldsymbol{\lambda}_0$
since the latter is not of complex type.  Viewing these as complexifications
of real representations on $\ker \dot{\mathbf{D}}_\tau^{\boldsymbol{\theta}}$
and $\coker \dot{\mathbf{D}}_\tau^{\boldsymbol{\theta}}$ as explained above,
it follows via the correspondence between real and complex irreducible
representations outlined in \S\ref{sec:realComplex} that all the irreducible
real subrepresentations are isomorphic to~$\boldsymbol{\theta}_0$.

The main difference if $\boldsymbol{\theta}_0$ is of complex type is that
$Y \oplus iY \subset \CC^G$ is no longer an isotypic component 
for~$\boldsymbol{\rho}_\CC$, but is instead the direct sum of two isotypic
components related to each other by complex conjugation
$$
Y \oplus iY = X \oplus \widebar{X} \subset \CC^G,
$$
corresponding to some complex irreducible representation
$\boldsymbol{\lambda}_0 : G \to \Aut_\CC(V)$ and its non-isomorphic dual
$\boldsymbol{\lambda}^*_0 : G \to \Aut_\CC(V^*)$.  Writing 
$\boldsymbol{\lambda} : G \to \Aut_\CC(X)$ and
$\bar{\boldsymbol{\lambda}} : G \to \Aut_\CC(\widebar{X})$ for the restriction
of $\boldsymbol{\rho}_\CC$ to these subspaces, we can then think of
$\dot{\mathbf{D}}_\tau^{\boldsymbol{\lambda} \oplus \bar{\boldsymbol{\lambda}}} =
\dot{\mathbf{D}}_\tau^{\boldsymbol{\lambda}} \oplus
\dot{\mathbf{D}}_\tau^{\bar{\boldsymbol{\lambda}}}$ as the complexification
of $\dot{\mathbf{D}}_\tau^{\boldsymbol{\theta}}$.  A repeat of the argument
above using the isomorphisms
$X \cong V \otimes_i V^*$ and $\widebar{X} \cong V^* \otimes_i V$ then gives
an $i$-complex-linear isomorphism
$$
\ker \dot{\mathbf{D}}_\tau^{\boldsymbol{\lambda} \oplus \bar{\boldsymbol{\lambda}}} \cong
(\ker \dot{\mathbf{D}}_\tau^{\boldsymbol{\lambda}_0} \otimes_i V^*) \oplus
(\ker \dot{\mathbf{D}}_\tau^{\boldsymbol{\lambda}_0^*} \otimes_i V),
$$
with $G$ acting via $\1 \otimes \boldsymbol{\lambda}^*_0$ on the first summand
and $\1 \otimes \boldsymbol{\lambda}_0$ on the second, and a similar
isomorphism for cokernels.  It follows that every irreducible subrepresentation in
either $\ker \dot{\mathbf{D}}_\tau^{\boldsymbol{\lambda} \oplus \bar{\boldsymbol{\lambda}}}$
or $\coker \dot{\mathbf{D}}_\tau^{\boldsymbol{\lambda} \oplus \bar{\boldsymbol{\lambda}}}$
is isomorphic to one of $\boldsymbol{\lambda}_0$ or $\boldsymbol{\lambda}_0^*$,
and the desired result for real subrepresentations again follows via the
correspondence between real and complex representations
in \S\ref{sec:realComplex}.
\end{proof}

Continuing in the setting of Lemma~\ref{lemma:ndimension}, let
$\KK = \End_G(W) \in \{\RR,\CC,\HH\}$ and write $k = \dim_\KK \ker \dot{\mathbf{D}}_\tau^{\boldsymbol{\theta}_0}$,
$c = \dim_\KK \coker \dot{\mathbf{D}}_\tau^{\boldsymbol{\theta}_0}$.  By
Lemma~\ref{lemma:regular}, $\boldsymbol{\theta} \cong \boldsymbol{\theta}_0^{\oplus m}$
with $m := \dim_\KK W$, so Lemma~\ref{lemma:kernelDecomp} gives
$\dim_\KK \ker \dot{\mathbf{D}}_\tau^{\boldsymbol{\theta}} = km$ and
$\dim_\KK \coker \dot{\mathbf{D}}_\tau^{\boldsymbol{\theta}} = cm$.
Lemma~\ref{lemma:ndimension} meanwhile decomposes the representation defined
by the $G$-action on $\ker \dot{\mathbf{D}}_\tau^{\boldsymbol{\theta}}$ as
$\boldsymbol{\theta}_0^{\oplus \ell}$ for some $\ell \ge 0$, so
$\ker \dot{\mathbf{D}}_\tau^{\boldsymbol{\theta}} \cong W^{\oplus \ell}$.
Comparing dimensions, we deduce $\ell = k$, and applying the
same argument to the cokernel then likewise identifies the representation defined by the
$G$-action on $\coker \dot{\mathbf{D}}_\tau^{\boldsymbol{\theta}}$
with~$\boldsymbol{\theta}_0^{\oplus c}$.  The following consequence is the 
origin of the codimension formula in Theorem~\ref{thm:submanifolds0}
(cf.~\ref{eqn:Schur2}).

\begin{cor}
\label{cor:codimension}
In the setting of Lemma~\ref{lemma:ndimension}, let
$\KK = \End_G(W)$.  Then the space of $G$-equivariant
real-linear maps $\ker \dot{\mathbf{D}}_\tau^{\boldsymbol{\theta}} \to
\coker \dot{\mathbf{D}}_\tau^{\boldsymbol{\theta}}$ satisfies
$$
\dim_\RR \Hom_G\big(\ker \dot{\mathbf{D}}_\tau^{\boldsymbol{\theta}},
\coker \dot{\mathbf{D}}_\tau^{\boldsymbol{\theta}}\big) =
\dim_\RR \KK \cdot \dim_\KK \ker \dot{\mathbf{D}}_\tau^{\boldsymbol{\theta}_0}
\cdot \dim_\KK \coker \dot{\mathbf{D}}_\tau^{\boldsymbol{\theta}_0}.
$$
\qed
\end{cor}

\subsection{Setting up the implicit function theorem}
\label{sec:IFT}

We assume throughout this section that $(\Theta,\dot{\Sigma}'',\pi,G,\rho,I,f)$
is the \emph{minimal} regular presentation of $\varphi : \Sigma' \to \Sigma$.
Suppose
$$
\boldsymbol{\theta}_i : G \to \Aut_\RR(W_i), \qquad i=1,\ldots,N
$$
is a complete list of pairwise non-isomorphic real irreducible representations
for~$G$, with
$$
\KK_i := \End_G(W_i), \quad\text{ and }\quad t_i := \dim_\RR \KK_i \in \{1,2,4\}.
$$
Recall that all of the data we have been considering depends smoothly
on a parameter $\tau$, which lives in a connected Banach manifold~$P$ as described
at the end of \S\ref{sec:regPres}.
Any $N$-tuples of nonnegative integers
$\mathbf{k} = (k_1,\ldots,k_N)$ and $\mathbf{c} = (c_1,\ldots,c_N)$
now determine subsets of this parameter space
$$
P(\mathbf{k},\mathbf{c}) := \Big\{ \tau \in P\ \Big|\ 
\text{$\dim_{\KK_i} \ker \dot{\mathbf{D}}_\tau^{\boldsymbol{\theta}_i} = k_i$ and
$\dim_{\KK_i} \coker \dot{\mathbf{D}}_\tau^{\boldsymbol{\theta}_i} = c_i$ for
all $i=1,\ldots,N$} \Big\}.
$$
Note that $P(\mathbf{k},\mathbf{c})$ is automatically empty unless
$k_i - c_i = \ind_{\KK_i} \dot{\mathbf{D}}_\tau^{\boldsymbol{\theta}_i}$ for all
$i=1,\ldots,N$, and these indices do not depend on the parameter~$\tau$.
Assuming this condition holds,
we would now like to present $P(\mathbf{k},\mathbf{c})$ 
locally as the zero-set of a smooth map to a finite-dimensional vector space, 
and to compute its derivative in a special case.

We start by translating the conditions defining $P(\mathbf{k},\mathbf{c})$ into
conditions on the pulled back operators
$\widehat{\varphi}_\tau^*\dot{\mathbf{D}}_\tau$ for a suitable family of
\emph{regular} covers
$\widehat{\varphi}_\tau : \widehat{\Sigma}_\tau \to \dot{\Sigma}_\tau$
with $\Aut(\widehat{\varphi}_\tau) = G$.
This can be defined by replacing the homomorphism $\rho : G \to S(I)$ with the
action of $G$ on itself by left multiplication, i.e.~let
$$
\widehat{\rho} : G \to S(G) : g \mapsto \widehat{\rho}_g, \qquad \widehat{\rho}_g(h) := gh,
$$
so that $(\Theta_\tau,\dot{\Sigma}'',\pi_\tau,G,\widehat{\rho},G,\Id)$ becomes a
minimal regular presentation for
$$
\widehat{\Sigma}_\tau := \left( \dot{\Sigma}''_\tau \times G \right) \Big/ G
\stackrel{\widehat{\varphi}_\tau}{\longrightarrow} \dot{\Sigma}_\tau :
[(z,g)] \mapsto \pi_\tau(z),
$$
or rather for the extension of this map to a branched cover of closed surfaces
as provided by
Lemma~\ref{lemma:closedUp}.  In keeping with our usual notational convention,
$\widehat{\Sigma}_\tau$ is a fixed smooth surface $\widehat{\Sigma}$ with a
fixed $G$-action by deck transformations but
a $\tau$-dependent family of conformal structures $\widehat{\jmath}_\tau =
\widehat{\varphi}_\tau^*j_\tau$, which are fixed on the cylindrical ends.

Denote the isotypic decomposition of the regular representation 
$\widehat{\boldsymbol{\rho}} : G \to \Aut_\RR(\RR^G)$ by
$$
\widehat{\boldsymbol{\rho}} = \bigoplus_{i=1}^N \widehat{\boldsymbol{\theta}}_i,
$$
where $\widehat{\boldsymbol{\theta}}_i \cong \boldsymbol{\theta}_i^{\oplus \ell_i}$
for integers $\ell_i$ which are strictly positive by
Lemma~\ref{lemma:regular}.  Then by Lemma~\ref{lemma:kernelDecomp},
\begin{equation*}
\begin{split}
\ker (\widehat{\varphi}_\tau^*\dot{\mathbf{D}}_\tau)  & \cong
\bigoplus_{i=1}^N \ker \dot{\mathbf{D}}_\tau^{\widehat{\boldsymbol{\theta}}_i}
\cong \bigoplus_{i=1}^N \left(\ker \dot{\mathbf{D}}_\tau^{\boldsymbol{\theta}_i}\right)^{\oplus \ell_i}, \\
\coker (\widehat{\varphi}_\tau^*\dot{\mathbf{D}}_\tau)  & \cong
\bigoplus_{i=1}^N \coker \dot{\mathbf{D}}_\tau^{\widehat{\boldsymbol{\theta}}_i}
\cong \bigoplus_{i=1}^N \left(\coker \dot{\mathbf{D}}_\tau^{\boldsymbol{\theta}_i}\right)^{\oplus \ell_i},
\end{split}
\end{equation*}
so $\tau \in P(\mathbf{k},\mathbf{c})$ implies
\begin{equation}
\label{eqn:kernelCondition}
\dim \ker (\widehat{\varphi}_\tau^*\dot{\mathbf{D}}_\tau) =
\sum_{i=1}^N t_i \ell_i k_i .
\end{equation}

\begin{lemma}
\label{lemma:localV}
Every $\sigma \in P(\mathbf{k},\mathbf{c})$ has a neighborhood
$\uU_\sigma \subset P$ such that $\uU_\sigma \cap P(\mathbf{k},\mathbf{c})$
is the set of
all $\tau \in \uU_\sigma$ for which \eqref{eqn:kernelCondition} holds.
\end{lemma}
\begin{proof}
Since all the operators $\dot{\mathbf{D}}_\tau^{\boldsymbol{\theta}_i}$ are Fredholm and
they depend continuously on~$\tau$, we can assume
$\dim \ker \dot{\mathbf{D}}_\tau^{\boldsymbol{\theta}_i} \le
\dim \ker \dot{\mathbf{D}}_{\sigma}^{\boldsymbol{\theta}_i}$ for all
$i=1,\ldots,N$ if $\tau$ is sufficiently close to~$\sigma$.
Thus \eqref{eqn:kernelCondition} can only be satisfied if none of these inequalities are
strict, which means $\tau \in P(\mathbf{k},\mathbf{c})$ since every
$\ell_i$ is positive.
\end{proof}

Recall from \S\ref{sec:CRpunctured} that the weighted Sobolev spaces
$W^{k,p,-\widehat{\varphi}^*\boldsymbol{\delta}}(\widehat{\varphi}_\tau^*\dot{E}_\tau)$
and $W^{k-1,p,-\widehat{\varphi}^*\boldsymbol{\delta}}(\widehat{\varphi}_\tau^*\dot{F}_\tau)$ are
defined in terms of fixed families of trivializations of $E_\tau$ near $\Theta_\tau$
and holomorphic cylindrical coordinates
which allow us to compute Sobolev norms on the cylindrical ends.
Given $\sigma \in P(\mathbf{k},\mathbf{c})$, choose a neighborhood
$\uU_{\sigma} \subset P$ that is diffeomorphic to a ball and
small enough to satisfy Lemma~\ref{lemma:localV}.  By assumption the bundles
$E_\tau$ depend smoothly on~$\tau$, which means there is a well-defined smooth
bundle $\widehat{E} \to P \times \Sigma$ with $\widehat{E}_{(\tau,z)} = (E_\tau)_z$.
Choosing a suitable connection on the latter, we can use parallel transport along
paths of the form $(\tau(t),\psi_{\tau(t)}(z)) \in \uU_\sigma \times \Sigma$ with $\tau(t)$ radiating
outward from~$\sigma$ to define a smooth family of
complex bundle isomorphisms
$$
\Psi_\tau : \psi_{\sigma}^*E_{\sigma} \to \psi_\tau^*E_\tau
$$
which respect these fixed trivializations near~$\Theta_\tau$ and satisfy
$\Psi_{\sigma} = \Id$.  These give rise to isomorphisms
$\dot{E}_{\sigma} \to \dot{E}_\tau$ covering the diffeomorphisms $\psi_\tau \circ \psi_{\sigma}^{-1} :
\dot{\Sigma}_{\sigma} \to \dot{\Sigma}_\tau$.
Notice that there are also natural real bundle isomorphisms
$$
d\psi_\tau : T\Sigma \to \psi_\tau^*T\Sigma,
$$
so that $d\psi_\tau \circ d\psi_{\sigma}^{-1}$ gives a family of isomorphisms
$T\dot{\Sigma}_{\sigma} \to T\dot{\Sigma}_\tau$
covering $\dot{\Sigma}_{\sigma} \stackrel{\psi_\tau \circ \psi_{\sigma}^{-1}}{\longrightarrow} \dot{\Sigma}_\tau$, 
and they respect the chosen holomorphic 
cylindrical coordinates on the ends.  These then induce 
smooth families of isomorphisms of complex bundles over~$\widehat{\Sigma}$,
$$
\widehat{\varphi}_{\sigma}^*\dot{E}_{\sigma} \to \widehat{\varphi}_\tau^*\dot{E}_\tau, \qquad
\widehat{\varphi}_{\sigma}^*\dot{F}_{\sigma} \to \widehat{\varphi}_\tau^*\dot{F}_\tau
$$
which again are the identity for $\tau=\sigma$ and are also
equivariant with respect to the natural $G$-action by bundle
isomorphisms covering deck transformations of~$\widehat{\Sigma}$.
Acting with these on sections produces
$\tau$-parametrized families of $G$-equivariant 
Banach space isomorphisms which we
shall also denote by~$\Psi_\tau$:
\begin{equation}
\label{eqn:Psi}
\begin{split}
W^{k,p,-\widehat{\varphi}^*\boldsymbol{\delta}}(\widehat{\varphi}_{\sigma}^*\dot{E}_{\sigma})
&\stackrel{\Psi_{\tau}}{\longrightarrow}
W^{k,p,-\widehat{\varphi}^*\boldsymbol{\delta}}(\widehat{\varphi}_{\tau}^*\dot{E}_{\tau}), \\
W^{k-1,p,-\widehat{\varphi}^*\boldsymbol{\delta}}(\widehat{\varphi}_{\sigma}^*\dot{F}_{\sigma})
&\stackrel{\Psi_{\tau}}{\longrightarrow}
W^{k-1,p,-\widehat{\varphi}^*\boldsymbol{\delta}}(\widehat{\varphi}_{\tau}^*\dot{F}_{\tau}).
\end{split}
\end{equation}
Here $\Psi_{\sigma} = \Id$.

We can now use these isomorphisms to define for $\tau \in \uU_{\sigma}$ a smooth 
family of $G$-equivariant Fredholm operators with fixed domain and target space,
\begin{equation}
\label{eqn:Dhat}
\widehat{\mathbf{D}}_{\tau} := \Psi_{\tau}^{-1} \circ 
\widehat{\varphi}_{\tau}^*\dot{\mathbf{D}}_{\tau} \circ
\Psi_{\tau} :
W^{k,p,-\widehat{\varphi}^*\boldsymbol{\delta}}(\widehat{\varphi}_{\sigma}^*\dot{E}_{\sigma})
\to 
W^{k-1,p,-\widehat{\varphi}^*\boldsymbol{\delta}}(\widehat{\varphi}_{\sigma}^*\dot{F}_{\sigma}),
\end{equation}
such that 
$$
\uU_{\sigma} \cap P(\mathbf{k},\mathbf{c}) = \left\{ \tau \in \uU_{\sigma}\ \Bigg|\ 
\dim \ker \widehat{\mathbf{D}}_{\tau} = \sum_{i=1}^N t_i \ell_i k_i \right\}.
$$
In order to present the latter as the zero-set of a smooth map, let us abbreviate
$$
\mathbf{X}_\sigma := W^{k,p,-\widehat{\varphi}^*\boldsymbol{\delta}}(\widehat{\varphi}_{\sigma}^*\dot{E}_{\sigma}), \qquad
\mathbf{Y}_\sigma := W^{k-1,p,-\widehat{\varphi}^*\boldsymbol{\delta}}(\widehat{\varphi}_{\sigma}^*\dot{F}_{\sigma}),
$$
so \eqref{eqn:Dhat} defines a smooth map
$$
\uU_{\sigma} \to \Lin_G(\mathbf{X}_\sigma,\mathbf{Y}_\sigma) : \tau \mapsto \widehat{\mathbf{D}}_{\tau},
$$
where $\Lin_G(\mathbf{X}_\sigma,\mathbf{Y}_\sigma)$ denotes the Banach space of bounded real-linear maps
$\mathbf{X}_\sigma \to \mathbf{Y}_\sigma$ that are $G$-equivariant. 
Since $\widehat{\mathbf{D}}_{{\sigma}} = \widehat{\varphi}_{\sigma}^*\dot{\mathbf{D}}_{\sigma}$
is Fredholm, we can choose a splitting
$$
\mathbf{X}_\sigma = \mathbf{V}_\sigma \oplus \ker (\widehat{\varphi}_{\sigma}^*\dot{\mathbf{D}}_{\sigma}),
$$
such that $\mathbf{V}_\sigma \subset \mathbf{X}_\sigma$ is a closed subspace and $\widehat{\mathbf{D}}_{\sigma}$
maps $\mathbf{V}_\sigma$ isomorphically to its image.  By Proposition~\ref{prop:formalAdjoint},
we can similarly split
$$
\mathbf{Y}_\sigma = \im (\widehat{\varphi}_{\sigma}^*\dot{\mathbf{D}}_{\sigma}) \oplus
\ker (\widehat{\varphi}_{\sigma}^*\dot{\mathbf{D}}_{\sigma}^*),
$$
where $\ker (\widehat{\varphi}_{\sigma}^*\dot{\mathbf{D}}_{\sigma}^*)$ is equivalently
the space of all sections in $W^{k-1,p,\widehat{\varphi}^*\boldsymbol{\delta}}(\widehat{\varphi}_{\sigma}^*\dot{F}_{\sigma})$
that are $L^2$-orthogonal to $\im (\widehat{\varphi}_{\sigma}^*\dot{\mathbf{D}}_{\sigma})$.
In terms of these splittings, $\widehat{\mathbf{D}}_{\tau}$ can be written
in block form
$$
\widehat{\mathbf{D}}_{\tau} = \begin{pmatrix}
\mathbf{D}_{\tau}^{11} & \mathbf{D}_{\tau}^{12} \\
\mathbf{D}_{\tau}^{21} & \mathbf{D}_{\tau}^{22}
\end{pmatrix},
$$
where after shrinking $\uU_{\sigma}$ if necessary, we can assume without loss of
generality that $\mathbf{D}_{\tau}^{11} : \mathbf{V}_\sigma \to 
\im (\widehat{\varphi}_{\sigma}^*\dot{\mathbf{D}}_{\sigma})$
is invertible for all $\tau \in \uU_{\sigma}$.  We can therefore define a map
\begin{equation}
\label{eqn:F}
\begin{split}
\mathbf{F}_\sigma  : \uU_{\sigma} &\to 
\Hom_G\big(\ker (\widehat{\varphi}_{\sigma}^*\dot{\mathbf{D}}_{\sigma}),
\ker (\widehat{\varphi}_{\sigma}^*\dot{\mathbf{D}}_{\sigma}^*)\big) \\
\tau &\mapsto \mathbf{D}_{\tau}^{22} - \mathbf{D}_{\tau}^{21} (\mathbf{D}_{\tau}^{11})^{-1} \mathbf{D}_{\tau}^{12}.
\end{split}
\end{equation}

\begin{lemma}
\label{lemma:zeroset}
A parameter $\tau \in \uU_\sigma$ belongs to $P(\mathbf{k},\mathbf{c})$
if and only if $\mathbf{F}_\sigma(\tau) = 0$.
\end{lemma}
\begin{proof}
Define for each $\tau \in \uU_\sigma$ the Banach space isomorphism
$$
\mathbf{T} = \begin{pmatrix}
\1 & -(\mathbf{D}_\tau^{11})^{-1}\mathbf{D}_\tau^{12} \\
0 & \1
\end{pmatrix} \in
\Lin(\mathbf{V}_\sigma \oplus \ker (\widehat{\varphi}_{\sigma}^*\dot{\mathbf{D}}_{\sigma})) =
\Lin(\mathbf{X}_\sigma).
$$
Then $\widehat{\mathbf{D}}_\tau \mathbf{T} =
\begin{pmatrix} \mathbf{D}_\tau^{11} & 0 \\ \mathbf{D}_\tau^{21} & \mathbf{F}_\sigma(\tau) \end{pmatrix}$,
and since $\mathbf{D}_\tau^{11}$ is invertible, 
$$
\ker \widehat{\mathbf{D}}_\tau \cong 
\ker (\widehat{\mathbf{D}}_\tau \mathbf{T})
= \{0\} \oplus \ker \mathbf{F}_\sigma(\tau) \cong \ker \mathbf{F}_\sigma(\tau).
$$
The latter can only have the same dimension as 
$\ker (\widehat{\varphi}_{\sigma}^*\dot{\mathbf{D}}_{\sigma})$ if
$\mathbf{F}_\sigma(\tau)$ vanishes.
\end{proof}

Observe that by Lemma~\ref{lemma:ndimension},
Corollary~\ref{cor:codimension} and Schur's lemma,
\begin{equation}
\label{eqn:Schur2}
\dim \Hom_G\big(\ker (\widehat{\varphi}_{\sigma}^*\dot{\mathbf{D}}_{\sigma}),
\ker (\widehat{\varphi}_{\sigma}^*\dot{\mathbf{D}}_{\sigma}^*)\big) =
\sum_{i=1}^N t_i k_i c_i.
\end{equation}
The lemma implies via the implicit function theorem that a neighborhood
of $\sigma$ in $P(\mathbf{k},\mathbf{c})$ is a smooth submanifold with
the same codimension that appears in Theorem~\ref{thm:submanifolds0} whenever
we can show that the linearization
$$
d\mathbf{F}_\sigma(\sigma) : T_\sigma P \to 
\Hom_G\big(\ker (\widehat{\varphi}_{\sigma}^*\dot{\mathbf{D}}_{\sigma}),
\ker (\widehat{\varphi}_{\sigma}^*\dot{\mathbf{D}}_{\sigma}^*)\big)
$$
is surjective.

We will need a precise formula for this linearization in the following special case.
Suppose we have a smooth path
$$
\gamma : (-\epsilon,\epsilon) \to P \qquad \text{ with 
$\gamma(0) = \sigma$ and $\dot{\gamma}(0) = Y \in T_\sigma P$}
$$
such that for all $\tau = \gamma(t)$:
\begin{enumerate}
\item $E_\tau = E_\sigma$ (i.e.~there is a canonical complex bundle isomorphism);
\item $\psi_\tau = \Id$;
\item $j_\tau = j_\sigma$.
\end{enumerate}
We are then free to choose the bundle isomorphisms $\Psi_\tau$ and consequently
the Banach space isomorphisms \eqref{eqn:Psi} to be the identity for
all $\tau = \gamma(t)$, so $\widehat{\mathbf{D}}_{\gamma(t)} = \widehat{\varphi}_\sigma^*\dot{\mathbf{D}}_{\gamma(t)}$,
where $\mathbf{D}_{\gamma(t)}$ is a smooth
family of Cauchy-Riemann operators on the fixed bundle
$E_\sigma \to \Sigma_\sigma$.
Differentiating this family gives a real-linear bundle map
$$
A_Y := \left. \p_t \mathbf{D}_{\gamma(t)}\right|_{t=0} \in
\Gamma(\Hom_\RR(E_\sigma,F_\sigma)),
$$
and we then find that
$$
\mathbf{L}(Y) := d\mathbf{F}_\sigma(\sigma) Y \in
\Hom_G\big(\ker (\widehat{\varphi}_{\sigma}^*\dot{\mathbf{D}}_{\sigma}),
\ker (\widehat{\varphi}_{\sigma}^*\dot{\mathbf{D}}_{\sigma}^*)\big)
$$
takes the form
\begin{equation}
\label{eqn:linearization}
\begin{split}
\mathbf{L}(Y) \eta = \pi\big( (\widehat{\varphi}_\sigma^*A_Y) \eta \big),
\end{split}
\end{equation}
where $\pi$ is the projection
$$
\mathbf{Y}_\sigma = \im (\widehat{\varphi}_{\sigma}^*\dot{\mathbf{D}}_{\sigma}) \oplus
\ker (\widehat{\varphi}_{\sigma}^*\dot{\mathbf{D}}_{\sigma}^*) 
\stackrel{\pi}{\longrightarrow} \ker (\widehat{\varphi}_{\sigma}^*\dot{\mathbf{D}}_{\sigma}^*).
$$
The local genericity result developed in \S\ref{sec:unique} below is
geared toward proving that operators such as $\mathbf{L}$ are surjective.

\section{Index computation}
\label{sec:index}

The goal of this section is to compute the Fredholm index of the twisted
Cauchy-Riemann type operators introduced in \S\ref{sec:twisted}.
We will use the notation of \S\ref{sec:prep} but dispense with the
parameter~$\tau$ since it is not important for the index computation,
hence $\varphi : (\Sigma',j') \to (\Sigma,j)$ is a fixed branched cover, and
$(\Theta,\dot{\Sigma}'',\pi,G,\rho,I,f)$ is a fixed regular presentation.
The complex vector bundles $E$ and $F$ with their restrictions $\dot{E}$ and
$\dot{F}$ to the punctured domain $\dot{\Sigma}$ are assumed to have rank
$$
m := \rank_\CC E \in \NN,
$$
and we assume
$$
\boldsymbol{\theta} : G \to \Aut_\RR(W)
$$
is a (not necessarily irreducible or faithful) representation of $G$ with
$$
n := \dim W \in \NN.
$$
The resulting twisted bundles over $\dot{\Sigma}$ can be written as
$$
\dot{E}^{\boldsymbol{\theta}} = \dot{E} \otimes_\RR W^{\boldsymbol{\theta}}, \qquad
\dot{F}^{\boldsymbol{\theta}} = \dot{F} \otimes_\RR W^{\boldsymbol{\theta}},
$$
in terms of the flat real vector bundle $W^{\boldsymbol{\theta}} :=
(\dot{\Sigma}'' \times W) / G \to \dot{\Sigma}$, and any Cauchy-Riemann type
operator $\mathbf{D} \in \CRR(E)$ then gives rise to the twisted operator
$$
\dot{\mathbf{D}}^{\boldsymbol{\theta}} : \Gamma(\dot{E}^{\boldsymbol{\theta}})
\to \Gamma(\dot{F}^{\boldsymbol{\theta}}).
$$
We need a bit more notation in order to state a formula for 
$\ind(\dot{\mathbf{D}}^{\boldsymbol{\theta}})$.  Recall that while the
deck transformations $G = \Aut(\pi)$ act on $\dot{\Sigma}''$ without fixed
points, their extensions to biholomorphic self-maps of $\Sigma''$ may fix
some of the punctures, so for each $w \in \Theta$ and
$\zeta \in \pi^{-1}(w) \subset \Theta'' := \pi^{-1}(\Theta)$, 
we can consider the stabilizer subgroup
$$
G_\zeta := \left\{ g \in G\ |\ g\zeta = \zeta \right\},
$$
which is necessarily cyclic.
Restricting $\boldsymbol{\theta}$ to $G_\zeta$ then defines a representation
$G_\zeta \to \Aut_\RR(W)$, which splits $W$ into $G_\zeta$-invariant subspaces
$W = W_\zeta \oplus W_\zeta'$ such that $G_\zeta$ acts on $W_\zeta$
trivially and on $W_\zeta'$ as a direct sum of nontrivial representations.
We define the number
$$
n_w := \dim W_\zeta' \in \{0,\ldots,n\}.
$$
As implied by the notation, this depends on $w \in \Theta$ but not on the
choice of preimage $\zeta \in \pi^{-1}(w)$: indeed, since $G$ acts transitively
on $\pi^{-1}(w)$, any two choices of $\zeta$ give rise to conjugate 
subgroups $G_\zeta$, and using orthonormality of characters, one can compute
$$
n_w = n - \dim W_\zeta = n - \frac{1}{|G_\zeta|} \sum_{g \in G_\zeta} \chi_{\boldsymbol{\theta}}(g),
$$
an expression which depends only on the congugacy class of~$G_\zeta$.

\begin{thm}
\label{thm:index}
Under the assumptions detailed above, the operator
$$
\dot{\mathbf{D}}^{\boldsymbol{\theta}} : W^{k,p,-\boldsymbol{\delta}}(\dot{E}^{\boldsymbol{\theta}})
\to W^{k-1,p,-\boldsymbol{\delta}}(\dot{F}^{\boldsymbol{\theta}})
$$
is Fredholm for any $k \in \NN$, $p \in (1,\infty)$ and negative exponential weights
$-\boldsymbol{\delta} = \{-\delta_w\}_{w \in \Theta}$ satisfying $0 < \delta_w < 2\pi/ |G|$
for all $w \in \Theta$.  Its index is
$$
\ind(\dot{\mathbf{D}}^{\boldsymbol{\theta}}) = n \cdot \ind(\mathbf{D}) - m \sum_{w \in \Theta} n_w.
$$
\end{thm}

The dimensions and indices in the above statement are all real, but note that if
$\boldsymbol{\theta}$ is irreducible with $\KK := \End_G(W) \in \{\CC,\HH\}$,
then the integers $n$ and
$n_w$ are automatically divisible by $t := \dim_\RR \KK \in \{2,4\}$, hence so 
is~$\ind(\dot{\mathbf{D}}^{\boldsymbol{\theta}})$.  Let us state the corollary
for the faithful case in terms of the $\KK$-linear index since it is most useful
in this form.

\begin{cor}[cf.~Lemma~\ref{lemma:index}]
\label{cor:faithful}
Assume $(\Theta,\dot{\Sigma}'',\pi,G,\rho,I,f)$ is the mimimal regular presentation,
and that $\boldsymbol{\theta}$ is faithful and irreducible with
$\End_G(W) \cong \KK \in \{\RR,\CC,\HH\}$.  Then
$$
\ind_\KK(\dot{\mathbf{D}}^{\boldsymbol{\theta}}) \le (\dim_\KK W) \cdot \ind_\RR(\mathbf{D}) - 
m |\Theta|,
$$
and if $\KK = \RR$, then
the inequality is strict unless all branch points of $\varphi$ have 
branching order~$2$.
\end{cor}
\begin{proof}
By Lemma~\ref{lemma:branchPoints}, the stabilizer subgroups $G_\zeta$ are
nontrivial for all $\zeta \in \Theta''$, and the conclusion about branch
points of order~$2$ will hold if and only if all of them are isomorphic to~$\ZZ_2$.
Now if $\boldsymbol{\theta}$ is faithful, it follows that all nontrivial
elements $g \in G_\zeta$ for $\zeta \in \Theta''$ also act nontrivially
on~$W$, hence the decomposition of $W$ into $G_\zeta$-invariant subspaces
contains at least a $1$-dimensional $\KK$-linear subspace on which $G_\zeta$
acts nontrivially, giving $n_w \ge \dim_\RR \KK$ for all $w \in \Theta$.
This implies the upper
bound, and in the case $\KK=\RR$,
it is an equality if and only if $n_w = 1$ for all $w \in \Theta$,
meaning each $G_\zeta$ acts on $W$ as the $(n-1)$-fold direct sum of the
trivial representation plus a real $1$-dimensional nontrivial representation,
which is required to be faithful.  But the only nontrivial faithful real
$1$-dimensional representation of any finite group is the nontrivial
representation of~$\ZZ_2$, hence $G_\zeta \cong \ZZ_2$.
\end{proof}

\begin{remark}
Doan and Walpuski have recently shown
that an index formula equivalent to that of Theorem~\ref{thm:index} can also be derived
from Kawasaki's orbifold Riemann-Roch theorem \cite{Kawasaki:RiemannRoch}.
From this perspective, branch points are regarded as orbifold
singularities instead of punctures; see
\cite{DoanWalpuski:BrillNoether}*{Appendix~2.B}.
\end{remark}

The remainder of this section is devoted to the proof of Theorem~\ref{thm:index},
which we shall break down into five steps.

\textbf{Step~1: Some notation}.\\
It will be convenient first to
complexify the representation.  We define $V := W \oplus iW$ and
the complex representation
$$
\boldsymbol{\lambda} : G \to \Aut_\CC(V)
$$
such that $\boldsymbol{\lambda}(g)|_W = \boldsymbol{\theta}(g)$ for all $g \in G$.
Note that for $w \in \Theta$ and $\zeta \in \pi^{-1}(w) \subset \Theta''$,
the trivial representation of $G_\zeta$ on $V$ is the complexification of
the trivial real representation on~$W$, so the splitting $W = W_\zeta \oplus
W_\zeta'$ explained above complexifies to a splitting $V = V_\zeta \oplus V_\zeta'$,
where $V_\zeta \subset V$ is the largest complex subspace on which $G_\zeta$
acts trivially, allowing us to write
$$
n_w = \dim_\CC V_\zeta' = n - \dim_\CC V_\zeta.
$$
The complexified representation now gives rise
to a complex flat bundle $V^{\boldsymbol{\lambda}} := (\dot{\Sigma}'' \times V) / G$,
corresponding twisted bundles
\begin{equation}
\label{eqn:EFcomplex}
\dot{E}^{\boldsymbol{\lambda}} := \dot{E} \otimes_\RR V^{\boldsymbol{\lambda}}, \qquad
\dot{F}^{\boldsymbol{\lambda}} := \dot{F} \otimes_\RR V^{\boldsymbol{\lambda}},
\end{equation}
and a twisted Cauchy-Riemann operator
$$
\dot{\mathbf{D}}^{\boldsymbol{\lambda}} : W^{k,p,-\boldsymbol{\delta}}(\dot{E}^{\boldsymbol{\lambda}})
\to W^{k-1,p,-\boldsymbol{\delta}}(\dot{F}^{\boldsymbol{\lambda}}).
$$
The following point is important to understand: the tensor products in \eqref{eqn:EFcomplex} are \emph{real},
thus $\dot{E}^{\boldsymbol{\lambda}}$ and $\dot{F}^{\boldsymbol{\lambda}}$
each inherit two complex structures $J$ and~$i$, where $J$ comes from the
complex structure of $E$ and $i$ from that of~$V$: they commute with
each other and are defined by
$$
J(\eta \otimes v) := J\eta \otimes v, \qquad i(\eta \otimes v) := \eta \otimes iv.
$$
In this sense, $\dot{\mathbf{D}}^{\boldsymbol{\lambda}}$ can be regarded as
the $i$-complex-linear extension of $\dot{\mathbf{D}}^{\boldsymbol{\theta}}$
to complexifications of the latter's domain and target space---this notion
of ``complexification'' ignores the fact that these spaces already have native complex
structures $J$ and treats them as \emph{real} vector spaces, which is appropriate
since $\dot{\mathbf{D}}^{\boldsymbol{\theta}}$ need not be $J$-complex linear.
We therefore obtain the relation
$$
\ind(\dot{\mathbf{D}}^{\boldsymbol{\theta}}) = \frac{1}{2}
\ind(\dot{\mathbf{D}}^{\boldsymbol{\lambda}}),
$$
and we shall compute $\ind(\dot{\mathbf{D}}^{\boldsymbol{\lambda}})$ by regarding
$\dot{\mathbf{D}}^{\boldsymbol{\lambda}}$ as a real-linear Cauchy-Riemann
type operator on the complex vector bundle $(\dot{E}^{\boldsymbol{\lambda}},J)$.
Since $\rank_\CC \dot{E}^{\boldsymbol{\lambda}} = \rank_\CC E \cdot \dim_\RR V = 2mn$,
the punctured Riemann-Roch formula from \cite{Schwarz}*{\S 3.3}
(or equivalently \cite{Wendl:SFT}*{Lecture~5}) gives
\begin{equation}
\label{eqn:RiemannRochReal}
\ind(\dot{\mathbf{D}}^{\boldsymbol{\lambda}}) = 2mn \cdot \chi(\dot{\Sigma}) + 
2 c_1^\Phi(\dot{E}^{\boldsymbol{\lambda}},J) + \sum_{w \in \Theta}
\muCZ^\Phi(\mathbf{A}^{\boldsymbol{\lambda}}_w - \delta_w),
\end{equation}
where $\Phi$ is an arbitrary choice of asymptotic trivialization,
and $\muCZ^\Phi(\mathbf{A}^{\boldsymbol{\lambda}}_w - \delta_w) \in \ZZ$ are Conley-Zehnder 
indices  that depend on certain asymptotic operators $\mathbf{A}^{\boldsymbol{\lambda}}_w$
to be discussed below and the exponential weight
$-\delta_w  \in (-2\pi/|G|,0)$ associated to each puncture $w \in \Theta$.
The main difficulty of the calculation is in choosing a suitable
asymptotic trivialization in which both
$c_1^\Phi(\dot{E}^{\boldsymbol{\lambda}},J)$ and 
$\muCZ^\Phi(\mathbf{A}^{\boldsymbol{\lambda}}_w - \delta_w)$ 
can be computed.

Denote
$$
d' := \deg(\pi) = |G|,
$$
and for each $w \in \Theta$ and $\zeta \in \pi^{-1}(w) \subset \Theta''$,
let
$$
k_\zeta \in \{1,\ldots,d'\}
$$
denote the branching order of $\pi$ at~$\zeta$, meaning $\pi$ is a
$k_\zeta$-to-$1$ map on a small punctured neighborhood of~$\zeta$.
We can then choose punctured neighborhoods $\uU_w \subset \dot{\Sigma}$ and 
$\uU_\zeta \subset \dot{\Sigma}''$ of 
$w$ and $\zeta$ respectively, with
holomorphic cylindrical coordinates $(s,t) \in [0,\infty) \times S^1$ on each
such that
$$
\pi(s,t) = (k_\zeta s, k_\zeta t)
$$
in coordinates on~$\uU_\zeta$.  In these coordinates, any $g \in G_\zeta$
necessarily preserves the end $\uU_\zeta$ and takes the form
$g(s,t) = (s , t + j/k_\zeta)$ for some $j \in \{0,\ldots,k_\zeta-1\}$.
This means that $G_\zeta$ is a cyclic group of order $k_\zeta$, and it
has a canonical generator $g_\zeta \in G_\zeta$ such that
$$
g_\zeta(s,t) = (s , t + 1/k_\zeta) \quad\text{ on $\uU_\zeta$}.
$$

In addition to the cylindrical coordinates, let us choose complex
trivializations of $E$ on each of the corresponding neighborhoods
of $\Theta$, thus giving an identification
\begin{equation}
\label{eqn:dotEend}
\dot{E}|_{\uU_w} = \left([0,\infty) \times S^1\right) \times E_w
\end{equation}
for each $w \in \Theta$.  For any choice $\zeta \in \pi^{-1}(w)
\subset \Theta''$, this also gives us an identification of
$\dot{E}^{\boldsymbol{\lambda}}|_{\uU_w}$ with
\begin{equation}
\label{eqn:dotEtensorEnd}
\left( \left( [0,\infty) \times S^1 \right) \times (E_w \otimes_\RR V ) 
\right)\Big/ G_\zeta,
\end{equation}
where the action of $G_\zeta = \ZZ_{k_\zeta}$ on
$\left( [0,\infty) \times S^1 \right) \times (E_w \otimes_\RR V )$
is determined by
$$
g_\zeta \cdot \big( (s,t) , \eta \otimes v \big) =
\big( (s , t + 1 / k_\zeta) , \eta \otimes \boldsymbol{\lambda}(g_\zeta) v \big).
$$
This picture can now easily be extended to the ``circle compactification''
of the punctured surface:
let $\widebar{\Sigma}$ and $\widebar{\Sigma}''$ denote the compact
surfaces with boundary obtained by replacing each cylindrical end
$[0,\infty) \times S^1$ in $\dot{\Sigma}$ and $\dot{\Sigma}''$ respectively
by the compact topological manifold
$[0,\infty] \times S^1$.  The connected components of $\p \widebar{\Sigma}$
and $\p \widebar{\Sigma}''$ are then in bijective correspondence with
the punctures $w \in \Theta$ or $\zeta \in \Theta''$ respectively, and 
the choice of cylindrical coordinates identifies each of these components
with~$S^1$.  We shall denote the boundary components accordingly by
$S^1_w, S^1_\zeta$ for $w \in \Theta$ or $\zeta \in \Theta''$, hence
$$
\p\widebar{\Sigma} = \bigsqcup_{w \in \Theta} S^1_w, \qquad
\p\widebar{\Sigma}'' = \bigsqcup_{\zeta \in \Theta''} S^1_\zeta.
$$
The covering map $\pi : \dot{\Sigma}'' \to \dot{\Sigma}$ now extends to a
continuous covering map
$$
\bar{\pi} : \widebar{\Sigma}'' \to \widebar{\Sigma}
$$
which restricts on the boundary components to
$$
\pi_\zeta := \bar{\pi}|_{S^1_\zeta} : S^1_\zeta \to S^1_{\pi(\zeta)} :
t \mapsto k_\zeta t,
$$
and each $g \in G$ also extends naturally to
a continuous deck transformation 
$\bar{g} : \widebar{\Sigma}'' \to \widebar{\Sigma}''$
of~$\bar{\pi}$, such that if $g(\zeta) = \zeta'$, then
$\bar{g}$ maps $S^1_\zeta \to S^1_{\zeta'}$ via the canonical 
diffeomorphism composed with a translation.  The identifications
\eqref{eqn:dotEend} and \eqref{eqn:dotEtensorEnd} then yield obvious extensions
of $\dot{E}$ and $\dot{E}^{\boldsymbol{\lambda}}$ as topological vector bundles
$$
\widebar{E} \to \widebar{\Sigma}, \qquad \widebar{E}^{\boldsymbol{\lambda}} \to \widebar{\Sigma},
$$
and we have
$$
\widebar{E}^{\boldsymbol{\lambda}} = \left( \bar{\pi}^*\widebar{E} \otimes_\RR V \right) \Big/ G.
$$

\textbf{Step~2: Asymptotic operators on the twisted bundle}.\\
With the essential notation in place, we can now
discuss asymptotic operators.  Recall that after choosing a suitable Hermitian 
inner product on $\dot{E}$ over the cylindrical ends, any Cauchy-Riemann type operator
$\dot{\mathbf{D}}$ on $\dot{E} \to \dot{\Sigma}$ with reasonable asymptotic behavior 
determines real-linear operators
$$
\mathbf{A}_w : \Gamma(\widebar{E}|_{S^1_w}) \to \Gamma(\widebar{E}|_{S^1_w}),
$$
for each $w \in \Theta$, see e.g.~\cite{Wendl:automatic}*{\S 2.1}.  These can 
be regarded as unbounded self-adjoint operators on $L^2(\widebar{E}|_{S^1_w})$ with
dense domain $H^1(\widebar{E}|_{S^1_w})$, and we say
$\mathbf{A}_w$ is \defin{nondegenerate} whenever its kernel is trivial,
in which case it determines a \defin{Conley-Zehnder index}
$$
\muCZ^\Phi(\mathbf{A}_w) \in \ZZ
$$
relative to any choice of complex trivialization $\Phi$ of $\widebar{E}|_{S^1_w}$.
In the case where $\dot{\mathbf{D}}$ is the restriction to $\dot{\Sigma}$ 
of some operator $\mathbf{D} \in \CRR(E)$ on $\Sigma$, the operators
$\mathbf{A}_w$ are very simple and were already computed in
\S\ref{sec:CRpunctured}: they are each the so-called \emph{trivial}
asymptotic operator
$$
\mathbf{A}_w = - J \p_t,
$$
where $\p_t$ is a well-defined differential operator on $\widebar{E}|_{S^1_w}$
since the fibers are all canonically identified with~$E_w$.
This operator is degenerate, but the introduction of negative
exponential weights $- \delta_w < 0$ identifies $\dot{\mathbf{D}}$ with
another Cauchy-Riemann type operator whose corresponding asymptotic
operators are $\mathbf{A}_w - \delta_w$, which are nondegenerate
for any $\delta_w > 0$ sufficiently small.

Denote by
$$
\mathbf{A}^{\boldsymbol{\lambda}}_w : \Gamma(\widebar{E}^{\boldsymbol{\lambda}}|_{S^1_w}) \to
\Gamma(\widebar{E}^{\boldsymbol{\lambda}}|_{S^1_w})
$$
the asymptotic operators associated to $\dot{\mathbf{D}}^{\boldsymbol{\lambda}}$ for each
$w \in \Theta$.
These are easiest to understand by considering the
pulled back Cauchy-Riemann operator
$$
\pi^*\dot{\mathbf{D}}^{\boldsymbol{\lambda}} : W^{1,p,-\pi^*\boldsymbol{\delta}}(
\pi^*\dot{E}^{\boldsymbol{\lambda}}) \to L^{p,-\pi^*\boldsymbol{\delta}}(\pi^*\dot{F}^{\boldsymbol{\lambda}}),
$$
whose asymptotic operators we will denote by
$$
\pi^*\mathbf{A}^{\boldsymbol{\lambda}}_\zeta : \Gamma\big((\bar{\pi}^*\widebar{E}^{\boldsymbol{\lambda}})|_{S^1_\zeta}\big) \to
\Gamma\big((\bar{\pi}^*\widebar{E}^{\boldsymbol{\lambda}})|_{S^1_\zeta}\big)
$$
for $\zeta \in \Theta''$.
The relation $\left(\pi^*\dot{\mathbf{D}}^{\boldsymbol{\lambda}} \right)
\left( \eta \circ \pi\right)
= \pi^* \left( \dot{\mathbf{D}}^{\boldsymbol{\lambda}} \eta \right)$ for sections
$\eta \in \Gamma(\dot{E}^{\boldsymbol{\lambda}})$ gives rise to the following relation
between asymptotic operators:
\begin{equation}
\label{eqn:asympPullback}
\left(\pi^*\mathbf{A}^{\boldsymbol{\lambda}}_\zeta \right) \left( f \circ \pi_\zeta \right)
= k_\zeta \cdot \left(\mathbf{A}^{\boldsymbol{\lambda}}_w f\right) \circ \pi_\zeta
\quad\text{ for $f \in \Gamma\big(\widebar{E}^{\boldsymbol{\lambda}}|_{S^1_w}\big)$ 
and $\zeta \in \pi^{-1}(w)$.}
\end{equation}
This can be proved via a local computation as in \S\ref{sec:CRpunctured}:
writing $\pi(s,t) = (ks,kt)$ in suitable holomorphic cylindrical coordinates
and $\dot{\mathbf{D}}^{\boldsymbol{\lambda}} \eta = (\dbar \eta + B\eta)(-ds + i\, dt)$
for some matrix-valued function $B(s,t)$ after a choice of trivialization for $\dot{E}^{\boldsymbol{\lambda}}$
over the end near~$w$, $\mathbf{A}^{\boldsymbol{\lambda}}_w$ is represented in this
trivialization by the operator $-i\p_t - B(\infty,t)$ by definition.
The corresponding trivialized formula for $\pi^*\dot{\mathbf{D}}^{\boldsymbol{\lambda}}$
then comes from
\begin{equation*}
\begin{split}
\pi^*\dot{\mathbf{D}}^{\boldsymbol{\lambda}}\left( \eta \circ \pi \right)(s,t) &=
\left.\pi^*\left( \dot{\mathbf{D}}^{\boldsymbol{\lambda}} \eta \right)\right|_{(s,t)}\\
&= \left( \dbar\eta(ks,kt) + B(ks,kt) \eta(ks,kt) \right)\left( -d(ks) + i\, d(kt) \right) \\
&= \left( \dbar + k\, B(ks,kt) \right) (\eta \circ \pi)(s,t) \cdot (-ds + i\, dt),
\end{split}
\end{equation*}
hence $\pi^*\dot{\mathbf{D}}^{\boldsymbol{\lambda}}$ appears in trivialized form as
the sum of $\dbar$ with the zeroth-order term $k B(ks,kt)$.
The trivialized formula for $\pi^*\mathbf{A}^{\boldsymbol{\lambda}}_\zeta$ is thus
$-i \p_t - k\, B(\infty,kt)$, which explains
the factor of $k_\zeta$ appearing in \eqref{eqn:asympPullback}.

For the following discussion, fix $w \in \Theta$ and $\zeta \in \pi^{-1}(w)$.
The definition of $\dot{\mathbf{D}}^{\boldsymbol{\lambda}}$ implies that
$\pi^*\dot{\mathbf{D}}^{\boldsymbol{\lambda}}$ acts on 
sections $\eta \otimes v \in \Gamma(\pi^*\dot{E} \otimes_\RR V)$
such that $(\pi^*\dot{\mathbf{D}}^{\boldsymbol{\lambda}}) (\eta \otimes v) =
\big[ \big(\pi^*\dot{\mathbf{D}}\big) \eta \big] \otimes v$ whenever 
$v : \dot{\Sigma}'' \to V$ is constant.  From this, one deduces that for
any section $f \otimes v \in 
\Gamma\big(\bar{\pi}^*\widebar{E} \otimes_\RR V|_{S^1_\zeta}\big)$
where $f$ is an arbitrary smooth map $S^1_\zeta \to E_w$ and
$v : S^1_\zeta \to V$ is constant, we have
\begin{equation}
\label{eqn:Apullback}
\pi^*\mathbf{A}^{\boldsymbol{\lambda}}_\zeta (f \otimes v) = - (J\, \p_t f) \otimes v.
\end{equation}
Now to write down a formula for $\mathbf{A}^{\boldsymbol{\lambda}}_w$, we can use the
natural identification of $\Gamma\big(\widebar{E}^{\boldsymbol{\lambda}}|_{S^1_w}\big)$ with
the space of $G_\zeta$-equivariant loops in $E_w \otimes_\RR V$,
$$
\Gamma\big(\widebar{E}^{\boldsymbol{\lambda}}|_{S^1_w}\big) = \left\{ F \in
C^\infty(S^1_\zeta , E_w \otimes_\RR V)\ \Big|\ F(t + 1/k_\zeta) =
g_\zeta \cdot F(t) \text{ for all $t \in S^1_\zeta$} \right\}.
$$
Acting on $G_\zeta$-equivariant loops~$F$, 
\eqref{eqn:asympPullback} and \eqref{eqn:Apullback} imply
\begin{equation}
\label{eqn:Amain}
\mathbf{A}^{\boldsymbol{\lambda}}_w F = - \frac{1}{k_\zeta} J\, \p_t F,
\end{equation}
where it is understood that $J \p_t$ acts on the tensor product by
taking $F = f \otimes v$ to $(J\, \p_t f) \otimes v$ whenever $v$ is
locally constant.

\textbf{Step~3: Trivializations and Conley-Zehnder indices}.\\
This is the step in which it is helpful to be working with the
complexification $\dot{\mathbf{D}}^{\boldsymbol{\lambda}}$ rather than
directly with~$\dot{\mathbf{D}}^{\boldsymbol{\theta}}$.
In order to choose a suitable trivialization $\Phi$ and
compute $\muCZ^\Phi(\mathbf{A}^{\boldsymbol{\lambda}}_w - \delta_w)$, we shall first split
$\mathbf{A}^{\boldsymbol{\lambda}}_w$ into a direct sum
of operators on $J$-complex line bundles.  Observe that
$\widebar{E}|_{S^1_w} = S^1_w \times E_w$ is already canonically
trivial, so any complex basis of $E_w$ gives a splitting
of $\mathbf{A}^{\boldsymbol{\lambda}}_w$ over an $m$-fold direct sum of isomorphic
$J$-complex bundles of rank $2n$,
$$
\widebar{E}^{\boldsymbol{\lambda}}|_{S^1_w} = \left( L^{\boldsymbol{\lambda}} \right)^{\oplus m},
$$
where
$$
L^{\boldsymbol{\lambda}} = S^1 \times \left( \CC \otimes_\RR V \right) \Big/ G_\zeta
$$
and the generator of $G_\zeta = \ZZ_{k_\zeta}$ acts by
$g_\zeta \cdot (t , f \otimes v ) = (t + 1 / k_\zeta, f \otimes
\boldsymbol{\lambda}(g_\zeta) v)$.  Note that $L^{\boldsymbol{\lambda}}$ carries two commuting complex
structures, $J$ and $i$, which act on the first and second factor
of the tensor product respectively.  Further: $V$ admits a complex basis
$(v_1,\ldots,v_n)$ consisting of eigenvectors of
$\boldsymbol{\lambda}(g_\zeta)$, and we can then define integers
$p_j \in \{0,\ldots,k_\zeta-1\}$ for $j=1,\ldots,n$ by
$$
\boldsymbol{\lambda}(g_\zeta) v_j = e^{2\pi i p_j / k_\zeta} v_j.
$$
Here we can identify $V_\zeta' \subset V$ as the subspace spanned by all $v_j$ such
that $p_j > 0$.
Identifying $V$ with $\CC^n$ via this eigenbasis yields a splitting
$$
L^{\boldsymbol{\lambda}} = L^{\boldsymbol{\lambda}}_1 \oplus \ldots \oplus L^{\boldsymbol{\lambda}}_n,
$$
where for $j=1,\ldots,n$,
$$
L^{\boldsymbol{\lambda}}_j := S^1 \times \left( \CC \otimes_\RR \CC \right) \Big/ \ZZ_{k_\zeta},
$$
with the generator $1 \in \ZZ_{k_\zeta}$ acting by 
$1 \cdot (t , f \otimes v) = (t + 1/ k_\zeta, f \otimes e^{2\pi i p_j / k_\zeta} v)$.
This bundle again carries the two commuting complex structures 
$J$ and $i$ acting on the first and second factors of the tensor product
respectively; it has complex rank~$2$ with respect to either one.  Finally, since 
$J$ acts $i$-complex-linearly on $\CC \otimes_\RR \CC$, we can find
eigenvectors $f_\pm \in \CC \otimes_\RR \CC$ such that $J f_\pm = \pm i f_\pm$,
so the splitting $\CC \otimes_\RR \CC = \CC f_+ \oplus \CC f_-$
gives a splitting of $J$- and $i$-complex vector bundles
$$
L^{\boldsymbol{\lambda}}_j = L^{\boldsymbol{\lambda}}_{j,+} \oplus L^{\boldsymbol{\lambda}}_{j,-},
$$
with
\begin{equation}
\label{eqn:Lrhopm}
L^{\boldsymbol{\lambda}}_{j,\pm} = (S^1 \times \CC) \Big/ \ZZ_{k_\zeta},
\end{equation}
where the generator $1 \in \ZZ_{k_\zeta}$ acts by
$1 \cdot (t, f) = (t + 1 / k_\zeta , e^{2\pi i p_j / k_\zeta} f)$.
Both $L^{\boldsymbol{\lambda}}_{j,+}$ and $L^{\boldsymbol{\lambda}}_{j,-}$ are complex line bundles over $S^1$,
carrying two complex structures $J$ and $i$, which satisfy
$J = i$ on $L^{\boldsymbol{\lambda}}_{j,+}$ but $J = -i$ on $L^{\boldsymbol{\lambda}}_{j,-}$.  This splitting
of bundles gives a splitting of $\mathbf{A}^{\boldsymbol{\lambda}}_w$ in the form
\begin{equation}
\label{eqn:bigSum}
\mathbf{A}^{\boldsymbol{\lambda}}_w = \left( \bigoplus_{j=1}^n \left(  
\mathbf{A}^{\boldsymbol{\lambda}}_{j,+} \oplus \mathbf{A}^{\boldsymbol{\lambda}}_{j,-}
\right)  \right)^{\oplus m},
\end{equation}
where for $j=1,\ldots,n$, $\mathbf{A}^{\boldsymbol{\lambda}}_{j,\pm}$ acts on
$$
\Gamma(L^{\boldsymbol{\lambda}}_{j,\pm}) = \left\{ f \in C^\infty(S^1 , \CC)\ \Big|\ 
f(t + 1 / k_\zeta) = e^{2\pi i p_j / k_\zeta} f(t) \text{ for all $t \in S^1$}
\right\}
$$
by
$$
\mathbf{A}^{\boldsymbol{\lambda}}_{j,\pm} f = \mp \frac{1}{k_\zeta} i\, \p_t f.
$$

Since $L^{\boldsymbol{\lambda}}_{j,\pm}$ are complex line bundles, 
$\muCZ^\Phi(\mathbf{A}^{\boldsymbol{\lambda}}_{j,\pm} - \delta_w)$ can be computed in terms of
winding numbers of eigenfunctions of $\mathbf{A}^{\boldsymbol{\lambda}}_{j,\pm}$, using
the relation proved in \cite{HWZ:props2}*{Theorem~3.10}.  In particular, if (as will turn
out to be true in our case) all eigenspaces
of $\mathbf{A}^{\boldsymbol{\lambda}}_{j,\pm}$ have real dimension~$2$, then
\begin{equation}
\label{eqn:CZwinding}
\muCZ^\Phi(\mathbf{A}^{\boldsymbol{\lambda}}_{j,\pm} - \delta_w) = 2\wind^\Phi(f_{j,\pm}) + 1,
\end{equation}
where $f_{j,\pm} \in \Gamma(L^{\boldsymbol{\lambda}}_{j,\pm})$ is any nontrivial eigenfunction
of $\mathbf{A}^{\boldsymbol{\lambda}}_{j,\pm} - \delta_w$ with the largest possible negative
eigenvalue.  A $\ZZ_{k_\zeta}$-equivariant
function $f : S^1 \to \CC$ satisfies $\mathbf{A}^{\boldsymbol{\lambda}}_{j,\pm} f = \lambda f$ 
if and only if it is a complex multiple of
\begin{equation}
\label{eqn:eigenfunctions}
f_\lambda(t) := e^{\pm i k_\zeta \lambda t}, \qquad
\lambda \mp \frac{2\pi p_j}{k_\zeta} \in 2\pi \ZZ.
\end{equation}
Observe that since $0 < \delta_w < 2\pi/d' \le 2\pi / k_\zeta$, every eigenvalue
$\lambda$ thus satisfies $\lambda - \delta_w \ne 0$; this proves that the
perturbed asymptotic operators $\mathbf{A}^{\boldsymbol{\lambda}}_{j,\pm}$ are
all nondegenerate and thus establishes the Fredholm property for
$\dot{\mathbf{D}}^{\boldsymbol{\lambda}}$.  Now to apply \eqref{eqn:CZwinding},
we need to find the unique eigenvalue $\lambda = 2\pi (\ell \pm p_j / k_\zeta)$ 
for $\ell \in \ZZ$ such that
$$
2\pi \left(\ell \pm \frac{p_j}{k_\zeta} \right) - \delta_w < 0 < 
2\pi \left[ (\ell +1 ) \pm \frac{p_j}{k_\zeta} \right] - \delta_w.
$$
Since $0 < \delta_w < 2\pi/d'$, this condition is equivalent to
$$
\ell \le \mp \frac{p_j}{k_\zeta} < \ell + 1,
$$
so choosing the appropriate $\ell \in \ZZ$ and plugging in \eqref{eqn:eigenfunctions}
leads to the formulas
\begin{equation}
\label{eqn:fjpm}
\begin{split}
f_{j,+}(t) &:= \begin{cases}
1 & \text{ if $p_j = 0$},\\
e^{- 2\pi i (k_\zeta - p_j) t} & \text{ if $p_j > 0$},
\end{cases} \\
f_{j,-}(t) &:= e^{2\pi i p_j t}.
\end{split}
\end{equation}
Let $\Phi_j^\pm$ for $j=1,\ldots,n$ denote a choice of $J$-complex 
trivializations of $L^{\boldsymbol{\lambda}}_{j,\pm}$ such that
$$
\wind^{\Phi_j^+}(f_{j,+}) = \wind^{\Phi_j^-}(f_{j,-}) = 0, \qquad
j=1,\ldots,n,
$$
and denote by $\Phi_w$ the resulting $J$-complex trivialization of
\begin{equation}
\label{eqn:lotsOfLineBundles} 
\widebar{E}^{\boldsymbol{\lambda}}\Big|_{S^1_w} = \left( \bigoplus_{j=1}^n \left(
L^{\boldsymbol{\lambda}}_{j,+} \oplus L^{\boldsymbol{\lambda}}_{j,-} \right) \right)^{\oplus m}.
\end{equation}
By \eqref{eqn:CZwinding}, we now have
$$
\muCZ^{\Phi_j^+}(\mathbf{A}^{\boldsymbol{\lambda}}_{j,+} - \delta_w) = 
\muCZ^{\Phi_j^-}(\mathbf{A}^{\boldsymbol{\lambda}}_{j,-} - \delta_w) = 1,
$$
and thus by \eqref{eqn:bigSum},
$\muCZ^{\Phi_w}(\mathbf{A}^{\boldsymbol{\lambda}}_w - \delta_w) = 2mn$.
Note that, a priori, this construction of $\Phi_w$ depends on an
arbitrary choice $\zeta \in \pi^{-1}(w)$, but the fact that
$\muCZ^{\Phi_w}(\mathbf{A}^{\boldsymbol{\lambda}}_w - \delta_w)$ turns out to
be independent of this choice tells us that $\Phi_w$ is uniquely
determined up to homotopy.  Performing this construction for
all punctures $w \in \Theta$, we will denote the resulting
asymptotic trivialization of $\dot{E}^{\boldsymbol{\lambda}}$ simply by~$\Phi$.
We've proved:

\begin{lemma}
\label{lemma:CZindexComp}
For the asymptotic trivialization $\Phi$ described above and each puncture
$w \in \Theta$, $\muCZ^{\Phi}(\mathbf{A}^{\boldsymbol{\lambda}}_w - \delta_w) = 2mn$.
\qed
\end{lemma}

\textbf{Step~4: The relative first Chern number}.\\
It remains to compute $c_1^\Phi(\dot{E}^{\boldsymbol{\lambda}},J)$.  Consider the
pullback $\pi^*\dot{E}^{\boldsymbol{\lambda}} = \pi^*\dot{E} \otimes_\RR V$.  The
first factor in this tensor product has a canonical homotopy class of
asymptotic trivializations, which we shall denote by~$\pi^*\Psi_0$, as it is
the pullback of an asymptotic trivialization~$\Psi_0$ for~$\dot{E}$,
satisfying $c_1^{\Psi_0}(\dot{E}) = c_1(E)$.
Moreover, the second factor is globally trivial, thus
$\pi^*\dot{E}^{\boldsymbol{\lambda}}$ carries a canonical asymptotic trivialization,
denoted by $\Psi$, such that
$$
c_1^{\Psi}(\pi^*\dot{E}^{\boldsymbol{\lambda}}) = \dim_\RR V \cdot 
c_1^{\pi^*\Psi_0}(\pi^*\dot{E}) = 2n \cdot 
\deg(\pi) \cdot c_1^{\Psi_0}(\dot{E}) = 2n d' \cdot c_1(E).
$$
If $\pi^*\Phi$ denotes the pullback of $\Phi$ to an asymptotic trivialization
of $\pi^*\dot{E}^{\boldsymbol{\lambda}}$, we then have
\begin{equation}
\label{eqn:c1tau1}
\begin{split}
c_1^\Phi(\dot{E}^{\boldsymbol{\lambda}}) &= \frac{1}{d'} c_1^{\pi^*\Phi}(\pi^*\dot{E}^{\boldsymbol{\lambda}})
= \frac{1}{d'} \left[ c_1^{\Psi}(\pi^*\dot{E}^{\boldsymbol{\lambda}}) + \deg^{\Psi}(\pi^*\Phi) \right] \\
&= 2n \cdot c_1(E) + \frac{1}{d'} \deg^{\Psi}(\pi^*\Phi),
\end{split}
\end{equation}
where $\deg^{\Psi}(\pi^*\Phi) \in \ZZ$ denotes the sum over all punctures 
$\zeta \in \Theta''$
of the degrees of the transition maps $S^1 \to \GL(2mn,\CC)$ that
change $\Psi$ to $\pi^*\Phi$.  We can compute the latter for each
$w \in \Theta$ and $\zeta \in \pi^{-1}(w) \subset \Theta''$ 
as a sum of winding numbers over a line bundle 
decomposition analogous to \eqref{eqn:lotsOfLineBundles}, namely
$$
\bar{\pi}^*\widebar{E}^{\boldsymbol{\lambda}}\big|_{S^1_\zeta} = 
\pi_\zeta^*\left(\widebar{E}^{\boldsymbol{\lambda}}\big|_{S^1_w}\right) =
\left( \bigoplus_{j=1}^n \left(
\pi_\zeta^*L^{\boldsymbol{\lambda}}_{j,+} \oplus \pi_\zeta^*L^{\boldsymbol{\lambda}}_{j,-} \right) 
\right)^{\oplus m},
$$
where pulling back \eqref{eqn:Lrhopm} via the projection
$\pi_\zeta : S^1 \to S^1 / \ZZ_{k_\zeta}$ gives the trivial line bundle
$$
\pi_\zeta^*L^{\boldsymbol{\lambda}}_{j,\pm} = S^1 \times \CC,
$$
with the pulled back trivialization $\pi_\zeta^*\Phi_j^\pm$ such that
the special eigenfunctions $f_{j,\pm}$ in \eqref{eqn:fjpm} have zero winding
as $t$ traverses~$S^1$.
The restriction $\Psi_\zeta$ of $\Psi$ to
$\bar{\pi}^*\widebar{E}^{\boldsymbol{\lambda}}\big|_{S^1_\zeta}$ is now the direct sum of
the standard trivializations on each of the factors $\pi_\zeta^*L^{\boldsymbol{\lambda}}_{j,\pm}$,
thus
\begin{equation}
\label{eqn:degtheta}
\deg^{\Psi_\zeta}(\pi_\zeta^*\Phi_w) = m \sum_{j=1}^n \left[ 
\wind_{S^1}(f_{j,+}) + \wind_{S^1}(f_{j,-}) 
\right].
\end{equation}
There is an important sublety here: recall that $J = \pm i$ on
$L^{\boldsymbol{\lambda}}_{j,\pm}$, hence the orientation induced by $J$ on $L^{\boldsymbol{\lambda}}_{j,-}$
is the \emph{opposite} of the obvious one, and the sign of
$\wind_{S^1}(f_{j,-})$ must be reversed accordingly, giving
\begin{equation*}
\begin{split}
\wind_{S^1}(f_{j,+}) &= \begin{cases}
0 & \text{ if $p_j = 0$},\\
p_j - k_\zeta & \text{ if $p_j > 0$},
\end{cases} \\
\wind_{S^1}(f_{j,-}) &= - p_j.
\end{split}
\end{equation*}
Plugging this into \eqref{eqn:degtheta}, we have
$$
\deg^{\Psi_\zeta}(\pi_\zeta^*\Phi_w) = 
m \sum_{j \in \{1,\ldots,n\},\  p_j \ne 0} (-k_\zeta)
= - m k_\zeta \dim_\CC V_\zeta'.
$$
Summing over all $\zeta \in \Theta''$ and plugging into \eqref{eqn:c1tau1}
then gives
$$
c_1^\Phi(\dot{E}^{\boldsymbol{\lambda}}) = 2n \cdot c_1(E) - \frac{m}{d'} 
\sum_{\zeta \in \Theta''} k_\zeta \dim_\CC V_\zeta'.
$$
Since $\dim_\CC V_\zeta' = n_w$ is independent
of $\zeta \in \pi^{-1}(w)$ for each $w \in \Theta$, and $\sum_{\zeta \in \pi^{-1}(w)} k_\zeta = d'$,
this implies:

\begin{lemma}
\label{lemma:c1comp}
$c_1^\Phi(\dot{E}^{\boldsymbol{\lambda}}) = 2n \cdot c_1(E) - m \sum_{w \in \Theta} n_w$.
\qed
\end{lemma}

\textbf{Step~5: Conclusion of the proof}.\\
Finally, we combine Lemmas~\ref{lemma:CZindexComp} and~\ref{lemma:c1comp}
and plug into \eqref{eqn:RiemannRochReal} to obtain
\begin{equation*}
\begin{split}
\ind(\dot{\mathbf{D}}^{\boldsymbol{\lambda}}) &= 2mn \cdot \chi(\dot{\Sigma}) + 
4 n \cdot c_1(E) - 2m \sum_{w \in \Theta} n_w + 2 m n |\Theta|  \\
&=
2 \left[ mn \cdot \chi(\Sigma) + 
2 n \cdot c_1(E) - m \sum_{w \in \Theta} n_w \right],
\end{split}
\end{equation*}
and thus
$$
\ind(\dot{\mathbf{D}}^{\boldsymbol{\theta}}) 
= n \left[ m \chi(\Sigma) + 2 c_1(E) \right] - m \sum_{w \in \Theta} n_w.
$$
The expression in brackets is $\ind(\mathbf{D})$,
so this completes the proof of Theorem~\ref{thm:index}.

\section{Petri's condition}
\label{sec:unique}

\subsection{The main local result}

Standard proofs of transversality results via the Sard-Smale theorem
(cf.~\cites{FloerHoferSalamon,McDuffSalamon:Jhol2}) typically require 
some kind of unique continuation lemma, which for $J$-holomorphic curves
usually means the similarity principle.  In this section we will establish
a local result about Cauchy-Riemann type operators that plays this role in 
the proof of Theorem~\ref{thm:submanifolds0}.  It combines the usual
unique continuation property with an additional ``quadratic'' local condition
that can be achieved under generic zeroth-order perturbations.

For any pair of smooth real vector bundles $E$ and $F$ over the same manifold~$M$,
one can define the \defin{Petri map}
$$
\Pi : \Gamma(E) \otimes \Gamma(F) \to \Gamma(E \otimes F),\qquad
\Pi(\eta \otimes \xi)(p) := \eta(p) \otimes \xi(p).
$$
Since we plan to discuss purely local conditions, let us amend this by
fixing a point $p \in M$ and considering the space of \emph{germs} of smooth sections at~$p$,
$$
\Gamma_p(E) := \Gamma(E)\big/{\sim},
$$
where $\eta,\eta' \in \Gamma(E)$ represent the same element of $\Gamma_p(E)$
if and only if they match on some neighborhood of~$p$.  The Petri map then
descends to a \defin{local Petri map} at $p$,
$$
\Pi : \Gamma_p(E) \otimes \Gamma_p(F) \to \Gamma_p(E \otimes F).
$$
It is easy to see that $\Pi$ is never injective, e.g.~its kernel contains
$f\eta \otimes \xi - \eta \otimes f\xi$ for any two sections
$\eta \in \Gamma(E)$, $\xi \in \Gamma(F)$ with a smooth
function $f : M \to \RR$.  It will sometimes become injective, however, if
the domain is restricted to certain spaces of solutions to linear PDEs.
To express this properly, let us assume $\mathbf{D} : \Gamma(E) \to \Gamma(F)$
is a linear partial differential operator with smooth coefficients, and
$\mathbf{D}^* : \Gamma(F) \to \Gamma(E)$ is its formal adjoint with respect
to a choice of bundle metrics on $E,F$ and volume form on~$M$.
For any point $p \in M$, both operators descend to linear maps on the spaces of
germs of smooth sections at~$p$, which we will denote by
$$
\mathbf{D}_p : \Gamma_p(E) \to \Gamma_p(F), \qquad \mathbf{D}_p^* : \Gamma_p(F) \to \Gamma_p(E).
$$
We will also assume $\mathbf{D}$ and $\mathbf{D}^*$ uniquely determine
(via extension or restriction) linear maps
$$
\mathbf{D} : \mathbf{X}(E) \to \mathbf{Y}(F), \qquad
\mathbf{D}^* : \mathbf{X}^*(F) \to \mathbf{Y}^*(E),
$$
where $\mathbf{X}(E)$, $\mathbf{Y}^*(E)$, $\mathbf{Y}(F)$ and $\mathbf{X}^*(F)$
are vector spaces of sections (or equivalence classes of sections defined almost everywhere)
of the respective bundles; in typical examples, these will be Sobolev spaces,
sometimes with exponential weight conditions if $M$ is a noncompact manifold
with cylindrical ends.  Let us add two conditions of a local nature, both of
which are satisfied for a wide class of elliptic operators, including those of
Cauchy-Riemann type:
\begin{itemize}
\item \textsc{(regularity)} Every section in $\ker\mathbf{D} \subset \mathbf{X}(E)$
or $\ker\mathbf{D}^* \subset \mathbf{X}^*(F)$ is smooth.
\item \textsc{(unique continuation at~$p$)} The maps $\ker\mathbf{D} \to \ker\mathbf{D}_p$
and $\ker\mathbf{D}^* \to \ker\mathbf{D}^*_p$ that send each section to its
germ at $p$ are injective.
\end{itemize}
The terminology in the following definition is adapted from the work of Doan 
and Walpuski \cite{DoanWalpuski:BrillNoether}, who
borrowed it in turn from algebraic geometry (see e.g.~\cite{ArbarelloCornalbaGriffithsHarris}).

\begin{defn}
\label{defn:Petri}
Suppose $\mathbf{D} : \mathbf{X}(E) \to \mathbf{Y}(F)$ is a differential operator
with formal adjoint $\mathbf{D}^* : \mathbf{X}^*(F) \to \mathbf{Y}^*(E)$ satisfying
the conditions specified above, and $p \in \uU \subset M$.
We say that $\mathbf{D}$ satisfies
\begin{enumerate}
\item \defin{Petri's condition}, if the restricted Petri map
$\displaystyle \ker \mathbf{D} \otimes \ker \mathbf{D}^* \stackrel{\Pi}{\longrightarrow} \Gamma(E \otimes F)$
is injective; 
\item \defin{Petri's condition over~$\uU$} if there is no nontrivial element
$t \in \ker\mathbf{D} \otimes \ker\mathbf{D}^*$ such that $\Pi(t) \in \Gamma(E \otimes F)$ vanishes
identically on~$\uU$;
\item the \defin{local Petri condition at $p$} if the map
$\displaystyle \ker \mathbf{D}_p \otimes \ker \mathbf{D}^*_p \stackrel{\Pi}{\longrightarrow} \Gamma_p(E \otimes F)$
is injective;
\item \defin{Petri's condition to infinite order at $p$} if there is no
nontrivial element $t \in \ker\mathbf{D}_p \otimes \ker\mathbf{D}^*_p$
such that $\Pi(t)$ has vanishing derivatives of all orders at~$p$.
\end{enumerate}
\end{defn}

Every condition on the list in Definition~\ref{defn:Petri} implies the
previous one; note that the implication $(3) \Rightarrow (2)$ in particular 
follows from our regularity and unique continuation assumptions.
The first two conditions are global in nature, as $\ker\mathbf{D}$ and
$\ker\mathbf{D}^*$ depend on the global properties of~$\mathbf{D}$, 
including the choice of domains $\mathbf{X}(E)$ and $\mathbf{X}^*(F)$.
These kernels will always be finite dimensional in the cases we consider, so
that it seems unsurprising (if non-obvious) that Petri's condition might hold.
In contrast, the third and fourth conditions are \emph{much} stronger and more
surprising because $\ker\mathbf{D}_p$ and $\ker\mathbf{D}_p^*$ are in general
infinite dimensional, but the local conditions are also more powerful,
e.g.~it will be extremely useful to observe that they are preserved under
pullbacks via branched covers of the base.

\begin{remark}
\label{remark:PetriIndep}
As defined above, the global versions of Petri's condition may 
in general depend not only on the operator $\mathbf{D}$ but also on the
auxiliary geometric data (bundle metrics and volume form) used to define~$\mathbf{D}^*$,
but the local conditions are independent of these choices.
Indeed, whenever $\mathbf{D}_1^*$ and $\mathbf{D}_2^*$ are two operators
arising as formal adjoints of $\mathbf{D}$ via different choices of the
geometric data, there is a smooth bundle automorphism $\Phi : F \to F$ that maps
local solutions of $\mathbf{D}_1^*\xi=0$ to local solutions of $\mathbf{D}_2^*\xi=0$,
so that $\1 \otimes \Phi : E \otimes F \to E \otimes F$ identifies the two
different versions of $\ker\Pi \subset \ker\mathbf{D}_p \otimes \ker\mathbf{D}^*_p$.
\end{remark}

\begin{remark}
\label{remark:openDensePetri}
It is clear from the definition that the set of points $p \in M$ at which the
local Petri condition is \emph{not} satisfied is open.  We will see in \S\ref{sec:PetriJ}
that Petri's condition to infinite order can sometimes be shown to hold at all points
in a dense subset of some region $\uU \subset M$, so it follows in this situation that the local
Petri condition also holds at \emph{all} points in~$\uU$.
\end{remark}

It should be emphasized that whenever we refer to the above definition,
we will be regarding all vector spaces as \emph{real} vector spaces so that
``$\otimes$'' means the real tensor product, even in cases where
$\mathbf{D}$ happens to be complex linear.  The only exception is 
Example~\ref{ex:complexPetri} below,
which is a digression from the main topic at hand.

\begin{example}
Elliptic operators over $1$-dimensional domains satisfy something much
stronger than the Petri condition to infinite order, because by local uniqueness of solutions to ODEs,
any linearly-independent set of local sections in $\ker\mathbf{D}$ 
or $\ker\mathbf{D}^*$ is also pointwise linearly independent.
For similar reasons, any Cauchy-Riemann type operator 
$\mathbf{D} : \Gamma(E) \to \Gamma(F)$ that splits over a direct sum of complex line bundles
with nonpositive first Chern numbers over a closed surface~$\Sigma$
must satisfy the \emph{global} Petri condition over arbitrary subsets $\uU \subset \Sigma$.
The reason for this is that on a line bundle $E \to \Sigma$ with $c_1(E) \le 0$,
the similarity principle guarantees that global solutions to
$\mathbf{D}\eta = 0$ are either trivial or nowhere vanishing, so that globally linearly-independent
sets of solutions are also linearly independent at every point.  This property
might not hold for the formal adjoint $\mathbf{D}^*$, but since solutions to
$\mathbf{D}^*\xi = 0$ satisfy unique continuation, any expression of the form
$\sum_{ij} c^{ij} \eta_i \otimes \xi_j$ with a nontrivial set of coefficients $c^{ij} \in \RR$
and linearly-independent sets $\{\eta_i \in \ker\mathbf{D}\}$
and $\{\xi_j \in \ker\mathbf{D}^*\}$ is still guaranteed to be nonzero at every point
outside a discrete subset.
Example~\ref{ex:realPetri} below shows however that the local Petri condition
in this situation is not always satisfied.
\end{example}

\begin{example}
\label{ex:complexPetri}
Complex-linear Cauchy-Riemann operators over a Riemann surface 
satisfy the complex version of Petri's condition to
infinite order at every point, i.e.~the definition above is satisfied if
real tensor products are replaced by complex tensor products.  
One can prove this by choosing holomorphic
trivializations and writing elements of $\ker\mathbf{D}$ and $\ker\mathbf{D}^*$
locally as Taylor series in $z$ or $\bar{z}$ respectively: it then turns
out that for any nontrivial $t \in \ker\mathbf{D} \otimes_\CC \ker\mathbf{D}^*$,
the Taylor series in $z$ and $\bar{z}$ for the resulting section of 
$E \otimes_\CC F$ at a given point is always nontrivial.
We omit the details since we will not need this fact.
\end{example}

\begin{example}
\label{ex:realPetri}
If we regard the standard Cauchy-Riemann operator $\mathbf{D}=\dbar$ on a 
trivial line bundle and its formal adjoint $\mathbf{D}^* = -\p$ as real-linear
operators, then they do not satisfy the local Petri condition at any point.
A local counterexample is given by
$$
1 \otimes i\bar{z} - i \otimes \bar{z} - z \otimes i + iz \otimes 1 \in
\ker \dbar \otimes_\RR \ker \p.
$$
It follows that the local Petri condition is also not satisfied by any
Cauchy-Riemann type operator that splits off a complex-linear summand.
\end{example}

\begin{example}
\label{ex:realPetri2}
Here is an example of a Cauchy-Riemann type operator that does not split off
any complex-linear summand but still fails to satisfy the local Petri condition:
take $E$ and $F$ to be the trivial complex line bundle over $\CC$, with standard bundle
metrics and the standard area form, and consider
$\mathbf{D} := \dbar + \kappa$, $\mathbf{D}^* = -\p + \kappa$, where
$\kappa : \CC \to \CC$ is complex conjugation.  Using coordinates $s + it \in \CC$,
one can associate to every $\lambda \in (-1,1)$ solutions $\eta_\lambda \in \ker\mathbf{D}$
and $\xi_\lambda \in \ker\mathbf{D}^*$ defined by\footnote{The inspiration for this
example comes from the asymptotic formulas in \cites{HWZ:props1,Siefring:asymptotics}:
in particular on the cylinder $\RR \times S^1$ with coordinates $(s,t)$, a translation-invariant Cauchy-Riemann
type equation $(\dbar + B(t)) \eta(s,t) = 0$ always has solutions of the form
$\eta(s,t) = e^{s \lambda} f(t)$, where $f$ is an eigenfunction of the asymptotic operator
$-i\p_t - B(t)$ with eigenvalue $\lambda \in \RR$.  In the asymptotic setting one
requires solutions to be periodic in~$t$, in which case the eigenvalue $\lambda$ can
only take a discrete set of values, but periodicity is not necessary in Example~\ref{ex:realPetri2},
and $\lambda$ can therefore be chosen much more freely.}
\begin{equation*}
\begin{split}
\eta_\lambda(s+it) &:= e^{\lambda s + \sqrt{1 - \lambda^2} t} \left( \sqrt{1-\lambda} + i \sqrt{1+\lambda} \right),\\
\xi_\lambda(s+it) &:= e^{-\lambda s - \sqrt{1-\lambda^2} t} \left( \sqrt{1-\lambda} - i \sqrt{1 + \lambda} \right).
\end{split}
\end{equation*}
Identifying the fibers $\CC$ with $\RR^2$ so that the fibers of $E \otimes_\RR F$
become the space of real $2$-by-$2$ matrices, the products $\Pi(\eta_\lambda \otimes \xi_\lambda)$
are now constant sections of $E \otimes_\RR F$:
$$
\Pi(\eta_\lambda \otimes \xi_\lambda)(s+it) = \begin{pmatrix}
1 - \lambda & -\sqrt{1 - \lambda^2} \\
\sqrt{1 - \lambda^2} & -1 - \lambda
\end{pmatrix}.
$$
Such products span the $3$-dimensional space of real matrices of the form
$\begin{pmatrix} a & b \\ -b & c \end{pmatrix}$,
thus any four such products must be linearly dependent, and the dependence
relation gives rise to nontrivial elements in $\ker \Pi$ by choosing four
distinct values of $\lambda \in (-1,1)$.
\end{example}

\begin{remark}
\label{remark:fail}
An earlier version of this paper 
(see Appendix~\ref{sec:catastropheThis})
claimed that every Cauchy-Riemann type
operator whose complex-antilinear part is invertible at a point $p$ satisfies
Petri's condition to infinite order at~$p$, but Example~\ref{ex:realPetri2} contradicts
that.
\end{remark}

The operators in Examples~\ref{ex:realPetri} and \ref{ex:realPetri2} are rather
special, and our main objective in this section is to prove that such
counterexamples cannot arise for \emph{generic} Cauchy-Riemann type operators.
To set up the result, assume now that $\Sigma$ is a Riemann surface with a 
Hermitian bundle metric
$\langle\ ,\ \rangle_\Sigma$ on~$T\Sigma$.  We will not require $\Sigma$ to be
compact since the discussion will be purely local, but fix a point $p \in \Sigma$
and an open neighborhood $\uU \subset \Sigma$ of $p$ with compact closure.
Fix also a complex vector bundle
$E \to \Sigma$ with a Hermitian bundle metric, let $F = \overline{\Hom}_\CC(T\Sigma,E)$,
and denote by $\CRR(E)$ the space of real-linear Cauchy-Riemann type operators
$\mathbf{D} : \Gamma(E) \to \Gamma(F)$.  We shall fix a specific $\mathbf{D}\fix \in \CRR(E)$
and define the space of all Cauchy-Riemann type operators $\mathbf{D}$
that match $\mathbf{D}\mathbf\fix$ outside of~$\uU$:
$$
\CRR(E\,;\,\uU,\mathbf{D}\fix) := \left\{ \mathbf{D} \in \CRR(E)\ \big|\ 
\mathbf{D} - \mathbf{D}\fix = 0 \text{ on $\Sigma \setminus \uU$} \right\}.
$$
This is an affine space over the Fr\'echet space of smooth sections of $\Hom_\RR(E,F)$
that vanish outside~$\uU$, so in particular it is a complete metric space.
For every $\mathbf{D} \in \CRR(E)$, $\mathbf{D}^*$
will denote the formal adjoint of $\mathbf{D}$ determined by the bundle
metrics on $E$ and~$\Sigma$.

For any $\eta \in \Gamma_p(E)$, we define the
\defin{vanishing order} of $\eta$ at $p$ by
$$
\ord(\eta;p) := \sup \left\{ k \in \{0\} \cup \NN\ \big|\ \text{all derivatives of
$\eta$ at $p$ up to order $k$ vanish} \right\}.
$$
For $t \in \Gamma_p(E) \otimes \Gamma_p(F)$, we will then say that $t$
\defin{vanishes to order $k$} if $t$ can be written as a finite sum
$t = \sum_j \eta_j \otimes \xi_j$ such that
$$
\ord(\eta_j;p) + \ord(\xi_j;p) \ge k \quad\text{ for every $j$.}
$$
The usual unique continuation results imply that for every $\mathbf{D} \in \CRR(E)$,
nontrivial local solutions to the equations $\mathbf{D}\eta = 0$ or
$\mathbf{D}^*\eta = 0$ satisfy $\ord(\eta;p) < \infty$ at every point.
One can easily prove from this that nontrivial elements 
$t \in \ker \mathbf{D}_p \otimes \ker \mathbf{D}_p^*$ also cannot vanish to
infinite order (see Proposition~\ref{prop:quadVanishing}).

The machinery developed in the next two subsections will prove:

\begin{thm}
\label{thm:Petri}
For every $\ell \in \NN$, there exists an integer $k \ge \ell$ and a
Baire subset 
$$
\CRR^{\ell,\reg}(E\,;\,\uU,\mathbf{D}\fix) \subset \CRR(E\,;\,\uU,\mathbf{D}\fix)
$$ 
with the following
significance: for every $\mathbf{D} \in \CRR^{\ell,\reg}(E\,;\,\uU,\mathbf{D}\fix)$,
if $\eta_1,\ldots,\eta_\ell \in \ker\mathbf{D}_p$ and
$\xi_1,\ldots,\xi_\ell \in \ker\mathbf{D}_p^*$ are $\ell$-tuples of local
solutions such that $t := \sum_{j=1}^\ell \eta_j \otimes \xi_j \in 
\Gamma_p(E) \otimes \Gamma_p(F)$ does not vanish to order~$\ell$, then
$\Pi(t) \in \Gamma_p(E \otimes F)$ does not vanish to order~$k$.
\end{thm}

In light of unique continuation, we now set
$$
\CRR^\reg(E\,;\,\uU,\mathbf{D}\fix) := \bigcap_{\ell \in \NN}
\CRR^{\ell,\reg}(E\,;\,\uU,\mathbf{D}\fix) \subset \CRR(E\,;\,\uU,\mathbf{D}\fix)
$$
and obtain:

\begin{cor}
\label{cor:Petri}
There exists a Baire subset
$$
\CRR^{\reg}(E\,;\,\uU,\mathbf{D}\fix) \subset \CRR(E\,;\,\uU,\mathbf{D}\fix)
$$
such that every $\mathbf{D} \in \CRR^{\reg}(E\,;\,\uU,\mathbf{D}\fix)$
satisfies Petri's condition to infinite order at the point $p \in \uU$.
\end{cor}

This result can be extended in various ways.  For instance, 
the regular set $\CRR^{\reg}(E\,;\,\uU,\mathbf{D}\fix)$ defined above depends
\textit{a priori} on the choice of a point $p \in \uU$, but one can also
find a Baire set of operators such that Petri's condition to infinite
order is satisfied simultaneously at every point in~$\uU$.  More generally, one can
consider smooth families of operators parametrized by a finite-dimensional
manifold and prove that for generic families, every operator in the
family satisfies these conditions.  In \S\ref{sec:PetriJ}, we will prove
that the normal Cauchy-Riemann operators of $J$-holomorphic curves can
all be assumed to satisfy Petri's condition to infinite order in regions 
where $J$ can be perturbed generically.  One of the advantages of
focusing on purely \emph{local} conditions is that once we establish this
result for somewhere injective curves, it carries over immediately to their
multiple covers, which will be a crucial ingredient in the proof of
Theorem~\ref{thm:submanifolds0}.

The aforementioned extensions of Corollary~\ref{cor:Petri} are all based
on the Sard-Smale theorem, but Theorem~\ref{thm:Petri} itself requires
(aside from unique continuation) only finite-dimensional analysis and linear
algebra.  Indeed, the conditions defining each of the spaces
$\CRR^{\ell,\reg}(E\,;\,\uU,\mathbf{D}\fix)$
in the statement of the theorem depend only on the $k$-jet of
$\mathbf{D} \in \CRR(E\,;\,\uU,\mathbf{D}\fix)$ at $p$ for some finite $k \in \NN$,
and this data varies in a finite-dimensional smooth manifold.  The
idea behind the proof is roughly to show that the set of jets of operators 
not satisfying the desired conditions lives in ``walls'' whose
codimensions can be assumed arbitrarily large
by making $k$ larger.  These walls are not submanifolds in general, but are
what we call ``$C^\infty$-subvarieties,'' whose local structure is nice
enough to apply Sard's theorem as if they were manifolds.  
(The necessary background on $C^\infty$-subvarieties is reviewed
in Appendix~\ref{sec:subvarieties}.)  The main
technical work behind the proof is then to estimate the ranks of certain
large matrices that determine the codimensions of these subvarieties.

The rest of this section will proceed as follows.
In \S\ref{sec:jetSpaces}, we introduce a general formalism for studying
differential operators via jet spaces at a point, and explain how results
such as Theorem~\ref{thm:Petri} can be reduced to a specific technical
lemma on estimating the ranks of certain finite-dimensional linear transformations.
We will then address this problem for Cauchy-Riemann operators in
\S\ref{sec:PetriCR}, leading to the proof of Theorem~\ref{thm:Petri}.
The extension to a result about normal Cauchy-Riemann operators of
holomorphic curves for generic $J$ will be stated and proved in~\S\ref{sec:PetriJ},
and \S\ref{sec:PetriLin} will then give an important application of Petri's
condition to global transversality problems as arising in
Theorem~\ref{thm:submanifolds0}.

\subsection{Jet space formalism}
\label{sec:jetSpaces}

The contents of this subsection are not specific to Cauchy-Riemann operators,
but may be relevant in principle to any linear partial differential operator
with smooth coefficients.

\subsubsection{Germs, jets, and the vanishing order filtration}
\label{sec:filtration}

Fix a smooth $n$-di\-men\-sion\-al manifold $M$ with a smooth vector bundle
$E \to M$ of real rank $m \in \NN$.  For a chosen
point $p \in M$, we continue to denote by
$$
\Gamma_p(E) := \Gamma(E)\big/{\sim}
$$
the vector space of germs of smooth sections of $E$ defined near~$p$.
This space has a natural filtration
\begin{equation}
\label{eqn:vanishingFiltration}
\Gamma_p(E) = \Gamma_p(E)^0 \supset \Gamma_p(E)^1 \supset \Gamma_p(E)^2 \supset \ldots,
\end{equation}
where for each $k \in \ZZ$ we define $\Gamma_p(E)^k \subset \Gamma_p(E)$ as the
space of germs of sections whose derivatives up to order $k-1$ at $p$
all vanish.  For $k \le 0$ this is a vacuous condition, hence
$\Gamma_p(E)^k = \Gamma_p(E)$.  For each $k \in \ZZ$ we define the space
of \defin{$k$-jets of sections} at $p$ by
$$
J^k_p E := \Gamma_p(E) \big/ \Gamma_p(E)^{k+1}.
$$
We will typically abuse notation by using a single symbol
such as $\eta$ to represent a section in $\Gamma(E)$, its germ in $\Gamma_p(E)$
and its $k$-jet in $J^k_p E$; when there is need for more clarity in the
notation, we will sometimes write the natural quotient projections as
$$
\Gamma(E) \text{ or } \Gamma_p(E) \stackrel{J^k_p}{\longrightarrow} J^k_p E,
$$
so that the $k$-jet of a section $\eta \in \Gamma(E)$ at $p$ can be denoted
by~$J^k_p \eta \in J^k_p E$.  The jet space inherits from
\eqref{eqn:vanishingFiltration} a finite filtration
\begin{equation}
\label{eqn:filtration}
J^k_p E = (J^k_p E)^0 \supset (J^k_p E)^1 \supset \ldots \supset
(J^k_p E)^{k} \supset (J^k_p E)^{k+1} = \{0\},
\end{equation}
where for each $\ell \le k$, $(J^k_p E)^{\ell+1}$ is the kernel of the
quotient projection $J^\ell_p : J^k_p E \to J^\ell_p E$.

There is an obvious isomorphism of $J^0_p E$ with the fiber
$E_p$, and the spaces $J^k_p E$ for $k < 0$ are all trivial.
If we choose local coordinates $(x_1,\ldots,x_n)$ for $M$ identifying
$p$ with $0 \in \RR^n$, together with a
trivialization of $E$ near~$p$, then $J^k_p E$ for each $k \in \ZZ$ becomes naturally identified
with the vector space of $\RR^m$-valued Taylor polynomials of degree at most~$k$,
\begin{equation}
\label{eqn:TaylorPoly}
\sum_{|\alpha| \le k} x^\alpha c_{\alpha}, \qquad c_{\alpha} \in \RR^m.
\end{equation}
The notation for the filtration above has been chosen so that under this
identification, $(J^k_p E)^\ell$ becomes
the space of Taylor polynomials of degree at most $k$ that are also
$O(|x|^\ell)$.

Given two vector spaces $V = V^0 \supset V^1 \supset V^2 \supset \ldots$
and $W = W^0 \supset W^1 \supset W^2 \supset \ldots$ with filtrations, we
will say in general that a linear map $T : V \to W$ \defin{preserves the
filtrations} if $T(V^n) \subset W^{n}$ for every
$n \ge 0$.

\subsubsection{Differential operators and formal adjoints}
\label{sec:diffOps}

Since we are mainly interested in Cauchy-Riemann type operators, for simplicity
we shall only consider differential operators of order~$1$ in the following
discussion, though the jet space formalism could easily be extended beyond this.

Given a second smooth vector bundle $F \to M$ of real rank $\ell \in \NN$ and
a first-order linear partial differential operator 
$\mathbf{D} : \Gamma(E) \to \Gamma(F)$ with smooth coefficients,
$\mathbf{D}$ descends to a map $\Gamma_p(E) \to \Gamma_p(F)$ that sends 
$\ker J_p^k \subset \Gamma_p(E)$ into $\ker J_p^{k-1} \subset \Gamma_p(F)$
for each $k \in \ZZ$, thus it also descends to a linear map
$$
\mathbf{D} : J^k_p E \to J^{k-1}_p F.
$$
Let us denote by
$$
\diff_p(E,F) \subset \Hom\big( \Gamma_p(E),\Gamma_p(F)\big)
$$
the vector space consisting of all germs at $p$ of linear differential operators 
$\Gamma(E) \to \Gamma(F)$ of order at most~$1$ with smooth coefficients.
The vector space of linear maps $J^k_p E \to J^{k-1}_p F$ that are induced by
operators in $\diff_p(E,F)$ will then be denoted by
$$
\diff^k_p(E,F) \subset \Hom\big( J^k_p E,J^{k-1}_p F \big),
$$
and we will again abuse notation by using a single symbol such as
$\mathbf{D}$ to denote a global differential operator $\Gamma(E) \to \Gamma(F)$,
its germ in $\diff_p(E,F)$, and the map in $\diff^k_p(E,F)$ that it determines.
Observe that $\diff^k_p(E,F)$ is a
finite-dimensional vector space isomorphic to the $(n+1)$-fold
product of $J^{k-1}_p\Hom(E,F)$: indeed, if we fix 
local coordinates $(x_1,\ldots,x_n)$ identifying a neighborhood
of $p$ with the $n$-disk $\DD^n_\epsilon$ of some radius $\epsilon > 0$, along with local trivializations
of $E$ and $F$ over the same neighborhood, then each $\mathbf{D} \in \diff_p(E,F)$
is represented by an operator $C^\infty(\DD^n_\epsilon,\RR^m) \to C^\infty(\DD^n_\epsilon,\RR^\ell)$ of the form
\begin{equation}
\label{eqn:DinCoords}
\mathbf{D} = \sum_{j=1}^n a_j \p_j + b
\end{equation}
for some smooth functions $a_1,\ldots,a_n,b : \DD^n_\epsilon \to \Hom(\RR^m,\RR^\ell)$.
For a given $\eta \in \Gamma(E)$, the $(k-1)$-jet of $\mathbf{D}\eta$ at $p$
is thus determined by the $(k-1)$-jets of the functions $a_1,\ldots,a_n,b$
at that point, and these are equivalent to bundle maps $E \to F$ defined near~$p$.

We will also consider a subset
$$
\diffaff_p(E,F) \subset \diff_p(E,F),
$$
which is assumed to have the property that for any given
$\mathbf{D} \in \diffaff_p(E,F)$, another operator $\mathbf{D}' \in\diff_p(E,F)$
satisfies
$$
\mathbf{D}' \in \diffaff_p(E,F)\quad\Leftrightarrow\quad
\mathbf{D}' = \mathbf{D} + A \text{ for some $A\in \Gamma_p(\Hom(E,F))$},
$$
i.e.~$\diffaff_p(E,F)$ is an affine space over $\Gamma_p(\Hom(E,F))$.
The space of maps $J^k_p E \to J^{k-1}_p F$ induced by operators
$\mathbf{D}\in \diffaff_p(E,F)$ then defines a subset
$$
\diffaff^k_p(E,F) \subset \diff^k_p(E,F),
$$
which is naturally an affine space over the finite-dimensional vector space $J^{k-1}_p\Hom(E,F)$.

In order to bring formal adjoints into this picture, we need to make choices
of bundle metrics for $E$ and $F$ and a volume form on $M$ near~$p$;
these choices will often be referred to collectively as the \defin{geometric data}.
It will be useful to fix geometric data once and for all at the point $p$
itself, while allowing it to vary at other points near~$p$.  Concretely, fix a
pair of inner products
$$
g_p = \langle\ ,\ \rangle_{E_p} \text{ on $E_p$}, \qquad
h_p = \langle\ ,\ \rangle_{F_p} \text{ on $F_p$},
$$
along with a nontrivial alternating multilinear $n$-form
$$
\mu_p \in \Lambda^n T^*_p M.
$$
Let us denote by $S^2 E^* \subset E^* \otimes E^* \to M$ the vector bundle of symmetric 
bilinear forms $E \oplus E \to \RR$.  The space of \defin{$k$-jets of bundle
metrics} on $E$ which match $g_p$ at $p$ is then
$$
J^k_p\big( \met(E)\big) := \left\{ g \in J^k_p\big( S^2 E^*\big) \ \Big|\ J^0_p g = g_p \right\},
$$
and it is naturally an affine space over the finite-dimensional vector space
$\ker J^0_p \subset J^k_p\big(S^2 E^*\big)$.  We similarly define the affine spaces
$$
J^k_p\big( \met(F)\big) := \left\{ h \in J^k_p\big( S^2 E^*\big) \ \Big|\ J^0_p h = h_p \right\}
$$
and
$$
J^k_p\big(\volume(M)\big) := \left\{ \mu \in J^k_p\big( \Lambda^n T^*M \big)\ \Big|\ J^0_p \mu = \mu_p \right\},
$$
which consist respectively of $k$-jets of bundle metrics on $F$ matching $h_p$ at $p$
and $k$-jets of volume forms on $M$ matching $\mu_p$ at~$p$.  We will again
abuse notation by using a single symbol such as $g$ or $\langle\ ,\ \rangle_E$
to denote a global bundle metric on~$E$ that matches $g_p$ at~$p$, or the
germ of such a metric near~$p$, or its $k$-jet in $J^k_p\big(\met(E)\big)$;
similar remarks apply to $J^k_p\big(\met(F)\big)$ and $J^k_p\big(\volume(M)\big)$.

Any choice of smooth bundle metrics $g = \langle\ ,\ \rangle_E$ on $E$
and $h = \langle\ ,\ \rangle_F$ on $F$ and a volume form
$\mu \in \Omega^n(M)$ assigns to each differential operator
$\mathbf{D} : \Gamma(E) \to \Gamma(F)$ a \defin{formal adjoint}
$\mathbf{D}^* : \Gamma(F) \to \Gamma(E)$ satisfying the relation
$$
\int_M \langle \xi,\mathbf{D} \eta \rangle_F \,\mu = 
\int_M \langle \mathbf{D}^*\xi , \eta \rangle_E \,\mu
\quad\text{ for all }\quad
\eta \in C_0^\infty(E),\ \xi \in C_0^\infty(F).
$$
Fix local coordinates and trivializations near $p$ to write 
$\mathbf{D}$ again in the form \eqref{eqn:DinCoords}.
The chosen bundle metrics and volume form can be written in terms of the
standard Euclidean inner product $\langle\ ,\ \rangle$ and volume form
$dx_1 \wedge \ldots \wedge dx_n$ as
$$
\langle\ ,\ \rangle_E = \langle\cdot,G\cdot\rangle, \quad
\langle\ ,\ \rangle_F = \langle\cdot,H\cdot\rangle, \quad
\mu = F\, dx_1 \wedge \ldots \wedge dx_n
$$
for some smooth functions $F : \DD^n_\epsilon \to \RR$, $G : \DD^n_\epsilon \to \End(\RR^m)$
and $H : \DD^n_\epsilon \to \End(\RR^\ell)$, where $F$ is everywhere nonzero and
$G$ and $H$ take values in the spaces of symmetric positive-definite matrices.
Note that the condition defining $\mathbf{D}^*$ does not change if the sign of
$\mu$ is reversed, so without loss of generality let us assume $F > 0$.
One can then compute a local formula for $\mathbf{D}^* : C^\infty(\DD^n_\epsilon,\RR^\ell)
\to C^\infty(\DD^n_\epsilon,\RR^m)$ as
\begin{equation}
\label{eqn:Dstar}
\mathbf{D}^* = -\sum_j (G^{-1} a_j^\transpose H) \p_j +
G^{-1} \left( b^\transpose H - \sum_j \left[ a_j^\transpose H \p_j(\ln F) + \p_j(a_j^\transpose H) \right] \right).
\end{equation}
We observe from this formula that the germ $\mathbf{D}^* \in \diff_p(F,E)$ at $p$ is determined
by the corresponding germs of the geometric data $g,h,\mu$ and
$\mathbf{D} \in \diff_p(E,F)$.  Moreover, if the first-order terms 
$a_j$ in~$\mathbf{D}$ are fixed, then for any
$\xi \in \Gamma(F)$, the $(k-1)$-jet of $\mathbf{D}^*\xi$ at $p$ is
determined by the $(k-1)$-jet of~$g$, the $k$-jets of $\mu$ and~$h$, and the
$(k-1)$-jet of the zeroth-order term $b$ in~$\mathbf{D}$.
It follows that the correspondence assigning to each $\mathbf{D} \in \diffaff_p(E,F)$
with germs of geometric data $g,h,\mu$ the germ of a formal adjoint $\mathbf{D}^* \in \diff_p(F,E)$
descends to a well-defined map
\begin{equation}
\label{eqn:smoothAdj}
\diffaff^k_p(E,F) \times J^{k-1}_p\big(\met(E)\big) \times J^k_p\big(\met(F)\big) \times J^k_p\big(\volume(M)\big) \stackrel{*}{\longrightarrow}
\diff^k_p(F,E).
\end{equation}
All the spaces involved in this map are finite-dimensional 
manifolds, and the map is smooth.

\subsubsection{Unique continuation in tensor products}

If $V = V^0 \supset V^1 \supset V^2 \supset \ldots$ and $W = W^0 \supset W^1 \supset W^2 \supset\ldots$
are two vector spaces with filtrations, then $V \otimes W$ inherits a
natural filtration
$$
V \otimes W = (V \otimes W)^0 \supset (V \otimes W)^1 \supset (V \otimes W)^2 \supset \ldots,
$$
where for each $n \ge 0$,
$$
(V \otimes W)^n := (V^0 \otimes W^n) + (V^1 \otimes W^{n-1}) + \ldots + (V^n \otimes W^0).
$$

\begin{lemma}
\label{lemma:linIndTensor}
Given two filtered vector spaces $V$ and~$W$, if $t \in (V \otimes W)^n$ is
nontrivial, then for some $r \in \NN$, $t$ can be written as
$$
t = \sum_{j=1}^r v_j \otimes w_j
$$
for two linearly-independent sets $v_1,\ldots,v_r \in V$ and
$w_1,\ldots,w_r \in W$ such that for all $j=1,\ldots,r$, we have
$$
v_j \in V^{k_j} \text{ and } w_j \in W^{\ell_j} \quad\text{ where }\quad
k_j + \ell_j = n.
$$
\end{lemma}
\begin{proof}
Suppose $t = \sum_{j=1}^r v_j \otimes w_j$ satisfies all of these conditions
except that the set $v_1,\ldots,v_r$ is linearly dependent, so there exist
constants $c_1,\ldots,c_r$ with $\sum_j c_j v_j = 0$ and not all of the $c_j$
are zero.  After reordering the set, we can assume without loss of generality
that $c_1 \ne 0$ and, for every $j=2,\ldots,r$ with $c_j \ne 0$,
$k_j \ge k_1$.  Writing $v_1 = \sum_{j=2}^r \frac{c_j}{c_1} v_j$ then gives
$$
t = \sum_{j=2}^r v_j \otimes \widehat{w}_j \quad\text{ where }\quad
\widehat{w}_j := w_j + \frac{c_j}{c_1} w_1.
$$
For each $j=2,\ldots,r$, we now have $\ell_j = n - k_j \le n - k_1 = \ell_1$,
thus $w_1 \in W^{\ell_1} \subset W^{\ell_j}$ and therefore
$\widehat{w}_j \in W^{\ell_j}$, hence the shortened sum also satisfies the
desired conditions.  One can apply a similar procedure to shorten the sum
if instead $w_1,\ldots,w_r$ is linearly dependent, and repeating this enough
times produces two sets that are both linearly independent.
\end{proof}

Let us say that a differential operator $\mathbf{D} \in \diff_p(E,F)$
has the \defin{strong unique continuation} property if there exists no
nontrivial solution $\eta \in \ker \mathbf{D}$ such that
$\eta \in \Gamma_p(E)^k$ for every $k \in \NN$.

\begin{prop}
\label{prop:quadVanishing}
If $\mathbf{D} \in \diff_p(E,F)$ and $\mathbf{D}^* \in \diff_p(F,E)$ both
have the strong unique continuation property, then there exists no
nontrivial element $t \in \ker\mathbf{D} \otimes \ker\mathbf{D}^*$ such that
$t \in \big( \Gamma_p(E) \otimes \Gamma_p(F) \big)^k$ for every $k \in \NN$.
\end{prop}
\begin{proof}
Given $t \in \ker\mathbf{D} \otimes \ker\mathbf{D}^*$ nonzero, there are
uniquely defined finite-dimensional subspaces $V \subset \ker\mathbf{D}$
and $W \subset \ker\mathbf{D}^*$ such that for any pair of linearly-independent
sets $\eta_1,\ldots,\eta_r \in \ker\mathbf{D}$ and $\xi_1,\ldots,\xi_r \in \ker\mathbf{D}^*$
with $t = \sum_j\eta_j \otimes \xi_j$,
$$
V = \Span\{\eta_1,\ldots,\eta_r\} \quad\text{ and }\quad
W = \Span\{\xi_1,\ldots,\xi_r\}.
$$
We claim there exists $k \in \NN$ such that no nontrivial $\eta \in V$
is in $\Gamma_p(E)^k$ and no nontrivial $\xi \in W$ is in $\Gamma_p(F)^k$.
Indeed, if there does not exist such a number for~$V$, then there exist
sequences $\eta_j \in V$ and $k_j \in \NN$ with $k_j \to \infty$ and
$\eta_j \in \Gamma_p(E)^{k_j}$ for every~$j$.  Since $V$ is finite dimensional, we can
normalize the $\eta_j$ and then find a convergent subsequence $\eta_j \to \eta_\infty \in V$
whose limit is nontrivial, but must also belong to $\bigcap_{k \in \NN} \Gamma_p(E)^k$,
giving a contradiction.  The same argument works for~$W$.

Now, fixing $k \in \NN$ as in the previous paragraph, suppose
$t \in (\Gamma_p(E) \otimes \Gamma_p(F))^{2k}$ and $t \ne 0$.
Lemma~\ref{lemma:linIndTensor} then
writes $t$ in the form $\sum_j \eta_j \otimes \xi_j$ where the $\eta_j$ and
$\xi_j$ are necessarily bases of $V$ and $W$ respectively, but they also
satisfy $\eta_j\in \Gamma_p(E)^{\ell_j}$ and $\xi_j \in \Gamma_p(F)^{m_j}$
with $\ell_j + m_j \ge 2k$ for each~$j$.  This implies either $\ell_j \ge k$ or
$m_j \ge k$ in each case, and is thus a contradiction.
\end{proof}

\subsubsection{Local rescaling}
\label{sec:rescaling}

Every differential operator is locally equivalent (up to choices of
coordinates and trivializations) to an arbitrarily small perturbation of an
operator with constant coefficients and no lower-order terms.
To make use of this observation, we shall from now on impose
the following additional condition on the 
affine space of local operators $\diffaff_p(E,F) \subset \diff_p(E,F)$:

\begin{assumption}
\label{ass:const}
There exists a choice of coordinates identifying a neighborhood $\uU \subset M$
of~$p$ with $\DD^n_\epsilon \subset \RR^n$ and $p$ with $0 \in \RR^n$, along with
local trivializations over~$\uU$, in which the
first-order coefficients $a_j : \DD^n_\epsilon \to \Hom(\RR^m,\RR^\ell)$ in
$\mathbf{D} = \sum_j a_j \p_j + b$ for each $\mathbf{D} \in \diffaff_p(E,F)$ are
constant.
\end{assumption}

Let us fix once and for all a neighborhood $\uU \subset M$ of $p$ with 
coordinates and trivializations for which the condition in
Assumption~\ref{ass:const} holds.  For every $\varepsilon \in [0,1]$,
we then associate to each $\mathbf{D} \in \diff_p(E,F)$ an
operator $\mathbf{D}_\varepsilon \in \diff_p(E,F)$
such that if $\mathbf{D}$ takes the form $\mathbf{D}\eta(x) = 
\sum_j a_j(x) \p_j \eta(x) + b(x)\eta(x)$ in the chosen coordinates and
trivializations, then $\mathbf{D}_\varepsilon$ is given by
$$
\mathbf{D}_\varepsilon \eta(x) = \sum_j a_j(\varepsilon x) \p_j\eta(x) + \varepsilon b(\varepsilon x) \eta(x),
\qquad
\varepsilon \in [0,1].
$$
We can similarly associate to each $\eta \in \Gamma_p(E)$ and $\xi \in \Gamma_p(F)$
germs of sections $\eta_\varepsilon \in \Gamma_p(E)$ and $\xi_\varepsilon \in \Gamma_p(F)$,
which in coordinates take the form
$$
\eta_\varepsilon(x) := \eta(\varepsilon x), \qquad \xi_\varepsilon(x) := \xi(\varepsilon x).
$$
We then have
$$
\mathbf{D}_\varepsilon \eta_\varepsilon = \varepsilon (\mathbf{D}\eta)_\varepsilon
$$
for every $\mathbf{D} \in \diff_p(E,F)$ and $\eta \in \Gamma_p(E)$.  Letting
these operators descend to jet spaces, we obtain 
for every $\mathbf{D} \in \diff^k_p(E,F)$ a smooth
$1$-parameter family of operators $\left\{ \mathbf{D}_\varepsilon \in \diff^k_p(E,F)\right\}_{\varepsilon \in [0,1]}$ and
linear maps
$$
J^k_p E \to J^k_p E : \eta \mapsto \eta_\varepsilon,
$$
which for $\varepsilon > 0$ are isomorphisms sending $\ker \mathbf{D} \subset J^k_p E$
onto $\ker \mathbf{D}_\varepsilon \subset J^k_p E$.

Next, fix geometric data consisting of bundle metrics
$g = \langle\ ,\ \rangle_E$ on $E$ and $h = \langle\ ,\ \rangle_F$ on $F$,
and a volume form $\mu$, such that all three match the fixed choices of
data $g_p$, $h_p$ and $\mu_p$ at~$p$.  
Using the same coordinates and trivializations over~$\uU$, 
we can write $g = \langle\cdot,G\cdot \rangle$, 
$h = \langle\cdot,H\cdot\rangle$ and $\mu = F\, dx_1 \wedge \ldots \wedge dx_n$,
and then define a smooth $1$-parameter family of geometric data
$g_\varepsilon,h_\varepsilon,\mu_\varepsilon$ for $\varepsilon \in [0,1]$ by replacing
the functions $G$, $H$ and $F$ with
$$
G_\varepsilon(x) := G(\varepsilon x), \qquad H_\varepsilon(x) := H(\varepsilon x), \qquad
F_\varepsilon(x) := F(\varepsilon x).
$$
Note that since $p$ corresponds to $0 \in \DD^n_\epsilon$ in coordinates, the modified
geometric data still matches the fixed choices $g_p,h_p,\mu_p$ at~$p$,
and we can then descend to jet spaces to obtain smooth $1$-parameter families
$$
g_\varepsilon \in J^{k-1}_p\big( \met(E)\big), \qquad
h_\varepsilon \in J^k_p\big(\met(F)\big),\qquad
\mu_\varepsilon \in J^k_p\big(\volume(M)\big)
$$
for $\varepsilon \in [0,1]$.
Now if $\mathbf{D}^* \in \diff_p(F,E)$ denotes the formal adjoint of 
$\mathbf{D} \in \diffaff_p(E,F)$ with respect to
the geometric data $g,h,\mu$ and $\mathbf{D}^*_\varepsilon \in \diff_p(F,E)$ is defined from
$\mathbf{D}^*$ via the same rescaling prescription as $\mathbf{D}_\varepsilon$ 
described above, then we
see from \eqref{eqn:Dstar} that $\mathbf{D}^*_\varepsilon$ is in fact the
formal adjoint of $\mathbf{D}_\varepsilon$ with respect to the data
$g_\varepsilon,h_\varepsilon,\mu_\varepsilon$.  Moreover, Assumption~\ref{ass:const}
implies that the map $\diff^k_p(E,F) \to \diff^k_p(E,F)$ induced by
$\diff_p(E,F) \to \diff_p(E,F) : \mathbf{D} \mapsto \mathbf{D}_\varepsilon$ preserves
$\diffaff^k_p(E,F)$, so we can now fit the smooth map \eqref{eqn:smoothAdj}
into the rows of a commutative diagram
$$
\begin{tikzcd}
\diffaff^k_p(E,F) \times J^{k-1}_p\big(\met(E)\big) \times J^k_p\big(\met(F)\big) \times J^k_p\big(\volume(M)\big) \ar[r,"*"] \ar[d,"\varepsilon"] & \diff^k_p(F,E) \ar[d,"\varepsilon"] \\
\diffaff^k_p(E,F) \times J^{k-1}_p\big(\met(E)\big) \times J^k_p\big(\met(F)\big) \times J^k_p\big(\volume(M)\big) \ar[r,"*"]                   & \diff^k_p(F,E),
\end{tikzcd}
$$
where the vertical maps abbreviated by ``$\varepsilon$'' are defined via the
corresondences $\mathbf{D} \mapsto \mathbf{D}_\varepsilon$,
$g \mapsto g_\varepsilon$, $h \mapsto h_\varepsilon$, $\mu \mapsto \mu_\varepsilon$
and $\mathbf{D}^* \mapsto \mathbf{D}^*_\varepsilon$.
The case $\varepsilon=0$ is special: since all $\mathbf{D} \in \diffaff_p(E,F)$ have
matching first-order terms and the geometric data $g,h,\mu$ all match at~$p$,
$\mathbf{D}_0$ and $\mathbf{D}_0^*$ are
\emph{uniquely-defined} operators that 
depend on the space $\diffaff_p(E,F)$ and the chosen inner products
$g_p$ and $h_p$, but
not otherwise on the specific choices of
operator $\mathbf{D} \in \diffaff_p(E,F)$ or volume form or bundle metrics.
Similarly, the volume form $\mu_0$ and bundle metrics 
$g_0$ and $h_0$ are fully determined by the fixed data $\mu_p$, $g_p$ and~$h_p$.

\subsubsection{Right-inverses}
\label{sec:rightInverse}

Henceforward we impose the following additional assumption.

\begin{assumption}
\label{ass:rightInverse}
The operators $\mathbf{D}_0 : J^k_p E \to J^{k-1}_p F$ and
$\mathbf{D}_0^* : J^k_p F \to J^{k-1}_p E$ obtained by the rescaling
procedure in \S\ref{sec:rescaling} are surjective.
\end{assumption}

\begin{remark}
It is not difficult to show that Assumption~\ref{ass:rightInverse}
is satisfied whenever the operators in $\diffaff_p(E,F)$ are elliptic.
For Cauchy-Riemann operators in particular, this is virtually obvious,
and we will write down explicit choices of right-inverses for
that case in \S\ref{sec:mainRank}.
\end{remark}

\begin{lemma}
\label{lemma:itsABundle}
Under Assumption~\ref{ass:rightInverse}, every
$\mathbf{D} \in \diffaff^k_p(E,F)$ is surjective,
and so is $\mathbf{D}^* \in \diff^k_p(F,E)$ for every
choice of geometric data
$g \in J^{k-1}_p\big(\met(E)\big)$, $h \in J^k_p\big(\met(F)\big)$ and
$\mu \in J^k_p\big(\volume(M)\big)$.
\end{lemma}
\begin{proof}
Since $\mathbf{D}_\varepsilon$ converges in $\Hom(J^k_p E,J^{k-1}_p F)$
to $\mathbf{D}_0$ as $\varepsilon \to 0$,
surjectivity of $\mathbf{D}_0$ implies for any given $\mathbf{D} \in \diff^k_p(E,F)$
that $\mathbf{D}_\varepsilon$ is also surjective for all $\varepsilon > 0$
sufficiently small.  The isomorphism $\ker \mathbf{D} \to \ker \mathbf{D}_\varepsilon$
induced by the correspondence $\eta \mapsto \eta_\varepsilon$ for all $\varepsilon > 0$
then implies that $\mathbf{D}$ is also surjective.  The same argument works
for the formal adjoints since $\mathbf{D}^*_\varepsilon \to \mathbf{D}^*_0$
as $\varepsilon \to 0$.
\end{proof}

Since we are working in finite-dimensional spaces, surjectivity 
allows us to choose right-inverses
$$
\mathbf{T}_0 : J^{k-1}_p F \to J^k_p E, \qquad
\mathbf{T}^*_0 : J^{k-1}_p E \to J^k_p F
$$
for $\mathbf{D}_0$ and $\mathbf{D}^*_0$ respectively.
We would now like to derive from these similar right-inverses for other 
operators that are close to
$\mathbf{D}_0$ and~$\mathbf{D}^*_0$, along with explicit isomorphisms between
the kernels of nearby operators.  To this end, consider an open neighborhood
$$
(g_0,h_0,\mu_0,\mathbf{D}_0) \in \uU \subset J^{k-1}_p\big(\met(E)\big) \times
J^k_p\big(\met(F)\big) \times J^k_p\big(\volume(M)\big) \times \diffaff^k_p(E,F),
$$
which we reserve the right to make smaller as necessary.  
Given $(g,h,\mu,\mathbf{D}) \in \uU$, we will as usual denote by $\mathbf{D}^*$
the formal adjoint of $\mathbf{D}$ with respect to the geometric data
$(g,h,\mu)$.  Since
$\mathbf{D}_0\mathbf{T}_0 = \1$ and $\mathbf{D}^*_0 \mathbf{T}^*_0 = \1$,
we can assume after shrinking $\uU$ that
for every $(g,h,\mu,\mathbf{D}) \in \uU$, the operators
$\mathbf{D}\mathbf{T}_0 : J^{k-1}_p F \to J^{k-1}_p F$ and
$\mathbf{D}^*\mathbf{T}^*_0 : J^{k-1}_p E \to J^{k-1}_p E$ are both close
enough to the identity to be invertible.  This gives rise to right-inverses
for $\mathbf{D}$ and $\mathbf{D}^*$, defined respectively by
$$
\mathbf{T} := \mathbf{T}_0(\mathbf{D}\mathbf{T}_0)^{-1} :
J^{k-1}_p F \to J^k_p E, \qquad
\mathbf{T}^* := \mathbf{T}^*_0(\mathbf{D}^*\mathbf{T}^*_0)^{-1} :
J^{k-1}_p E \to J^k_p F.
$$
Notice that $\mathbf{T}$ and $\mathbf{T}^*$ 
depend smoothly on $(g,h,\mu,\mathbf{D}) \in \uU$.

For a fixed $(g,h,\mu,\mathbf{D}) \in \uU$, arbitrary operators close to $\mathbf{D}$
in $\diffaff^k_p(E,F)$ have the form $\widehat{\mathbf{D}} := \mathbf{D} + A$
for $A \in J^{k-1}_p\Hom(E,F)$ small, and the formal adjoint $\widehat{\mathbf{D}}^*$
with respect to the geometric data $(g,h,\mu)$ is then $\mathbf{D}^* + A^*$,
where $A^* \in J^{k-1}_p \Hom(F,E)$ is the $(k-1)$-jet of the fiberwise
transpose (with respect to $g$ and~$h$) of a smooth bundle map 
$E \to F$ representing~$A$.  If $A$ is small enough,\footnote{We will not need
this detail, but it is often possible to choose $\mathbf{T}_0$ and $\mathbf{T}^*_0$
so that they have degree $+1$ with respect to the vanishing-order filtration,
in which case the operators $A\mathbf{T}$,
$\mathbf{T}A$, $A^*\mathbf{T}^*$ and $\mathbf{T}^*A^*$ also have this property 
and are therefore nilpotent.  It follows in this case
that all infinite series appearing in this discussion are actually finite sums,
so $A$ does not really need to be small.} then we can use the
same trick again to write down right-inverses of $\widehat{\mathbf{D}}$ and
$\widehat{\mathbf{D}}^*$ in the form
\begin{equation*}
\begin{split}
\widehat{\mathbf{T}} &:= \mathbf{T}(\widehat{\mathbf{D}}\mathbf{T})^{-1} =
\mathbf{T} \left( \1 + A\mathbf{T} \right)^{-1} = 
\mathbf{T} \sum_{j=0}^\infty (-1)^j (A\mathbf{T})^j,\\
\widehat{\mathbf{T}}^* &:= \mathbf{T}^*(\widehat{\mathbf{D}}^*\mathbf{T}^*)^{-1} =
\mathbf{T}^* \left( \1 + A^*\mathbf{T}^* \right)^{-1} =
\mathbf{T}^* \sum_{j=0}^\infty (-1)^j (A^* \mathbf{T}^*)^j.
\end{split}
\end{equation*}
Shrinking the size of $A$ further if necessary, we can then define isomorphisms
\begin{equation*}
\begin{split}
\Psi_{(\mathbf{D},A)} &:= \1 - \widehat{\mathbf{T}} A = \sum_{j=0}^\infty (-1)^j (\mathbf{T} A)^j :
J^k_p E \to J^k_p E,\\
\Psi^*_{(\mathbf{D},A)} &:= \1 - \widehat{\mathbf{T}}^*A^* = \sum_{j=0}^\infty (-1)^j (\mathbf{T}^* A^*)^j :
J^k_p F \to J^k_p F,
\end{split}
\end{equation*}
which satisfy
$$
\widehat{\mathbf{D}} \Psi_{(\mathbf{D},A)} = \widehat{\mathbf{D}} - A = \mathbf{D} \quad\text{ and }\quad
\widehat{\mathbf{D}}^* \Psi^*_{(\mathbf{D},A)} = \widehat{\mathbf{D}}^* - A^* = \mathbf{D}^*,
$$
so they restrict to isomorphisms $\ker\mathbf{D} \stackrel{\Psi_{(\mathbf{D},A)}}{\longrightarrow} \ker\widehat{\mathbf{D}}$
and $\ker\mathbf{D}^* \stackrel{\Psi^*_{(\mathbf{D},A)}}{\longrightarrow} \ker\widehat{\mathbf{D}}^*$ 
respectively.  The operators $\Psi_{(\mathbf{D},A)}$ and $\Psi^*_{(\mathbf{D},A)}$
depend smoothly on both $(g,h,\mu,\mathbf{D}) \in \uU$ and $A \in J^{k-1}_p\Hom(E,F)$.

\subsubsection{The universal Petri moduli space}
\label{sec:universalPetri}

We now consider the subset
$$
\vV^k \subset J^{k-1}_p\big(\met(E)\big) \times J^k_p\big(\met(F)\big) \times
J^k_p\big(\volume(M)\big) \times \diffaff^k_p(E,F) \times \left( J^k_p E \otimes J^k_p F \right)
$$
consisting of all tuples $(g,h,\mu,\mathbf{D},t)$ such that
$$
t \in \ker \mathbf{D} \otimes \ker \mathbf{D}^* \subset J^k_p E \otimes J^k_p F,
$$
where it should be understood that $\mathbf{D}^*$ is the formal adjoint
of $\mathbf{D}$ with respect to the geometric data $g,h,\mu$.
In light of Assumption~\ref{ass:rightInverse} and Lemma~\ref{lemma:itsABundle},
the obvious projection endows $\vV^k$ with a
natural vector bundle structure
$$
\vV^k \to J^{k-1}_p\big(\met(E)\big) \times J^k_p\big(\met(F)\big) \times
J^k_p\big(\volume(M)\big) \times \diffaff^k_p(E,F),
$$
whose fiber over $(g,h,\mu,\mathbf{D})$ is
$\ker\mathbf{D} \otimes \ker \mathbf{D}^*$.  
We will prefer to think of $\vV^k$ rather as a \emph{family} of vector
bundles over the space of operators $\diffaff^k_p(E,F)$, parametrized by the
space of geometric data $(g,h,\mu) \in J^{k-1}_p\big(\met(E)\big) \times J^k_p\big(\met(F)\big) \times
J^k_p\big(\volume(M)\big)$.  Thus for each
$(g,h,\mu)$, denote
$$
\vV^k(g,h,\mu) := \left\{ (\mathbf{D},t)\ \Big|\ (g,h,\mu,\mathbf{D},t) \in \vV^k \right\}.
$$

It will be useful to amend these definitions in two ways.  Given a pair of
real vector spaces $V$ and~$W$, let us say that an element $t \in V \otimes W$
has \defin{rank $r$} if $t = \sum_{j=1}^r v_j \otimes w_j$ for two
linearly-independent sets $v_1,\ldots,v_r \in V$ and $w_1,\ldots,w_r \in W$.
Note that if $V$ is finite dimensional, then the rank of $t \in V \otimes W$
under the canonical isomorphism $V \otimes W \cong \Hom(V^*,W)$ is just the
rank of the corresponding linear map $V^* \to W$.  As a consequence, whenever
$V$ and $W$ are both finite dimensional, the set of elements of rank $r \in\NN$
in $V \otimes W$ is a smooth submanifold whose codimension is the dimension
of $\Hom(\ker T,\coker T)$ for a linear map $T : V^* \to W$ of rank~$r$, giving
\begin{equation}
\label{eqn:codimRank}
\begin{split}
\dim \left\{ t \in V \otimes W\ \big|\ \rank t = r\right\} &=
\dim V \cdot \dim W -
\left( \dim V - r \right) \cdot \left( \dim W - r \right) \\
&= r(\dim V + \dim W) - r^2.
\end{split}
\end{equation}
With this understood, we can define for each $r \in \NN$ a smooth submanifold
$$
\vV^k_r := \left\{ (g,h,\mu,\mathbf{D},t) \in \vV^k\ \Big|\ \rank t = r \right\},
$$
which is foliated by the smooth family of smooth submanifolds
$$
\vV^k_r(g,h,\mu) := \left\{ (\mathbf{D},t) \in \vV^k(g,h,\mu)\ \Big|\ \rank t = r \right\}
$$
parametrized by the space of geometric data
$(g,h,\mu) \in J^{k-1}_p\big(\met(E)\big) \times J^k_p\big(\met(F)\big) \times
J^k_p\big(\volume(M)\big)$.
Finally, recalling the filtration by vanishing orders in 
\S\ref{sec:filtration}, we define for each $\ell \in \{1,\ldots,k\}$ the open subset
$$
\vV^k_{r,\ell} := \left\{ (g,h,\mu,\mathbf{D},t) \in \vV^k_r\ \big|\ t \not\in \big( J^k_p E \otimes J^k_p F \big)^\ell \right\},
$$
which is likewise foliated by a smooth family of submanifolds
$$
\vV^k_{r,\ell}(g,h,\mu) := \left\{ (\mathbf{D},t) \in \vV^k_r(g,h,\mu)\ \big|\ 
t \not\in \big( J^k_p E \otimes J^k_p F\big)^\ell \right\}.
$$
parametrized by the geometric data $(g,h,\mu)$.

The Petri map $\Pi : \Gamma_p(E) \otimes \Gamma_p(F) \to \Gamma_p(E \otimes F)$ descends for each $k \in \ZZ$
to a linear map
$$
\Pi^k : J^k_p E \otimes J^k_p F \to J^k_p (E \otimes F)
$$
that preserves the filtration by vanishing orders.  Since the projection map
$\vV^k_{r,\ell}(g,h,\mu) \to J^k_p E \otimes J^k_p F$ sending $(g,h,\mu,\mathbf{D},t)$ to $t$
is smooth and also depends smoothly on the geometric data $(g,h,\mu)$,
$\Pi^k$ gives rise to a smooth family of smooth maps
\begin{equation}
\label{eqn:PetriMapkrl}
\Pi^k_{r,\ell} : \vV^k_{r,\ell}(g,h,\mu) \to J^k_p(E \otimes F) :
(\mathbf{D},t) \mapsto \Pi^k(t),
\end{equation}
whose zero-set we shall denote by
$$
\Petrimod^k_{r,\ell}(g,h,\mu) := (\Pi^k_{r,\ell})^{-1}(0) = \left\{ (\mathbf{D},t) \in \vV^k_{r,\ell}(g,h,\mu) \ \Big|\ \Pi^k(t) = 0 \right\}.
$$
This is the so-called \emph{universal Petri moduli space}.
Our main goal is to prove under suitable assumptions that
it is a $C^\infty$-subvariety in $\vV^k_{r,\ell}(g,h,\mu)$ 
and to establish an effective lower bound $R \in\NN$ on its codimension.
Once this is done, Sard's theorem (see Appendix~\ref{sec:subvarieties})
will imply that for almost every $\mathbf{D} \in \diffaff^k_p(E,F)$, the space
$$
\Petrimod^k_{r,\ell}(g,h,\mu,\mathbf{D}) := \left\{ t\ \Big|\ (\mathbf{D},t) \in \Petrimod^k_{r,\ell}(g,h,\mu) \right\}
$$
is a $C^\infty$-subvariety of codimension at least $R$ in the manifold
$$
\vV^k_{r,\ell}(g,h,\mu,\mathbf{D}) := \left\{ t \in \ker\mathbf{D} \otimes \ker\mathbf{D}^*\ \Big|\ \rank t = r,\ t \not\in (J^k_p E \otimes J^k_p F)^\ell \right\}.
$$
If the codimension $R$ is large enough, this will imply that $\Petrimod^k_{r,\ell}(g,h,\mu,\mathbf{D})$ is empty.

Denote the linearization of the map \eqref{eqn:PetriMapkrl} at the point
$(\mathbf{D},t) \in \Petrimod^k_{r,\ell}(g,h,\mu)$ by
$$
d_2\Pi^k_{r,\ell}(g,h,\mu,\mathbf{D},t) : T_{(\mathbf{D},t)} \vV^k_{r,\ell}(g,h,\mu) \to J^k_p(E \otimes F),
$$
where the subscript in ``$d_2$'' is
meant to emphasize that this is a partial derivatve---we differentiate with
respect to $(\mathbf{D},t)$ while holding $(g,h,\mu)$ constant.
Estimating the rank of $d_2\Pi^k_{r,\ell}$ requires
being able to write down a sufficiently large space of tangent vectors
in $T_{(\mathbf{D},t)} \vV^k_{r,\ell}(g,h,\mu)$.  Suppose that
$(g,h,\mu,\mathbf{D})$ belongs to the neighborhood $\uU$ of
$(g_0,h_0,\mu_0,\mathbf{D}_0)$ chosen in \S\ref{sec:rightInverse}, so we have
right-inverses $\mathbf{T},\mathbf{T}^*$ and isomorphisms
$\Psi_{(\mathbf{D},A)},\Psi^*_{(\mathbf{D},A)}$ that depend smoothly on
$(g,h,\mu,\mathbf{D}) \in \uU$ and a small zeroth-order perturbation 
$A \in J^{k-1}_p\Hom(E,F)$.  We can use this to associate to every
$A \in J^{k-1}_p\Hom(E,F)$ and $t \in \ker\mathbf{D} \otimes \ker\mathbf{D}^*$
a smooth path
$$
(-\delta,\delta) \to \vV^k(g,h,\mu) : s \mapsto (\mathbf{D} + s A , (\Psi_{(\mathbf{D},sA)} \otimes \Psi^*_{(\mathbf{D},sA)}) t)
$$
which passes through $(\mathbf{D},t)$ at $s=0$.
Observe that if $t = \sum_{j=1}^r \eta_j \otimes \xi_j$ for two
linearly-independent sets $\eta_1,\ldots,\eta_r \in J^k_p E$ and 
$\xi_1,\ldots,\xi_r \in J^k_p F$, then $\Psi_{(\mathbf{D},sA)}$ and
$\Psi^*_{(\mathbf{D},sA)}$ map these to linearly-independent sets
when $s$ is close enough to~$0$, since both operators are then close
to the identity.  It follows that if $(\mathbf{D},t) \in \Petrimod^k_{r,\ell}(g,h,\mu)$,
then the path above is in $\vV^k_{r,\ell}(g,h,\mu)$ for $\delta > 0$
sufficiently small.
Differentiating it at $s=0$, then feeding the resulting tangent
vector into $d_2\Pi^k_{r,\ell}(g,h,\mu,\mathbf{D},t)$ and
multiplying the result by $-1$ for cosmetic purposes, we obtain the linear map
\begin{equation*}
\begin{split}
\mathbf{L}(g,h,\mu,\mathbf{D},t) : J^{k-1}_p\Hom(E,F) &\to J^k_p(E \otimes F),\\
A &\mapsto \Pi^k \circ (\mathbf{T} A \otimes \1 + \1 \otimes \mathbf{T}^* A^*)(t).
\end{split}
\end{equation*}
This depends smoothly on the data $(g,h,\mu,\mathbf{D},t)$ and is
well defined whenever $(g,h,\mu,\mathbf{D})$ is sufficiently close
to $(g_0,h_0,\mu_0,\mathbf{D}_0)$.  The rank of this operator is clearly
less than or equal to that of $d_2\Pi^k_{r,\ell}(g,h,\mu,\mathbf{D},t)$.
We shall abbreviate the special case
\begin{equation}
\label{eqn:thisIsTheKey}
\mathbf{L}_t := \mathbf{L}(g_0,h_0,\mu_0,\mathbf{D}_0,t) :
J^{k-1}_p\Hom(E,F) \to J^k_p(E \otimes F)
\end{equation}
for $t \in \ker\mathbf{D}_0 \otimes \ker\mathbf{D}_0^*$, as this will
turn out to be the only case that matters in practice.  In fact,
we can now use the rescaling trick from \S\ref{sec:rescaling} to reduce
the local analysis of the space $\Petrimod^k_{r,\ell}(g,h,\mu)$ to the problem of
estimating the rank of~$\mathbf{L}_t$.

For every $\varepsilon \in (0,1]$ and $q \in \ZZ$ and every 
choice of the geometric data $(g,h,\mu)$, one can define a diffeomorphism
\begin{equation}
\label{eqn:PhiEps}
\Phi_\varepsilon : \vV^k_{r,\ell}(g,h,\mu) \stackrel{\cong}{\longrightarrow} \vV^k_{r,\ell}(g_\varepsilon,h_\varepsilon,\mu_\varepsilon) :
(\mathbf{D},t) \mapsto \left( \mathbf{D}_\varepsilon,t_\varepsilon \right),
\end{equation}
where the map $\ker \mathbf{D} \otimes \ker\mathbf{D}^* \to 
\ker \mathbf{D}_\varepsilon \otimes \ker \mathbf{D}^*_\varepsilon :
t \mapsto t_\varepsilon$ is defined via
\begin{equation}
\label{eqn:PhiEps2}
\eta \otimes \xi \mapsto \frac{1}{\varepsilon^q} \eta_\varepsilon \otimes \xi_\varepsilon.
\end{equation}
The scaling factor $\varepsilon^q$ here is not strictly necessary, but has
been added for use in the proof of Lemma~\ref{lemma:thisIsTheKey} below.
We see that $\Phi_\varepsilon$ maps
$\Petrimod^k_{r,\ell}(g,h,\mu)$ bijectively onto $\Petrimod^k_{r,\ell}(g_\varepsilon,h_\varepsilon,\mu_\varepsilon)$
for each $\varepsilon \in (0,1]$.
This map is not defined for $\varepsilon=0$, but the data
$g_\varepsilon$, $h_\varepsilon$, $\mu_\varepsilon$, $\mathbf{D}_\varepsilon$
and $\mathbf{D}^*_\varepsilon$ do have
well-defined limits as $\varepsilon \to 0$; in particular,
$\mathbf{D}_0$ and $\mathbf{D}^*_0$ are both operators with constant coefficients
and no zeroth-order term in our chosen local coordinates and trivializations.  
The following definition is highly dependent on this
choice of coordinates, but so is the map~$\Phi_\varepsilon$; there will be no problem
as long as the same choices are used for both.

\begin{defn}
We will say that an element of $J^k_p E$ or
$J^k_p F$ is \defin{homogeneous of degree~$d$} if, under the natural identifications
of these spaces with spaces of Taylor polynomials determined by the chosen
coordinates and trivializations from Assumption~\ref{ass:const},
it is represented by a homogeneous polynomial of degree~$d$.
Similarly, we will call an element $t = \sum_j \eta_j \otimes \xi_j \in
J^k_p E \otimes J^k_p F$ \defin{homogeneous of degree~$d$} if for every~$j$,
the elements $\eta_j \in J^k_p E$ and $\xi_j \in J^k_p F$ are homogeneous 
with degrees adding up to~$d$.
\end{defn}

\begin{remark}
The homogeneous elements $t \in J^k_p E \otimes J^k_p F$ of degree $q$ are
precisely those which are fixed under the map \eqref{eqn:PhiEps2}
for every $\varepsilon > 0$.
\end{remark}

\begin{lemma}
\label{lemma:thisIsTheKey}
Suppose that for every homogeneous element $t \in \ker \Pi^k \subset J^k_p E \otimes J^k_p F$ of degree
less than~$\ell$ that also belongs to $\ker \mathbf{D}_0 \otimes \ker \mathbf{D}^*_0$,
the linear map $\mathbf{L}_t : J^{k-1}_p\Hom(E,F) \to J^k_p(E \otimes F)$
has rank at least $R \in \NN$.  Then for every $r \in \NN$,
$\Petrimod^k_{r,\ell}(g,h,\mu)$ is a
$C^\infty$-subvariety of codimension at least $R$ in $\vV^k_{r,\ell}(g,h,\mu)$.
\end{lemma}
\begin{proof}
Suppose $(\mathbf{D},t) \in \Petrimod^k_{r,\ell}(g,h,\mu)$ and let
$q \in \{0,\ldots,\ell-1\}$ denote the largest integer such that
$t \in (J^k_p E \otimes J^k_p F)^q$. Use this value of $q$ to define the
scaling factor in \eqref{eqn:PhiEps2} for the definition of the
diffeomorphisms $\Phi_\varepsilon$ in \eqref{eqn:PhiEps}.  Identifying
$k$-jets with Taylor polynomials as in \eqref{eqn:TaylorPoly}, we can
write $t$ as a finite sum $\sum_j \eta_j \otimes \xi_j$, where for
each individual value of~$j$, $\eta_j \in \ker\mathbf{D}$ is a polynomial
of degree at most $k$ with lowest-order term of degree $u_j \ge 0$, $\xi_j \in \ker\mathbf{D}^*$
is likewise a polynomial of degree at most $k$ with lowest-order term of
degree $v_j \ge 0$, and $u_j + v_j \ge q$, with equality $u_j + v_j = q$
in at least one case.  It follows that $t_\varepsilon \in J^k_p E \otimes J^k_p F$
converges as $\varepsilon \to 0$ to a nontrivial homogenous element
$t_0 \in \ker\mathbf{D}_0 \otimes \ker\mathbf{D}^*_0 \subset J^k_p E \otimes J^k_p F$
of degree~$q < \ell$, and $\Pi^k(t_0) = 0$ since $\Pi^k(t_\varepsilon) = \Pi^k(t) = 0$
for every~$\varepsilon > 0$.  As a consequence,
$(g_\varepsilon,h_\varepsilon,\mu_\varepsilon,\mathbf{D}_\varepsilon,t_\varepsilon) \in
\vV^k$ converges as $\varepsilon \to 0$ to
$(g_0,h_0,\mu_0,\mathbf{D}_0,t_0) \in\vV^k$.
Since $\mathbf{L}_{t_0}$ has rank
at least~$R$ by the hypothesis of the lemma, it follows for all
$\varepsilon > 0$ sufficiently small that
$$
\rank d_2\Pi^k_{r,\ell}(g_\varepsilon,h_\varepsilon,\mu_\varepsilon,\mathbf{D}_\varepsilon,t_\varepsilon) \ge
\rank \mathbf{L}(g_\varepsilon,h_\varepsilon,\mu_\varepsilon,\mathbf{D}_\varepsilon,t_\varepsilon) \ge R.
$$
Fix $\varepsilon > 0$ in this range.
Then an arbitrary element $(\mathbf{D}',t') \in \vV^k_{r,\ell}(g,h,\mu)$ in some 
small neighborhood of $(\mathbf{D},t)$ belongs to $\Petrimod^k_{r,\ell}(g,h,\mu)$
if and only if $\Pi^k \circ \Phi_\varepsilon(\mathbf{D}',t') = 0$.
Since $\Phi_\varepsilon$ is a diffeomorphism,
the linearization of $\Pi^k \circ \Phi_\varepsilon : \vV^k_{r,\ell}(g,h,\mu) \to J^k_p(E \otimes F)$
at $(\mathbf{D},t)$ has the same image as the operator
$d_2\Pi^k_{r,\ell}(g_\varepsilon,h_\varepsilon,\mu_\varepsilon,\mathbf{D}_\varepsilon,t_\varepsilon)$,
and thus has rank at least~$R$.
\end{proof}

\subsection{Application to Cauchy-Riemann operators}
\label{sec:PetriCR}

We shall now apply Lemma~\ref{lemma:thisIsTheKey} for
the specific case of Cauchy-Riemann type operators.  For the rest of this
section, assume $M$ is a Riemann surface $(\Sigma,j)$, $E$ is a complex
vector bundle of complex rank~$m \in \NN$, $F = \overline{\Hom}_\CC(T\Sigma,E)$,
and $\diffaff_p(E,F)$ is the space of germs of real-linear Cauchy-Riemann
type operators on~$E$ near~$p \in \Sigma$.  This space of operators
satisfies Assumption~\ref{ass:const} since one can always choose
trivializations and coordinates in which every $\mathbf{D} \in \diffaff_p(E,F)$
is a zeroth-order perturbation of $\dbar := \p_s + i\p_t$.
To define formal adjoints, we assume $g = \langle\ ,\ \rangle_E$ is the
real part of a Hermitian bundle metric on~$E$, $\mu$ is the area form
on $\Sigma$ determined by a Hermitian bundle metric
$\langle\ ,\ \rangle_\Sigma$ on~$T\Sigma$, and $h = \langle\ ,\ \rangle_F$
is the real part of the Hermitian bundle metric determined on $F$
via the natural isomorphism $F \cong T\Sigma \otimes_\CC E$.  

\begin{remark}
It is
important to keep in mind that the operators $\mathbf{D} \in \diffaff_p(E,F)$
are in general real- and not complex-linear, thus throughout this section,
the symbols $\Hom(V,W)$ and $V \otimes W$ will always refer to \emph{real}-linear
maps and real tensor products unless otherwise noted, even in cases where
$V$ and $W$ are both complex.  We will use the notation 
$\Hom_\CC(V,W)$ and $V \otimes_\CC W$ to specify the complex analogues of
these operations.
\end{remark}

\subsubsection{A digression on real and complex tensor products}
\label{sec:tensor}

Suppose $V$ and $W$ are complex vector spaces, and let $\widebar{W}$
denote the complex conjugate of~$W$, i.e.~it is the same real vector
space, but with a sign inserted in the definition of its complex structure.
There is then a canonical complex-antilinear isomorphism $W \to \widebar{W}$
defined by the identity map, and we shall denote it by
$$
W \to \widebar{W} : w \mapsto \bar{w}.
$$
The spaces $V \otimes_\CC W$ and $V \otimes_\CC \widebar{W}$ are both
quotients of the real tensor product $V \otimes W$, e.g.~we obtain
$V \otimes_\CC W$ from $V \otimes W$ by introducing the equivalence relation
$iv \otimes w \sim v \otimes iw$, and for $V \otimes_\CC \widebar{W}$ the
relation is instead $iv \otimes w \sim -v \otimes iw$.  If the resulting
quotient projections are denoted by $\pi_+ : V \otimes W \to V \otimes_\CC W$
and $\pi_- : V \otimes W \to V \otimes_\CC \widebar{W}$, then we obtain
an isomorphism
$$
(\pi_+,\pi_-) : V \otimes W  \stackrel{\cong}{\longrightarrow} (V \otimes_\CC W) \oplus (V \otimes_\CC \widebar{W}).
$$
This discussion carries over verbatim to a pair of complex vector bundles
$E$ and $F$ over~$\Sigma$, giving a bundle isomorphism
$$
(\pi_+,\pi_-) : E \otimes F \to (E \otimes_\CC F) \oplus (E \otimes_\CC \widebar{F}).
$$
The Petri map then fits into a commutative diagram
\begin{equation}
\label{eqn:PetriCommute}
\begin{tikzcd}[column sep=large]
\Gamma(E) \otimes \Gamma(F) \ar[r,"\Pi"] \ar[d,"\cong"] & \Gamma(E \otimes F) \ar[d,"\cong"] \\
\left( \Gamma(E) \otimes_\CC \Gamma(F)\right) \oplus \left(\Gamma(E) \otimes_\CC \Gamma(\widebar{F})\right) \ar[r,"\Pi_\CC \oplus \Pi_\CC"] & \Gamma\big( (E \otimes_\CC F) \oplus (E \otimes_\CC \widebar{F}) \big),
\end{tikzcd}
\end{equation}
where $\Pi_\CC : \Gamma(E) \otimes_\CC \Gamma(F) \to \Gamma(E \otimes_\CC F)$
denotes the obvious \emph{complex-linear Petri map} that is defined for any
two complex vector bundles.
Suppose in particular that $E$ and $F$ are line bundles and we have 
chosen complex trivializations for both over some region~$\uU$.
The bundle $\widebar{F}$ inherits from this a trivialization
over $\uU$ such that the canonical map $F \to \widebar{F}$ looks like complex conjugation,
and $E \otimes_\CC F$ and $E \otimes_\CC \widebar{F}$ likewise inherit
natural trivializations.
The diagram now allows us to identify the real Petri map with
\begin{equation}
\label{eqn:PetriComplex}
\begin{split}
C^\infty(\uU,\CC) \otimes C^\infty(\uU,\CC) &\to
C^\infty(\uU,\CC) \oplus C^\infty(\uU,\CC),\\
f \otimes g &\mapsto (fg,f\bar{g}).
\end{split}
\end{equation}

\subsubsection{The main rank estimate}
\label{sec:mainRank}

Fix a holomorphic coordinate chart near $p \in \Sigma$ and a corresponding complex
local trivialization of $E$ such that the Hermitian bundle metrics on $T\Sigma$ and
$E$ both match the standard Hermitian inner product
at~$p$.  The bundle $F$ naturally inherits from these choices a local trivialization in which its
Hermitian bundle metric also appears standard at~$p$.  These choices
identify elements of $J^k_p E$ with polynomials in $z$ and $\bar{z}$,
$$
\sum_{j+\ell \le k} z^j \bar{z}^\ell c_{j,\ell}, \qquad c_{j,\ell} \in \CC^m,
$$
hence
\begin{equation}
\label{eqn:dimJet}
\dim_\CC J^k_p E = \dim_\CC J^k_p F = m \left( 1 + 2 + \ldots + (k+1)\right) =
\frac{m(k+1)(k+2)}{2}.
\end{equation}
Every $\mathbf{D} \in \diffaff_p(E,F)$ is now identified with an
operator of the form
$$
\mathbf{D} = \dbar + A : C^\infty(\DD_\epsilon,\CC^m) \to C^\infty(\DD_\epsilon,\CC^m),
$$
where $\dbar = \p_s + i\p_t$ and $A : \DD_\epsilon \to \End_\RR(\CC^m)$.
The operator $\mathbf{D}_0$ obtained by rescaling as in
\S\ref{sec:rescaling} is then simply
$$
\mathbf{D}_0 = \dbar := \p_s + i\p_t = 2\frac{\p}{\p\bar{z}},
$$
and since the rescaled bundle
metrics $g_0$, $h_0$ and area form $\mu_0$ are all standard in these
coordinates, the formal adjoint of $\mathbf{D}_0$ with respect to this
geometric data is
$$
\mathbf{D}_0^* = -\p = -(\p_s - i \p_t) = - 2 \frac{\p}{\p z}.
$$
We can therefore choose right-inverses
$\mathbf{T}_0 : J^{k-1}_p F \to J^k_p E$ and
$\mathbf{T}^*_0 : J^{k-1}_p E \to J^k_p F$ that are uniquely determined
in coordinates by the conditions
\begin{equation}
\label{eqn:T0}
\mathbf{T}_0 \left( z^j \bar{z}^\ell c\right) := \frac{1}{2(\ell+1)} z^j \bar{z}^{\ell+1} c, \qquad
0 \le j+\ell \le k-1, \quad c \in \CC^m,
\end{equation}
and
\begin{equation}
\label{eqn:T0star}
\mathbf{T}^*_0 \left( z^j \bar{z}^\ell c\right) := - \frac{1}{2(j+1)} z^{j+1} \bar{z}^\ell c, \qquad
0 \le j+\ell \le k-1, \quad c \in \CC^m.
\end{equation}
These choices determine the maps $\mathbf{L}_t : J^{k-1}_p\Hom(E,F) \to J^k_p(E \otimes F)$
in \eqref{eqn:thisIsTheKey}.  Observe now that the domain of this operator has a natural splitting
$$
J^{k-1}_p\Hom(E,F) = J^{k-1}_p\Hom_\CC(E,F) \oplus J^{k-1}_p\overline{\Hom}_\CC(E,F).
$$
If we were to restrict to complex-linear zeroth-order terms $A \in \Gamma(\Hom_\CC(E,F))$,
then the perturbed operators $\mathbf{D} = \mathbf{D}_0 + A$ would always be 
equivalent to $\mathbf{D}_0$ under changes of trivialization, killing any hope that
$\Petrimod^k_{r,\ell}(g_0,h_0,\mu_0,\mathbf{D})$ might be a smaller space than $\Petrimod^k_{r,\ell}(g_0,h_0,\mu_0,\mathbf{D}_0)$.
For this reason, we shall restrict $\mathbf{L}_t$ to the complementary subspace
consisting of $(k-1)$-jets of \emph{antilinear} perturbations.
Having done this, the following additional detail becomes relevant:
for $A \in J^{k-1}_p\overline{\Hom}_\CC(E,F)$ and $t = \sum_j \eta_j \otimes \xi_j \in \ker\mathbf{D}_0 \otimes \ker\mathbf{D}_0^*$,
the commutative diagram \eqref{eqn:PetriCommute} implies
\begin{equation*}
\begin{split}
\pi_- \circ \mathbf{L}_t(A) &= \Pi^k_\CC \circ \pi_- \left( \sum_j \left( \mathbf{T}_0 A \eta_j \otimes \xi_j + \eta_j \otimes \mathbf{T}_0^*A^*\xi_j\right)\right) 
\\
&= \Pi^k_\CC \sum_j\left( \mathbf{T}_0 A \eta_j \otimes_\CC \bar{\xi}_j + 
\eta_j \otimes_\CC \overline{\mathbf{T}_0^*A^*\xi_j} \right),
\end{split}
\end{equation*}
where $\Pi^k_\CC$ denotes the map induced on $k$-jets by the complex Petri map~$\Pi_\CC$.
Since $\mathbf{T}_0$ and $\mathbf{T}_0^*$ are complex linear while 
$A$ and $\xi_j \mapsto \bar{\xi}_j$ are antilinear, the expression 
on the right hand side is the result of applying some real-linear map to
$\pi_+(t) = \sum_j \eta_j \otimes_\CC \xi_j$; the point here is that
real-linear operators of the form $\phi \otimes \psi$ are well defined on 
the complex tensor product whenever $\phi$ and $\psi$ are either both
complex linear or both complex antilinear.  But as mentioned in
Example~\ref{ex:complexPetri}, $\mathbf{D}_0$ satisfies the complex Petri
condition, so the fact that $\Pi^k(t)=0$ implies that $\Pi^k_\CC \circ \pi_+(t) = 0$
and thus $\pi_+(t)=0$, so that the expression vanishes automatically.
We conclude from this discussion that all interesting information in
$\mathbf{L}_t$ is carried by the map
\begin{equation}
\label{eqn:Lhat}
\widehat{\mathbf{L}}_t := \left.\pi_+ \circ \mathbf{L}_t\right|_{J^{k-1}_p\overline{\Hom}_\CC(E,F)} :
J^{k-1}_p\overline{\Hom}_\CC(E,F) \to J^k_p(E \otimes_\CC F).
\end{equation}
Clearly the rank of $\widehat{\textbf{L}}_t$ gives a lower bound for the rank
of~$\mathbf{L}_t$.  The workhorse result behind
Theorem~\ref{thm:Petri} is now the following: 

\begin{prop}
\label{prop:workhorse}
For every $\ell \in \NN$, there exists a constant $C_\ell > 0$ that depends on
$\ell$ but not on~$k$, such that for all
$t \in \ker \Pi^k \subset \ker\mathbf{D}_0 \otimes \ker\mathbf{D}^*_0$ that are homogeneous
elements of degree less than $\ell$ in $J^k_p E \otimes J^k_p F$, the operator
$\widehat{\mathbf{L}}_t : J^{k-1}_p\overline{\Hom}_\CC(E,F) \to J^k_p(E \otimes_\CC F)$ satisfies
$$
\rank \widehat{\mathbf{L}}_t \ge C_\ell k^2.
$$
\end{prop}

\begin{lemma}
If Proposition~\ref{prop:workhorse} holds in the case $\rank_\CC E = 1$,
then it holds in general.
\end{lemma}
\begin{proof}
For $\rank_\CC E = m \in \NN$,
the chosen trivializations furnish local splittings $E = E_1 \oplus \ldots \oplus E_m$ 
and $F = F_1 \oplus \ldots \oplus F_m$ that are
respected by $\mathbf{D}_0$ and $\mathbf{D}_0^*$, i.e.~both are $m$-fold
direct sums of identical operators given by $\dbar$ or $-\p$ respectively.
Their chosen right-inverses $\mathbf{T}_0$ and $\mathbf{T}_0^*$ also
respect these splittings.
Let us denote the resulting splittings of the kernels by
$\ker\mathbf{D}_0 = K_1 \oplus \ldots \oplus K_m$ and
$\ker\mathbf{D}_0^* = L_1 \oplus \ldots \oplus L_m$, so that
$\ker\mathbf{D}_0 \otimes \ker\mathbf{D}^*_0$ splits into $m^2$ identical factors
of the form $K_i \otimes L_j$.  Similarly, $J^k_p(E \otimes F)$ splits into
$m^2$ identical factors of the form $J^k_p(E_i \otimes F_j)$, and the
Petri map $\Pi^k : J^k_p E \otimes J^k_p F \to J^k_p(E \otimes F)$ sends
$J^k_p E_i \otimes J^k_p F_j$ to $J^k_p (E_i \otimes F_j)$ for every $i$ and~$j$.
A homogeneous element $t \in \ker \Pi^k \subset \ker\mathbf{D}_0 \otimes \ker\mathbf{D}_0^*$
of degree $q < \ell$ is now defined by its $m^2$ components
$t_{ij} \in \ker \Pi^k \cap ( K_i \otimes L_j)$, at least one of which must be a nontrivial
homogeneous element of degree~$q$; call this component~$t_{uv}$.  Now consider
the restriction of $\widehat{\mathbf{L}}_t$ to the subspace
$$
J^{k-1}_p\overline{\Hom}_\CC(E_u,F_v) \subset J^{k-1}_p\overline{\Hom}_\CC(E,F),
$$
defined as the $(k-1)$-jets of bundle maps $A : E \to F$ that annihilate
$E_i$ for all $i \ne u$ and have image in~$F_v$.  Since the bundle metrics
$g_0$ and $h_0$ are standard in our chosen trivializations, $A^*$ then
belongs to the corresponding subspace $J^{k-1}_p\overline{\Hom}_\CC(F_v,E_u)
\subset J^{k-1}_p\overline{\Hom}_\CC(F,E)$.  Composing our restriction
of $\widehat{\mathbf{L}}_t$ with the natural projection
$J^k_p(E \otimes F) \to J^k_p(E_u \otimes F_v)$ then produces an operator
$J^{k-1}_p\overline{\Hom}_\CC(E_u,F_v) \to J^k_p(E_u \otimes_\CC F_v)$ that
matches the rank~$1$ case of $\widehat{\mathbf{L}}_t$, and its rank gives a
lower bound for the rank of~$\widehat{\mathbf{L}}_t$.
\end{proof}

The remainder of this subsection is devoted to proving the $\rank_\CC E = 1$
case of Proposition~\ref{prop:workhorse}.  

We shall write everything in the chosen
coordinates and trivializations so that elements of
$J^k_p E$, $J^k_p F$ and $J^k_p(E \otimes_\CC F)$ are now identified with 
complex-valued polynomials of degree
at most $k$ in the variables $z$ and~$\bar{z}$.  The holomorphic
polynomials form $\ker \mathbf{D}_0$, while the antiholomorphic polynomials
form $\ker\mathbf{D}^*_0$.  Using \eqref{eqn:PetriComplex} to compute the
kernel of the Petri map, it turns out that arbitrary elements of
$\ker \Pi^k \subset J^k_p E \otimes J^k_p F$ now take the form\footnote{This seems
a good moment to remind the reader that all tensor products in this section
are \emph{real} tensor products unless the symbol ``$\otimes_\CC$'' is used.}
$$
t = \sum_{j,n=0}^k \left[ a_{jn}\left( z^j \otimes \bar{z}^n + iz^j \otimes i\bar{z}^n \right)
 + b_{jn} \left( i z^j \otimes \bar{z}^n - z^j \otimes i\bar{z}^n \right) \right]  + R,
$$
where $a_{jn},b_{jn} \in \RR$ are real coefficients subject to the condition
$\sum_{j+n=q} a_{jn} = \sum_{j+n=q} b_{jn} = 0$ for every $q=0,\ldots,k$,
and $R$ is an arbitrary sum of homogeneous elements that have degrees
greater than $k$ and therefore vanish automatically under~$\Pi^k$.
For Proposition~\ref{prop:workhorse} we are interested only in homogeneous
elements of some degree less than~$\ell$, so let us fix an integer $q \le \ell$ and write
$$
t = \sum_{j=0}^{q-1} \left[ a_j \left( z^j \otimes \bar{z}^{q-1-j} + i z^j \otimes i \bar{z}^{q-1-j} \right)
 + b_j \left( i z^j \otimes \bar{z}^{q-1-j} - z^j \otimes i \bar{z}^{q-1-j} \right) \right],
$$
where $a_j,b_j \in\RR$ are now subject to the conditions
$\sum_{j=0}^{q-1} a_j = \sum_{j=0}^{q-1} b_j = 0$ and we explicitly assume
that at least one of these coefficients is nonzero.
The action of an antilinear bundle map $A \in \Gamma(\overline{\Hom}_\CC(E,F))$ 
on a section $\eta \in\Gamma(E)$ can be written in trivializations as
$$
(A\eta)(z) := \alpha(z) \widebar{\eta(z)}
$$
for some complex-valued function~$\alpha$, thus the 
map $A : J^k_p E \to J^{k-1}_p F$ can be written as
$$
A \eta = \sum_{u+v \le k-1} \alpha_{uv} z^u \bar{z}^v \bar{\eta}
$$
for some coefficients $\alpha_{uv} \in \CC$.  The transpose 
$A^* : J^k_p F \to J^{k-1}_p E$ is given by exactly the same formula---here 
we are taking transposes of the $1$-by-$1$ matrices $\alpha_{uv}$ and
thus leaving them unchanged, as the antilinearity of $A$ makes the
transpose the appropriate transformation here instead of the Hermitian
adjoint.  With this data in place and the explicit formulas given in
\eqref{eqn:T0} and \eqref{eqn:T0star} for $\mathbf{T}_0$ and $\mathbf{T}^*_0$,
we now obtain an explicit formula for $\widehat{\mathbf{L}}_t(A) \in J^k_p(E \otimes_\CC F)$ as
$$
\widehat{\mathbf{L}}_t(A) = \sum_{j=0}^{q-1} \sum_{u+v \le k-q}
\left( \frac{\bar{c}_j \alpha_{uv}}{v+j+1} z^u \bar{z}^{v+q} -
\frac{c_j \alpha_{uv}}{u+q-j} z^{u+q}\bar{z}^v \right),
$$
where we have defined
$$
c_j := a_j + i b_j \in \CC \quad\text{ for }\quad j=0,\ldots,q-1.
$$
Two immediate remarks are in order: first, the second summation in this
formula stops at $k-q$ instead of $k-1$ because all terms in $A$ with
degree larger than $k-q$ produce terms in $\widehat{\mathbf{L}}_t(A)$ that
have degree greater than $k$ and thus vanish in $J^k_p(E \otimes_\CC F)$.
Along the same lines, we notice that whenever $A$ is given by a homogeneous
polynomial of degree $n$, $\widehat{\mathbf{L}}_t(A)$ is likewise
homogeneous with degree~$n+q$, indicating a natural splitting of the
map $\widehat{\mathbf{L}}_t : J^{k-1}_p\overline{\Hom}_\CC(E,F) \to
J^k_p(E \otimes_\CC F)$ into factors
$$
\widehat{\mathbf{L}}_t = \widehat{\mathbf{L}}_t^{(0)} \oplus \ldots \oplus \widehat{\mathbf{L}}_t^{(k-q)},
$$
where for each $n=0,\ldots,k-q$, $\widehat{\mathbf{L}}_t^{(n)}$ is defined
on the space of homogeneous degree~$n$ polynomials in $J^{k-1}_p\overline{\Hom}_\CC(E,F)$.
(Strictly speaking, there are additional factors defined on homogeneous polynomials
of higher degree, but we will ignore them because they are trivial.)

For each individual $n \in \{0,\ldots,k-q\}$, the map $\widehat{\mathbf{L}}_t^{(n)}$
takes the form
$$
\widehat{\mathbf{L}}_t^{(n)}(A) = \sum_{u+v=n} \alpha_{uv} \cdot \left[ 
\left( \sum_{j=0}^{q-1} \frac{\bar{c}_j}{v+j+1} \right) z^u \bar{z}^{v+q} -
\left( \sum_{j=0}^{q-1} \frac{c_j}{u+q-j} \right) z^{u+q} \bar{z}^v \right].
$$
To simplify this expression, we can write $\mathbf{c} = (c_0,\ldots,c_{q-1}) \in \CC^q$
as a column vector and define for integers $u,v \ge 0$ the complex numbers
$$
\theta_v := \begin{pmatrix} \frac{1}{v+1} & \cdots & \frac{1}{v+q} \end{pmatrix} \bar{\mathbf{c}}
\quad\text{ and }\quad
\kappa_u := \begin{pmatrix} \frac{1}{u+q} & \cdots & \frac{1}{u+1} \end{pmatrix} \mathbf{c},
$$
so that now
$$
\widehat{\mathbf{L}}_t^{(n)}(A) = \sum_{u+v=n} \alpha_{uv} \cdot 
\left( \theta_{v} z^u \bar{z}^{v+q} - \kappa_u z^{u+q} \bar{z}^v \right).
$$
If we now identify the homogeneous degree~$n$ part of $A$ with the vector 
in $\CC^{n+1}$ given by $(\alpha_{n,0},\alpha_{n-1,1},\ldots,\alpha_{0,n})$,
and use the monomials 
$$
z^{n+q}, z^{n+q-1}\bar{z}, z^{n+q-2}\bar{z}^2,\ldots, z \bar{z}^{n+q-1}, \bar{z}^{n+q}
$$
as a complex basis for the homogeneous degree $n+q$ part of $J^k_p(E \otimes_\CC F)$,
then $\widehat{\mathbf{L}}_t^{(n)}$ is represented by the $(n+q+1)$-by-$(n+1)$ complex matrix
\begin{equation}
\label{eqn:Lmatrix}
\widehat{\mathbf{L}}_t^{(n)} =
\begin{pmatrix}
-\kappa_n &               &          &           \\
\vdots    & -\kappa_{n-1} &          &           \\
\theta_0  &  \vdots       & \ddots   &           \\
          & \theta_1      &          & -\kappa_0 \\
          &               & \ddots   &  \vdots   \\
	  &               &          & \theta_n   \\
\end{pmatrix}.
\end{equation}
In this matrix, all entries not written explicitly are understood to be~$0$.

\begin{lemma}
\label{lemma:invertible}
For any set of distinct positive integers $i_1,\ldots,i_q$, the matrix
$$
\begin{pmatrix}
\dfrac{\strut 1}{\strut i_1 + q} & \cdots & \dfrac{\strut 1}{\strut i_1 + 1} \\
\vdots                           & \ddots & \vdots                           \\
\dfrac{\strut 1}{\strut i_q + q} & \cdots & \dfrac{\strut 1}{\strut i_q + 1} 
\end{pmatrix}
$$
is invertible.
\end{lemma}
\begin{proof}
This follows from the well-known formula for so-called \emph{Cauchy determinants},
$$
\det 
\begin{pmatrix}
\dfrac{\strut 1}{\strut z_1+w_1} & \cdots & \dfrac{\strut 1}{\strut z_1 + w_q} \\
\vdots            & \ddots & \vdots              \\
\dfrac{\strut 1}{\strut z_q+w_1} & \cdots & \dfrac{\strut 1}{\strut z_q + w_q}
\end{pmatrix} = 
\frac{\displaystyle \prod_{i=1}^q \prod_{j=i}^{i-1} (z_i - z_j) (w_i - w_j)}{\displaystyle \prod_{i,j=1}^q (z_i + w_j)},
$$
see e.g.~\cite{PolyaSzego2}*{pp.~92 and~279}.
\end{proof}

Since at least one of the coefficients $a_j$ or $b_j$ is nonzero,
the vector $\mathbf{c} \in \CC^q$ cannot be annihilated by $q$ linearly
independent vectors, so we conclude:

\begin{cor}
In the matrix \eqref{eqn:Lmatrix}, at most $q-1$ of the entries
$\kappa_0,\ldots,\kappa_n$ can be zero.
\qed
\end{cor}

This result implies that at most $q-1$ columns of the matrix \eqref{eqn:Lmatrix}
need to be eliminated in order to produce a matrix whose columns are all
linearly independent, hence if $n \ge q-1$, we have
$$
\rank_\CC \widehat{\mathbf{L}}_t^{(n)} \ge n - (q - 1).
$$
If $k \ge 2q$, then summing this estimate for $n=q,\ldots,k-q$ gives
$$
\rank_\CC \widehat{\mathbf{L}}_t \ge 1 + 2 + \ldots + k - 2q + 1 =
\frac{1}{2}(k-2q+1)(k-2q+2),
$$
and thus 
$$
\rank \widehat{\mathbf{L}}_t \ge (k-2q+1)(k-2q+2) \ge (k-2\ell+1)(k-2\ell+2)
$$
whenever $k \ge 2\ell$.  This estimate might not be satisfied for
$k$ underneath this threshold, but since that is only finitely many cases,
we can now just choose a constant $C_\ell > 0$ small enough to achieve
$C_\ell k^2 \le \rank\widehat{\mathbf{L}}_t$ for those cases and
$C_\ell k^2 \le (k-2\ell+1)(k-2\ell+2)$ for all $k \ge 2\ell$.
With this, the proof of Proposition~\ref{prop:workhorse} is complete.

\subsubsection{Proof of Theorem~\ref{thm:Petri}}

Consider the $C^\infty$-subvarieties $\Petrimod^k_{r,\ell}(g,h,\mu) \subset
\vV^k_{r,\ell}(g,h,\mu)$ from \S\ref{sec:universalPetri}
in the specific setting of local Cauchy-Riemann type
operators $\mathbf{D} \in \diffaff_p(E,F)$ with $\rank_\CC E = m$.
For any given operator $\mathbf{D} \in \diffaff^k_p(E,F)$, 
we know from Lemma~\ref{lemma:itsABundle} that
$\mathbf{D} : J^k_p E \to J^{k-1}_p F$ and
$\mathbf{D}^* : J^k_p F \to J^{k-1}_p E$ are both surjective, thus
\eqref{eqn:dimJet} gives
\begin{equation*}
\begin{split}
\dim \ker \mathbf{D} = \dim \ker \mathbf{D}^* &= \dim J^k_p E - \dim J^{k-1}_p F\\
&= m(k+1)(k+2) - mk(k+1) = 2m(k+1),
\end{split}
\end{equation*}
and plugging this into \eqref{eqn:codimRank},
\begin{equation}
\label{eqn:dimVkr}
\dim \vV^k_{r,\ell}(g,h,\mu,\mathbf{D}) = 4 r m(k+1) - r^2.
\end{equation}
Next, combining Proposition~\ref{prop:workhorse} with Lemma~\ref{lemma:thisIsTheKey}
gives:

\begin{prop}
\label{prop:PetriMain}
For every $\ell \in \NN$, there exists a constant $C_\ell > 0$ such that
for all integers $k \ge \ell$ and all $r \in \NN$,
$\Petrimod^k_{r,\ell}(g,h,\mu) \subset \vV^k_{r,\ell}(g,h,\mu)$
is a $C^\infty$-subvariety of codimension at least $C_\ell k^2$.
\qed
\end{prop}

Sard's theorem (see Propsition~\ref{prop:SardSmaleVariety}) now provides
a Baire subset 
$$
\diffaff^{k,\reg}_p(E,F\,;\,r,\ell) \subset \diffaff^k_p(E,F)
$$
such that for all $\mathbf{D} \in \diffaff^{k,\reg}_p(E,F\,;\,r,\ell)$,
$\Petrimod^k_{r,\ell}(g,h,\mu,\mathbf{D})$ is a $C^\infty$-subvariety in
$\vV^k_{r,\ell}(g,h,\mu,\mathbf{D})$ of codimension at least $C_\ell k^2$.
Since this codimension grows quadratically with $k$ while the dimension of
$\vV^k_{r,\ell}(g,h,\mu,\mathbf{D})$ grows only linearly, we conclude
that for any fixed $r,\ell \in \NN$, the space $\Petrimod^k_{r,\ell}(g,h,\mu,\mathbf{D})$
is empty for all $k$ sufficiently large.

To conclude the proof of Theorem~\ref{thm:Petri}, we choose for each
$\ell \in \NN$ some $k \ge \ell$ large enough so that 
$\Petrimod^k_{\ell,\ell}(g,h,\mu,\mathbf{D}) = \emptyset$ for every
$\mathbf{D} \in \diffaff^{k,\reg}_p(E,F\,;\,\ell,\ell)$, and then define
$\CRR^{\ell,\reg}(E\,;\,\uU,\mathbf{D}\fix)$ to be the set of all
operators in $\CRR(E\,;\,\uU,\mathbf{D}\fix)$ whose $k$-jets at $p$
belong to $\diffaff^{k,\reg}_p(E,F\,;\,\ell,\ell)$.

\subsection{Petri's condition is satisfied for generic~$J$}
\label{sec:PetriJ}

We now return to the setting of \S\ref{sec:mainIdea} and
consider the moduli space
$\mM_{g}(A,J)$ of unparametrized closed $J$-holomorphic curves
$u : (\Sigma,j) \to (M,J)$ of genus $g \ge 0$ homologous to $A \in H_2(M)$ in a
symplectic manifold $(M,\omega)$ of dimension $2n \ge 4$ with
$J \in \jJ(M,\omega\,;\,\uU,\Jfix)$.  Here $\uU \subset M$ is an open
subset with compact closure, $\Jfix$ is a fixed compatible almost complex
structure, and all $J \in \jJ(M,\omega\,;\,\uU,\Jfix)$ are assumed to
match $\Jfix$ outside of~$\uU$.

\begin{thm}
\label{thm:PetriJ}
There exists a Baire subset $\jJ^\reg \subset \jJ(M,\omega\,;\,\uU,\Jfix)$
such that for all $J \in \jJ^\reg$ and every $u \in \mM_g(A,J)$
with parametrization $u : (\Sigma,j) \to (M,J)$, the normal Cauchy-Riemann
operator $\mathbf{D}_u^N \in \CRR(N_u)$ satisfies Petri's condition to infinite order
on an open and dense set of points in $u^{-1}(\uU)$.  In particular,
$\mathbf{D}_u^N$ satisfies the local Petri condition at every point
in~$u^{-1}(\uU)$ (cf.~Remark~\ref{remark:openDensePetri}).
\end{thm}

We will deduce Theorem~\ref{thm:PetriJ} from the results of the previous
subsection after showing essentially that the natural map from the
universal moduli space of simple holomorphic curves with one marked point to the space of $k$-jets of 
normal Cauchy-Riemann operators at the marked point is always a submersion.  Up to some technical details
still to be addressed, the next lemma implies this.
Recall that a point $z \in \Sigma$ in the domain of a smooth map
$v : \Sigma \to M$ is called an \defin{injective point} if $dv(z) :
T_z\Sigma \to T_{v(z)}M$ is injective and $\{z\} = v^{-1}(v(z))$.
For a simple $J$-holomorphic curve, the complement of the set of injective 
points is a discrete set.

\begin{lemma}
\label{lemma:operatorPert}
Assume $J \in \jJ(M,\omega\,;\,\uU,\Jfix)$, and $v : (\Sigma,j) \to (M,J)$ is a simple
$J$-holomorphic curve with generalized normal bundle $N_v \subset v^*TM$
defined as the $\omega$-symplectic complement of the generalized tangent
bundle $T_v \subset v^*TM$.
Given any $A \in \Omega^{0,1}(\Sigma,\End_\RR(N_v))$ with support
contained in the set of injective
points in $v^{-1}(\uU)$, there exists a smooth family of almost complex structures 
$$
\{ J_\tau \in \jJ(M,\omega\,;\,\uU,\Jfix) \}_{\tau \in (-\epsilon,\epsilon)}
$$
such that $J_0 = J$, $J_\tau(v(z)) = J(v(z))$ for all $\tau$ and~$z$,
and the resulting family of normal Cauchy-Riemann operators
$\mathbf{D}_{v,\tau}^N \in \CRR(N_v)$ for $v$ defined with respect to $J_\tau$ satisfies
$$
\left. \p_\tau \mathbf{D}_{v,\tau}^N \eta\right|_{\tau=0} = \pi_N \circ \nabla_{\eta} Y \circ Tv \circ j = A \eta
$$
for $\eta \in \Gamma(N_v)$, where $Y := \left.\p_\tau J_\tau\right|_{\tau=0} \in \Gamma(\overline{\End}_\CC(TM,J))$,
$\nabla$ is any connection on~$M$, and $\pi_N : v^*TM \to N_v$ denotes the projection along~$T_v$.
\end{lemma}
\begin{proof}
If $\{J_\tau\}$ is any smooth path in $\jJ(M,\omega\,;\,\uU,\Jfix)$ with $J_0 = J$,
$J_\tau(v) \equiv J(v)$ for all $\tau$ and $Y := \p_\tau J|_{\tau=0}$, then $Y(v) \equiv 0$,
hence $\nabla Y$ is well defined along $v$ independently of any connection.
For $\eta \in \Gamma(N_v)$, let us write $\nabla_\eta Y$ in block form as
\begin{equation}
\label{eqn:tangentNormal}
\nabla_\eta Y = \begin{pmatrix}
\nabla^T_\eta Y & \nabla^{TN}_\eta Y \\
\nabla^{NT}_\eta Y & \nabla^N_\eta Y
\end{pmatrix} \in \Gamma(\overline{\End}_\CC(v^*TM,J))
\end{equation}
with respect to the tangent-normal decomposition $v^*TM = T_v \oplus N_v$.
Since $N_v$ is the $\omega$-symplectic orthogonal
complement of~$T_v$, the fact that $J_\tau$ is always $\omega$-compatible then 
translates into conditions that constrain $\nabla^T_\eta Y$ and
$\nabla^N_\eta Y$ separately and another condition that determines
$\nabla^{TN}_\eta Y$ in terms of $\nabla^{NT}_\eta Y$, namely
$$
\omega( (\nabla^{NT}_\eta Y) v , w) + \omega(v , (\nabla^{TN}_\eta Y) w) = 0
$$
for all $(v,w) \in T_v \oplus N_v$.  This means that $\omega$-compatibility
does not prevent us from freely choosing $\nabla^{NT}_\eta Y$ so long as
we (1)~do not mind $\nabla^{TN}_\eta Y$ being determined by this choice, and
(2)~do this only in regions where $v$ has no double points, so that the
splitting of $TM$ into $T_v \oplus N_v$ is unambiguous.  Now using the
definition of the normal Cauchy-Riemann operator, one computes that for
any $\eta \in \Gamma(N_v)$,
$$
\left. \p_\tau \mathbf{D}_{v,\tau}^N \eta\right|_{\tau=0} = \nabla^{NT}_\eta Y \circ Tv \circ j.
$$
On a region where $v$ has neither critical points nor double points and its
image lies in the perturbation domain~$\uU$, we can therefore choose the
normal derivatives of $Y$ along $v$ to make the above expression match~$A$.
\end{proof}

To prove Theorem~\ref{thm:PetriJ}, we will use
the Floer $C_\varepsilon$-topology (cf.~\cite{Floer:action}*{\S 5})
to define spaces of perturbed data.  
Given any $\Jref \in \jJ(M,\omega\,;\,\uU,\Jfix)$, we define
$$
T_{\Jref}\jJ(J,\omega\,;\,\uU,\Jfix) \subset \Gamma(\overline{\End}_\CC(TM,\Jref))
$$
as the space of smooth $\Jref$-antilinear bundle maps $Y$
that vanish outside~$\uU$ and satisfy $\omega(\cdot,Y\cdot) + \omega(Y\cdot,\cdot) \equiv 0$;
intuitively, this is the tangent space at $\Jref$ to the smooth
Fr\'echet manifold $\jJ(J,\omega\,;\,\uU,\Jfix)$.  There is a natural embedding
\begin{equation}
\label{eqn:expJ}
Y \mapsto J_Y := \left( \1 + \frac{1}{2} \Jref Y \right) \Jref \left(\1 + \frac{1}{2} \Jref Y \right)^{-1}
\end{equation}
which takes a $C^0$-small neighborhood of $0$ in $T_{\Jref}\jJ(J,\omega\,;\,\uU,\Jfix)$
homeomorphically to a neighborhood of $\Jref$ in $\jJ(M,\omega\,;\,\uU,\Jfix)$.
Now choose a Riemannian metric on $M$ in order to define the
$C^\nu$-norms on $\Gamma(\overline{\End}_\CC(TM,\Jref))$ for each integer $\nu \ge 0$,
fix a sequence of positive numbers $\varepsilon_\nu \to 0$, and define the
$C_\varepsilon$-norm
\begin{equation}
\label{eqn:FloerNorm}
\| Y \|_{C_\varepsilon} := \sum_{\nu=0}^\infty \varepsilon_\nu
\| Y \|_{C^\nu}
\end{equation}
for $Y \in \Gamma(\overline{\End}_\CC(TM,\Jref))$.
Fixing any $\delta > 0$ sufficiently small, this gives rise to
a smooth, separable and metrizable Banach manifold
$$
\jJ_\varepsilon := \left\{ J_Y\ \big|\ 
\text{$Y \in T_{\Jref}\jJ(J,\omega\,;\,\uU,\Jfix)$, $\| Y \|_{C_\varepsilon} < \infty$ and $\|Y\|_{C^0} < \delta$}
\right\}
$$
which embeds continuously into $\jJ(J,\omega\,;\,\uU,\Jfix)$ and contains
arbitrarily $C^\infty$-small perturbations of~$\Jref$.
Note that since $\uU \subset M$ has compact closure, the equivalence classes of 
the individual $C^\nu$-norms are each independent of auxiliary choices such as connections or
local trivializations, but the equivalence class of the $C_\varepsilon$-norm 
may in fact depend on these choices.  This is immaterial, as the choice of the
sequence $\{\epsilon_\nu\}_{\nu=0}^\infty$ carries no geometric meaning in itself;
what is important is rather that the space of sections of class $C_\varepsilon$ can always be
enlarged by making $\varepsilon_\nu$ converge to $0$ faster.  To say this more
precisely, let us endow the set
$$
\seqs := \left\{ \text{sequences } \varepsilon = \{\varepsilon_\nu\}_{\nu=0}^\infty \ \Big|\ 
\varepsilon_\nu > 0 \text{ for all $\nu$, and }
\lim_{\nu \to \infty} \varepsilon_\nu = 0 \right\}
$$
with a pre-order $\prec$ defined by
$$
\varepsilon \prec \varepsilon' \qquad \Longleftrightarrow \qquad
\limsup_{\nu \to\infty} \frac{\varepsilon_\nu}{\varepsilon_\nu'} < \infty.
$$
\begin{defn}
Given a statement $S(\varepsilon)$ dependent on a choice of $\varepsilon \in \seqs$,
we will say that $S(\varepsilon)$ holds for all $\varepsilon \in \seqs$ 
\defin{with sufficiently rapid decay} if there exists $\varepsilon_0 \in \seqs$ such
that $S(\varepsilon)$ holds for all $\varepsilon \prec \varepsilon_0$.
\end{defn}

\begin{lemma}
\label{lemma:CepsFreedom}
The $C_\varepsilon$-norms on sections
$Y \in T_{\Jref}\jJ(J,\omega\,;\,\uU,\Jfix)$ have the following properties:
\begin{enumerate}
\item If $\varepsilon \prec \varepsilon'$ in $\seqs$, then there exists a constant
$c > 0$ such that $\|Y\|_{C_\varepsilon} \le c \|Y\|_{C_{\varepsilon'}}$
for all~$Y$.
\item For any given $Y$, $\|Y\|_{C_\varepsilon} < \infty$ for all $\varepsilon \in \seqs$
with sufficiently rapid decay.
\item
Every countable subset of $\seqs$ has a lower bound in~$\seqs$ with respect
to the pre-order~$\prec$.
\end{enumerate}
\end{lemma}
\begin{proof}
Property~(1) follows easily from the observation that $\varepsilon \prec \varepsilon'$
if and only if there exist constants $C > 0$ and $\nu_0 \in \NN$ such that
$\varepsilon_\nu \le C \varepsilon_\nu'$ for all $\nu > \nu_0$.
To prove~(2), observe that any nontrivial smooth section $Y$ vanishing outside of $\uU$
is of class $C_\varepsilon$ for $\varepsilon_\nu := 1 / \left( 2^\nu \cdot \|Y\|_{C^\nu} \right)$,
then apply~(1).
Finally, $\varepsilon \in \seqs$ is a lower bound for the countable subset
$\{\varepsilon^{(1)},\varepsilon^{(2)},\varepsilon^{(3)},\ldots \} \subset \seqs$ whenever
$\varepsilon_\nu \le \min\left\{\varepsilon_\nu^{(1)},\ldots,\varepsilon_\nu^{(\nu)}\right\}$
for every~$\nu$.
\end{proof}

Let us discuss the geometric data to be used in formulating the local
Petri condition for a holomorphic curve.
Given $J \in \jJ(M,\omega\,;\,\uU,\Jfix)$, the complex vector bundle $(TM,J)$ carries a natural
Hermitian metric whose real part is $g_J := \omega(\cdot,J\cdot)$.
If $u : (\Sigma,j) \to (M,J)$ is $J$-holomorphic and is immersed
at the point $\zeta \in \Sigma$, then $g_J$ can be pulled back to define
a Riemannian metric on $\Sigma$ near $\zeta$ in the conformal class of~$j$, thus
giving rise to an area form $\mu_u$ on $\Sigma$ and compatible bundle metrics
$g_u$ on $N_u$ and $h_u$ on $\overline{\Hom}_\CC(T\Sigma,N_u)$ near~$\zeta$, where for
concreteness we are also free to assume $N_u \subset u^*TM$ is the
$g_J$-orthogonal complement of~$T_u$.  In order to avoid ambiguity, we shall 
assume in the following that $\mathbf{D}_u^N$ and $(\mathbf{D}_u^N)^*$ are
defined via these specific choices of geometric data near any given immersed point $\zeta \in \Sigma$;
note that this would not be a valid global definition for $(\mathbf{D}^N_u)^*$
since the pulled back metric on $\Sigma$ becomes singular at critical points,
but this will not matter since we only intend to study finite jets 
of $(\mathbf{D}^N_u)^*$ at a specific immersed point.  Recall from Remark~\ref{remark:PetriIndep}
that Petri's condition does not depend on choices of geometric data.
Moreover, while the global topological type of $N_u$ may change (because
the number of critical points may change) as $u$ moves about in its moduli 
space, the germs of $\mathbf{D}^N_u$ and $(\mathbf{D}^N_u)^*$ at an immersed
point can still be assumed to depend smoothly on~$u$.

Let us denote by
$$
\mM^*_{g,1}(A,J) \subset \mM_{g,1}(A,J)
$$
the open subset consisting of simple curves with one marked point such
that the marked point is an injective point with image in~$\uU$.
We will abuse notation and write elements of $\mM^*_{g,1}(A,J)$
as $(u,\zeta)$, where $u : (\Sigma,j) \to (M,J)$ is a specific parametrization
and $\zeta \in \Sigma$ is the marked point.  Using the notation of
\S\ref{sec:universalPetri}, we then define for each $k,r,\ell \in \NN$
with $\ell \le k$ the space
$$
\widehat{\mM}^{k,r,\ell}_{g,1}(A,J) := \left\{ (u,\zeta,t)\ \Big|\ 
(u,\zeta) \in \mM^*_{g,1}(A,J),\ t \in \vV^k_{r,\ell}(g_u,h_u,\mu_u,\mathbf{D}^N_u) \right\}
$$
where $g_u,h_u,\mu_u$ are the specific choices of geometric data determined by
$u$ and $g_J$ as described in the
previous paragraph.  The extra term $t$ is an element
in the tensor product of the $k$-jet versions of $\ker\mathbf{D}^N_u$ and
$\ker(\mathbf{D}^N_u)^*$ at~$\zeta$, having rank $r$ and not vanishing to order~$\ell$.
We will be interested especially in the subset
$$
\mM^{k,r,\ell}_{g,1}(A,J) := \left\{ (u,\zeta,t) \in \widehat{\mM}^{k,r,\ell}_{g,1}(A,J)\ \Big|\ 
\Pi^k(t)=0 \right\}.
$$
To understand the structure of these spaces, we define corresponding universal
moduli spaces:
\begin{equation*}
\begin{split}
\univ^*_{g,1}(A,\jJ_\varepsilon) &:= \left\{ (u,\zeta,J)\ \Big|\ J \in \jJ_\varepsilon,\ (u,\zeta) \in \mM^*_{g,1}(A,J) \right\},\\
\widehat{\univ}^{k,r,\ell}_{g,1}(A,\jJ_\varepsilon) &:= \left\{ (u,\zeta,t,J) \ \Big|\ J \in \jJ_\varepsilon,\ (u,\zeta,t) \in \widehat{\mM}^{k,r,\ell}_{g,1}(A,J) \right\},\\
\univ^{k,r,\ell}_{g,1}(A,\jJ_\varepsilon) &:= \left\{ (u,\zeta,t,J)\ \Big|\ J \in \jJ_\varepsilon,\ (u,\zeta,t) \in \mM^{k,r,\ell}_{g,1}(A,J) \right\}.
\end{split}
\end{equation*}
We shall always choose $\varepsilon \in \seqs$ to have sufficiently rapid decay so that,
by standard arguments as in \cite{McDuffSalamon:Jhol2},
$\univ^*_{g,1}(A,\jJ_\varepsilon)$ is a
smooth, metrizable and separable Banach manifold such that the projection
$\univ^*_{g,1}(A,\jJ_\varepsilon) \to \jJ_\varepsilon : (u,\zeta,J) \mapsto J$ is a smooth
Fredholm map whose index is the virtual dimension of
$\mM^*_{g,1}(A,J)$.  It follows that the same is true for
$\widehat{\univ}^{k,r,\ell}_{g,1}(A,\jJ_\varepsilon)$, as the additional
$k$-jet data $t$ varies in a smooth finite-dimensional manifold that
depends smoothly on the $k$-jet of the operator $\mathbf{D}^N_u$
at the immersed point~$\zeta$, and this in turn depends smoothly on
$(u,\zeta,J) \in \univ^*_{g,1}(A,\jJ_\varepsilon)$.

It will be convenient to impose an extra condition defining an open subset
of $\univ^{k,r,\ell}_{g,1}(A,\jJ_\varepsilon)$.  For each $\ell \in\NN$, let
$C_\ell > 0$ denote the constant furnished by Proposition~\ref{prop:PetriMain}
in \S\ref{sec:PetriCR}, with the roles of the bundles $E,F$ and point $p$ in that
subsection played by $N_u$, $\overline{\Hom}_\CC(T\Sigma,N_u)$ and $\zeta \in \Sigma$ 
respectively.

\begin{defn}
\label{defn:epsReg}
Given $J \in \jJ(M,\omega\,;\,\uU,\Jfix)$ and $\varepsilon \in \seqs$, we will say
that an element $(u,\zeta,t) \in \mM^{k,r,\ell}_{g,1}(A,J)$ is
\defin{$\varepsilon$-regular} if $J \in \jJ_\varepsilon$ and 
$(u,\zeta,t,J)$ has a neighborhood $\oO \subset \widehat{\univ}^{k,r,\ell}_{g,1}(A,\jJ_\varepsilon)$
such that $\oO \cap \univ^{k,r,\ell}_{g,1}(A,\jJ_\varepsilon)$ is a
$C^\infty$-subvariety of $\widehat{\univ}^{k,r,\ell}_{g,1}(A,\jJ_\varepsilon)$ with
codimension at least $C_\ell k^2$.
\end{defn}
Note that $\varepsilon$-regularity is an open condition by construction,
i.e.~the set of tuples $(u,\zeta,t,J) \in \univ^{k,r,\ell}_{g,1}(A,\jJ_\varepsilon)$
such that $(u,\zeta,t)$ is $\varepsilon$-regular is open.  
The important consequence of Lemma~\ref{lemma:operatorPert} will
be that it is generally also nonempty.

\begin{lemma}
\label{lemma:epsReg}
Any given $(u,\zeta,t) \in \mM^{k,r,\ell}_{g,1}(A,\Jref)$
is $\varepsilon$-regular for all $\varepsilon \in \seqs$ with
sufficiently rapid decay.
\end{lemma}
\begin{proof}
Observe first that $\Jref \in \jJ_\varepsilon$ for every $\varepsilon \in \seqs$.
Now given $(u,\zeta,t) \in \mM^{k,r,\ell}_{g,1}(A,\Jref)$,
define the Fr\'echet space
$$
\yY_0 := \left\{ Y \in T_{\Jref}\jJ(M,\omega\,;\,\uU,\Jfix)\ \Big|\ Y|_{u(\Sigma)} \equiv 0 \right\}
$$
and for each $\varepsilon \in \seqs$ the Banach space
$$
\yY_\varepsilon := \left\{ Y \in T_{\Jref}\jJ(M,\omega\,;\,\uU,\Jfix)\ \Big|\ Y|_{u(\Sigma)} \equiv 0 \text{ and } \|Y\|_{C_\varepsilon} < \infty \right\},
$$
where the latter is regarded as a closed subspace of $T_{\Jref}\jJ_\varepsilon$ with the
$C_\varepsilon$-topology.  Abbreviating $E := N_u$ and $F := \overline{\Hom}_\CC(T\Sigma,N_u)$,
Lemma~\ref{lemma:operatorPert} provides a surjective linear map
$$
\Psi_0 : \yY_0 \to J^{k-1}_\zeta(\Hom_\RR(E,F)) : Y \mapsto J^{k-1}_\zeta(A_Y),
$$
where $A_Y$ denotes (the germ near $\zeta$ of) the zeroth-order term
determined by $Y$ according to the formula $A_Y\eta = \pi_N \circ \nabla_\eta Y \circ Tu \circ j$.
Since the target space of $\Psi_0$ is finite dimensional, 
Lemma~\ref{lemma:CepsFreedom} implies that it remains surjective when
restricted to the subspace $\yY_\varepsilon$ for $\varepsilon \in \seqs$
with sufficiently rapid decay.
Each $Y \in \yY_\varepsilon$ now gives rise to a $1$-parameter family
of almost complex structures
$J_\tau := J_{\tau Y} \in \jJ_\varepsilon$ defined via \eqref{eqn:expJ},
which match $\Jref$ along $u$ and satisfy $J_0=\Jref$.  This defines
a smooth family $(u,\zeta,J_\tau) \in \univ^*_{g,1}(A,\jJ_\varepsilon)$
that deforms the normal Cauchy-Riemann operator of $u$ in the direction
of $A_Y$ but leaves the geometric data along $u$ unchanged.
It follows that the linearization at $(u,\zeta,t,\Jref)$
of the natural projection map\footnote{Strictly speaking,
the definition of $\vV^k_{r,\ell}$ in this context depends on the germs near $\zeta \in \Sigma$
of the vector bundles $N_u$ and $\overline{\Hom}_\CC(T\Sigma,N_u)$,
which vary as $(u,\zeta,t,J)$ moves in $\widehat{\univ}^{k,r,\ell}_{g,1}(A,\jJ_\varepsilon)$,
so for the purposes of \eqref{eqn:fraglicheProjektion}, $\vV^k_{r,\ell}$ should
be replaced with a suitable fiber bundle over $\widehat{\univ}^{k,r,\ell}_{g,1}(A,\jJ_\varepsilon)$,
of which the map in \eqref{eqn:fraglicheProjektion} is a section.
This detail makes little difference for the present argument, however,
since the family $(u,\zeta,J_\tau) \in \univ^*_{g,1}(A,\jJ_\varepsilon)$
involves a fixed curve with a fixed marked point and $J_\tau|_{\im(u)}$
also fixed.}
\begin{equation}
\label{eqn:fraglicheProjektion}
\widehat{\univ}^{k,r,\ell}_{g,1}(A,\jJ_\varepsilon) \to \vV^k_{r,\ell} :
(u,\zeta,t,J) \mapsto (g_u,h_u,\mu_u,\mathbf{D}^N_u,t)
\end{equation}
is surjective onto $T_{(\mathbf{D}^u_N,t)}\vV^k_{r,\ell}(g_u,h_u,\mu_u)$, and the
result then follows from Proposition~\ref{prop:PetriMain}.
\end{proof}

Applying the Sard-Smale theorem to the projection
$\widehat{\univ}^{k,r,\ell}_{g,1}(A,\jJ_\varepsilon) \to \jJ_\varepsilon :
(u,\zeta,t,J) \mapsto J$
as in Proposition~\ref{prop:SardSmaleVariety}, we can associate to each
$\varepsilon \in \seqs$ and each set of positive integers $k,r,\ell$ with $k \ge \ell$ a Baire subset
$$
\jJ_\varepsilon^\reg(k,r,\ell) \subset \jJ_\varepsilon
$$
such that for all $J \in \jJ_\varepsilon^\reg(k,r,\ell)$,
$\widehat{\mM}^{k,r,\ell}_{g,1}(A,J)$ is a smooth finite-dimensional manifold
and the open set of $\varepsilon$-regular elements in
$$
\mM^{k,r,\ell}_{g,1}(A,J) \subset \widehat{\mM}^{k,r,\ell}_{g,1}(A,J)
$$
is a $C^\infty$-subvariety of codimension at least $C_\ell k^2$.
The dimension of $\widehat{\mM}^{k,r,\ell}_{g,1}(A,J)$ is the Fredholm index of the
projection $\widehat{\univ}^{k,r,\ell}_{g,1}(A,\jJ_\varepsilon) \to \jJ_\varepsilon$, which is larger
than that of $\univ^*_{g,1}(A,\jJ_\varepsilon) \to \jJ_\varepsilon$ by
$\dim \vV^k_{r,\ell}(g_u,h_u,\mu_u,\mathbf{D}_u^N)$.  Plugging in \eqref{eqn:dimVkr},
this gives
$$
\dim \widehat{\mM}^{k,r,\ell}_{g,1}(A,J) = 
\virdim \mM_{g,1}(A,J) + 4r(n-1)(k+1) - r^2.
$$
This number grows linearly with~$k$, while the codimension $C_\ell k^2$ grows
quadratically, thus for any fixed $r,\ell,g,A$, the integer
\begin{equation}
\label{eqn:virdimPetri}
\virdim \mM^{k,r,\ell}_{g,1}(A,J) := \virdim \mM_{g,1}(A,J) + 4r(n-1)(k+1) - r^2 - C_\ell k^2
\end{equation}
becomes negative for all $k \in \NN$ sufficiently large.
Taking the countable intersection of the Baire sets
$\jJ^\reg_\varepsilon(k,r,\ell)$ for all $k,r,\ell,g,A$, we obtain:

\begin{cor}
\label{cor:itsEmpty!}
For every $\varepsilon \in \seqs$, 
there exists a Baire subset $\jJ^\reg_\varepsilon \subset \jJ_\varepsilon$
such that for all $J \in \jJ^\reg_\varepsilon$ and any given
$g \ge 0$, $A \in H_2(M)$ and $r,\ell \in \NN$, the set of
$\varepsilon$-regular elements in
$\mM^{k,r,\ell}_{g,1}(A,J)$ is empty whenever $k$ is large enough for the
integer in \eqref{eqn:virdimPetri} to be negative.
\qed
\end{cor}

For the proof of Theorem~\ref{thm:PetriJ}, we will use a variation on a
popular trick due to Taubes, presenting the desired set $\jJ^\reg \subset \jJ(M,\omega\,;\,\uU,\Jfix)$
as the intersection of an explicit countable collection of open and dense subsets.
This depends on the ability to decompose the relevant moduli space into
a countable union of compact subsets, and as preparation, the following lemma gives
a way of doing this for the moduli space of complex structures.  Given
a smooth oriented surface $\Sigma$, we let $\jJ(\Sigma)$ denote the space of smooth
complex structures on $\Sigma$ compatible with the orientation, with its
natural $C^\infty$-topology.  For integers $g,m \ge 0$, $\mM_{g,m}$ will denote the (uncompactified)
moduli space of Riemann surfaces with genus $g$ and $m$ marked points;
recall that elements of the latter are equivalence classes of tuples 
$(\Sigma,j,\Theta)$ where $(\Sigma,j)$ is a Riemann surface of genus $g$ and
$\Theta \subset \Sigma$ is an ordered set of $m$ points.

\begin{lemma}
\label{lemma:TaubesDM}
Given integers $g,m \ge 0$, fix a closed surface $\Sigma$ of genus~$g$
and an ordered set of $m$ points $\Theta = \{\zeta_1,\ldots,\zeta_m\} \in \Sigma_g$.  Then there exists
a nested sequence of compact subsets 
$$
\jJ^1(\Sigma,\Theta) \subset \jJ^2(\Sigma,\Theta) \subset \jJ^3(\Sigma,\Theta) \subset \ldots \subset \jJ(\Sigma)
$$
such that every element of $\mM_{g,m}$ has a representative $(\Sigma,j,\Theta)$
for some $j \in \jJ^K(\Sigma,\Theta)$, $K \in \NN$.
\end{lemma}
\begin{proof}
Let $\pi : \jJ(\Sigma) \to \mM_{g,m} : j \mapsto [(\Sigma,j,\Theta)]$ denote
the natural projection.  Choose for each $j \in \jJ(\Sigma)$ a smooth 
slice $\tT_j \subset \jJ(\Sigma)$ through $j$ for the natural action of 
$\Diff_0(\Sigma,\Theta)$ on $\jJ(\Sigma)$, i.e.~$\tT_j$ locally parametrizes
the Teichm\"uller space of $(\Sigma,\Theta)$ near~$j$.  Since
Teichm\"uller space is finite dimensional, $\tT_j$ contains a compact
neighborhood $\vV_j \subset \tT_j$ of~$j$, and the image of $\vV_j$ under
$\pi$ is then a neighborhood of $[(\Sigma,j,\Theta)]$ in~$\mM_{g,m}$.
Since the latter is second countable, we can then find a sequence
$j_1,j_2,j_3,\ldots \in \jJ(\Sigma)$ such that $\bigcup_{i \in \NN} \pi(\vV_{j_i}) =
\mM_{g,m}$.  Set $\jJ^K(\Sigma,\Theta) := \vV_{j_1} \cup \ldots \cup \vV_{j_K}$.
\end{proof}

\begin{proof}[Proof of Theorem~\ref{thm:PetriJ}]
For the following definition, we fix a model surface $\Sigma_g$ of genus~$g$ and a point
$\zeta \in \Sigma_g$, along with Riemannian metrics on $\Sigma_g$
and $M$, denoting the various induced distance functions by 
$\dist(\cdot,\cdot)$.  The Levi-Civit\`a connection then
induces connections on the bundles $E = N_u$ and $F = \overline{\Hom}_\CC(T\Sigma_g,N_u)$
appearing below, which can be used in defining metrics on the jet spaces
$J^k_\zeta E$ and $J^k_\zeta F$.
For each $K,\ell \in \NN$, fix an
integer $k := k(K,\ell) \ge \ell$ large enough so that
\begin{equation}
\label{eqn:negativeDim}
\virdim \mM^{k,r,\ell}_{g,1}(A,J) < 0 \quad\text{ for all }\quad
r \in \{1,\ldots,K\}.
\end{equation}
With this choice in place, we define
$$
\nN^{K}(J) \subset \bigcup_{r=1}^{K} \mM^{k,r,\ell}_{g,1}(A,J)
$$
as a set of elements $(u,\zeta,t)$ satisfying quantitative versions of the
various conditions defining the spaces $\mM^{k,r,\ell}_{g,1}(A,J)$.
Concretely, we require every
element of $\nN^{K}(J)$ to be representable as a curve
$u : (\Sigma_g,j) \to (M,J)$ with marked point $\zeta \in \Sigma_g$ and
$t \in \vV^{k}_{r,\ell}(g_u,h_u,\mu_u,\mathbf{D}_u^N)$ with $|t|=1$ such that:
\begin{enumerate}
\item
\textsl{Domains do not degenerate:}
$j$ belongs to the compact set $\jJ^K(\Sigma_g,\{\zeta\})$ from
Lemma~\ref{lemma:TaubesDM}.
\item
\textsl{Bubbles do not form:}
$\sup_{z \in \Sigma_g} |du(z)| \le K$.
\item
\textsl{The marked point does not escape:}
$\dist(u(\zeta),M \setminus \uU) \ge 1/K$.
\item
\textsl{The marked point remains an injective point:}
$$
|du(\zeta)| \ge \frac{1}{K} \quad\text{ and }\quad
\inf_{z \in \Sigma_g \setminus \{\zeta\}} \frac{\dist(u(\zeta),u(z))}{\dist(\zeta,z)} \ge \frac{1}{K}.
$$
\item
\textsl{The rank of $t$ does not blow up:}
$\rank t \le K$.
\item
\textsl{The vanishing order of $t$ does not increase:}
Writing $E=N_u$ and $F=\overline{\Hom}_\CC(T\Sigma_g,E)$,
the distance of $t \in J^k_\zeta E \otimes J^k_\zeta F$ from
the subspace $(J^k_\zeta E \otimes J^k_\zeta F)^\ell$ is at least $1/K$.
\end{enumerate}
Now let
$$
\jJ^{\reg,K} := \left\{ J \in \jJ(M,\omega\,;\,\uU,\Jfix)\ \Big|\ 
\nN^{K}(J) = \emptyset \right\}.
$$
To see that $\jJ^{\reg,K}$ is open, suppose the contrary: then there exist
sequences $J_\nu \in \jJ(M,\omega\,;\,\uU,\Jfix)$ and
$(u_\nu,\zeta,t_\nu) \in \nN^{K}(J_\nu)$ with $J_\nu \to J \in \jJ^{\reg,K}$
as $\nu \to \infty$.
Assuming the parametrizations $u_\nu : (\Sigma_g,j_\nu) \to (M,J_\nu)$
satisfy all of the conditions listed above, elliptic regularity 
combined with the compactness of $\jJ^K(\Sigma_g,\{\zeta\})$ 
and the condition $|t_\nu| = 1$ then gives a
subsequence converging to an element of $\nN^{K}(J)$, 
which is a contradiction.

We claim that $\jJ^{\reg,K}$ is also dense.  To see this, recall that
the reference structure $\Jref$ in the definition of $\jJ_\varepsilon$ was
arbitrary, so it will suffice to prove that for some $\varepsilon \in \seqs$,
$\jJ_\varepsilon$ contains arbitrarily $C_\varepsilon$-small perturbations
of $\Jref$ that are in~$\jJ^{\reg,K}$.  The argument of the previous paragraph
shows that $\nN^K(\Jref)$ is compact, so since
$\varepsilon$-regularity is an open condition, 
Lemma~\ref{lemma:epsReg} implies after taking a lower bound 
for finitely many choices of $\varepsilon \in \seqs$ that every element of
$\nN^K(\Jref)$ is $\varepsilon$-regular, and so therefore is everything
in some open neighborhood of $\nN^K(\Jref) \times \{\Jref\}$
in $\univ^{k,r,\ell}_{g,1}(A,\jJ_\varepsilon)$.
Since $\jJ^\reg_\varepsilon \subset \jJ_\varepsilon$ is a Baire subset,
we can choose a sequence $J_\nu \in \jJ^\reg_\varepsilon$ with $J_\nu \to \Jref$,
and we claim that $J_\nu \in \jJ^{\reg,K}$ for all $\nu$ sufficiently large.
If not, then after restricting to a subsequence, there exists a sequence
$(u_\nu,\zeta_\nu,t_\nu) \in \nN^K(J_\nu)$ which converges 
by the compactness argument in the previous paragraph to
an element of $\nN^K(\Jref)$, implying that $(u_\nu,\zeta_\nu,t_\nu)$
is $\varepsilon$-regular for $\nu$ large.  In light of the assumption
$\virdim \mM^{k,r,\ell}_{g,1}(A,J) < 0$, this contradicts
Corollary~\ref{cor:itsEmpty!}.

The space
$$
\jJ^\reg := \bigcap_{K \in \NN} \jJ^{\reg,K} \subset \jJ(M,\omega\,;\,\uU,\Jfix)
$$
is now a Baire subset.  If $J \in \jJ^\reg$ and there exists a
simple $J$-holomorphic curve $u : (\Sigma,j) \to (M,J)$ of genus $g$
with an injective point $\zeta \in u^{-1}(\uU) \subset \Sigma$ at which
Petri's condition is not satisfied to infinite order, then we can 
define $u$ as an element of $\mM^*_{g,1}(A,J)$ by calling $\zeta$ the
marked point.  Since nontrivial elements 
$t \in \ker\mathbf{D}^N_u \otimes \ker(\mathbf{D}^N_u)^*$ have finite
rank and cannot vanish to infinite order at any point, we can then 
normalize $t$ and thus find
an element $(u,\zeta,t) \in \nN^{K}(J)$ for $K$ sufficiently large,
which is a contradiction.  This proves that for $J \in \jJ^\reg$, all
simple curves $v : (\Sigma,j) \to (M,J)$ satisfy Petri's condition to infinite order at every injective
point in $v^{-1}(\uU)$, which is an open and dense subset of~$v^{-1}(\uU)$.
It follows that the condition is also satisfied for all multiple covers
$u = v \circ \varphi$ at points in $u^{-1}(\uU) = \varphi^{-1}(v^{-1}(\uU))$
that are not branch points and are preimages of injective points; that is
likewise an open and dense subset of $u^{-1}(\uU)$.
\end{proof}

\begin{remark}
\label{remark:families}
The proof above would work equally well to find generic \emph{families}
of almost complex structures depending on finitely many parameters such that
Petri's condition is always satisfied.
The key point is that for the parametric moduli spaces analogous to
$\mM^{k,r,\ell}_{g,1}(A,J)$ and $\widehat{\mM}^{k,r,\ell}_{g,1}(A,J)$,
the codimension of the former in the latter
grows quadratically with~$k$, while the dimension of the larger space
grows only linearly, so that the space analogous to
$\mM^{k,r,\ell}_{g,1}(A,J)$ will always turn out to be empty for
generic choices if $k$ is made sufficiently large, no matter how many extra dimensions are
added to the original moduli space by introducing parameters.
The extension to families is important for the bifurcation theory
discussed in \S\ref{sec:bifurcations}.
\end{remark}

\subsection{A global application}
\label{sec:PetriLin}

We now give an application of Petri's condition which will be
crucial for the proof of Theorem~\ref{thm:submanifolds0}.  The setting is as
follows: assume $E$ and $F$ are smooth real vector bundles over a smooth
(not necessarily compact)
manifold~$M$, with chosen bundle metrics $\langle\ ,\ \rangle_E$,
$\langle\ ,\ \rangle_F$ and a chosen volume from $\mu$ on~$M$ which are
used to define $L^2$-pairings
$$
\langle \eta,\eta' \rangle_{L^2} := \int_M \langle \eta,\eta'\rangle_E \,\mu,
\qquad
\langle \xi,\xi' \rangle_{L^2} := \int_M \langle \xi,\xi'\rangle_F \,\mu
$$
for $\eta,\eta' \in \Gamma(E)$ and $\xi,\xi' \in \Gamma(F)$.
The product $\langle \eta,\eta' \rangle_{L^2}$ is well defined for two
(not necessarily smooth or compactly supported) sections 
$\eta,\eta'$ of $E$ whenever the function $\langle \eta,\eta' \rangle_E$ 
belongs to $L^1(M,\mu)$, and in this case we will say they are
\defin{$L^2$-orthogonal} if $\langle \eta,\eta' \rangle_{L^2}=0$;
an analogous definition applies for sections of~$F$.
Consider a linear partial differential operator
$\mathbf{D} : \Gamma(E) \to \Gamma(F)$ and its formal adjoint
$\mathbf{D}^*: \Gamma(F) \to\Gamma(E)$ defined via $\langle\xi,\mathbf{D}\eta\rangle_{L^2}
= \langle \mathbf{D}^*\xi,\eta \rangle_{L^2}$ for all smooth sections
$\eta,\xi$ with compact support.  We will consider the extensions of both of
these operators to certain Banach space completions,
$$
\mathbf{D} : \mathbf{X}(E) \to \mathbf{Y}(F),\qquad
\mathbf{D}^* : \mathbf{X}^*(F) \to \mathbf{Y}^*(E),
$$
where $\mathbf{X}(E)$ and $\mathbf{Y}^*(E)$ are Banach
spaces of sections of~$E$ in some regularity class defined almost everywhere,
while $\mathbf{Y}(F)$ and $\mathbf{X}^*(F)$ are likewise Banach spaces of sections
of~$F$.  In this functional-analytic setting, we impose the following assumptions:
\begin{enumerate}
\item
$\mathbf{D}$ and $\mathbf{D}^*$ are Fredholm operators whose kernels consist 
only of smooth sections;
\item
$\ker\mathbf{D}^* \subset \mathbf{Y}(F)$, and the $L^2$-product 
$\langle \xi,\xi' \rangle_{L^2}$ is well defined whenever
$\xi \in \mathbf{Y}(F)$ and $\xi' \in \ker\mathbf{D}^*$, so in particular
it is well defined whenever both are in~$\ker\mathbf{D}^*$;
\item
$\mathbf{Y}(F) = \im \mathbf{D} \oplus \ker \mathbf{D}^*$,
where the two factors in this splitting are closed $L^2$-orthogonal subspaces.
\end{enumerate}
We shall denote the natural projection resulting from the third assumption by
$$
\pi : \mathbf{Y}(F) \to \ker\mathbf{D}^*.
$$

\begin{remark}
\label{remark:satisfy}
In the setting of \S\ref{sec:CRpunctured}, the assumptions above are 
satisfied for a Cauchy-Riemann type operator $\dot{\mathbf{D}} :
\Gamma(\dot{E}) \to \Gamma(\dot{F})$ over a
punctured Riemann surface $\dot{\Sigma}$, using the weighted
Sobolev spaces
$\mathbf{X}(\dot{E}) := W^{k,p,-\boldsymbol{\delta}}(\dot{E})$ and
$\mathbf{Y}(\dot{F}) := W^{k-1,p,-\boldsymbol{\delta}}(\dot{F})$
for $k \in \NN$, $p \in (1,\infty)$ and exponential weights
$\boldsymbol{\delta} = \{ \delta_w > 0 \}_{w \in \Theta}$; recall that
$\dot{\mathbf{D}}$ is 
Fredholm if all $\delta_w$ are chosen to be sufficiently small.
For the formal adjoint $\dot{\mathbf{D}}^*$, we then define
$\mathbf{X}^*(\dot{F}) := W^{k,p,\boldsymbol{\delta}}(\dot{F})$ and
$\mathbf{Y}^*(\dot{E}) := W^{k-1,p,\boldsymbol{\delta}}(\dot{E})$, so that
Proposition~\ref{prop:formalAdjoint} provides the necessary splitting
of $\mathbf{Y}(\dot{F})$.
\end{remark}

\begin{lemma}
\label{lemma:surjective}
Given the assumptions above,
suppose $\uU \subset M$ is an open subset such that $\mathbf{D}$
satisfies Petri's condition over~$\uU$.
Assume moreover that
$V \subset \Gamma(\Hom(E,F))$ is a linear subspace satisfying
the following conditions:
\begin{enumerate}
\item $\Phi\eta \in \mathbf{Y}(F)$ for all $\Phi \in V$ and $\eta \in \ker\mathbf{D}$.
\item There exists a dense subset $\Delta \subset \uU$ with the following property:
for every 
$z \in \Delta$ and $\Phi_0 \in \Hom(E_z,F_z)$,
there exists a $\Phi \in \Gamma(\Hom(E,F))$ satisfying
$\Phi(z) = \Phi_0$ such that for every neighborhood $\uU' \subset \uU$ of~$z$,
$\beta\Phi \in V$ for some smooth function $\beta : M \to [0,1]$ with
compact support in $\uU'$ satisfying $\beta(z)=1$.
\end{enumerate}
Then the linear map
$\mathbf{L} : V \to \Hom(\ker\mathbf{D},\ker\mathbf{D}^*)$
defined by $\mathbf{L}(\Phi) \eta = \pi(\Phi\eta)$ is surjective.
\end{lemma}
\begin{proof}
Fix bases $\eta_1,\ldots,\eta_m \in \ker\mathbf{D}$ and
$\xi_1,\ldots,\xi_n \in \ker\mathbf{D}^*$.  Since
$\im\mathbf{D} = \ker \pi$ is $L^2$-orthogonal to $\ker\mathbf{D}^*$, we then have
$$
\langle \mathbf{L}(\Phi) \eta_i , \xi_j \rangle_{L^2} =
\langle \Phi \eta_i , \xi_j \rangle_{L^2} \quad \text{ for all }\quad
i=1,\ldots,m,\ j=1,\ldots,n,
$$
and these matrix elements determine 
$\mathbf{L}(\Phi) : \ker\mathbf{D} \to \ker\mathbf{D}^*$.
Now if $\mathbf{L}$ is not surjective, there exists a nontrivial linear map
$\Psi : \ker\mathbf{D} \to \ker\mathbf{D}^*$ which is ``orthogonal''
to every $\mathbf{L}(\Phi)$ in the sense that its matrix elements
$\Psi^{ij} := \langle \Psi \eta_i , \xi_j \rangle_{L^2} \in \RR$ satisfy
$$
\sum_{i,j} \Psi^{ij} \langle \Phi \eta_i , \xi_j \rangle_{L^2} = 0
$$
for every $\Phi \in V$.  We can rewrite this as
$$
0 = \sum_{i,j} \Psi^{ij} \int_{\uU} \langle \Phi\eta_i , \xi_j  \rangle_F \, \mu 
 = \int_{\uU} \langle\ ,\ \rangle_F \circ (\Phi \otimes \1) \circ
\bigg( \sum_{i,j} \Psi^{ij}\, \eta_i \otimes \xi_j \bigg) \, \mu ,
$$
where $\sum_{i,j} \Psi^{ij} \eta_i \otimes \xi_j$ is regarded as a section
of $E \otimes F$.
Since the $\Psi^{ij}$ are not all zero, this section is
the image of a nontrivial element of $\ker\mathbf{D} \otimes \ker\mathbf{D}^*$
under the Petri map, so by assumption, it does not vanish identically on~$\uU$.
Now choose a point $z \in \Delta$ at which this section is nonzero.
Lemma~\ref{lemma:tensor} below provides a 
linear map $\Phi_0 : E_z \to F_z$ such that the integrand is positive
near~$z$ for any $\Phi \in V$ satisfying $\Phi(z) = \Phi_0$, and we can then
make the entire integral positive after multiplying $\Phi$ by smooth bump 
functions with sufficiently small support.
\end{proof}

We used:

\begin{lemma}
\label{lemma:tensor}
Suppose $V$ and $W$ are real finite-dimensional vector spaces,
$\langle\ ,\ \rangle : W \otimes W \to \RR$ is an inner product on~$W$, 
and $T \in V \otimes W$ is nonzero.  Then there exists a linear map
$\Phi : V \to W$ such that $\langle\ ,\ \rangle \circ (\Phi \otimes \1)(T) > 0$.
\end{lemma}
\begin{proof}
Choosing a basis $v_1,\ldots,v_n$ of~$V$, we have $T = \sum_{j=1}^n v_j \otimes w_j$
for unique vectors $w_1,\ldots,w_n \in W$, which do not all vanish since $T \ne 0$.
Choosing $\Phi : V \to W$ such that $\Phi(v_j) = w_j$ for all $j$ then gives
$\langle\ ,\ \rangle \circ (\Phi \otimes \1)(T) = \sum_j \langle w_j,w_j \rangle > 0$.
\end{proof}

\section{Proof of the stratification theorem}
\label{sec:SardSmale}

We are now in a position to prove Theorem~\ref{thm:submanifolds0}.
The main idea behind the proof is standard, though some details are less so: 
we will write down a universal moduli space with a projection to a suitable
Banach manifold of perturbed data whose regular values
have the property stated in the theorem.  The hard part is of course to
prove that the universal moduli space is a smooth Banach manifold---this
follows from the implicit function theorem after proving that some version
of the operator defined in \eqref{eqn:linearization} is surjective, and
that is where the results of the previous section on Petri's condition are
needed.

Fix $\Jref \in \jJ(M,\omega\,;\,\uU,\Jfix)$ and consider again the space
$\jJ_\varepsilon$ of Floer $C_\varepsilon$-small perturbations of
$\Jref$ as constructed in \S\ref{sec:PetriJ} via a choice of decaying
positive sequence $\varepsilon = \{\varepsilon_\nu\}_{\nu=0}^\infty \in \seqs$.
For each of the choices of data in the statement of Theorem~\ref{thm:submanifolds0},
we define a universal moduli space
$$
\univ^d(\jJ_\varepsilon\,;\,\ell_1,\ldots,\ell_m)
$$
consisting of pairs $(u,J)$ with $J \in \jJ_\varepsilon$ and $u$ belonging to
the isosymmetric stratum
$$
\mM^d(J\,;\,\ell_1,\ldots,\ell_m) := \mM^d_{\mathbf{b}}(\mM_{g,m}(A,J\,;\,\ell_1,\ldots,\ell_m)).
$$
We shall denote elements of $\mM^d(J\,;\,\ell_1,\ldots,\ell_m)$
by $u=v\circ \varphi$, where we have chosen parametrizations of the
underlying simple curve $v : (\Sigma,j) \to (M,J)$ and the $d$-fold branched
cover $\varphi : (\Sigma',j') \to (\Sigma,j)$.
Recall from \S\ref{sec:prep} that for every such element
$u = v \circ \varphi$, there is a unique isomorphism class of
minimal regular presentations for~$\varphi$, giving rise to a regular cover
$$
\widehat{\varphi} : (\widehat{\Sigma},\widehat{\jmath}) \to (\dot{\Sigma},j)
$$
with automorphism group $G := \Aut(\widehat{\varphi})$,
where $\dot{\Sigma}$ is the punctured surface obtained from $\Sigma$ by
removing the critical values of~$\varphi$.  We can then consider the
$J$-holomorphic curve $\widehat{u} := v \circ \widehat{\varphi} :
(\widehat{\Sigma},\widehat{\jmath}) \to (M,J)$ and its normal Cauchy-Riemann
operator $\dot{\mathbf{D}}^N_{\widehat{u}}$, defined as in \S\ref{sec:CRpunctured}
on a Sobolev space of sections of $E := N_{\widehat{u}}$ 
over the punctured domain $\widehat{\Sigma}$
with negative exponential weights close to zero.  Recall that its formal adjoint
$(\dot{\mathbf{D}}^N_{\widehat{u}})^*$ is defined on a similar Sobolev
space of sections of $F := \overline{\Hom}_\CC(T\widehat{\Sigma},N_{\widehat{u}})$, 
but with corresponding positive exponential weights.  The notation associating
to each $(u = v \circ \varphi,J) \in \univ^d(\jJ_\varepsilon\,;\,\ell_1,\ldots,\ell_m)$
a regular covering map $\widehat{\varphi}$ of potentially larger degree
and corresponding $J$-holomorphic curve $\widehat{u} = v \circ \widehat{\varphi}$
will be used consistently in the following.

\begin{defn}
\label{defn:universalWalls}
Given integers $k,c \ge 0$ and an almost complex structure $J$, we define the subset
$$
\mM^d(J\,;\,\ell_1,\ldots,\ell_m\,;\,k,c) := 
\Big\{ u \in \mM^d(J\,;\,\ell_1,\ldots,\ell_m) \ \Big|\ 
\dim \ker \dot{\mathbf{D}}_{\widehat{u}}^N = k 
\text{ and }\dim \coker \dot{\mathbf{D}}_{\widehat{u}}^N = c \Big\}.
$$
This gives rise to a universal moduli space
$$
\univ^d(\jJ_\varepsilon\,;\,\ell_1,\ldots,\ell_m\,;\,k,c) \subset
\univ^d(\jJ_\varepsilon\,;\,\ell_1,\ldots,\ell_m)
$$
consisting of all pairs $(u,J)$ such that $J \in \jJ_\varepsilon$ and
$u \in \mM^d(J\,;\,\ell_1,\ldots,\ell_m\,;\,k,c)$.
\end{defn}

By the results of \S\ref{sec:IFT}, in particular Lemma~\ref{lemma:localV},
the connected components of the subsets $\mM^d(J\,;\,\ell_1,\ldots,\ell_m\,;\,k,c)$
for individual values of $k$ and $c$ are precisely the walls 
described in Theorem~\ref{thm:submanifolds0} (see also Remark~\ref{remark:whatisawall}).
We would thus be able to apply the standard Sard-Smale argument toward
a proof of Theorem~\ref{thm:submanifolds0} if we could show that
$\univ^d(\jJ_\varepsilon\,;\,\ell_1,\ldots,\ell_m\,;\,k,c) \subset
\univ^d(\jJ_\varepsilon\,;\,\ell_1,\ldots,\ell_m)$ is a smooth Banach submanifold
of the correct finite codimension on each component.  What we will actually
show is that this is true for a certain open subset
of $\univ^d(\jJ_\varepsilon\,;\,\ell_1,\ldots,\ell_m\,;\,k,c)$,
which suffices due to the genericity of Petri's condition.

\begin{defn}
\label{defn:PetriReg}
An element $u = v \circ \varphi \in \mM^d(J\,;\,\ell_1,\ldots,\ell_m\,;\,k,c)$
will be called \defin{Petri regular} if for the regular covering map
$\widehat{\varphi}$ and corresponding $J$-holomorphic curve $\widehat{u} = v \circ \widehat{\varphi}$
described above, the operator $\dot{\mathbf{D}}_{\widehat{u}}^N$ satisfies Petri's condition
over~$\widehat{u}^{-1}(\uU)$.  We will denote the set of Petri regular curves by
$$
\mM^d_\Pi(J\,;\,\ell_1,\ldots,\ell_m\,;\,k,c) \subset
\mM^d(J\,;\,\ell_1,\ldots,\ell_m\,;\,k,c),
$$
and define the corresponding universal moduli space
$$
\univ^d_\Pi(\jJ_\varepsilon\,;\,\ell_1,\ldots,\ell_m\,;\,k,c) \subset
\univ^d(\jJ_\varepsilon\,;\,\ell_1,\ldots,\ell_m\,;\,k,c)
$$
to be the set of pairs $(u,J) \in \univ^d(\jJ_\varepsilon\,;\,\ell_1,\ldots,\ell_m\,;\,k,c)$
such that $u$ belongs to the moduli
space $\mM^d_\Pi(J\,;\,\ell_1,\ldots,\ell_m\,;\,k,c)$.
\end{defn}

\begin{remark}
The condition defining $\mM^d_\Pi(J\,;\,\ell_1,\ldots,\ell_m\,;\,k,c)$
is clearly satisfied by any curve $u = v \circ \varphi$
for which
$\mathbf{D}_v^N$ satisfies the local Petri condition on $v^{-1}(\uU)$,
thus by Theorem~\ref{thm:PetriJ}, there is a Baire subset in
$\jJ(M,\omega\,;\,\uU,\Jfix)$ for which
$\mM^d_\Pi(J\,;\,\ell_1,\ldots,\ell_m\,;\,k,c) = 
\mM^d(J\,;\,\ell_1,\ldots,\ell_m\,;\,k,c)$.
\end{remark}

The next several results are aimed at proving that for suitable choices
of the sequence~$\varepsilon$, $\univ^d_\Pi(\jJ_\varepsilon\,;\,\ell_1,\ldots,\ell_m\,;\,k,c)$
is a finite-codimensional Banach submanifold of
$\univ^d(\jJ_\varepsilon\,;\,\ell_1,\ldots,\ell_m)$.

\begin{lemma}
\label{lemma:contin}
For $\varepsilon \in \seqs$ with sufficiently rapid decay,
$\univ^d(\jJ_\varepsilon\,;\,\ell_1,\ldots,\ell_m)$ carries a smooth
Banach manifold structure such that every
$(u_0=v_0 \circ \varphi_0,J_0) \in 
\univ^d(\jJ_\varepsilon\,;\,\ell_1,\ldots,\ell_m)$ admits a neighborhood
$\vV \subset \univ^d(\jJ_\varepsilon\,;\,\ell_1,\ldots,\ell_m)$ 
with a smooth family of vector bundle isomorphisms
$$
v_0^*TM \stackrel{\cong}{\longrightarrow} v^*TM,
\quad\text{ for }\quad (u=v \circ \varphi,J) \in \vV
$$
mapping $N_{v_0}$ isomorphically to~$N_{v}$.
\end{lemma}
\begin{proof}
For each $(u_0=v_0 \circ \varphi_0,J_0) \in \univ^d(\jJ_\varepsilon\,;\,\ell_1,\ldots,\ell_m)$,
the underlying simple curve $v_0 : (\Sigma,j_0) \to (M,J)$ 
lives in the universal moduli space $\univ^*(\jJ_\varepsilon)$
defined in Appendix~\ref{sec:critical}, more specifically in the
subset
$$
\univ^*(\jJ_\varepsilon\,;\,\ell_1,\ldots,\ell_m) \subset \univ^*(\jJ_\varepsilon)
$$
of this space defined by the condition that the $i$th marked point should
have critical order $\ell_i$ and curves are immersed everywhere else.
If $\varepsilon$ has sufficiently rapid decay, then $\univ^*(\jJ_\varepsilon)$
is a smooth Banach manifold, and $\univ^*(\jJ_\varepsilon\,;\,\ell_1,\ldots,\ell_m)$
is an open subset of the space
$\widehat{\univ}^*(\jJ_\varepsilon\,;\,\ell_1,\ldots,\ell_m) \subset \univ^*(\jJ_\varepsilon)$,
which is shown in Lemma~\ref{lemma:submersion} to be a smooth finite-codimensional
submanifold of~$\univ^*(\jJ_\varepsilon)$.  In particular, we can identify
an open neighborhood of the element
$(v_0,J_0)$ in $\univ^*(\jJ_\varepsilon\,;\,\ell_1,\ldots,\ell_m)$
with a smooth finite-codimensional submanifold 
$$
X_\varepsilon \subset \dbar^{-1}(0) \subset \tT \times \bB \times \jJ_\varepsilon
$$
of the zero-set of the nonlinear Cauchy-Riemann operator~$\dbar$, where
$\tT$ denotes a Teichm\"uller slice through $j_0$ in the space of complex
structures on~$\Sigma$, and $\bB$ is a suitable Banach manifold
of maps $v : \Sigma \to M$.

We claim that there exists a neighborhood $\vV_0 \subset X_\varepsilon$
of $(j_0,v_0,J_0)$ that parametrizes a smooth family of bundle isomorphisms
$v_0^*TM \to v^*TM$ sending $N_{v_0}$ to~$N_v$.  Note that this would be
clearly false if we did not impose the critical point constraints on~$v$,
as e.g.~$v_0$ might then have critical points while $v$ is immersed, in
which case $N_{v_0}$ and $N_v$ would have different topological types.
Assuming $N_v \subset v^*TM$
is always defined as the symplectic orthogonal complement of $T_v \subset v^*TM$
with $T_v := \im dv$ away from critical points,
let us recall from \cite{Wendl:automatic} how the latter is defined at
critical points.  We have a smooth family of bundles
$v^*TM$ carrying linearized Cauchy-Riemann operators $\mathbf{D}_v$, whose
complex-linear parts $\mathbf{D}_v^\CC$ define a smooth family of holomorphic
structures on $v^*TM$.  The crucial observation is then 
that $dv \in \Gamma(\Hom_\CC(T\Sigma,v^*TM))$ is always a holomorphic
section with respect to the holomorphic bundle structures on $v^*TM$ and~$T\Sigma$,
so choosing a smooth family of holomorphic trivializations and holomorphic
coordinates near the $i$th marked point, each $dv$ is represented by some
holomorphic function of the form
$$
f^{(i)}_v : \DD \to \CC^m, \qquad  f^{(i)}(z) = z^{\ell_i} g^{(i)}_v(z),
$$
where $g^{(i)}_v : \DD \to \CC^m$ is another family of holomorphic functions
which depend smoothly on $(j,v,J) \in X_\varepsilon$ but also are nonzero at~$0$.
The main point here is that the critical orders $\ell_i$ do not vary with~$v$.
The span of $g^{(i)}_v(0)$ thus defines the fibers of $T_v$ near each critical point,
so we deduce smooth dependence of $T_v$ on $(j,v,J) \in X_\varepsilon$, and therefore
also of~$N_v$.

We can parametrize
a neighborhood of $\varphi_0$ in $\mM^d_{\mathbf{b}}(j_0)$ as explained in
Examples~\ref{ex:Teichmueller} and~\ref{ex:Teichmueller3}, meaning that if
$\Theta = \{w_1,\ldots,w_r\} \subset \Sigma$ is the set of critical values of
$\varphi_0$, we choose a smooth family of diffeomorphisms
$\psi_\tau : \Sigma \to \Sigma$ parametrized by $\tau \in B^{2r}$ which are
holomorphic near $\Theta$ and supported on a slightly larger neighborhood
of $\Theta$ such that $\psi_0 = \Id$ and
$$
B^{2r} \to \Sigma^{\times r} : \tau \mapsto (\psi_\tau(w_1),\ldots,\psi_\tau(w_r))
$$
is a diffeomorphism onto an open set.  The neighborhood of $(u_0,J_0)$ in
the space 
$\univ^d(\jJ_\varepsilon\,;\,\ell_1,\ldots,\ell_m)$ can now be identified
with $B^{2r} \times X_\varepsilon$ by associating to each $(\tau,(j,v,J)) \in B^{2r} \times X_\varepsilon$
the curve $v \circ (\psi_\tau \circ \varphi_0)$, making
$\univ^d(\jJ_\varepsilon\,;\,\ell_1,\ldots,\ell_m)$ a smooth fiber bundle
over $\univ^*(\jJ_\varepsilon\,;\,\ell_1,\ldots,\ell_m)$.
\end{proof}

\begin{lemma}
The subset $\univ^d_\Pi(\jJ_\varepsilon\,;\,\ell_1,\ldots,\ell_m\,;\,k,c)
\subset \univ^d(\jJ_\varepsilon\,;\,\ell_1,\ldots,\ell_m\,;\,k,c)$ is open.
\end{lemma}
\begin{proof}
Lemma~\ref{lemma:contin} implies that the operators
$\dot{\mathbf{D}}^N_{\widehat{u}}$ and $(\dot{\mathbf{D}}^N_{\widehat{u}})^*$ can
both be understood as varying continuously with $(u,J) \in
\univ^d(\jJ_\varepsilon\,;\,\ell_1,\ldots,\ell_m)$, and 
the dimensions of their kernels are locally constant
as long as $(u,J)$ moves only in the subset
$\univ^d(\jJ_\varepsilon\,;\,\ell_1,\ldots,\ell_m\,;\,k,c)$.
It follows that the family of Petri maps defined on
$\ker \mathbf{D}_{\widehat{u}}^N \otimes \ker (\mathbf{D}_{\widehat{u}}^N)^*$ and
then restricted to $\widehat{u}^{-1}(\uU)$ depends continuously on
$(u,J) \in \univ^d(\jJ_\varepsilon\,;\,\ell_1,\ldots,\ell_m\,;\,k,c)$,
and since their domains are finite dimensional, the injectivity of these
maps is an open condition.
\end{proof}

Following Example~\ref{ex:Teichmueller3},
the smooth family of operators $\mathbf{D}_v^N$ parametrized by
$\univ^d(\jJ_\varepsilon\,;\,\ell_1,\ldots,\ell_m)$ can now be fit into
the general picture from \S\ref{sec:prep} of a parametrized
family of bundles with Cauchy-Riemann operators.
In particular, we choose the parameter space $P$ to be the local model of
$\univ^d(\jJ_\varepsilon\,;\,\ell_1,\ldots,\ell_m)$ near $(u_0,J_0)$ described in the
proof of Lemma~\ref{lemma:contin} above,
$$
P := B^{2r} \times X_\varepsilon \subset B^{2r} \times \dbar^{-1}(0) \subset
B^{2r} \times (\tT \times \bB \times \jJ_\varepsilon),
$$
and in the notation of \S\ref{sec:prep}, associate to each
$\tau = (\sigma,(j,v,J)) \in P$ the data
$$
\psi_\tau := \psi_\sigma,\qquad
j_\tau := j,\qquad
(E_\tau,J_\tau) := (N_v,J),\qquad
\mathbf{D}_\tau := \mathbf{D}_v^N.
$$
If $(u_0,J_0) \in
\univ^d(\jJ_\varepsilon\,;\,\ell_1,\ldots,\ell_m\,;\,k,c)$,
then using the setup in
\S\ref{sec:IFT}, we now find a smooth map
\begin{equation}
\label{eqn:bigF}
\mathbf{F}_\varepsilon : B^{2r} \times X_\varepsilon \to \Hom_G\big( \ker \dot{\mathbf{D}}_{\widehat{u}_0}^N,
\ker (\dot{\mathbf{D}}_{\widehat{u}_0}^N)^* \big)
\end{equation}
whose zero-set is a neighborhood of $(u_0,J_0)$ in
$\univ^d(\jJ_\varepsilon\,;\,\ell_1,\ldots,\ell_m\,;\,k,c)$.

\begin{defn}
\label{defn:epsReg2}
We will say that $(u_0,J_0) \in
\univ^d(\jJ_\varepsilon\,;\,\ell_1,\ldots,\ell_m\,;\,k,c)$
is \defin{$\varepsilon$-regular} if $\varepsilon\in \seqs$ has sufficiently rapid decay
to satisfy the conclusions of Lemma~\ref{lemma:contin} and, additionally,
the linearization of the map \eqref{eqn:bigF} at $(0,(j_0,v_0,J_0))$ is surjective.
Given $J \in \jJ(M,\omega\,;\,\uU,\Jfix)$ and $\varepsilon \in \seqs$, an element $u$
in the space $\mM^d(J\,;\,\ell_1,\ldots,\ell_m\,;\,k,c)$ will similarly
be called \defin{$\varepsilon$-regular} if $J \in \jJ_\varepsilon$ and
$(u,J)$ is $\varepsilon$-regular.
\end{defn}

In analogy with Definition~\ref{defn:epsReg}, $\varepsilon$-regularity for an element
$(u_0,J_0) \in \univ^d(\jJ_\varepsilon\,;\,\ell_1,\ldots,\ell_m\,;\,k,c)$
just means that a neighborhood of $(u_0,J_0)$ in this space
is a smooth Banach submanifold with the ``correct'' finite codimension 
in $\univ^d(\jJ_\varepsilon\,;\,\ell_1,\ldots,\ell_m)$.
It could be phrased alternatively as the condition that $(u_0,J_0)$ is a
transverse intersection of the map $(u,J) \mapsto \dot{\mathbf{D}}^N_{\widehat{u}}$
from $\univ^d(\jJ_\varepsilon\,;\,\ell_1,\ldots,\ell_m)$ to the
relevant space of $G$-equivariant Fredholm operators with 
the finite-codimensional submanifold $\{ \mathbf{T}\ |\ \dim \ker\mathbf{T}
= \dim \ker \dot{\mathbf{D}}^N_{\widehat{u}_0} \}$; 
expressed in this way, $\varepsilon$-regularity is clearly an open condition
and is independent of the choices involved (except of course for the choice
of $\varepsilon \in \seqs$).

Let us define the analogous condition for moduli spaces with fixed~$J$.
Note that if the simple curve $v_0$ is regular for the \emph{constrained}
moduli space $\mM_{g,m}(A,J_0\,;\,\ell_1,\ldots,\ell_m)$ as defined in
Appendix~\ref{sec:critical}, then the set
$$
X(J_0) := \left\{ (j,v,J_0) \in X_\varepsilon\ \big|\ j \in \tT,\ v \in \bB \right\} 
\subset \dbar_{J_0}^{-1}(0)
$$
is independent of $\varepsilon \in \seqs$ and is a smooth finite-dimensional submanifold pa\-ra\-met\-riz\-ing a neighborhood of
$v_0$ in $\mM_{g,m}(A,J_0\,;\,\ell_1,\ldots,\ell_m)$.  A neighborhood of
$u_0$ in $\mM^d(J\,;\,\ell_1,\ldots,\ell_m)$
is then parametrized by the submanifold $B^{2r} \times X(J_0) \subset B^{2r} \times X_\varepsilon$.
We will say that $u_0 \in \mM_{g,m}(A,J_0\,;\,\ell_1,\ldots,\ell_m\,;\,k,c)$ is \defin{regular in its stratum} if regularity of
$v_0$ in the sense above holds and, additionally, the restricted linearization
$$
\begin{tikzcd}[column sep=huge]
T_{(0,(j_0,v_0,\varphi_0))}(B^{2r} \times X(J_0)) \ar[r,"{d\mathbf{F}_\varepsilon(0,(j_0,v_0,\varphi_0))}"] &
\Hom_G\big( \ker \dot{\mathbf{D}}_{\widehat{u}_0}^N,
\ker (\dot{\mathbf{D}}_{\widehat{u}_0}^N)^* \big)
\end{tikzcd}
$$
is surjective.  This can also be rephrased as a transverse intersection
condition in the space of Fredholm operators, and is thus open and independent
of choices (including $\varepsilon$).
Our goal is to show that all curves satisfy this condition
for generic~$J$.

\begin{lemma}
\label{lemma:epsReg2}
If $u=v \circ \varphi \in \mM^d(\Jref,;\,\ell_1,\ldots,\ell_m\,;\,k,c)$
is Petri regular, then it is $\varepsilon$-regular for all $\varepsilon \in \seqs$ with sufficiently rapid decay.
\end{lemma}
\begin{proof}
Clearly $(u,\Jref) \in \univ^d_\Pi(\jJ_\varepsilon\,;\,\ell_1,\ldots,\ell_m\,;\,k,c)$
for every $\varepsilon \in \seqs$, and we shall assume $\varepsilon$ has sufficiently rapid
decay so that $\univ^d(\jJ_\varepsilon\,;\,\ell_1,\ldots,\ell_m)$ is a smooth
Banach manifold.
By Lemma~\ref{lemma:operatorPert}, there is a large space of smooth
perturbations $Y \in T_{\Jref}\jJ(M,\omega\,;\,\uU,\Jfix)$ that give rise 
via \eqref{eqn:expJ} to smooth $1$-parameter families 
$J_\tau := J_{\tau Y} \in \jJ(M,\omega\,;\,\uU,\Jfix)$ for which
$v$ remains $J_\tau$-holomorphic, and the normal Cauchy-Riemann operator
$\mathbf{D}_v^N$ is perturbed in the direction of an arbitrary smooth zeroth-order term
$A_Y$ with support in $v^{-1}(\uU)$ away from the discrete set of critical and double points of~$v$.
Such a perturbation defines a tangent vector $(0,Y) \in T_{(u,\Jref)} \univ^d(\jJ_\varepsilon\,;\,\ell_1,\ldots,\ell_m)$
whenever $\varepsilon$ has sufficiently rapid decay for $Y$ to be of class~$C_\varepsilon$.
Assuming this for the moment,
the resulting perturbation to $\dot{\mathbf{D}}_{\widehat{u}}^N$ is
$$
\dot{\mathbf{D}}_{\widehat{u}}^N \leadsto \dot{\mathbf{D}}_{\widehat{u}}^N + \widehat{\varphi}^*A_Y,
$$
hence differentiating $\mathbf{F}_\varepsilon$ in the direction $(0,Y)$
produces a $G$-equivariant 
linear map $\mathbf{L}(Y) : \ker \dot{\mathbf{D}}_{\widehat{u}}^N \to
\ker (\dot{\mathbf{D}}_{\widehat{u}}^N)^*$ given by \eqref{eqn:linearization}, namely
$$
\mathbf{L}(Y) \eta = \pi\big( (\widehat{\varphi}^*A_Y ) \eta \big),
$$
in terms of the projection
$$
\pi : W^{k-1,p,-\widehat{\varphi}^*\boldsymbol{\delta}}(N_{\widehat{u}}) = 
\im (\dot{\mathbf{D}}_{\widehat{u}}^N) \oplus 
\ker (\dot{\mathbf{D}}_{\widehat{u}}^N)^* \to 
\ker (\dot{\mathbf{D}}_{\widehat{u}}^N)^*.
$$
We claim that $Y$ can be chosen to make $\mathbf{L}(Y)$ equal to any given element
$$
\Psi \in \Hom_G\big( \ker \dot{\mathbf{D}}_{\widehat{u}}^N,
\ker (\dot{\mathbf{D}}_{\widehat{u}}^N)^* \big).
$$
Indeed, let us abbreviate $E = N_{v}$ and $F = \overline{\Hom}_\CC(T\Sigma,N_{v})$,
and let $\Delta \subset v^{-1}(\uU) \subset \Sigma$ denote the set of injective
points of $v$ that are not critical values of $\widehat{\varphi}$ and have image in~$\uU$;
these form an open and dense subset of $v^{-1}(\uU)$.
Since $\dot{\mathbf{D}}_{\widehat{u}}^N$ satisfies Petri's condition
over $\widehat{u}^{-1}(\uU)$, 
Lemma~\ref{lemma:surjective} then provides for any given $\Psi$
a section $\widehat{A} \in
\Gamma(\Hom_\RR(\widehat{\varphi}^*\dot{E},\widehat{\varphi}^*\dot{F}))$ 
with compact support in the open and dense subset
$\widehat{\varphi}^{-1}(\Delta) \subset \widehat{u}^{-1}(\uU)$ such that
$$
\langle \xi , \widehat{A}\eta \rangle_{L^2} = \langle \xi , \Psi \eta \rangle_{L^2}
$$
for all $\xi \in \ker (\dot{\mathbf{D}}_{\widehat{u}}^N)^*$
and $\eta \in \ker \dot{\mathbf{D}}_{\widehat{u}}^N$.
Note that we are free to assume the $L^2$-product is invariant under the action
of $G$ via deck transformations.  Then since $\Psi$ is $G$-equivariant,
we also have for every $g \in G$,
\begin{equation*}
\langle \xi,(g\widehat{A}) \eta \rangle_{L^2} = 
\langle g^{-1} \xi , \widehat{A} (g^{-1} \eta) \rangle_{L^2} =
\langle g^{-1} \xi, \Psi (g^{-1} \eta) \rangle_{L^2} =
\langle g^{-1} \xi, g^{-1} (\Psi \eta) \rangle_{L^2} 
= \langle \xi,\Psi \eta \rangle_{L^2},
\end{equation*}
implying that the symmetrization
$\widehat{A}_G := \frac{1}{|G|} \sum_{g \in G} g\widehat{A}$ also satisfies
$$
\langle \xi , \widehat{A}_G \eta \rangle_{L^2} = \langle \xi, \Psi\eta \rangle_{L^2}
$$
for all $\xi,\eta$.  But the $G$-invariance of $\widehat{A}_G$ implies
$\widehat{A}_G = \widehat{\varphi}^*A$ for some
$A \in \Gamma(\Hom_\RR(\dot{E},\dot{F}))$ with compact support in~$\Delta$,
hence $A = A_Y$ for some $Y \in T_{\Jref}\jJ(M,\omega\,;\,\uU,\Jfix)$, and this proves the claim.
We can now choose any finite collection of perturbations
$Y_1,\ldots,Y_N \in T_{\Jref}\jJ(M,\omega\,;\,\uU,\Jfix)$ 
such that the $\mathbf{L}(Y_i)$ span $\Hom_G\big( \ker \dot{\mathbf{D}}_{\widehat{u}}^N,
\ker (\dot{\mathbf{D}}_{\widehat{u}}^N)^* \big)$,
and choose $\varepsilon \in \seqs$ so that all of them are of class~$C_\varepsilon$.
\end{proof}

By the implicit function theorem, the open set of $\varepsilon$-regular elements in
$$
\univ^d(\jJ_\varepsilon\,;\,\ell_1,\ldots,\ell_m\,;\,k,c)
\subset \univ^d(\jJ_\varepsilon\,;\,\ell_1,\ldots,\ell_m)
$$
is a smooth Banach submanifold whose codimension near any given
element $(u,J)$ is given by the formula
in \eqref{eqn:Schur2}, and thus matches $\codim(u)$ as specified
by Definition~\ref{defn:codimCurve}.
We can then apply the Sard-Smale theorem to the projection
$$
\univ^d(\jJ_\varepsilon\,;\,\ell_1,\ldots,\ell_m\,;\,k,c) \to
\jJ_\varepsilon : (u,J) \mapsto J
$$
and thus find a Baire subset $\jJ_\varepsilon^\reg \subset \jJ_\varepsilon$
such that for all $J \in \jJ_\varepsilon^\reg$, all 
$\varepsilon$-regular elements of
$\mM^d(J\,;\,\ell_1,\ldots,\ell_m\,;\,k,c)$
are regular in their stratum. 

To turn this into a Baire subset of $\jJ(M,\omega\,;\,\uU,\Jfix)$ and
drop the $\varepsilon$-regularity condition,
we now apply another variation on the Taubes trick that was used in
the proof of Theorem~\ref{thm:PetriJ}, i.e.~we exhaust the moduli space 
of Petri regular curves by a countable collection of compact subsets 
$$
\nN^K(J) \subset 
\mM^d_\Pi(J\,;\,\ell_1,\ldots,\ell_m\,;\,k,c),
\qquad K \in \NN,
$$
in order to define open and dense subsets
of $\jJ(M,\omega\,;\,\uU,\Jfix)$ whose intersection has the desired properties.  
As in \S\ref{sec:theBigIdea},
let $h \ge 0$ denote the genus of $d$-fold branched covers of a genus~$g$
surface as determined by the branching data $\mathbf{b}$ and the Riemann-Hurwitz
formula.  We shall again write $\textbf{b} = (\mathbf{b}_1,\ldots,\mathbf{b}_r)$
for some $r \ge 0$, where each individual $\mathbf{b}_i$ is a
tuple $(b_i^1,\ldots,b_i^{q_i})$ of natural numbers satisfying
$\sum_{j=1}^{q_i} b_i^j = d$.  
Now fix a closed model surface
$\Sigma_g$ of genus~$g$ along with an ordered set of distinct points
$\Theta = (x_1,\ldots,x_m)$ in~$\Sigma_g$ and a continuous function
$F_g : \Sigma_g \to [0,\infty)$ that is positive on $\Sigma \setminus \Theta$
and, using local complex coordinates $z$ to identify a neighborhood of
each $x_j$ with $\DD \subset \CC$ so that $x_j$ becomes $0 \in \DD$, satisfies
$$
F_g(z) = |z|^{\ell_j} \text{ near } x_j,\qquad j=1,\ldots,m.
$$
Similarly, fix a closed model surface $\Sigma_h$ of genus~$h$, an
ordered set of distinct points
$$
\Theta' = (\zeta_1^1,\ldots,\zeta_1^{q_1},\ldots,\zeta_r^1,\ldots,\zeta_r^{q_r})
$$
in~$\Sigma_h$, and a continuous function $F_h : \Sigma_h \to [0,\infty)$
that is positive on $\Sigma_h \setminus \Theta'$ and takes the form
$$
F_h(z) = |z|^{b_i^j-1} \text{ near } \zeta_i^j, \qquad j=1,\ldots,q_i,\ i=1,\ldots,r
$$
in suitable local coordinates.
We also make arbitrary choices of Riemannian metrics on $\Sigma_g$,
$\Sigma_h$ and $M$ so as to define the various distance functions $\dist(\ ,\ )$
and norms referred to below.
We then define $\nN^K(J)$ to consist of every element in 
$\mM^d(J\,;\,\ell_1,\ldots,\ell_m\,;\,k,c)$
that admits
a representative of the form $u = v \circ \varphi : (\Sigma_h,j') \to
(M,J)$, with $v : (\Sigma_g,j) \to (M,J)$ simple and $\varphi : (\Sigma_h,j')
\to (\Sigma,j)$ a $d$-fold holomorphic branched cover, such that
$v$ is critical of order $\ell_i$ at $x_i$ for $i=1,\ldots,m$ and
$\varphi$ has branching order $b_i^j$ at $\zeta_i^j$ for $j=1,\ldots,q_i$
and $i=1,\ldots,r$, and the following quantitative conditions are also satisfied:
\begin{enumerate}
\item \label{item:domain}
\textsl{Domains do not degenerate:} Using the compact sets of complex
structures provided by Lemma~\ref{lemma:TaubesDM},
$j \in \jJ^K(\Sigma_g,\Theta)$ and $j' \in \jJ^K(\Sigma_h,\Theta')$.
\item \label{item:bubbles}
\textsl{Bubbles do not form:}
$\sup_{z \in \Sigma_g} |dv(z)| \le K$ and $\sup_{z \in \Sigma_h} |d\varphi(z)| \le K$.
\item \label{item:injective}
\textsl{Injective points do not disappear:}
There exists a point $\zeta \in \Sigma_g$ such that
$$
|dv(\zeta)| \ge \frac{1}{K},\quad \inf_{z \in\Sigma_g \setminus \{\zeta\}}
\frac{\dist(v(\zeta),v(z))}{\dist(\zeta,z)} \ge \frac{1}{K}, \quad\text{ and }\quad
\dist(v(\zeta),M \setminus \uU) \ge \frac{1}{K}.
$$
\item \label{item:critical}
\textsl{Critical orders do not increase:}
$$
\inf_{z \in \Sigma_g \setminus\Theta} \frac{|dv(z)|}{F_g(z)} \ge \frac{1}{K} \quad\text{ and }\quad
\inf_{z \in \Sigma_h \setminus\Theta'} \frac{|d\varphi(z)|}{F_h(z)} \ge \frac{1}{K}.
$$
\item \label{item:collide}
\textsl{Images of branch points do not collide:}
There exist distinct points $w_i = \varphi(\zeta_i^1) = \ldots = \varphi(\zeta_i^{q_i}) \in \Sigma_g$
for $i=1,\ldots,r$ such that
$$
\dist(w_i,w_j) \ge \frac{1}{K} \quad\text{ for all }\quad i,j=1,\ldots,r \text{ with } i \ne j.
$$
\item \label{item:kernel}
\textsl{Kernels do not get larger:}
Writing $\dot{E} := N_{\widehat{u}}$ and $\dot{F} := \overline{\Hom}_\CC(T\widehat{\Sigma},N_{\widehat{u}})$
for the canonically defined regular cover $\widehat{u} : \widehat{\Sigma} \to M$ of~$v$,
the operator $\dot{\mathbf{D}}_{\widehat{u}}^N : W^{k,p,-\boldsymbol{\delta}}(\dot{E}) \to
W^{k-1,p,-\boldsymbol{\delta}}(\dot{F})$ satisfies
$$
\left\| \dot{\mathbf{D}}_{\widehat{u}}^N \eta \right\|_{W^{k-1,p,-\boldsymbol{\delta}}} 
\ge \frac{1}{K} \inf_{\xi \in \ker\dot{\mathbf{D}}_{\widehat{u}}^N} \| \eta - \xi \|_{W^{k,p,-\boldsymbol{\delta}}}
\quad\text{ for all }\quad \eta \in W^{k,p,-\boldsymbol{\delta}}(N_{\widehat{u}}).
$$
\item \label{item:Petri}
\textsl{Curves remain Petri regular:}
For the regular cover~$\widehat{u}$, the Petri map 
$\Pi : \ker \dot{\mathbf{D}}_{\widehat{u}}^N \otimes \ker (\dot{\mathbf{D}}_{\widehat{u}}^N)^* \to
\Gamma(\dot{E} \otimes \dot{F})$ satisfies the estimate
$$
\|\Pi(t)\|_{C^0(\widehat{\Sigma}^K)} \ge \frac{1}{K} \|t\|,
$$
where $\widehat{\Sigma}^K := \left\{ z \in \widehat{\Sigma} \ \big|\ 
\dist(\widehat{u}(z) , M \setminus \uU) \ge 1/K \right\}$ and the
norm on the tensor product
$\ker \dot{\mathbf{D}}_{\widehat{u}}^N \otimes \ker (\dot{\mathbf{D}}_{\widehat{u}}^N)^*$
is defined via any norms on $\ker \dot{\mathbf{D}}_{\widehat{u}}^N$
and $\ker (\dot{\mathbf{D}}_{\widehat{u}}^N)^*$ that vary continuously
with $u \in \mM^d(J\,;\,\ell_1,\ldots,\ell_m\,;\,k,c)$.
\end{enumerate}
Clearly every element of $\mM^d_\Pi(J\,;\,\ell_1,\ldots,\ell_m\,;\,k,c)$
belongs to some $\nN^K(J)$ for $K \in \NN$ sufficiently large.  Now define
$$
\jJ^{\reg,K} \subset \jJ(M,\omega\,;\,\uU,\Jfix)
$$
via the property that $J \in \jJ^{\reg,K}$ if and only if every element of
$\nN^K(J)$ is regular in its stratum.  

We claim that $\jJ^{\reg,K}$ is open in $\jJ(M,\omega\,;\,\uU,\Jfix)$.
Indeed, suppose $J_\nu \in \jJ(M,\omega\,;\,\uU,\Jfix)$ is
a sequence converging to $J \in \jJ^{\reg,K}$ as $\nu \to \infty$
such that for every~$\nu$, there exists a curve $u_\nu \in \nN^K(J_\nu)$
that is not regular in its stratum.  Given parametrizations
$u_\nu = v_\nu \circ \varphi_\nu$ with $v_\nu : (\Sigma_g,j_\nu) \to (M,J_\nu)$
and $\varphi_\nu : (\Sigma_h,j_\nu') \to (\Sigma_g,j_\nu)$ satisfying the
conditions above, conditions~\ref{item:domain} and~\ref{item:bubbles} imply
via standard elliptic regularity arguments that there are
$C^\infty$-convergent subsequences $v_\nu \to v$, $j_\nu \to j$,
$\varphi_\nu \to \varphi$ and $j_\nu'\to j'$, so that $u_\nu$ itself
converges to the composition of a $J$-holomorphic curve 
$v : (\Sigma_g,j) \to (M,J)$ and another $d$-fold holomorphic branched
cover $\varphi : (\Sigma_h,j') \to (\Sigma,j)$.  Since all conditions
in the definition of $\nN^K(J)$ are closed, they are also satisfied for
the limit~$u$.  Condition~\ref{item:injective}
then guarantees that $v$ has an injective point mapped into~$\uU$, 
conditions~\ref{item:critical} and~\ref{item:collide} ensure that both $v$
and $\varphi$ satisfy the given constraints on critical orders and branching
data, and condition~\ref{item:kernel} implies via Lemma~\ref{lemma:FredGood}
below that
$\dim \ker\dot{\mathbf{D}}_{\widehat{u}}^N = \dim \ker\dot{\mathbf{D}}_{\widehat{u}_\nu}^N$.
It follows that $u \in \mM^d(J\,;\,\ell_1,\ldots,\ell_m\,;\,k,c)$,
thus $u$ also belongs to~$\nN^K(J)$ and must therefore be regular in its
stratum.  Regularity must then also hold for $u_\nu$ with $\nu$ sufficiently
large, since it is an open condition, and this is a contradiction.

The use of condition~\ref{item:kernel} in the above argument depends on
interpreting it in terms of
the injective map induced by $\dot{\mathbf{D}}_{\widehat{u}}^N$ on the quotient of its
domain by its kernel, and then feeding this into the following 
functional-analytic lemma:

\begin{lemma}
\label{lemma:FredGood}
Suppose $X$ and $Y$ are Banach spaces, $\mathbf{T}_n : X \to Y$ is a sequence
of Fredholm operators converging to a Fredholm operator $\mathbf{T} : X \to Y$,
and there exists a constant $c > 0$ such that
$$
\|\mathbf{T}_n x \|_Y \ge c \|\pi_n x\|_{X/\ker\mathbf{T}_n},
$$
where $\pi_n : X \to X / \ker\mathbf{T}_n$ is the quotient projection.
Then $\dim \ker\mathbf{T}_n = \dim \ker\mathbf{T}$ for all $n$ sufficiently large.
\end{lemma}
\begin{proof}
One can use the same trick as in the proof of Lemma~\ref{lemma:zeroset} to
find a sequence of Banach space isomorphisms $\boldsymbol{\Phi}_n : X \to X$
converging to $\1$ such that $\ker\mathbf{T}_n \subset \boldsymbol{\Phi}_n(\ker\mathbf{T})$
for every $n$ sufficiently large.  Then if $\dim \ker\mathbf{T}_n < \dim\ker\mathbf{T}$
for all~$n$, we can find a bounded sequence $x_n \in \boldsymbol{\Phi}_n(\ker\mathbf{T})$
such that the norm of $\pi_n(x_n)$ in $X / \ker\mathbf{T}_n$ is bounded
away from zero.  Equivalently, $x_n = \boldsymbol{\Phi}_n(v_n)$ for a
bounded sequence $v_n \in \ker\mathbf{T}$, which then has a subsequence
convergent to some $v_\infty \in \ker\mathbf{T}$ since $\dim\ker\mathbf{T} < \infty$, implying a corresponding
subsequence $x_n \to x_\infty$ and thus $\mathbf{T}_n x_n \to 0$.  The latter
contradicts the estimate in the hypothesis.
\end{proof}

We claim that $\jJ^{\reg,K}$ is also dense in $\jJ(M,\omega\,;\,\uU,\Jfix)$.
Since the reference structure $\Jref \in \jJ(M,\omega\,;\,\uU,\Jfix)$ can
be chosen arbitrarily, it suffices to find some $\varepsilon \in \seqs$ and
a sequence $J_\nu \in \jJ^{\reg,K}$ such that $J_\nu \to \Jref$ in the
$C_\varepsilon$-topology.  The argument used above for openness shows that
$\nN^K(\Jref)$ is compact, and condition~\ref{item:Petri} implies that
every curve in $\nN^K(\Jref)$ is Petri regular,
so by Lemma~\ref{lemma:epsReg2}, one can choose a lower bound for a finite
set of choices $\varepsilon \in \seqs$ and thus assume that every curve
in $\nN^K(\Jref)$ is $\varepsilon$-regular.  Now pick a sequence
$J_\nu \in \jJ_\varepsilon^\reg$ with $J_\nu \to \Jref$, and arguing by
contradiction, suppose $J_\nu \not\in \jJ^{\reg,K}$, meaning there exists
a sequence $u_\nu \in \nN^K(J_\nu)$ such that each $u_\nu$ is not regular
in its stratum.  After passing to a subsequence, the previous compactness
argument shows that $u_\nu$ converges to some $u \in \nN^K(\Jref)$, implying
that $u_\nu$ is $\varepsilon$-regular for all $\nu$ sufficiently large.
That contradicts the definition of $\jJ_\varepsilon^\reg$ and thus proves
the claim.

To conclude, $\bigcap_{K \in \NN} \jJ^{\reg,K}$ is now a Baire subset
of $\jJ(M,\omega\,;\,\uU,\Jfix)$ containing almost complex structures $J$
such that every Petri regular curve in
$\mM^d(J\,;\,\ell_1,\ldots,\ell_m\,;\,k,c)$ is
regular in its stratum.  By Theorem~\ref{thm:PetriJ}, we can intersect
this with another Baire subset in order to assume that every curve in
$\mM^d(J\,;\,\ell_1,\ldots,\ell_m\,;\,k,c)$ is Petri
regular.  The resulting Baire subset depends on the choices of data
$d$, $\mathbf{b}$, $G$, $g$, $m$, $A$, $\ell_1,\ldots,\ell_m$, $k$,
but since there are only countably many such choices, a further countable
intersection of Baire subsets now produces a Baire subset of almost complex
structures for which the result of Theorem~\ref{thm:submanifolds0} holds.
The proof of Theorem~\ref{thm:submanifolds0} is thus complete.

\section{Super-rigidity in dimension four}
\label{sec:dimensionFour}

We now prove the $4$-dimensional case of Theorem~\ref{thm:super},
using intersection-theoretic arguments that are essentially unrelated to
the rest of the paper.  Throughout this section, assume $(M,J)$ is
an almost complex manifold with
$$
\dim M = 4.
$$

The genus zero case is an ``automatic'' phenomenon, i.e.~it does not require
any genericity condition except for ensuring that the index~$0$ simple curve
is immersed:

\begin{prop}
\label{prop:genusZero}
Every simple immersed $J$-holomorphic sphere $v : (S^2,i) \to (M,J)$ of
index~$0$ in an almost complex $4$-manifold is super-rigid.
\end{prop}
\begin{proof}
Assume $\varphi : (\Sigma',j') \to (S^2,i)$ is a $d$-fold branched cover and 
$u = v \circ \varphi$.
Since $v$ is immersed, the Riemann-Roch formula implies
$$
0 = \ind(v) = \ind \mathbf{D}_v^N = \chi(S^2) + 2 c_1(N_v),
$$
hence $c_1(N_v) = -1$.  Then $c_1(N_u) = c_1(\varphi^*N_v) = -d$, so
if $\eta \in \ker \mathbf{D}_u^N$ is nontrivial, its algebraic count of zeroes
is negative, violating the similarity principle.
\end{proof}

For the genus one case, we use a variant of the ``magic trick'' proposed by
Hutchings \cite{Hutchings:magic} in the context of Embedded Contact Homology.

\begin{prop}
\label{prop:genusOne}
A simple immersed $J$-holomorphic torus $v : (\TT^2,j) \to (M,J)$ of
index~$0$ in an almost complex $4$-manifold is super-rigid if and only if 
all its unbranched covers are Fredholm regular.
\end{prop}
\begin{proof}
We will assume for most of the proof that $v : (\Sigma,j) \to (M,J)$ has
unspecified genus $g \ge 1$.  Since $v$ is immersed
with index~$0$, it is regular if and only if its normal Cauchy-Riemann operator
$\mathbf{D}_v^N$ is injective, so given this and the assumption that the
same holds for all unbranched covers $u = v \circ \varphi$,
we need to show that
$\mathbf{D}_u^N$ is also injective for $u = v \circ \varphi$ 
where $\varphi : (\Sigma',j') \to (\Sigma,j)$ is \emph{any} holomorphic branched cover.
We will prove this by induction on the degree $d := \deg(\varphi)$, thus
assume it is true for all covers up to degree $d-1$.  Note that
since $\ind(v) = 0$, we have
\begin{equation}
\label{eqn:c1Condition}
\ind \mathbf{D}_v^N = \chi(\Sigma) + 2 c_1(N_v) = 0,
\end{equation}
and if $\varphi$ has branch points, then $\Sigma'$ has genus $g' > 1$ by the Riemann-Hurwitz
formula.

By the construction in the proof of Proposition~\ref{prop:isolated}, one can 
endow the total space of the normal bundle $\pi : N_v \to \Sigma$ with an almost 
complex structure $J_N$ such that $J_N$-holomorphic curves 
$u_\eta : (S,i) \to (N_v,J_N)$ correspond to sections 
$\eta \in \ker \mathbf{D}_{v \circ \psi}^N$ along holomorphic branched covers
$\psi = \pi \circ u_\eta : (S,i) \to (\Sigma,j)$.  If $\ker \mathbf{D}_u^N$
contains a nontrivial element $\eta$, the inductive hypothesis implies that
the corresponding $J_N$-holomorphic curve $u_\eta$ is somewhere injective.
We can view $v$ itself as a $J_N$-holomorphic embedding into~$N_v$, and
$u_\eta$ is homologous to its $d$-fold cover, so applying the adjunction 
formula to both $u_\eta$ and~$v$ as $J_N$-holomorphic curves in~$N_v$,
\begin{equation*}
\begin{split}
u_\eta \bullet u_\eta &= 2 \delta(u_\eta) + c_1(u_\eta^*T N_v) - \chi(\Sigma')
= 2 \delta(u_\eta) + d \cdot c_1(v^*T N_v) - \chi(\Sigma') \\
&= d^2 (v \bullet v) = d^2 \cdot c_1(N_v) = d^2 \cdot c_1(v^*T N_v) -
d^2 \cdot \chi(\Sigma),
\end{split}
\end{equation*}
where $\delta(u_\eta) \ge 0$ denotes the algebraic count of double points
and critical points of~$u_\eta$.  Solving for $\delta(u_\eta)$ and plugging
in \eqref{eqn:c1Condition} to compute $c_1(v^*TN_v) = \chi(\Sigma) + c_1(N_v)
= \frac{1}{2} \chi(\Sigma) = 1 - g$, we have
\begin{equation*}
\begin{split}
2 \delta(u_\eta) &= d(d-1) \cdot c_1(v^*TN_v) - d^2 \cdot \chi(\Sigma) + \chi(\Sigma') \\
&= d (d-1) (1 - g) - 2 d^2 (1-g) + 2 - 2 g' = d(d+1)(g-1) - 2(g'-1)
\end{split}
\end{equation*}
Plugging in $g=1$ and the fact that $g' > 1$, this gives a contradiction
since $\delta(u_\eta)$ cannot be negative.
\end{proof}

\begin{remark}
\label{remark:bifurcationsFour}
In the spirit of \S\ref{sec:bifurcations}, the two results above show that
the story of super-rigidity and bifurcations is simpler in dimension
four.  In the genus zero case bifurcations can be avoided altogether:
since having a critical point is a codimension~$2$ condition
(see Appendix~\ref{sec:critical}), index~$0$ simple curves for generic
$1$-parameter families of almost complex structures can be assumed immersed,
and therefore super-rigid by Prop.~\ref{prop:genusZero}.
This is no longer true in the genus one case since regularity of some
unbranched cover might fail under a generic homotopy, producing the
birth-death or degree-doubling bifurcations in \cite{Taubes:counting},
but Prop.~\ref{prop:genusOne} implies that this
is the only danger---the only bifurcations that can happen involve
unbranched covers with $g'=1$ and $d \in \{1,2\}$, and they are already described 
in \cite{Taubes:counting}.
\end{remark}

\appendix

\section{Moduli spaces with prescribed orders of critical points}
\label{sec:critical}

The proposition below is well known to experts, but a proof of it is difficult
to find in the literature, so we will sketch one here.

Fix a symplectic manifold $(M,\omega)$ of dimension~$2n$, $n \in \NN$,
and suppose $J \in \jJ(M,\omega)$.  Recall that if $(\Sigma,j)$ is a connected
Riemann surface and $u : (\Sigma,j) \to
(M,J)$ is a nonconstant $J$-holomorphic curve with a critical point
$du(z) = 0$, then the critical point is isolated and has a well-defined 
positive \defin{order},
$$
\ord(du;z) \in \NN,
$$
characterized by the property that $\ord(du;z) = \ell$ if $z$ is a zero
of order $\ell$ for the section $du \in \Gamma(\Hom_\CC(T\Sigma,u^*TM))$, where the 
latter is viewed as a holomorphic section with respect to a natural holomorphic
bundle structure on $u^*TM$ determined by the linearized Cauchy-Riemann operator,
see e.g.~\cite{Wendl:automatic}*{\S 3.3}.  
When $(\Sigma,j)$ is closed, we
denote the resulting algebraic count of critical points by
$$
Z(du) := \sum_{\{z \in \Sigma\ |\ du(z) = 0\}} \ord(du;z) \ge 0,
$$
and note that it vanishes if and only if $u$ is immersed.  Given integers
$g,m \ge 0$, a homology class $A \in H_2(M)$ and a tuple of positive 
integers $(\ell_1,\ldots,\ell_m)$, let
$$
\mM_{g,m}(A,J\,;\, \ell_1,\ldots,\ell_m) \subset \mM_{g,m}(A,J)
$$
denote the following subset of the moduli space of unparametrized $J$-ho\-lo\-mor\-phic
curves homologous to $A$ with genus~$g$ and $m$ marked points: 
a map $u : (\Sigma,j) \to (M,J)$ with marked points $\zeta_1,\ldots,\zeta_m \in \Sigma$
representing an element of $\mM_{g,m}(A,J)$ belongs to 
$\mM_{g,m}(A,J\,;\, \ell_1,\ldots,\ell_m)$ if and only if it is critical
at all marked points,
$$
\ord(du;\zeta_j) = \ell_j \quad\text{ for } \quad j=1,\ldots,m,
$$
and it is immersed everywhere else.

\begin{prop}
\label{prop:critical}
Fix an open subset $\uU \subset M$ with compact closure and a compatible
almost complex structure $\Jfix \in \jJ(M,\omega)$.  There exists a Baire
subset
$$
\jJ^\reg \subset \jJ(M,\omega\,;\,\uU,\Jfix)
$$
such that for all $J \in \jJ^\reg$ and all
$g,m \ge 0$, $A \in H_2(M)$ and $(\ell_1,\ldots,\ell_m) \in \NN^m$, 
the open subset of $\mM_{g,m}(A,J\,;\, \ell_1,\ldots,\ell_m)$ consisting of
somewhere injective curves that pass through $\uU$ is a smooth manifold with
dimension equal to its virtual dimension, where
$$
\virdim \mM_{g,m}(A,J\,;\, \ell_1,\ldots,\ell_m) = \virdim \mM_g(A,J) -
\sum_{i=1}^m (2 n \ell_i - 2).
$$
\end{prop}
\begin{cor}
For generic compatible $J$ in any closed symplectic $2n$-manifold, all
closed, connected and somewhere injective $J$-holomorphic curves $u$ 
with $m \ge 0$ critical points satisfy $\ind(u) \ge 2n Z(du) - 2m$.
\end{cor}

One well-known consequence of this result is that for generic~$J$,
somewhere injective index~$0$ curves in almost complex manifolds of dimension
at least four are always immersed.  Another proof of this is given in
\cite{OhZhu:embedding}, though it is analytically somewhat more complicated
than the one given below.

It will suffice to prove that the same statement as in Prop.~\ref{prop:critical}
holds for the slightly larger moduli space
$$
\widehat{\mM}_{g,m}(A,J\,;\, \ell_1,\ldots,\ell_m)
$$
characterized by the condition $\ord(du;\zeta_j) \ge \ell_j$ for all
$j=1,\ldots,m$ without requiring $u$ to be immersed outside the marked points.
Indeed, $\mM_{g,m}(A,J\,;\, \ell_1,\ldots,\ell_m) \subset
\widehat{\mM}_{g,m}(A,J\,;\, \ell_1,\ldots,\ell_m)$ is an open subset.
We shall borrow from Zehmisch \cite{Zehmisch:jets} the notion of
\emph{holomorphic jets}: given a point $p$ in an almost complex manifold
$(M,J)$ and an integer $r > 0$, a \defin{holomorphic $r$-jet} at $p$ is
an equivalence class of $J$-holomorphic curves
$$
u : (\DD_\epsilon,i) \to (M,J)
$$
with $u(0) = p$, where $(\DD_\epsilon,i)$ denotes the $\epsilon$-disk in~$\CC$,
and two curves are considered equivalent if their partial derivatives at $0$ 
match up to order~$r$.  The nonlinear Cauchy-Riemann equation implies that
the holomorphic $r$-jet represented by $u$ is determined by the holomorphic
part of its Taylor polynomial of degree~$r$ (see \cite{Wendl:lecturesV2}*{Prop.~2.99}),
and moreover, every holomorphic Taylor polynomial of degree~$r$ is realizable
as the $r$-jet of a local $J$-holomorphic curve (\cite{Wendl:lecturesV2}*{Theorem~2.100}).
Thus the space of all holomorphic $r$-jets at $p$ is a real $2rn$-dimensional
vector space, and the union of these spaces for all $p \in M$ forms a smooth
manifold 
$$
\Jethol^r(M)
$$
of real dimension $2n(r+1)$.

We shall analyze the local structure of $\widehat{\mM}_{g,m}(A,J\,;\, \ell_1,\ldots,\ell_m)$
following a minor modification of the scheme outlined in
\cite{Wendl:lecturesV2}*{Chapter~4}.  For simplicity, we shall
assume in this exposition that $2g + m \ge 3$, so that we only need to
deal with \emph{stable} marked Riemann surfaces. (For the finitely many
non-stable cases, see Remark~\ref{remark:nonstable}.)
Given $(\Sigma,j_0,\Theta,u_0)$
representing an element of $\widehat{\mM}_{g,m}(A,J\,;\, \ell_1,\ldots,\ell_m)$,
with marked points $\Theta := (\zeta_1,\ldots,\zeta_m)$,
choose a \defin{Teichm\"uller slice} through~$j_0$: this means a smooth 
$(6g-6+2m)$-dimensional family
$\tT$ of complex structures on $\Sigma$ that includes $j_0$ and parametrizes
a neighborhood of $[j_0]$ in the Teichm\"uller space of complex structures
modulo diffeomorphisms that are homotopic to the identity and fix~$\Theta$.
The tangent space $T_{j_0} \tT$ is also required
to define a closed complement of the image of the canonical Cauchy-Riemann
operator on $T\Sigma$ restricted to the space of vector fields vanishing at~$\Theta$,
cf.~\cite{Wendl:lecturesV2}*{Definition~4.29}.  Moreover, we can arrange
for $\tT$ to have the following two properties (cf.~\cite{Wendl:automatic}*{Lemmas~3.3 and~3.4}):
\begin{itemize}
\item $\tT$ is invariant under the action of the group $\Aut(\Sigma,j_0,\Theta)$
of biholomorphic maps on $(\Sigma,j_0)$ fixing~$\Theta$;
\item There exists a neighborhood of $\Theta$ on which every $j \in \tT$ matches $j_0$.
\end{itemize}
Now let $r := \max\{\ell_1,\ldots,\ell_m\}$, and choose any $k \in \NN$ and $p \in (1,\infty)$
such that
\begin{equation}
\label{eqn:SobolevCondition}
(k-r)p > 2,
\end{equation}
so the Sobolev embedding theorem implies that functions of class $W^{k,p}$
on $\Sigma$ are also in~$C^r$.  We define the Banach manifold
$$
\bB := W^{k,p}(\Sigma,M)
$$
and smooth Banach space bundle $\eE \to \tT \times \bB$ with fibers
$$
\eE_{(j,u)} := W^{k-1,p}\big(\overline{\Hom}_\CC((T\Sigma,j),(u^*TM,J))\big),
$$
so that
$$
\dbar_J : \tT \times \bB \to \eE : (j,u) \mapsto Tu + J \circ Tu \circ j
$$
defines a smooth section.  We say that $(\Sigma,j_0,\Theta,u_0)$ is
\defin{Fredholm regular} if the linearization 
$$
D\dbar_J(j_0,u_0) : T_{j_0}\tT \oplus W^{k,p}(u_0^*TM) \to
W^{k-1,p}\big(\overline{\Hom}_\CC((T\Sigma,j_0),(u_0^*TM,J))\big)
$$
of this section at $(j_0,u_0)$ is surjective, in which case a neighborhood of
$(j_0,u_0)$ in $\dbar_J^{-1}(0)$ is a smooth finite-dimensional manifold,
and its quotient by the natural action of $\Aut(\Sigma,j_0,\Theta)$ can
be identified naturally with a neighborhood of $[(\Sigma,j_0,\Theta,u_0)]$
in $\mM_{g,m}(A,J)$.  To incorporate the critical point condition,
fix holomorphic coordinates identifying a neighborhood of each marked
point $\zeta_j$ with the standard unit disk $(\DD,i)$; note that this can
be done for all $j \in \tT$ at once since they are assumed to match $j_0$
near~$\Theta$.  Then since $\bB$ has a continuous inclusion into
$C^r(\Sigma,M)$, there is a well-defined and smooth\footnote{The smoothness
of $\ev$ is clear because it is the restriction to $\dbar_J^{-1}(0)$ of a
map $\bB \to \Jethol^{\ell_1}(M) \times \ldots \times \Jethol^{\ell_m}(M)$
which in the natural Banach manifold charts provided by
\cite{Eliasson} looks like a linear map
evaluating derivatives of functions at the fixed points $\Theta \subset \Sigma$.
This works because we are choosing to represent elements of
$\mM_{g,m}(A,J)$ by maps with marked
points at fixed positions; of course there is no actual constraint on the
movement of the marked points, but this freedom is seen in our setup by
varying $j$ in $\tT$ instead of varying the points $\zeta_1,\ldots,\zeta_m$.
This is a notable difference from the setup in \cite{OhZhu:embedding}.}
\emph{jet evaluation map}
$$
\ev : \dbar_J^{-1}(0) \to \Jethol^{\ell_1}(M) \times \ldots \times \Jethol^{\ell_m}(M),
$$
whose $i$th factor for $i=1,\ldots,m$ is the holomorphic $\ell_i$-jet
represented by $u$ in its parametrization by $(\DD,i)$ at~$\zeta_i$.
We will say that $(\Sigma,j_0,\Theta,u_0)$ is \defin{regular for the
constrained moduli space} $\widehat{\mM}_{g,m}(A,J\,;\, \ell_1,\ldots,\ell_m)$ if
it is Fredholm regular and the jet evaluation map is transverse to the
submanifold 
$$
Z \subset \Jethol^{\ell_1}(M) \times \ldots \times \Jethol^{\ell_m}(M)
$$
consisting of $m$-tuples of jets of constant maps.
Note that this condition does not depend on the chosen holomorphic coordinates
near the marked points, as it is equivalent to the condition that $u$ should
have vanishing derivatives up to order $\ell_i$ at $\zeta_i$ for each $i=1,\ldots,m$.
Whenever the regularity condition is satisfied, $\ev^{-1}(Z) \subset
\dbar_J^{-1}(0)$ inherits the structure of a smooth submanifold with
real codimension $2n\sum_i \ell_i$, so $\widehat{\mM}_{g,m}(A,J\,;\, \ell_1,\ldots,\ell_m)$ 
in general becomes an orbifold near $[(\Sigma,j_0,\Theta,u_0)]$, with
\begin{equation*}
\begin{split}
\dim \widehat{\mM}_{g,m}(A,J\,;\, \ell_1,\ldots,\ell_m) &=
\dim \mM_{g,m}(A,J) - 2n \sum_i \ell_i \\
&= \dim \mM_g(A,J) + 2m - 2n\sum_i \ell_i \\
&= \dim \mM_g(A,J) - \sum_{i=1}^m (2n \ell_i - 2).
\end{split}
\end{equation*}

To prove that the constrained regularity condition can be achieved generically,
fix $\Jref \in \jJ(M,\omega\,;\,\uU,\Jfix)$ and a suitable sequence of positive
numbers $\varepsilon_\nu \to 0$, and consider a Banach manifold
$\jJ_\varepsilon$ of almost complex structures in $\jJ(M,\omega\,;\,\uU,\Jfix)$
that are $C_\varepsilon$-close to~$\Jref$ (cf.~\S\ref{sec:PetriJ}).
This gives rise to two universal moduli spaces,
$$
\univ^*(\jJ_\varepsilon) := \left\{ (u,J) \ \big|\ 
\text{$J \in \jJ_\varepsilon$ and $u \in \mM^*_{g,m}(A,J)$} \right\}
$$
and
$$
\widehat{\univ}^*(\jJ_\varepsilon\,;\,\ell_1,\ldots,\ell_m) := \left\{ (u,J) \ \big|\ 
\text{$J \in \jJ_\varepsilon$ and $u \in \widehat{\mM}^*_{g,m}(A,J\,;\,\ell_1,\ldots,\ell_m)$} \right\},
$$
where we abbreviate by
\begin{equation*}
\begin{split}
\mM^*_{g,m}(A,J) &\subset \mM_{g,m}(A,J),\\
\widehat{\mM}^*_{g,m}(A,J\,;\,\ell_1,\ldots,\ell_m) &\subset \widehat{\mM}_{g,m}(A,J\,;\,\ell_1,\ldots,\ell_m)
\end{split}
\end{equation*}
the subspaces defined via the condition that $u$ be somewhere injective and
pass through~$\uU$.  As is well known, $\univ^*(\jJ_\varepsilon)$ is a 
separable and metrizable smooth Banach manifold
if $\varepsilon_\nu$ converges to $0$ fast enough, and for
$[(\Sigma,j_0,\Theta,u_0)] \in \mM^*_{g,m}(A,J_0)$, a neighborhood of
$(u_0,J_0)$ in $\univ^*(\jJ_\varepsilon)$ can be identified with the zero-set
of a smooth section
$$
\dbar : \tT \times \bB \times \jJ_\varepsilon \to \eE : (j,u,J) \mapsto
\dbar_J(u),
$$
where $\eE$ now denotes the Banach space bundle with fibers
$$
\eE_{(j,u,J)} = W^{k-1,p}\big(\overline{\Hom}_\CC((T\Sigma,j),(u^*TM,J))\big).
$$
The tangent space $T_{(u_0,J_0)} \univ^*(\jJ_\varepsilon)$ is the kernel of
the surjective operator
\begin{equation*}
\begin{split}
\mathbf{L} := D\dbar(j_0,u_0,J_0) : T_{j_0} \tT \oplus W^{k,p}(u_0^*TM) \oplus T_{J_0}\jJ_\varepsilon &\to
W^{k-1,p}(\overline{\Hom}_\CC(T\Sigma,u_0^*TM)) \\
(y,\eta,Y) &\mapsto J_0 \circ Tu_0 \circ y + \mathbf{D}_{u_0} \eta + Y \circ Tu_0 \circ j_0,
\end{split}
\end{equation*}
where $\mathbf{D}_{u_0}$ is the linearized Cauchy-Riemann operator associated
to $u_0 : (\Sigma,j_0) \to (M,J_0)$.  We can again define the smooth jet
evaluation map
\begin{equation}
\label{eqn:jetEval}
\ev : \dbar^{-1}(0) \to \Jethol^{\ell_1}(M) \times \ldots \times \Jethol^{\ell_m}(M)
\end{equation}
and identify a neighborhood of $(u_0,J_0)$ in
$\widehat{\univ}^*(\jJ_\varepsilon\,;\,\ell_1,\ldots,\ell_m)$ with $\ev^{-1}(Z)$.
The main technical ingredient behind Proposition~\ref{prop:critical} is now the
following.

\begin{lemma}
\label{lemma:submersion}
If $\varepsilon_\nu \to 0$ fast enough, then the jet evaluation map \eqref{eqn:jetEval} is a submersion.
\end{lemma}
\begin{proof}
We need to show that for any $X \in T_{\ev(u_0)}(\Jethol^{\ell_1}(M) \times \ldots \times \Jethol^{\ell_m}(M))$,
there exists an element $(y,\eta,Y) \in \ker \mathbf{L}$ with
$$
d\ev(u_0) \eta = X.
$$
Let us first observe that this problem can be solved locally near the marked
points: in fact, there exists a smooth section $\eta \in \Gamma(u_0^*TM)$ with
$$
\mathbf{D}_{u_0} \eta = 0 \text{ near $\Theta$} \quad \text{ and } \quad
d\ev(u_0) \eta = X.
$$
This follows from the local existence theorem for $J$-holomorphic curves with
prescribed holomorphic derivatives at a point, cf.~\cite{Wendl:lecturesV2}*{Theorem~2.100}.
More precisely, choose a smooth path $\gamma = (\gamma_1,\ldots,\gamma_m) : (-\delta,\delta) \to
\Jethol^{\ell_1}(M) \times \ldots \times \Jethol^{\ell_m}(M)$ with
$\gamma(0) = \ev(u_0)$ and $\dot{\gamma}(0) = X$.  Then the local existence
theorem provides for each $i=1,\ldots,m$ a smooth family of $J$-holomorphic curves
$u_\tau^{(i)} : \DD_\epsilon \to M$ defined on sufficiently small disks $\DD_\epsilon \subset \CC$
such that the holomorphic $\ell_i$-jet represented by $u_\tau^{(i)}$ is $\gamma_i(\tau)$
for each~$\tau$.  The desired section 
$\eta \in \Gamma(u_0^*TM)$ can now be constructed by writing it in our
chosen holomorphic coordinates near each marked point $\zeta_i$ as
$\p_\tau u_\tau^{(i)}|_{\tau=0}$ and then extending it arbitrarily outside these
neighborhoods.

Given $\eta$ as above, we aim now to find a pair
$(\xi,Y) \in W^{k,p}(u_0^*TM) \oplus T_{J_0}\jJ_\varepsilon$ such that
$$
\mathbf{L}(0,\eta+\xi,Y) = \mathbf{L}(0,\xi,Y) + \mathbf{D}_{u_0}\eta = 0
\quad \text{ and }\quad
d\ev(u_0)\xi = 0,
$$
in which case $(0,\eta+\xi,Y) \in T_{(u_0,J_0)} \univ^*(\jJ_\varepsilon)$
and $d\ev(u_0,J_0) (0,\eta+\xi,Y) = X$.
We will use the weighted Sobolev spaces described in \S\ref{sec:CRpunctured}.
Let $\dot{\Sigma} := \Sigma \setminus \Theta$, and assume without loss of
generality that $u_0^{-1}(\uU) \subset \Sigma$ is disjoint from~$\Theta$;
this can be achieved at the cost of shrinking $\uU$ and therefore the
space of perturbations~$\jJ_\varepsilon$.  As a consequence, $Y \circ Tu_0 \circ j_0$
now has compact support in $\dot{\Sigma}$ for any $Y \in T_{J_0} \jJ_\varepsilon$.
Using the fixed holomorphic coordinates on neighborhoods of marked points
$\zeta_i \in \Theta$, we can identify them biholomorphically with half-cylinders
$[0,\infty) \times S^1$ and fix trivializations of $u_0^*TW$ on these
neighborhoods to define weighted Sobolev norms and a bounded linear map
$$
\dot{\mathbf{D}}_{u_0} : W^{k,p,\delta}(u_0^*TM|_{\dot{\Sigma}}) \to
W^{k-1,p,\delta}\big(\overline{\Hom}_\CC(T\dot{\Sigma},u_0^*TM)|_{\dot{\Sigma}}\big),
$$
where sections $\eta$ of class $W^{k,p,\delta}$ are required to satisfy
$e^{\delta s} \eta \in W^{k,p}([0,\infty) \times S^1)$ when expressed in the chosen trivialization
and holomorphic coordinates $(s,t) \in [0,\infty) \times S^1$ on each
cylindrical end near~$\Theta$.  As explained in \S\ref{sec:CRpunctured},
$\dot{\mathbf{D}}_{u_0}$ is asymptotic to the trivial asymptotic operator at
each puncture and is thus Fredholm for any $\delta \in \RR \setminus 2\pi\ZZ$.
We claim that whenever this condition is satisfied, the linear map
\begin{equation*}
\begin{split}
\mathbf{L}_\delta : W^{k,p,\delta}(u_0^*TM|_{\dot{\Sigma}}) \oplus
T_{J_0} \jJ_\varepsilon &\to
W^{k-1,p,\delta}\big(\overline{\Hom}_\CC(T\dot{\Sigma},u_0^*TM)|_{\dot{\Sigma}}\big) \\
(\xi,Y) &\mapsto \dot{\mathbf{D}}_{u_0} \xi + Y \circ Tu_0 \circ j_0
\end{split}
\end{equation*}
is surjective.  The proof is more or less standard: we start with the case $k=1$
and note that since
$\dot{\mathbf{D}}_{u_0}$ is Fredholm, $\mathbf{L}_\delta$ has closed range,
so it is not surjective if and only if there exists a nontrivial section
$\lambda \in (L^{p,\delta})^* = L^{q,-\delta}$ for $1/p + 1/q = 1$ which is
$L^2$-orthogonal to the images of both $\eta \mapsto \dot{\mathbf{D}}_{u_0}\eta$
and $Y \mapsto Y \circ Tu_0 \circ j_0$.  Since $u_0$ has an injective point
$z_0 \in \dot{\Sigma}$ with $u(z_0) \in \uU$, the latter implies that
$\lambda$ vanishes near~$z_0$; this depends on $\varepsilon_\nu$ converging to $0$
fast enough for $T_{J_0}\jJ_\varepsilon$ to contain an abundance of bump functions
with arbitrarily small support.  The former implies in turn that $\lambda$ is a weak
solution to the formal adjoint equation $\dot{\mathbf{D}}_{u_0}^*\lambda = 0$
and is therefore smooth with isolated zeroes, giving a contradiction.
The case of general $k \in \NN$ follows from this via elliptic regularity,
namely Lemma~\ref{lemma:regularity}.

With this claim in place, we observe that $- \mathbf{D}_{u_0}\eta$ vanishes
near $\Theta$ and thus restricts to $\dot{\Sigma}$ as a section of class
$W^{k-1,p,\delta}$ for any $\delta > 0$, thus we can find $\xi \in W^{k,p,\delta}(u_0^*TM|_{\dot{\Sigma}})$
and $Y \in T_{J_0}\jJ_\varepsilon$ such that
$$
\mathbf{L}(0,\xi,Y) = -\mathbf{D}_{u_0}\eta \quad \text{ on } \quad \dot{\Sigma}.
$$
Since $Y$ has compact support in $\dot{\Sigma}$ and $\mathbf{D}_{u_0}\eta = 0$
near~$\Theta$, this equation implies $\mathbf{D}_{u_0} \xi = 0$ near~$\Theta$.
The continuous inclusion $W^{k,p,\delta} \hookrightarrow C^0$ implies
that $\xi$ also has a continuous extension over $\Sigma$ that vanishes 
on~$\Theta$; moreover, since \eqref{eqn:SobolevCondition} implies
a continuous inclusion $W^{k,p} \hookrightarrow C^1$, $\xi$ has a bounded
first derivative on the cylindrical ends, implying via a short computation
that for $1 < q < 2$, the $L^q$-norm of its derivative on punctured
disk-like neighborhoods of $\Theta$ is finite.  It follows that the extension
of $\xi$ over the punctures is in $W^{1,q}$ on~$\Sigma$, and elliptic
regularity then implies that it is smooth everywhere.  Finally, the exponential
weight condition implies that in each holomorphic coordinate system identifying
the neighborhood of a marked point $\zeta_i \in \Theta$ with $\DD$ such that
$\zeta_i$ is at the origin, we have
$$
|\xi(z)| \le c |z|^{\delta/2\pi}
$$
for some constant $c > 0$.
But the choice of $\delta > 0$ in this discussion was arbitrary, so choosing
it large enough, we can arrange for $\xi$ to have vanishing derivatives
of arbitrarily large finite order at~$\Theta$, proving
$d\ev(u_0) \xi = 0$.
\end{proof}

The lemma implies that $\widehat{\univ}^*(\jJ_\varepsilon\,;\,\ell_1,\ldots,\ell_m)$ is
a separable and metrizable smooth Banach manifold, so we can now apply the
Sard-Smale theorem to the projection 
$$
\widehat{\univ}^*(\jJ_\varepsilon\,;\,\ell_1,\ldots,\ell_m)
\to \jJ_\varepsilon : (u,J) \mapsto J,
$$
giving a Baire subset
of $\jJ_\varepsilon$ for which $\widehat{\mM}^*_{g,m}(A,J\,;\,\ell_1,\ldots,\ell_m)$ is
a manifold of the correct dimension, and the countable intersection of these
subsets for all $g$, $m$, $A$ and $(\ell_1,\ldots,\ell_m)$ is again
comeager in~$\jJ_\varepsilon$, proving that there is a $C^\infty$-dense 
subset of $\jJ(M,\omega\,;\,\uU,\Jfix)$ for which the statement of the
theorem holds.  To turn this into a Baire subset of $\jJ(M,\omega\,;\,\uU,\Jfix)$,
one can use the standard Taubes trick (see e.g.~\cite{Wendl:lecturesV2}*{\S 4.4.2}): 
present $\widehat{\mM}^*_{g,m}(A,J\,;\,\ell_1,\ldots,\ell_m)$ as a countable union of compact
subsets, and associate to each one a set of regular almost complex structures,
which is open by construction and dense due to the argument above, so its
intersection is comeager.

\begin{remark}
Lemma~\ref{lemma:submersion} implies that for generic $J$, the jet evaluation
map can be made transverse to any given submanifold, hence this method can
be used to understand any moduli space of holomorphic curves with marked 
points satisfying conditions on their derivatives, e.g.~the incidence/tangency
conditions studied by Cieliebak-Mohnke \cites{CieliebakMohnke:transversality,
CieliebakMohnke:Audin} or McDuff-Siegel \cite{McDuffSiegel:tangency}.
\end{remark}

\begin{remark}
\label{remark:nonstable}
The assumption $2g + m \ge 3$ misses only four special cases, and for these
the discussion above is modified as follows:
\begin{enumerate}
\item The automorphism group $\Aut(\Sigma,j_0,\Theta)$ is not finite, but
is instead a nontrivial Lie group;
\item The usual formula $\dim \tT = 6g - 6 + 2m$ for the dimension of Teichm\"uller space
is wrong.
\end{enumerate}
In fact, these two differences cancel each other out in the sense that
$$
\dim \tT - \dim \Aut(\Sigma,j_0,\Theta) = 6g - 6 + 2m,
$$
which is why the stated formulas for the virtual dimensions of the moduli spaces $\mM_{g,m}(A,J)$
and $\mM_{g,m}(A,J\,;\,\ell_1,\ldots,\ell_m)$ remain correct in these
non-stable cases.  In the cases with genus zero, Teichm\"uller space is
trivial and there is thus no need to include a Teichm\"uller slice in the
argument; the only difference is then the fact that 
dividing $\dbar_J^{-1}(0)$ by $\Aut(\Sigma,j_0,\Theta)$ changes its dimension.
There is no need to discuss the non-stable genus one case here since that case also
has $m=0$, and thus does not involve critical point constraints.
\end{remark}

\section{Super-rigid curves are isolated}
\label{sec:isolated}

In this appendix we prove the following precise version of the statement
that the multiple covers of a super-rigid curve form an open and closed
subset of the ambient moduli space.

\begin{prop}
\label{prop:isolated}
Suppose $(M,J_k)$ is a sequence of almost complex manifolds with $J_k \to J_\infty$
in $C^\infty$ on some compact subset containing a super-rigid
$J_\infty$-holomorphic curve $u_\infty : ({\Sigma},j_\infty) \to
(M,J_\infty)$.  Then for sufficiently large~$k$, there exists a
sequence of $J_k$-holomorphic
curves $u_k : ({\Sigma},j_k) \to (M,J_k)$ with $j_k \to j_\infty$
and $u_k \to u_\infty$ in $C^\infty$,
and if $v_k$ is any sequence of smooth closed $J_k$-holomorphic
curves Gromov-convergent to a stable nodal $J_\infty$-holomorphic curve
with image contained in $u_\infty(\Sigma)$, then for all $k$
sufficiently large, every $v_k$ is either a biholomorphic reparametrization
or a multiple cover of~$u_k$.
\end{prop}

Note that this statement belongs to the almost complex category and makes no
reference to any symplectic structure.  Other than that detail, a nearly
identical statement has been proved before by Zinger, see 
\cite{Zinger:comparison}*{Prop.~3.2}.  The proof given below is essentially
the same and is included mainly for the sake of completeness;
it just requires the extra step of introducing an auxiliary
symplectic structure in order to use Gromov's compactness theorem.
Recall from \S\ref{sec:defn} that if $u \in \mM_g(A,J)$ and 
$d \ge 1$ and $h \ge 0$ are integers, we denote by
$$
\widebar{\mM}_h(d ; u) \subset \widebar{\mM}_h(dA,J)
$$
the moduli space of all stable nodal $d$-fold covers of $u$ with
arithmetic genus~$h$.

Suppose $J_k \to J_\infty$ is a $C^\infty$-convergent sequence of almost complex 
structures on a manifold $M$, and 
$[(\Sigma,j_\infty,u_\infty)] \in \mM_g(A,J_\infty)$
is a super-rigid curve.  Then $u_\infty$ is Fredholm
regular with index~$0$, so the implicit function theorem implies the
existence of curves $u_k  : (\Sigma,j_k) \to (M,J_k)$ for sufficiently
large~$k$ such that $j_k \to j_\infty$ and $u_k \to u_\infty$ in
$C^\infty$; these curves are unique up to biholomorphic reparametrization,
and are also simple and immersed for sufficiently large~$k$.
Assume $v_k \in \mM_h(dA,J_k)$
is a sequence of $J_k$-holomorphic curves converging to a nodal
cover $\tilde{u} \in \widebar{\mM}_h(d ; u_\infty)$ for some $d > 0$.
We will show that if the curves $v_k$ are not covers of $u_k$ for all
sufficiently large~$k$, then rescaling the normal fibers near $u_k$ as
$k \to \infty$ gives rise to a nontrivial section in the kernel of the 
normal Cauchy-Riemann operator on some cover of~$u_\infty$, contradicting
super-rigidity.

Choose a convergent sequence of $J_k$-invariant Riemannian metrics 
and corresponding Levi-Civita connections~$\nabla^k$.
Since the maps $u_k$ are immersed, we can define $J_k$-invariant 
normal bundles $N_{u_k} \to \Sigma$ 
as the orthogonal complements of $\im du_k$.  These are all isomorphic
as real vector bundles, so we can identify them all with the real bundle
$N := N_{u_\infty} \subset u_\infty^*TM$ carrying
a sequence of complex structures
$$
(N,J_k) \stackrel{\pi}{\to} (\Sigma,j_k),
$$
and then use the sequence of exponential maps determined by $\nabla^k$
to define a $C^\infty$-convergent sequence of immersions
$$
\Psi_k : \nN(\Sigma) \to M
$$
of some fixed neighborhood $\nN(\Sigma) \subset N$ of the zero section
$\Sigma \subset N$ onto some neighborhood of $u_k(\Sigma)$, 
such that $\Psi_k|_{\Sigma} = u_k$.
Let $\widehat{J}_k = \Psi_k^*J_k$ for $k=1,2,3,\ldots,\infty$,
so that for $k$ sufficiently large, the curves $v_k$ can be
identified with $\widehat{J}_k$-holomorphic curves in the total space of~$N$,
and each $u_k$ is identified with the zero section.

Let $\pi_N : u_\infty^*TM \to N$ denote the normal projection, so that
$\widehat{\nabla} := \pi_N \circ \nabla^\infty$ induces a 
connection on $N \to \Sigma$ 
(as a \emph{real} vector bundle), and thus defines a splitting
into horizontal and vertical subbundles
$$
TN = HN \oplus VN.
$$
This splitting is invariant under the diffeomorphisms on $N$ defined by real
scalar multiplication.  For $z \in \Sigma$ and $\eta \in N_z$,
the fibers in the splitting admit canonical identifications
$$
H_{(z,\eta)} N = T_z\Sigma, \qquad V_{(z,\eta)} N = N_z,
$$
and we can write $\widehat{J}_k$ with respect to the splitting as
$$
\widehat{J}_k(z,\eta) = \begin{pmatrix}
\alpha_k(z,\eta) & \beta_k(z,\eta) \\
\gamma_k(z,\eta) & \delta_k(z,\eta)
\end{pmatrix},
$$
for some smoothly varying linear maps $\alpha_k(z,\eta) : T_z\Sigma \to T_z\Sigma$,
$\beta_k(z,\eta) : N_z \to T_z\Sigma$ and so forth.
Since $u_k : (\Sigma,j_k) \to (M,J_k)$ is $J_k$-holomorphic 
and the fibers of $N_{u_k}$ are $J_k$-invariant
along~$u_k$, we have
$$
\alpha_k(z,0) = j_k(z), \quad \delta_k(z,0) = J_k(u_k(z)), \quad 
\beta_k(z,0) = 0, \quad \gamma_k(z,0) = 0.
$$
Now for any constant $r > 0$, the diffeomorphism
$$
\Phi_r : N \to N : (z,\eta) \mapsto (z,r\eta)
$$
transforms $\widehat{J}_k$ to 
$$
\widehat{J}_k^r(z,\eta) := \Phi_r^*\widehat{J}_k|_{(z,\eta)} = \begin{pmatrix}
\alpha_k(z,r\eta) & r \beta_k(z,r\eta) \\
\frac{1}{r} \gamma_k(z,r\eta) & \delta_k(z,r\eta)
\end{pmatrix},
$$
so given any positive sequence $r_k \to 0$, the sequence
$\widehat{J}_k^{r_k}$ converges in $C^\infty$ on compact subsets
of $N$~to
\begin{equation}
\label{eqn:J0}
\widehat{J}_\infty^0(z,\eta) := \begin{pmatrix}
j_\infty(z) & 0 \\
d \gamma_\infty(z,0) \eta & J_\infty(u_\infty(z))
\end{pmatrix}.
\end{equation}

\begin{lemma}
\label{lemma:tame}
A neighborhood of $\Sigma$ in $N$ admits a
symplectic form $\omega$ that tames $\widehat{J}_\infty^0$.
\end{lemma}
\begin{proof}
We use a variation on Thurston's method for constructing symplectic forms on fibrations
(cf.~\cite{McDuffSalamon:ST3}*{Theorem~6.1.4}).  For any open subset
$\uU \subset \Sigma$, let $\Lambda(\uU)$ denote the space of smooth $1$-forms
$\lambda$ on $\pi^{-1}(\uU)$ satisfying the following conditions:
\begin{enumerate}
\renewcommand{\labelenumi}{(\roman{enumi})}
\item At any point $(z,0) \in \uU \subset N|_\uU$ in the zero section,
$$
\lambda|_{(z,0)} = 0 \quad\text{ and }\quad
d\lambda|_{T_z\Sigma \times N_z} = 0;
$$
\item The restriction of $d\lambda$ to fibers in $\pi^{-1}(\uU)$ defines a
symplectic vector bundle structure on $N|_{\uU}$ taming~$J_\infty$.
\end{enumerate}
We observe that $\Lambda(\uU)$ is nonempty whenever there exists a
complex trivialization of $(N,J_\infty)$ over~$\uU$, and moreover,
it is \emph{$C^\infty$-convex} in the sense that if $\lambda_0, \lambda_1
\in \Lambda(\uU)$, then
$$
(\psi \circ \pi) \lambda_1 + (1 - \psi \circ \pi) \lambda_0 \in
\Lambda(\uU)
$$
for every smooth function $\psi : \uU \to [0,1]$.
It follows that an element of $\Lambda(\Sigma)$ can be constructed by
patching together local constructions via a partition of unity.

Now given $\lambda \in \Lambda(\Sigma)$, choose an area
form $\sigma$ on $\Sigma$ taming $j_\infty$.  Then for a sufficiently
large constant $K > 0$,
$$
\omega := K \pi^*\sigma + d\lambda
$$
is a closed $2$-form that tames $\widehat{J}_\infty^0$ at $\Sigma$ and hence 
also in a neighborhood of~$\Sigma$.
\end{proof}
\begin{remark}
The above proof did not use any special properties of $\widehat{J}_\infty^0$
except that the zero section is pseudoholomorphic and the normal fibers
along the zero section are also complex.  The same argument shows that
for any embedded closed $J$-holomorphic curve in any almost complex
manifold $(M,J)$, a neighborhood of the curve admits a symplectic form
that tames~$J$.
\end{remark}

\begin{lemma}
\label{lemma:normalSections}
Suppose $\psi : \widetilde{\Sigma} \to \Sigma$ is a smooth map,
$\tilde{\jmath}$ is a complex structure on $\widetilde{\Sigma}$, and
$\xi \in \Gamma(\psi^*N)$ is a smooth section along~$\psi$.  Then
the map $z \mapsto \xi(z)$ from $\widetilde{\Sigma}$ into the
total space of $N$
is a pseudoholomorphic map $(\widetilde{\Sigma},\tilde{\jmath}) \to
(N,\widehat{J}_\infty^0)$ if and only if
$\psi : (\widetilde{\Sigma},\tilde{\jmath}) \to (\Sigma,j_\infty)$ is
holomorphic and $\xi \in \ker \mathbf{D}_{u_\infty \circ \psi}^N$.
\end{lemma}
\begin{proof}
Denote by $v : \widetilde{\Sigma} \to N$ the smooth map into the
total space of $N$ defined by $v(z) := \xi(z) \in N_{\psi(z)} \subset N$.
Then using \eqref{eqn:J0}, the equation 
$Tv + \widehat{J}_\infty^0 \circ Tv \circ \tilde{\jmath} = 0$
translates into the two equations
$$
d\psi(z) + j_\infty(\psi(z)) \circ d\psi(z) \circ \tilde{\jmath}(z) = 0,
$$
and
$$
\widehat{\nabla} \eta(z) + J_\infty(u_\infty(\psi(z))) \circ
\widehat{\nabla} \eta(z) \circ \tilde{\jmath} +
\left[ d\gamma_\infty(\psi(z),0) \eta(z) \right] d\psi(z) \circ \tilde{\jmath} = 0
$$
for $z \in \widetilde{\Sigma}$.
The first equation says that $\psi : (\widetilde{\Sigma},\tilde{\jmath}) \to
(\Sigma,j_\infty)$ is holomorphic, and under this assumption, the second
matches $\mathbf{D}_{u_\infty \circ \psi}^N \eta = 0$ after observing
$$
[d\gamma_\infty(\psi,0)\eta] \circ d\psi \circ \tilde{\jmath} = 
\pi_N \circ (\nabla_\eta J_\infty) \circ T(u_\infty \circ \psi) \circ \tilde{\jmath}.
$$
\end{proof}

We now prove Proposition~\ref{prop:isolated} as
follows.  Arguing by contradiction, assume after taking a subsequence
that the curves $v_k : (\widetilde{\Sigma},\tilde{\jmath}_k) \to
(M,J_k)$ are not covers of $u_k$ for any~$k$ as $k \to \infty$.
Choose a symplectic form $\omega$ near the zero section in
$N = N_{u_\infty}$ as given by Lemma~\ref{lemma:tame}, and choose
$\delta > 0$ such that
$\omega$ tames $\widehat{J}_\infty^0$ on $\{ \eta \in N\ |\ |\eta| < 2\delta \}$.
Writing $v_k(z) = \xi_k(\psi_k(z))$ for sequences
$\psi_k : \widetilde{\Sigma} \to \Sigma$ and $\xi_k \in \Gamma(\psi_k^*N)$,
we have
$$
r_k := \frac{1}{\delta} \max_{z \in \widetilde{\Sigma}} | \xi_k(z) | > 0
$$
and $r_k \to 0$ by assumption.  Then
$$
w_k := \Phi_{r_k}^{-1} \circ v_k : (\widetilde{\Sigma},\tilde{\jmath}_k) \to
(N , \widehat{J}_k^{r_k})
$$
is a sequence of smooth pseudoholomorphic curves in {a compact subset of} 
the neighborhood
$\{ \eta \in N\ |\ |\eta| < 2\delta \}$, which can be written as
$w_k(z) = \eta_k(\psi_k(z))$ where $\eta_k = \frac{1}{r_k} \xi_k$
satisfies
\begin{equation}
\label{eqn:nonzero}
\max_{z \in \widetilde{\Sigma}} | \eta_k(z)| = \delta.
\end{equation}
Note that since $v_k$ converges to a nodal curve in $\widebar{\mM}_h(d ; u_\infty)$,
we can also assume the maps $\psi_k : \widetilde{\Sigma} \to \Sigma$ have
fixed degree~$d$.  Then
since $\widehat{J}_k^{r_k} \to \widehat{J}_\infty^0$ and the latter is
tamed by $\omega$ in the region under consideration, Gromov compactness
applies to $w_k$ and yields a subsequence convergent to a stable nodal
curve $w_\infty \in \widebar{\mM}_h(d[\Sigma],\widehat{J}_\infty^0)$.
By Lemma~\ref{lemma:normalSections}, each {smooth} component $w$ of
$w_\infty$ has the form $w(z) = \eta(\psi(z))$ where
$\psi : (\widetilde{\Sigma},\tilde{\jmath}) \to (\Sigma,j_\infty)$ is
holomorphic and $\mathbf{D}_{u_\infty \circ \psi}^N \eta = 0$.
We claim there must be at least one such component for which
$\deg(\psi) > 0$ and $\eta \ne 0$.   Indeed, \eqref{eqn:nonzero}
implies that there is at least one component with $\eta \ne 0$.
If every such component also satisfies $\deg(\psi) = 0$, then $\eta$
is a nonzero constant on this component, as the normal operator
$\mathbf{D}_{u_\infty \circ \psi}^N$ is simply the standard Cauchy-Riemann
operator on a trivial bundle when $\psi$ is constant.  But
since $\deg(\psi_k) = d > 0$, any component with $\deg(\psi) = 0$
is necessarily connected by a chain of nodes to another component
with $\deg(\psi) > 0$, and on this component, $\eta$ is nonzero
at the nodal point.  This implies the existence of a nontrivial
element $\eta \in \ker \mathbf{D}_{u_\infty \circ \psi}^N$ for some
positive degree holomorphic cover~$\psi$, and thus violates
super-rigidity.  The proof of Proposition~\ref{prop:isolated} is
complete.

\section{The Sard-Smale theorem for $C^\infty$-subvarieties}
\label{sec:subvarieties}

The proof of Petri's condition in \S\ref{sec:unique} requires
a version of the Sard-Smale theorem for objects
that are not Banach manifolds but are almost as nice
in some analytically quantifiable sense.  The results in
this appendix are easy consequences of standard results
in the analysis of smooth Banach manifolds,
but expresed in a slightly more general framework.

Suppose $X$ is a smooth Banach manifold and $Y \subset X$ is a
subset.  Given $k \in \NN$, we will say that $Y$ is a 
\defin{$C^\infty$-subvariety of codimension at least~$k$} if for every $x \in Y$,
there exists a neighborhood $\uU \subset X$ of $x$, a finite-dimensional
vector space $V$ and a smooth map $f : \uU \to V$ such that:
\begin{enumerate}
\item $Y \cap \uU = f^{-1}(0)$;
\item $\rk df(x) \ge k$.
\end{enumerate}

\begin{prop}
\label{prop:localSubmanifold}
If $Y \subset X$ is a $C^\infty$-subvariety of codimension at least~$k$, then for every
$x \in Y$, there exists a smooth Banach submanifold $\widetilde{Y} \subset X$ 
of codimension~$k$ such that a neighborhood of $x$ in $Y$ is contained 
in~$\widetilde{Y}$.
\end{prop}
\begin{proof}
Given $x \in Y$, we have $Y \cap \uU = f^{-1}(0)$ for some open neighborhood
$x \in \uU \subset X$ and smooth map $f : \uU \to V$, with $V$ a 
finite-dimensional vector space and $\dim \im df(x) \ge k$.  Then we can
choose a linear map $\Lambda : V \to \RR^k$ whose restriction to
$\im df(x) \subset V$ is surjective onto~$\RR^k$, hence
$\Lambda \circ df(x) : T_x X \to \RR^k$ is surjective.  Define
$\widetilde{Y} \subset X$ to be a neighborhood of $x$ in
$(\Lambda \circ f)^{-1}(0)$.  The implicit function theorem implies that this
is a Banach submanifold of codimension~$k$ if the neighborhood is taken
sufficiently small.
\end{proof}

The discussion so far makes sense under a very unrestrictive definition of
the term ``Banach manifold,'' e.g.~in \cite{Lang:geometry}, such objects
need not even be Hausdorff.  In practice, of course, the Banach manifolds
one encounters in applications are typically at least metrizable
(hence Hausdorff and paracompact) and separable.  The latter is the
condition required for the Sard-Smale theorem \cite{Smale:Sard}.
We will need the following standard bit of general topology:

\begin{lemma}
\label{lemma:countable}
If $X$ is a paracompact and separable topological space, then every open
cover of $X$ has a countable subcover.  \qed
\end{lemma}

The following is the main result of this appendix.  The proof of Theorem~\ref{thm:Petri}
uses the special case in which all manifolds are finite dimensional,
so the Fredholm assumption is automatic and only the finite-dimensional
version of Sard's theorem is needed.  The infinite-dimensional version with
the Sard-Smale theorem is required for the proof of Theorem~\ref{thm:PetriJ}.

\begin{prop}
\label{prop:SardSmaleVariety}
Assume $\univ$ and $Z$ are separable and metrizable smooth Banach manifolds,
$\pi : \univ \to Z$ is a smooth Fredholm map, and $X \subset \univ$ is a
$C^\infty$-subvariety of codimension at least~$k \in \NN$.  For each $z \in Z$,
denote
$$
\mM(z) := \pi^{-1}(z) \subset \univ, \qquad X(z) := X \cap \mM(z) \subset \mM(z),
$$
and let $Z^\reg_\pi \subset Z$
denote the Baire subset consisting of regular values of~$\pi$.  Then there
exists a further Baire subset $Z^\reg_X \subset Z$ such that for all
$z \in Z^\reg_\pi \cap Z^\reg_X$, $X(z)$ is a $C^\infty$-subvariety
of codimension at least $k$ in~$\mM(z)$.
\end{prop}
\begin{proof} 
Suppose $x \in X$, so by assumption, there exists a neighborhood
$$
x \in \univ_x \subset \univ,
$$
a finite-dimensional vector space $V_x$ and a smooth map
$f_x : \univ_x \to V_x$ such that
$f_x^{-1}(0) = X \cap \univ_x$ and $\rank df_x(x) \ge k$.  After possibly
shrinking $\univ_x$ to a smaller neighborhood of~$x$,
we can use the argument
in the proof of Proposition~\ref{prop:localSubmanifold} to find a
linear map $\Lambda_x : V_x \to \RR^k$ such that $0 \in \RR^k$ is a
regular value of $\Lambda_x \circ f_x : \univ_x \to \RR^k$ and
$$
\widetilde{\univ}_x := (\Lambda_x \circ f_x)^{-1}(0) \subset \univ
$$
is a smooth Banach submanifold of codimension~$k$ containing
$X \cap \univ_x$.

Since $\univ$ is metrizable and separable, $X$ also has
both of these properties, thus Lemma~\ref{lemma:countable} implies that
we can find a sequence $\{x_n \}_{n=1}^\infty$ of points in $X$
such that every $x \in X$ lies in at least one of the neighborhoods~$\univ_{x_n}$.
Let $Z^\reg_n \subset Z$ denote the set of regular
values of the projection
$$
\widetilde{\univ}_{x_n} \stackrel{\pi}{\longrightarrow} Z,
$$
The latter is a smooth Fredholm map since $\widetilde{\univ}_{x_n}$ is a
smooth finite-codimensional submanifold of~$\univ$.  
The Sard-Smale theorem thus implies that
$Z^\reg_n \subset Z$ is a Baire subset, and consequently,
$$
Z^\reg_X := \bigcap_{n=1}^\infty Z^\reg_n \subset Z
$$
is also a Baire subset.

Now for any $z \in Z^\reg_X \cap Z^\reg_\pi$ and $x \in X(z)$, pick
$n \in \NN$ such that $x \in \univ_{x_n}$, and consider the restricted map
$$
g_n : \mM(z) \cap \univ_{x_n} \to V_{x_n} : x \mapsto f_{x_n}(x),
$$
whose zero-set is a neighborhood of $x$ in~$X(z)$.
Regularity and the implicit
function theorem imply that $\widetilde{\univ}_{x_n} \subset \univ$ and
$\mM(z) \subset \univ$ are transverse submanifolds,
so that $0$ is also a regular value
of $\Lambda_{x_n} \circ g_n : \mM(z) \cap \univ_{x_n} \to \RR^k$.
It follows that $\Lambda_{x_n} \circ d g_n(x) : T_x \mM(z) \to \RR^k$ is 
surjective, and thus $\rank d g_n(x) \ge k$.
\end{proof}

The results of this discussion combine to yield the following useful
consequence:

\begin{cor}
In the setting of Proposition~\ref{prop:SardSmaleVariety}, if
the smooth Fredholm map $\pi : \univ \to Z$ satisfies $\ind d\pi(x) < k$ for all
$x \in \univ$, then $X(z)$ is empty for generic $z \in Z$.
\qed
\end{cor}

\section{History of errors}
\label{sec:catastrophe}

This appendix has been added (at the
suggestion of an anonymous referee) in the interest of transparency: its
purpose is to clarify more precisely what went wrong with previous attempts
to prove Theorem~\ref{thm:super}, and how those attempts are related to the
proof in this paper.  There were at least two claims of proofs of super-rigidity
that were publicized and then withdrawn before I ever started thinking about
the problem, but since it is not my place to comment on those, I will only
discuss the attempts that I have been involved in.

\subsection{Analytic perturbation theory}

The original version of \cite{GerigWendl} was a preprint under a different
title \cite{GerigWendl:v1}, which claimed a proof of Theorem~\ref{thm:super}
(also in dimension four) for embedded index~$0$ curves that are fully contained in
the perturbation domain $\uU \subset M$.  The ideas behind that argument were
almost totally disjoint from those of the present paper, excepting the
superficial feature that both derive originally from (separate) ideas
developed in Taubes's work on the Gromov invariant.
The literature on the Gromov invariant contains two quite different methods
to prove transversality for the doubly covered tori that must be counted:
one (from \cite{Taubes:counting}) is based on a splitting
of Cauchy-Riemann type operators with respect to irreducible representations,
and gives rise to dimension-counting arguments that
provided the original inspiration for this paper.  The other, from 
\cite{Taubes:SWtoGr}*{Proof of Prop.~7.1, Step~7}, is in some respects more novel: it is based on
a Weitzenb\"ock formula for Cauchy-Riemann type operators and analytic
perturbation theory.  In the setting of \cite{Taubes:SWtoGr}, where one needs
to prove that a $\ZZ_2$-equivariant index~$0$ Cauchy-Riemann type operator $\mathbf{D} :
\Gamma(E) \to \Gamma(F)$ on
a trivial line bundle $E \to \TT^2$ can 
always be perturbed equivariantly to one that is
invertible, these two ingredients combine in the following way:
\begin{enumerate}
\item The Weitzenb\"ock formula implies that for any complex-antilinear
bundle isomorphism $A : E \to F$, the deformed operator 
$\mathbf{D}_\tau := \mathbf{D} + \tau A$ is invertible for all $\tau \gg 0$.
\item Since the deformed operators $\mathbf{D}_\tau$ depend analytically on
the parameter $\tau \in \RR$, analytic perturbation theory as in \cite{Kato}
implies that the set $\{ \tau \in \RR\ |\ \text{$\mathbf{D}_\tau$ is not invertible} \}$
is either $\RR$ or is discrete.  The first possibility has already been ruled out
via the Weitzenb\"ock formula, so it follows that $\mathbf{D}_\tau$ is invertible
for all $\tau \ne 0$ in some neighborhood of~$0$.
\end{enumerate}
This technique has the appealing feature that it does not care how symmetric
the perturbation term $A \in \Gamma(\overline{\Hom}_\CC(E,F))$ is, thus it can work equally well for
simple holomorphic curves and multiple covers.  The preprint \cite{GerigWendl:v1}
was motivated by the insight that both parts of the argument can be made to
work somewhat more generally: the operator $\mathbf{D}$ can have negative
index if we talk about \emph{injectivity} of $\mathbf{D}_\tau$ instead of
invertibility, and $E$ can also be a higher-rank bundle if $A$ is required
to satisfy an extra condition which, for topological reasons, can be assumed
without loss of generality.  Applying the argument to normal Cauchy-Riemann
operators of branched covers then produces the following result:

\begin{lemma}[\cite{GerigWendl}]
\label{lemma:whatWeProved}
Suppose $\dim M \ge 4$, $J \in \jJ(M,\omega\,;\,\uU,\Jfix)$,
$v : (\Sigma,j) \to (M,J)$ is an embedded closed $J$-holomorphic curve of index~$0$
with image contained in~$\uU$,
and $u = v \circ \varphi$ where $\varphi : (\widetilde{\Sigma},\widetilde{\jmath}) \to (\Sigma,j)$
is a holomorphic branched cover of degree $d \in \NN$ between closed connected
Riemann surfaces.  Then there exists a smooth $1$-parameter family
$\{J_\tau \in \jJ(M,\omega\,;\,\uU,\Jfix)\}_{\tau \in (-\epsilon,\epsilon)}$
such that $J_0=J$, $v$ and $u$ are $J_\tau$-holomorphic for every~$\tau$, and
the resulting normal Cauchy-Riemann operators $\mathbf{D}_{u,\tau}^N$
for $u$ with respect to $J_\tau$ are injective for all $\tau \ne 0$.
\qed
\end{lemma}

A proof of generic super-rigidity would follow via relatively straightforward
topological arguments if one instead had the following stronger statement:\footnote{The question
mark in the statement indicates that I do not know whether Lemma~\ref{lemma:whatWeShouldHaveProved}
is true, and I do not have a strong enough opinion about it to call it a conjecture.}

\begin{lemmq}
\label{lemma:whatWeShouldHaveProved}
In the setting of Lemma~\ref{lemma:whatWeProved}, the family
of almost complex structures
$\{J_\tau \in \jJ(M,\omega\,;\,\uU,\Jfix) \}_{\tau \in (-\epsilon,\epsilon)}$ can be chosen so that for some neighborhood $\oO(\varphi)$ of
$\varphi$ in the moduli space of $d$-fold holomorphic branched covers,
the normal Cauchy-Riemann operators $\mathbf{D}_{v \circ \varphi',\tau}^N$
are injective for all $\tau \ne 0$ and $\varphi' \in \oO(\varphi)$.
\end{lemmq}

Unfortunately, Lemma~\ref{lemma:whatWeProved} does not imply
Lemma~\ref{lemma:whatWeShouldHaveProved}, as analytic perturbation theory gives 
no obvious way to control
the size of the range of parameter values $\tau \in (-\epsilon,\epsilon) \setminus \{0\}$
for which injectivity is guaranteed as $\varphi$ varies in the moduli space
of branched covers.  This detail was overlooked in \cite{GerigWendl:v1};
the crucial gap in our argument was pointed out by Ionel and Parker.
What can still be salvaged from Lemma~\ref{lemma:whatWeProved}, and eventually
appeared as the main result of the published paper \cite{GerigWendl}, is a
result similar to Theorem~\ref{thm:unbranched} about transversality for \emph{unbranched} covers: 
in the unbranched case there
is no distinction between Lemmas~\ref{lemma:whatWeProved} and~\ref{lemma:whatWeShouldHaveProved}
because the moduli space that $\varphi$ lives in is discrete.

I currently believe the proof of Theorem~\ref{thm:super} originally attempted
in \cite{GerigWendl:v1} to be unsalvageable.  There are also strong philosophical
arguments for preferring the approach of the present paper over analytic
perturbation theory: notably, the use of the Weitzenb\"ock formula requires a
more \emph{global} class of perturbations ($u$ must be contained in the 
perturbation domain $\uU \subset M$ rather than merely intersecting it), and the whole
strategy seems completely unsuitable 
for studying the wall-crossing phenomena mentioned in \S\ref{sec:bifurcations}.
On the other hand, the Weitzenb\"ock argument (minus analytic perturbation theory)
has been usefully exploited by other authors in certain special settings
where geometric information removes the need to assume $\tau \gg 0$;
see \cites{LeeParker:structure,IonelParker:GV}.

\subsection{Earlier versions of the present paper}
\label{sec:catastropheThis}

The main ideas behind the proofs of Theorems~\ref{thm:super}--\ref{thm:submanifolds0}
have changed very little since the first version of this paper appeared on the arXiv, but
one important technical detail has changed a lot: the proof that generic Cauchy-Riemann
type operators satisfy Petri's condition.\footnote{The term
``Petri's condition'' did not appear in the first three versions of this paper
on the arXiv, but the same notion was there under the label of
``unique continuation for tensor products'' and has sometimes also been
advertised as ``quadratic unique continuation''.  The current terminology
was introduced by Doan and Walpuski \cite{DoanWalpuski:BrillNoether}
after the first version of this paper appeared.}  The intuition from the
beginning had been that Petri's condition was the main analytical lemma needed for the
proof of Theorem~\ref{thm:submanifolds0} (on which Theorems~\ref{thm:super}--\ref{thm:branched} all depend), 
and that it should hold due to
unique continuation except for some special class of non-generic Cauchy-Riemann
type operators.  Up to version~3 on the arXiv \cite{Wendl:superV3}, a much
more naive approach to this lemma was taken, in which the word ``generic''
was given a precise characterization:

\begin{lemmf}[\cite{Wendl:superV3}*{Corollary~5.2 and Lemma~3.11}]
\label{lemma:falseLemma}
Suppose $E,F \to \Sigma$ are complex vector bundles and $\mathbf{D} : \Gamma(E) \to \Gamma(F)$ is a Cauchy-Riemann
type operator such that the bundle map $\mathbf{D}^{0,1} \in\Gamma(\overline{\Hom}_\CC(E,F))$
given by the complex-antilinear part of $\mathbf{D}$ defines an invertible map
$E_z \to F_z$ at some point $z \in \Sigma$.  Then $\mathbf{D}$ satisfies Petri's condition
to infinite order at~$z$.
\end{lemmf}

It is relatively easy to show (see \cite{Wendl:superV3}*{Lemma~6.2})
that the hypothesis on invertibility of
complex-antilinear parts is generic, i.e.~all normal Cauchy-Riemann operators
of $J$-holomorphic curves satisfy it for generic (and necessarily non-integrable)~$J$.  The benefit of this
condition is that it forces $\ker\mathbf{D} \subset \Gamma(E)$ and 
$\ker\mathbf{D}^* \subset \Gamma(F)$ to be
totally real subspaces, meaning that any real-linearly independent set of
vectors in one of these spaces is also complex-linearly independent.
The original reason to believe in Lemma~\ref{lemma:falseLemma} was
the elementary observation mentioned in Example~\ref{ex:complexPetri} that for complex-linear Cauchy-Riemann type operators,
which can always be expressed locally as the standard one, the complex
version of Petri's condition (involving complex tensor products) does hold to
infinite order at every point;
a proof of this may be found on page~48
of \cite{Wendl:superV3}.  Lemma~\ref{lemma:falseLemma} was thus an attempt to
fit real-linear Cauchy-Riemann type operators into a complex-linear context
with the aid of the totally real hypothesis.
The proof was destroyed by a careless mistake in linear algebra:
Equations~(5.3) and~(5.4) in \cite{Wendl:superV3} define certain functions
$\eta^\nu_\alpha$ and $\xi^\mu_\beta$ that are meant to be in $\ker\mathbf{D}$
and $\ker\mathbf{D}^*$ respectively because they are linear combinations 
of functions in those spaces, but in fact, the coefficients in those linear
combinations are complex rather than real, while $\mathbf{D}$ and $\mathbf{D}^*$
are only real-linear.  Similarly, the claim in the final paragraph
of that proof that certain linear combinations $\sum_i c^{ij} \xi_i$ and
$\sum_j c^{ij}\eta_j$ satisfy linear Cauchy-Riemann or anti-Cauchy-Riemann
equations does not hold, again because the coefficients $c^{ij}$ are complex instead
of real.  These errors were noticed by Doan and Walpuski while working on
their own alternative exposition of the super-rigidity proof \cite{DoanWalpuski:BrillNoether}.
Example~\ref{ex:realPetri2} was found later, showing that
Lemma~\ref{lemma:falseLemma} is in fact false.

After Lemma~\ref{lemma:falseLemma} fell apart, the intuition remained that
the failure of the local Petri condition for a Cauchy-Riemann type operator
should be overdetermined in some sense, and the jet space approach in
the current \S\ref{sec:unique} was then developed to make this intuition
precise.  Lemma~\ref{lemma:falseLemma} has now been replaced by
Corollary~\ref{cor:Petri}, whose proof is completely different from what was
attempted in \cite{Wendl:superV3}, and has an additional advantage over the
earlier approach in that the jet space formalism can potentially be
applied to more general classes of operators beyond Cauchy-Riemann
(\S\ref{sec:jetSpaces} has been written with this in mind).
A more detailed informal discussion of
the fix may be found in the blog post \cite{Wendl:blogPetri}.

For completeness, I should mention a somewhat serious but non-fatal error
that was also pointed out by Doan and Walpuski but corrected between arXiv versions 2 and~3 of this paper.
The definition of the walls appearing in Theorem~\ref{thm:submanifolds0}
was slightly wrong in earlier versions, because it was overlooked that
in the splitting of the normal Cauchy-Riemann operator $\mathbf{D}_u^N$ into
summands $\dot{\mathbf{D}}_{u,\boldsymbol{\theta}_i}^N$ corresponding to
irreducible representations~$\boldsymbol{\theta}_i$, the kernels and cokernels
of these summands are always modules over the equivariant endomorphism algebra
($\RR$, $\CC$ or $\HH$) of~$\boldsymbol{\theta}_i$, and this structure must
be respected in talking about their dimensions.  The result was a mistake in
\cite{Wendl:superV2}*{Theorem~D} that was hard to spot, because the statement
looked the same as in the current version, but its meaning was different.
The source of the problem was
an erroneous representation-theoretic dimension calculuation in \cite{Wendl:superV2}*{Corollary~3.23},
which was stated without proof.
A corrected version of that result appears in this version as
Corollary~\ref{cor:codimension}, with a proof given in the preceding paragraph.

\end{document}